\documentclass[11pt,reqno]{amsart}

\newif\ifdraftnotes
\draftnotesfalse

\usepackage[margin=1.15in]{geometry}
\usepackage[T1]{fontenc}
\usepackage[utf8]{inputenc}
\usepackage{lmodern}
\usepackage{microtype}
\usepackage{amsmath,amssymb,amsthm,mathtools,mathrsfs,bm}
\mathtoolsset{showonlyrefs=true, showmanualtags=true}
\numberwithin{equation}{section}
\allowdisplaybreaks
\usepackage{graphicx}
\usepackage{booktabs,longtable,array,tabularx,multirow}
\usepackage{enumitem}
\usepackage{xcolor}
\usepackage{comment}
\usepackage[most]{tcolorbox}
\usepackage{tikz}
\usetikzlibrary{arrows.meta,positioning,calc,fit}

\setlist[enumerate]{leftmargin=2.2em}
\setlist[itemize]{leftmargin=1.9em}

\usepackage[hidelinks]{hyperref}

\hypersetup{
  pdfauthor={Bonan Chen, Xiang Fang, Yajun Zhang},
  pdftitle={A Joint Functional Central Limit Theorem for Lp-Norms of Random Littlewood Polynomials},
  pdfsubject={Joint exponent and logarithmic-time fluctuations of random Littlewood polynomials},
  pdfkeywords={Littlewood polynomial, Lp norm, functional central limit theorem, logarithmic time, Laguerre polynomial, Rademacher chaos}
}
\usepackage{lineno}

\usepackage{hyperref} 
\newtheorem{theorem}{Theorem}[section]
\newtheorem{proposition}[theorem]{Proposition}
\newtheorem{lemma}[theorem]{Lemma}
\newtheorem{corollary}[theorem]{Corollary}

\theoremstyle{definition}
\newtheorem{definition}[theorem]{Definition}

\theoremstyle{remark}
\newtheorem{remark}[theorem]{Remark}

\newcommand{\E}{\mathbb{E}}
\newcommand{\Prob}{\mathbb{P}}
\newcommand{\R}{\mathbb{R}}
\newcommand{\C}{\mathbb{C}}

\newcommand{\Z}{\mathbb{Z}}
\newcommand{\Torus}{\mathbb{T}}
\newcommand{\dd}{\mathop{}\!\mathrm{d}}
\newcommand{\ind}{\mathbf{1}}

\newcommand{\poly}{\mathrm{poly}}
\newcommand{\fl}{\mathrm{fl}}
\newcommand{\ec}{\mathrm{ec}}
\newcommand{\dc}{\mathrm{dc}}
\newcommand{\ren}{\mathrm{ren}}

\newcommand{\LagProj}{\Pi^{\mathrm{Lag}}}
\newcommand{\LagGen}{\mathcal{L}_{\mathrm{Lag}}}
\newcommand{\LagNum}{\mathcal{N}_{\mathrm{Lag}}}

\DeclareMathOperator{\Cov}{Cov}
\DeclareMathOperator{\Var}{Var}
\DeclareMathOperator{\supp}{supp}

\newcommand{\boldparagraph}[1]{\par\addvspace{0.5em}
  \paragraph{\textbf{#1}}
}

\ifdraftnotes
  
  \newcommand{\DraftNote}[1]{%
    \par\smallskip\noindent
    \begin{tcolorbox}[
      colback=yellow!7!white,
      colframe=yellow!50!black,
      boxrule=0.4pt,
      arc=1mm,
      left=1mm,right=1mm,top=1mm,bottom=1mm]
      \small\raggedright\textbf{Draft note.} #1
    \end{tcolorbox}%
    \par\smallskip}
  \newcommand{\StatementPlaceholder}[1]{%
    \textit{Skeleton statement.} #1}
\else
  \excludecomment{sectioncontract}
  \newcommand{\DraftNote}[1]{}
  \newcommand{\StatementPlaceholder}[1]{%
    \PackageError{joint-pt-publication}{Unreplaced theorem placeholder}{Replace every placeholder before setting draftnotesfalse.}}
\fi

\makeatletter
\def\@tocline#1#2#3#4#5#6#7{\relax
  \ifnum #1>\c@tocdepth %
  \else
    \par \addpenalty\@secpenalty\addvspace{#2}%
    \begingroup \hyphenpenalty\@M
    \@ifempty{#4}{%
      \@tempdima\csname r@tocindent\number#1\endcsname\relax
    }{%
      \@tempdima#4\relax
    }%
    \parindent\z@ \leftskip#3\relax \advance\leftskip\@tempdima\relax
    \rightskip\@pnumwidth plus4em \parfillskip-\@pnumwidth
    #5\leavevmode\hskip-\@tempdima
      \ifcase #1
       \or\or \hskip 2.5em \or \hskip 2em \else \hskip 3em \fi%
      #6\nobreak\relax
    \hfill\hbox to\@pnumwidth{\@tocpagenum{#7}}\par%
    \nobreak
    \endgroup
  \fi}
\makeatother

\title[Joint functional CLT for Littlewood $L^p$-norms]
{A Joint Functional Central Limit Theorem for the $L^p$-Norm Profile of Random Littlewood Polynomials}

\author[B. Chen]{Bonan Chen}
\address{
School of Mathematical Sciences,
Soochow University,
Suzhou, 215006, P. R. China
}
\email{bnchen@suda.edu.cn}

\author[X. Fang]{Xiang Fang}
\address{
Department of Applied Mathematics,
National Yang Ming Chiao Tung University,
Hsinchu, 300093, Taiwan
}
\email{xfang@nycu.edu.tw}

\author[Y. Zhang]{Yajun Zhang}
\address{School of Mathematical Sciences,
Dalian University of Technology,
Dalian, 116024, P. R. China}
\email{zhang-yajun@mail.dlut.edu.cn}


\date{}

\begin{document}
\begin{abstract}
For a random Littlewood polynomial
\[
q_n(z)=\sum_{j=0}^{n-1}\varepsilon_jz^j,
\]
with independent Rademacher coefficients \((\varepsilon_j)\), the
normalized angular intensity
\[
S_n(\theta)=\frac{|q_n(e^{i\theta})|^2}{n}
\]
has total mass one with respect to normalized arc measure \(\mu\) on
\(\mathbb T\). Its positive moments
\[
T_n^{(p)}
=
\int_{\mathbb T}S_n(\theta)^{p/2}\,d\mu(\theta)
=
\frac{\|q_n\|_p^p}{n^{p/2}}
\]
form the moment profile of the intensity distribution as a function of \(p\). We prove a
joint functional central limit theorem for its fluctuations: for \(p\)
in a fixed compact interval \(K\subset(0,\infty)\) and logarithmic time
\[
n=\lfloor Ne^t\rfloor,\qquad t\in[-T,T],
\]
the centered values
\[
\sqrt n\bigl(T_n^{(p)}-\Gamma(1+p/2)\bigr)
\]
have piecewise linear interpolations that converge to a Gaussian field
\(G(p,t)\) in \(C(K\times[-T,T])\).  The corresponding step
processes converge in \(D([-T,T];C(K))\) with the Skorokhod
\(J_1\)-topology.
The limiting field admits the representation
\[
G(p,t)=\sum_{a=2}^{\infty}c_a(p)G^{(a)}(t),
\]
where \(c_a(p)\) are the coefficients of \(x^{p/2}\) in the ordinary
Laguerre expansion, with polynomials normalized by \(L_a(0)=1\).
The \(G^{(a)}\) are independent stationary Gaussian modes with covariance
\[
\operatorname{Cov}(G^{(a)}(s),G^{(a)}(t))
=
2V_a e^{-(a-\frac12)|t-s|}.
\]
Here \(V_a\) is the squared \(L^2\)-norm of the density of the sum of
\(a\) independent uniform random variables on \([0,1]\).
The exponent changes the deterministic weights of these modes, while
logarithmic time describes their dependence across nested lengths.
We also obtain the functional limit theorem for the normalized
\(L^p\)-norms
\[
\frac{\|q_n\|_p}{\sqrt n}.
\]
At \(p=4\), we recover the known scalar autocorrelation fluctuation law
at fixed logarithmic time and identify the dependence of these
fluctuations across nested lengths. The theorem also gives the joint covariance of the
positive-exponent norm profile.
Finite-mode convergence follows from the multivariate functional de Jong
theorem after model-specific variance, influence, and fourth-moment
calculations.  The additional maximal estimates control the infinite
Laguerre tail and its exponent derivative uniformly over nested lengths.
\end{abstract}

\maketitle
\setcounter{tocdepth}{2}
\newpage
\tableofcontents
\newpage
\section{Introduction}
\label{sec:introduction}

\subsection{Littlewood polynomials and their \texorpdfstring{\(L^p\)}{Lp}-profiles}
\label{sec:intro-lp-profiles}

Littlewood polynomials have coefficients in
\(\{-1,1\}\). They form a classical family at the intersection of harmonic analysis,
probability, combinatorics, and analytic number theory;
see \cite{Montgomery2017}.
These polynomials are related to the design of binary sequences
with low autocorrelation for radar and communications
\cite{BrestBoskovic2021,DingTang2010}.
Problems concerning such sequences also appear in
computer science \cite{BoskovicBrest2017,GCF2009} and statistical 
physics \cite{Bernasconi1987,PackebuschMertens2016}.

A central goal in the study of Littlewood polynomials is to understand how the discrete structure of binary sequences shapes the global behavior of Fourier series and the waves they represent.
Their \(L^p\)-norms provide fundamental ways to quantify the resulting flatness and concentration of the waves. Indeed, understanding the delicate cancellation among the Fourier terms of Littlewood polynomials is a classical problem. In the past a focus is on the celebrated Littlewood's flatness problem which concerns how flat an extreme Littlewood polynomial can be; see \cite{Littlewood1966}, \cite{BalisterBollobasMorrisSahasrabudheTiba2020}.
In the present work, we seek to understand how typical Littlewood polynomials fluctuate, via a functional central limit theorem, with respect to two parameters varying at the same time: the length of the signal and the exponent of the measurement.
This is a genuine second-order probability problem.

In this paper we consider the random polynomials
\begin{equation}
q_n(z)=\sum_{j=0}^{n-1}\varepsilon_jz^j,
\qquad
\mathbb P(\varepsilon_j=1)=\mathbb P(\varepsilon_j=-1)=\frac12,
\qquad n\ge1,
\label{eq:intro-qn}
\end{equation}
where \((\varepsilon_j)_{j\ge0}\) are independent.  The same infinite
sign sequence is used for every \(n\), so increasing the length extends
the polynomial rather than resampling its coefficients.  The present
paper proves a joint functional central limit theorem for their
normalized \(L^p\)-norms, with both the exponent and the polynomial
length as parameters.

We identify \(\mathbb T\) with \(\mathbb R/(2\pi\mathbb Z)\) and write
\(\mathrm d\mu(\theta)=\mathrm d\theta/(2\pi)\).  For \(p>0\), let
\(\|q_n\|_p=(\int_{\mathbb T}|q_n(e^{i\theta})|^p\,
\mathrm d\mu(\theta))^{1/p}\), with the usual quasinorm interpretation
when \(0<p<1\).  Parseval's identity gives \(\|q_n\|_2^2=n\) for every
choice of signs.  We therefore introduce the normalized intensity and
its power integrals:
\begin{align}
S_n(\theta)
&=\frac{|q_n(e^{i\theta})|^2}{n},
\label{eq:intro-Sn}\\
T_n^{(p)}
&=\frac{\|q_n\|_p^p}{n^{p/2}}
=\int_{\mathbb T}S_n(\theta)^{p/2}\,\mathrm d\mu(\theta).
\label{eq:intro-Tnp}
\end{align}
Thus \(\int_{\mathbb T}S_n\,\mathrm d\mu=T_n^{(2)}=1\).  The total
intensity is fixed; what varies with the signs is its distribution
across the circle.

For each realization, the function \(p\mapsto T_n^{(p)}\) records the
positive moments of the intensity distribution obtained by
sampling the angle according to \(\mu\).  We call it the moment
profile.  The corresponding norm profile is the family of normalized
norms viewed as a function of the exponent.  It contains the same
information as these power integrals.  Both profiles depend only on
the distribution of the intensity values, not on their
locations on the circle.  Increasing \(p\) gives greater weight to the larger values,
while the function \(x^{p/2}\) changes from concave to convex as \(p\)
crosses \(2\).  Since the mean intensity is one, Jensen's inequality
gives \(T_n^{(p)}\le1\) for \(0<p<2\) and \(T_n^{(p)}\ge1\) for
\(p>2\). This poses essential difficulty in the analysis of small $p$. Indeed, varying \(p\) probes different features of one and the
same normalized intensity distribution.  Our joint fluctuation problem
asks how these measurements vary together.

The quartic norm has a concrete interpretation in terms of the
correlations of a binary sequence.  For \(n\ge2\), define the aperiodic
autocorrelations
\(C_{n,k}=\sum_{j=0}^{n-1-k}\varepsilon_j\varepsilon_{j+k}\),
\(1\le k<n\).  The Fourier coefficients of
\(|q_n(e^{i\theta})|^2\) at the nonzero frequencies \(\pm k\) are
\(C_{n,k}\).  A second application of Parseval's identity therefore gives
\begin{equation}
T_n^{(4)}-1
=\int_{\mathbb T}\bigl(S_n(\theta)-1\bigr)^2\,\mathrm d\mu(\theta)
=\frac{2}{n^2}\sum_{k=1}^{n-1}C_{n,k}^2.
\label{eq:intro-quartic-autocorrelation}
\end{equation}
The merit factor of the binary sequence is
\begin{equation}
\operatorname{MF}(q_n)
=\frac{n^2}{2\sum_{k=1}^{n-1}C_{n,k}^2}
=\frac{1}{T_n^{(4)}-1}.
\label{eq:intro-merit-factor}
\end{equation}
Thus \(T_n^{(4)}-1\) measures both the squared \(L^2\)-deviation of
\(S_n\) from the constant intensity \(1\) and the normalized sum of
squared autocorrelations.  Maximizing the merit factor, equivalently
minimizing the quartic norm at fixed length, is a classical problem
in binary-sequence design \cite{JedwabKatzSchmidt2013}.
 Up to normalization, the same sum of squared
autocorrelations is the energy minimized in Bernasconi's spin model
\cite{Bernasconi1987}.
In short, the significance of \(p=4\) is established by an exact correlation identity;
our work, on the other hand, illustrates how to pass from autocorrelation energy at \(p=4\) to the full intensity distribution as \(p\) varies.

Norm questions are also tied to deterministic flatness.  One asks for
Littlewood polynomials \(P_n\) of length \(n\) satisfying
\[
c\sqrt n\le |P_n(e^{i\theta})|\le C\sqrt n
\qquad(\theta\in\mathbb T),
\]
where \(0<c<C<\infty\) are independent of \(n\).  Balister,
Bollob\'as, Morris, Sahasrabudhe, and Tiba proved that such polynomials
exist for every sufficiently large \(n\)
\cite{BalisterBollobasMorrisSahasrabudheTiba2020}.
These are uniform two-sided bounds, not a requirement that the
normalized modulus converge uniformly to \(1\).  The random problem
considered here concerns typical coefficient choices rather than the
existence of such deterministic constructions.

More generally, one may study nonlinear angular integrals
\(\int_{\mathbb T}f(S_n)\,\mathrm d\mu\); the norm statistics correspond
to \(f(x)=x^{p/2}\).  Our focus is the norm profile itself, with the
exponent treated as a parameter and with the 
coefficient sequence followed through a multiplicative window of
lengths, meaning lengths between fixed positive multiples of a reference
length.  The finite-mode part of the argument applies more generally
to arbitrary fixed polynomial tests, but passing to the full
moment profile requires uniform control of the infinite
Laguerre tail and of its derivative with respect to the exponent; see
Section~\ref{sec:intro-proof-outline}.

\subsection{Previous results and the joint fluctuation problem}
\label{sec:intro-first-order-to-fluctuations}

The first-order theory identifies a deterministic limit for each
normalized \(L^p\)-statistic.  Borwein and Lockhart
\cite{BorweinLockhart2001} proved that, for random Littlewood
polynomials and every fixed \(p>0\),
\begin{equation}
\mathbb E T_n^{(p)}\longrightarrow
m_p:=\Gamma\!\left(1+\frac p2\right),
\label{eq:intro-first-order-mean}
\end{equation}
and that \(T_n^{(p)}\to m_p\) in probability.  They also obtained
\(\mathbb E\|q_n\|_p/\sqrt n\to m_p^{1/p}\).  The constant \(m_p\) is
the \(p/2\)-moment of an exponential random variable of mean one.  Related triangular-array
extensions of the moment-convergence argument were obtained by Cuny and
Weber \cite{CunyWeber2006}.  Choi and Erd\'elyi
\cite[Theorem~1.6]{ChoiErdelyi2014} obtained all positive moments of the normalized
norms: for every fixed \(p,\alpha>0\),
\[
\mathbb E\left(\frac{\|q_n\|_p}{\sqrt n}\right)^\alpha
\longrightarrow m_p^{\alpha/p}.
\]
In particular, this already implies convergence of the first two
moments of \(T_n^{(p)}\) with respect to the signs, without specifying the
rate or a fluctuation law.

Almost-sure convergence at the quartic exponent was proved by Duan,
Fang, and Zhan \cite{DuanFangZhan2024}.  The companion manuscript \cite{ChangChenDuanFangZhangAS} establishes, on one event
of probability one and simultaneously for every \(p>0\),
\begin{equation}
T_n^{(p)}\longrightarrow m_p
\quad\text{almost surely},
\qquad
\frac{\|q_n\|_p}{\sqrt n}\longrightarrow m_p^{1/p}
\quad\text{almost surely}.
\label{eq:intro-first-order-as}
\end{equation}
The companion also proves, for each fixed \(p>0\), the quantitative
bounds
\[
\mathbb E T_n^{(p)}=m_p+O_p(n^{-1}),
\qquad
\operatorname{Var}(T_n^{(p)})=O_p(n^{-1}).
\]
Here the subscript \(p\) records dependence of the implied constants on
the exponent.
The companion's smoothing, small-ball, Stein comparison, and Gaussian anisotropy
arguments are shared ingredients of the expectation analysis reproduced
in Appendix~\ref{app:expectation-envelope}.  There we verify a common
bound uniform over compact exponent intervals, including intervals
crossing \(p=2\).  First-order convergence and these moment bounds do
not determine the limiting fluctuation law, its mixed-exponent
covariance, or its evolution across nested lengths.

For the quartic statistic, the fluctuation law is already known from
work on binary sequences.  By
\eqref{eq:intro-quartic-autocorrelation}, the autocorrelation demerit
factor is \(T_n^{(4)}-1\).  The mean calculation of Sarwate
\cite{Sarwate1984}, equivalently the quartic calculation of Newman and
Byrnes \cite{NewmanByrnes1990}, and Jedwab's exact variance formula
\cite{Jedwab2024} give
\begin{equation}
\mathbb E T_n^{(4)}=2-\frac1n,
\qquad
\operatorname{Var}(T_n^{(4)})
=\frac{16}{3n}+O(n^{-2}).
\label{eq:intro-known-quartic-moments}
\end{equation}
Katz and Ramirez \cite{KatzRamirezMoments2025} developed a
combinatorial method for higher central moments of the demerit factor.
Their subsequent theorem
\cite[Theorem~1.6]{KatzRamirezLimiting2025} proves that every fixed
standardized moment converges to the corresponding moment of the
standard normal distribution.  Since the normal distribution is
determined by its moments, this theorem and
\eqref{eq:intro-known-quartic-moments} imply
\begin{equation}
\sqrt n\bigl(T_n^{(4)}-2\bigr)
\ \Longrightarrow\
N\!\left(0,\frac{16}{3}\right).
\label{eq:intro-known-quartic-clt}
\end{equation}
Thus the scalar Gaussian limit at the quartic exponent is an established result.  The joint
theorem below places this marginal within the fluctuating norm profile
and describes its dependence across nested lengths.

A related limit theory concerns nonlinear functions of real random
trigonometric polynomials.  The classical Salem--Zygmund theorem
\cite{SalemZygmund1954} concerns their value distributions under
uniform sampling of the angular variable; see \cite{AngstPoly2021}
for almost-sure and functional extensions.  Angst and Poly
\cite{AngstPoly2023} prove a \(\sqrt n\)-scale central limit theorem for
angular integrals of a fixed nonlinear test function, under their
coefficient and test-function hypotheses.  Their real trigonometric
model uses independent cosine and sine coefficient sequences.  In the
present model, the real and imaginary parts of \(q_n\) are generated by
the same signs, and the observable is a function of their squared
modulus.  The Hermite expansion, Newton--Girard treatment of repeated indices,
reduction to homogeneous sums, and tail approximation in
\cite{AngstPoly2023} are methodological antecedents of the present
proof. Their higher-order homogeneous sums are handled using the
universality theorem of Nourdin, Peccati, and Reinert
\cite{NourdinPeccatiReinert2010}.  In the Gaussian
setting, Sodin and Tsirelson \cite[Sections~2.2--2.3]{SodinTsirelson2004}
use radial Wick expansions and diagram identities for nonlinear
functionals of complex Gaussian processes.  Their normalized radial
Wick polynomials are, up to sign, the Laguerre modes used here.  These
Gaussian identities underlie our comparison calculations; the
non-Gaussian error bounds and the maximum over nested lengths are
verified separately.

There is also a close connection with nonlinear periodogram
statistics: \(S_n/(2\pi)\) is the usual periodogram of the observed
Rademacher sequence.  Fa\"y and Soulier \cite{FaySoulier2001} prove
Gaussian limits for weighted sums of nonlinear functions of
periodogram values at Fourier frequencies, under hypotheses on the
coefficient law and the test functions; their results also include,
under additional distributional assumptions, a functional central
limit theorem for the empirical distribution of these values.  Fa\"y,
Moulines, and Soulier \cite{FayMoulinesSoulier2002} obtain related
results for stationary linear processes, allowing nonlinearities that
may be singular or unbounded.

McElroy and Holan \cite{McElroyHolan2009} obtain joint Gaussian
limits for finite vectors of different integer powers of the
periodogram, summed over Fourier frequencies; these sums must be
distinguished from the continuous angular integrals considered here.
Kl\"uppelberg and Mikosch \cite{KlueppelbergMikosch1996} establish
Gaussian functional limits indexed by the observed prefix and a
frequency cutoff for integrated periodograms, a statistic linear
in the periodogram.

Earlier work on nonlinear spectral statistics includes Taniguchi's
\cite{Taniguchi1980} estimation of integrals of functions of the spectral
density of a stationary Gaussian process and von Sachs's
\cite[Theorem~3.1]{vonSachs1994} central limit theorem for nonlinear
integrals of Gaussian periodograms with an asymptotically vanishing
data taper. These results have their own assumptions on the spectral
density, nonlinearity, and taper.

For continuous squared-periodogram integrals, Deo and Chen
\cite{DeoChen2000} derive a Gaussian limit for Gaussian white noise and
extend their distributional analysis to stationary Gaussian series.
They also show that replacing a continuous integral by its sum over
Fourier frequencies can change the asymptotic variance. Thus their work
provides both a continuous-integral fluctuation theorem and a reason to
distinguish angular integration from Fourier-frequency sampling.\enlargethispage{2pt}

Continuous weighted integrals of integer powers
of tapered periodograms are studied by Sakhno \cite{Sakhno2014}, who
states joint asymptotic normality results for finite vectors involving
different integer powers, under hypotheses on the spectral densities.
This is a relevant continuous-integral antecedent.  For
constant weights, the rectangular taper, and independent Rademacher
data, the covariance displayed in \cite[Theorem~1(2)]{Sakhno2014},
used in its Theorem~5, gives zero at the squared-periodogram statistic;
this disagrees with the positive quartic variance in
\eqref{eq:intro-known-quartic-moments}.\footnote{Indeed,
\(f=(2\pi)^{-1}\), \(f_4=-2(2\pi)^{-3}\), and the covariance
bracket at \(k=l=2\) equals
\(2(2\pi)^{-3}-2(2\pi)^{-3}=0\).  For \(n=2T+1\),
\(T_n^{(4)}=2\pi J_{2,T}(1)\).  This is a check of the displayed
specialization, not a report of a published correction.}
We therefore do not use that covariance formula as an established
Rademacher predecessor.  The nonlinear periodogram literature already
contains functional results.  The present theorem combines continuous
angular integration of powers of a Rademacher polynomial,
uniformity over compact positive-exponent intervals, and a common
nested sequence of polynomial lengths.

We study the moment profile over a multiplicative window of lengths.  Let \(K=[p_-,p_+]\) be a compact interval in
\((0,\infty)\), let \(T>0\), and set \(I_T=[-T,T]\) and
\(n_N(t)=\lfloor Ne^t\rfloor\).  The field of interest is
\begin{equation}
(p,t)\longmapsto
\sqrt{n_N(t)}\,
\bigl(T_{n_N(t)}^{(p)}-m_p\bigr),
\qquad (p,t)\in K\times I_T.
\label{eq:intro-joint-fluctuation-question}
\end{equation}
The parameter \(p\) changes the moment of the intensity distribution,
whereas \(t\) changes the length of the initial coefficient block.
Logarithmic time keeps ratios of lengths on a fixed scale as
\(N\to\infty\).  Throughout the window, the polynomials are built from
the same signs.

A scalar central limit theorem, even when available for every fixed
exponent, does not determine the mixed limiting law of
\eqref{eq:intro-joint-fluctuation-question}.  Nor does convergence for
finitely many exponents and lengths provide the uniform control of
oscillations needed for a functional limit.  The theorem below gives
the jointly indexed Gaussian limit together with path-space
convergence.  Its covariance is described by a common family of
independent stationary Gaussian modes, and its quartic marginal at
each fixed time agrees with \eqref{eq:intro-known-quartic-clt}.  For
example, it gives a limit law for the largest scaled absolute deviation
of the normalized norms over the whole exponent and length range; see
\eqref{eq:intro-simultaneous-norm-maximum}.  The precise path-space
limits for the power integrals and the normalized norms follow.

\subsection{Main results}
\label{sec:intro-main-results}
\label{subsec:intro-main-results}

Let \(K=[p_-,p_+]\Subset(0,\infty)\), with \(p_-<p_+\), let
\(T>0\), and put \(I_T=[-T,T]\).  For \(N\ge1\), set
\[
n_N(t)=\lfloor Ne^t\rfloor,\qquad
n_-=\lfloor Ne^{-T}\rfloor,\qquad
n_+=\lfloor Ne^T\rfloor,
\]
taking \(N\) sufficiently large that \(n_-\ge1\).  With
\(m_p=\Gamma(1+p/2)\), define
\begin{align}
Y_n^{\mathrm{dc}}(p)
&=\sqrt n\bigl(T_n^{(p)}-m_p\bigr),
\label{eq:intro-centered-power-values}\\
\mathcal Z_N^{\mathrm{dc},\mathrm{fl}}(t)(p)
&=Y_{n_N(t)}^{\mathrm{dc}}(p),
\qquad (p,t)\in K\times I_T.
\label{eq:intro-power-floor-field}
\end{align}
Let \(\mathcal Z_N^{\mathrm{dc},\mathrm{poly}}\) denote the canonical
logarithmic-time polygonal interpolation of the \(C(K)\)-valued data
\((Y_n^{\mathrm{dc}})_{n_-\le n\le n_+}\), with the endpoint
convention of Definition~\ref{def:canonical-polygonal-interpolation}.
We view the polygonal field in \(C(K\times I_T)\), equipped with the
supremum norm, and the floor field in \(D(I_T;C(K))\), equipped with
the Skorokhod \(J_1\)-topology.

For every integer \(a\ge2\), define
\begin{align}
c_a(p)
&=\Gamma\!\left(1+\frac p2\right)(-1)^a\binom{p/2}{a},
\label{eq:intro-ca}\\
f_a
&=\underbrace{\mathbf1_{[0,1]}*\cdots*\mathbf1_{[0,1]}}_{a\text{ factors}},
\qquad
V_a=\int_0^a f_a(u)^2\,\mathrm du,
\label{eq:intro-Va}
\end{align}
where \(\binom za=z(z-1)\cdots(z-a+1)/a!\) and \(*\) denotes
convolution on \(\mathbb R\).  The coefficients \(c_a(p)\) are the
Laguerre coefficients of \(x^{p/2}\); see
Section~\ref{sec:spectral-framework}.  Let
\((G^{(a)})_{a\ge2}\) be independent centered stationary Gaussian
processes on \(\mathbb R\), with continuous paths and covariance
\begin{equation}
\operatorname{Cov}\!\left(G^{(a)}(s),G^{(a)}(t)\right)
=2V_a e^{-(a-\frac12)|t-s|}.
\label{eq:intro-mode-covariance}
\end{equation}
The associated Gaussian series is
\begin{equation}
G(p,t)=\sum_{a=2}^{\infty}c_a(p)G^{(a)}(t).
\label{eq:intro-G-mode}
\end{equation}
The following is the deterministically centered form of
Theorem~\ref{thm:main-joint-field}.

\begin{theorem}[Joint functional central limit theorem]
\label{thm:intro-joint-field}
For every \(K\) and \(T\) as above, the series
\eqref{eq:intro-G-mode} converges in
\[
L^2\!\left(\Omega;C^1(K;C(I_T))\right)
\]
to a centered Gaussian field \(G\).  In particular, the limit may be
chosen in \(C^1(K;C(I_T))\).  As \(N\to\infty\),
\begin{align}
\mathcal Z_N^{\mathrm{dc},\mathrm{poly}}
&\ \Longrightarrow\ G
&&\text{in }C(K\times I_T),
\label{eq:intro-main-polygonal}\\
\mathcal Z_N^{\mathrm{dc},\mathrm{fl}}
&\ \Longrightarrow\ \bigl(t\mapsto G(\cdot,t)\bigr)
&&\text{in }D(I_T;C(K)).
\label{eq:intro-main-floor}
\end{align}
For \(p,q\in K\) and \(s,t\in I_T\),
\begin{equation}
\operatorname{Cov}\!\left(G(p,s),G(q,t)\right)
=2\sum_{a=2}^{\infty}
V_a c_a(p)c_a(q)e^{-(a-\frac12)|t-s|}.
\label{eq:intro-main-covariance}
\end{equation}
The covariance series converges absolutely and uniformly on
\(K\times K\times I_T\times I_T\).
\end{theorem}

The \(C^1\)-statement concerns the limiting field.  The prelimit
convergence in \eqref{eq:intro-main-polygonal} and
\eqref{eq:intro-main-floor} is asserted only in the displayed value
topologies, not in a \(C^1\)-topology.

An expectation-centered version has the same limit.  Indeed,
Proposition~\ref{prop:uniform-deterministic-expectation} proves
\begin{equation}
\sup_{p\in K}\left|\mathbb E T_n^{(p)}-m_p\right|
\le \frac{C_K}{n}.
\label{eq:intro-centering-envelope}
\end{equation}
Since \(n\asymp_T N\) throughout the length window, replacing \(m_p\)
by \(\mathbb E T_n^{(p)}\) changes the \(\sqrt n\)-scaled fields by
\(O_{K,T}(N^{-1/2})\) uniformly over the window.

For the normalized norms, set
\[
r_p=m_p^{1/p},
\qquad
W_n(p)
=\sqrt n\left[\frac{\|q_n\|_p}{\sqrt n}-r_p\right],
\]
and define
\begin{equation}
\mathcal W_N^{\mathrm{fl}}(t)(p)=W_{n_N(t)}(p).
\label{eq:intro-norm-field}
\end{equation}
Let \(\mathcal W_N^{\mathrm{poly}}\) be the same canonical
logarithmic-time interpolation applied to the values \(W_n\).
Thus the transformation
\[
T_n^{(p)}\longmapsto
\bigl(T_n^{(p)}\bigr)^{1/p}
=\frac{\|q_n\|_p}{\sqrt n}
\]
is performed at each integer length before interpolation.
We give the corresponding
functional limit theorem below.

\begin{corollary}[Normalized-\(L^p\) functional central limit theorem]
\label{cor:intro-normalized-norm}
Let
\begin{equation}
H(p,t)=\mathfrak a_pG(p,t),
\qquad
\mathfrak a_p=\frac1p\,m_p^{1/p-1}.
\label{eq:Hpt-def}
\end{equation}
Then, as \(N\to\infty\),
\begin{align}
\mathcal W_N^{\mathrm{poly}}
&\ \Longrightarrow\ H
&&\text{in }C(K\times I_T),
\label{eq:intro-norm-polygonal}\\
\mathcal W_N^{\mathrm{fl}}
&\ \Longrightarrow\ \bigl(t\mapsto H(\cdot,t)\bigr)
&&\text{in }D(I_T;C(K)).
\label{eq:intro-norm-floor}
\end{align}
The field \(H\) is centered Gaussian and may be chosen in
\(C^1(K;C(I_T))\).  Its covariance is obtained by multiplying the
right-hand side of \eqref{eq:intro-main-covariance} by
\(\mathfrak a_p\mathfrak a_q\).
\end{corollary}
\subsection{The two parameters and the Gaussian limit}
\label{sec:intro-gaussian-structure}

The first parameter, that is, the logarithmic time coordinate,  reflects the nesting of the coefficient
blocks.  For lengths \(m\le n\), the covariance calculation in
Theorem~\ref{thm:nested-top-chaos-covariance} gives the factor
\((m/n)^{a-1/2}\) for the leading contribution of the \(a\)-th
Laguerre mode.  For \(s\le t\),
\[
\left(\frac{n_N(s)}{n_N(t)}\right)^{a-\frac12}
\longrightarrow e^{-(a-\frac12)(t-s)}.
\]
Thus multiplicative changes of length become additive time increments,
and the limiting covariance depends only on the time difference.  This
explains stationarity in logarithmic time.  
For each \(a\geq2\), the
covariance kernel
\[
2V_a e^{-(a-\frac12)|t-s|}
\]
is that of a stationary Ornstein--Uhlenbeck process with decay rate
\(a-1/2\). Thus \(G^{(a)}\) has the corresponding stationary Gaussian
law, and lower-index modes retain correlations over longer intervals of
logarithmic time.

The second parameter, that is, the exponent coordinate, has a different role.  In
\eqref{eq:intro-G-mode}, varying \(p\) changes the deterministic
weights \(c_a(p)\), not the underlying family of random processes.
Different moments of the same intensity distribution therefore have
limits built from common modes.  In the covariance formula, the
product \(c_a(p)c_a(q)\) describes how two exponents weight mode
\(a\), while \(e^{-(a-\frac12)|t-s|}\) describes its persistence
as the length changes.  The constant and linear Laguerre terms do not
fluctuate, because \(L_0(x)=1\), \(L_1(x)=1-x\), and
\(\int_{\mathbb T}S_n\,\mathrm d\mu=1\).  Moreover,
\eqref{eq:intro-ca} gives \(c_a(2k)=0\) for \(a>k\).
Every positive even exponent thus involves only finitely many modes.
For a positive exponent that is not an even integer, all the
coefficients \(c_a(p)\), \(a\ge2\), are nonzero, and infinitely
many modes contribute.

The quadratic and quartic cases give exact illustrations whenever
these exponents belong to \(K\).  At \(p=2\), Parseval's identity
makes the normalized statistic deterministic.  At \(p=4\), only
\(c_2(4)=2\) remains among the fluctuating modes.  Hence
\[
G(2,t)=H(2,t)=0,
\qquad
G(4,t)=2G^{(2)}(t).
\]
Since \(V_2=2/3\), the quartic covariance is
\begin{equation}
\operatorname{Cov}\!\left(G(4,s),G(4,t)\right)
=\frac{16}{3}e^{-\frac32|t-s|}.
\label{eq:intro-quartic-process-covariance}
\end{equation}
Its equal-time variance agrees with
\eqref{eq:intro-known-quartic-clt}, while the exponential factor
describes the dependence across nested lengths.  The quartic process
is therefore a single-mode specialization of the full field.

The joint theorem also determines how different power statistics,
and hence the corresponding normalized norms, fluctuate together.
Whenever \(4\in K\), only the second mode contributes to the
covariance between \(G(p,\cdot)\) and the quartic process, for every
\(p\in K\).  Using
\(c_2(p)=m_p p(p-2)/8\) in
\eqref{eq:intro-main-covariance} gives
\begin{equation}
\operatorname{Cov}\!\left(G(p,s),G(4,t)\right)
=\frac{m_p p(p-2)}{3}e^{-\frac32|t-s|}.
\label{eq:intro-quartic-mixed-covariance}
\end{equation}
This covariance is negative for \(0<p<2\) and positive for
\(p>2\).  The same sign distinction holds for the normalized-norm
limits, since the multipliers \(\mathfrak a_p\) are positive.
Also, \(V_2>0\) and the formula for \(c_2(p)\) show that
\(\operatorname{Var}(G(p,t))>0\) for every
\(p\in K\setminus\{2\}\).  These covariance relations give
information about the joint norm profile that separate scalar limit
laws do not determine.

Functional convergence also controls statistics that depend on the
whole exponent--length window.  For example,
Corollary~\ref{cor:intro-normalized-norm} yields
\begin{equation}
\max_{n_-\le n\le n_+}\ \sup_{p\in K}
\sqrt n\left|\frac{\|q_n\|_p}{\sqrt n}-r_p\right|
\ \Longrightarrow\
\sup_{(p,t)\in K\times I_T}|H(p,t)|.
\label{eq:intro-simultaneous-norm-maximum}
\end{equation}
The left-hand side is exactly the supremum norm of
\(\mathcal W_N^{\mathrm{poly}}\): the interpolation attains every
grid value, and each intermediate value is a convex combination of
two adjacent values.  The continuous mapping theorem therefore gives
the assertion.  It describes the largest absolute deviation from the
first-order norm profile over the entire exponent--length window,
rather than at a prescribed exponent or length.

\subsection{Main difficulties and ideas of the proof}
\label{sec:intro-proof-outline}

We prove the joint theorem by approximating the power function
\(h_p(x)=x^{p/2}\) by a finite sum of Laguerre polynomials.
We first identify the joint limit for this finite approximation, then
show that the discarded terms are negligible at the fluctuation scale,
simultaneously over exponents and polynomial lengths. To express the
approximation error, we first use exact expectation centering and set
\begin{equation}
X_n(f)=\sqrt n\left[
\int_{\mathbb T}f(S_n(\theta))\,\mathrm d\mu(\theta)
-\mathbb E\int_{\mathbb T}f(S_n(\theta))\,\mathrm d\mu(\theta)
\right].
\label{eq:intro-proof-Xn}
\end{equation}
Let \((L_a)_{a\ge0}\) be the Laguerre basis orthonormal for
\(e^{-x}\,\mathrm dx\), the measure in the limiting intensity law.
We use \(c_a(p)\) for the coefficients of \(h_p\), including
\(a=0,1\). For a fixed integer \(d\ge2\), put
\begin{equation}
R_{p,d}(x)=h_p(x)-\sum_{a=0}^d c_a(p)L_a(x).
\label{eq:intro-proof-remainder}
\end{equation}
Since \(L_0(x)=1\), \(L_1(x)=1-x\), and the angular mean of
\(S_n\) is one, \(X_n(L_0)=X_n(L_1)=0\). Linearity therefore gives
\[
X_n(h_p)=\sum_{a=2}^d c_a(p)X_n(L_a)+X_n(R_{p,d}).
\]
The finite sum and the remainder require different arguments.
Convergence of the Laguerre expansion in the weighted \(L^2\)-space
does not by itself control the centered random integrals at the
\(\sqrt n\) scale or their maximum over lengths.

\boldparagraph{\textit{Finite-mode limits.}}
We first determine the joint process limit when only finitely many
Laguerre terms are retained. The Laguerre polynomials are useful because
cancellation among lower-degree terms leaves a leading sum of products
of \(2a\) distinct signs in \(X_n(L_a)\). The error in this
homogeneous multilinear approximation is negligible over the length
window. Counting the contributing index tuples determines the variance
and covariance across nested lengths.

For a fixed finite family, we apply the multivariate functional de Jong
theorem of D\"obler, Kasprzak, and Peccati
\cite[Theorem~1.4]{DoblerKasprzakPeccati2022}. We verify its hypotheses
by computing the variances of initial coefficient blocks, bounding the
influence of individual signs, and controlling the fourth moments.
The normalization and these checks are given in
Section~\ref{subsec:finite-laguerre-vector-convergence} and
Appendix~\ref{app:homogeneous-chaos}.

Let \(\mathcal Z_N^{\mathrm{ec},\mathrm{poly}}\) denote the polygonal
interpolation of \(X_n(h_p)\), and let
\(\mathcal Z_{N,d}^{\mathrm{poly}}\) denote that of the finite sum
above. For every fixed \(d\ge2\), we obtain
\(\mathcal Z_{N,d}^{\mathrm{poly}}\Longrightarrow G_d\) in
\(C(K\times I_T)\), where
\(G_d(p,t)=\sum_{a=2}^d c_a(p)G^{(a)}(t)\); see
\eqref{eq:finite-mode-joint-polygonal-convergence}.
This identifies the finite Gaussian approximation. To recover the
original field, we must now control the terms omitted from it.

\boldparagraph{\textit{Uniform approximation over lengths.}}
We next ask whether the discarded terms remain small simultaneously
over the whole length window. There are two obstacles. The fixed-degree
constants grow too rapidly to sum against the algebraically decaying
coefficients at non-even exponents. Also, a bound at one length does
not by itself control a maximum over order-\(N\) lengths.

We keep \(d\) fixed and divide its remainder at an auxiliary,
slowly increasing cutoff \(\Lambda_N\). For the degrees
\(d<a\le\Lambda_N\), Theorem~\ref{thm:coefficient-uniform-low-block}
controls the relevant portions of dyadic groups \(A<a\le2A\) through
the weighted sum of squared
coefficients \(\sum_{A<a\le2A}a|u_a|^2\), uniformly in those
coefficients, within the range \(1\le A\le\Lambda_N\).
For the remaining contribution, the Laguerre semigroup
\(P_\tau L_a=e^{-a\tau}L_a\) damps large-degree coefficients.
We write \(h_p-P_{1/\Lambda_N}h_p\) as a sum of differences of
smoothed functions, called Abel blocks. A separate comparison of
\(P_{1/\Lambda_N}h_p\) with the sharp Laguerre truncation at \(\Lambda_N\) returns
to the required remainder.

Before splitting the intensity range, we remove the \(L_1\) component
of each smooth block and subtract the resulting value at zero.
These operations leave its centered statistic unchanged and give
functions vanishing at zero. Their order matters because multiplication
by an intensity cutoff need not preserve the original cancellations;
see Section~\ref{subsec:process-equivalent-normalization}.
We then estimate the pieces at small, intermediate, and large intensity
values separately. Near zero, the estimates must account for possible
singular derivatives of the original power functions. At large values,
the weight \(e^{-x}\) in the Laguerre energy integrals supplies decay.

For each localized test \(u\), control of the length maximum uses
the fact that nearby lengths share most of their signs. For \(m<n\),
conditioning on the first \(m\) signs splits
\(X_n(u)-X_m(u)\) into the fluctuation due to the added signs and
the change in conditional average. Fourier
projection controls the leading new-sign contribution; the nonlinear
terms are estimated separately. Gaussian comparison identifies the
conditional-average change, up to a controlled error, with
\(X_m(\sqrt{n/m}\,P_{\log(n/m)}u-u)\).
For the fluctuating degrees \(a\ge2\), the multiplier is
\((m/n)^{a-1/2}-1\), which is
small when \(a(n-m)/n\) is small; see
\eqref{eq:direct-conditional-mean-identity} and
\eqref{eq:direct-semigroup-Laguerre-expansion}.
We therefore use finer time grids for larger degrees and sum the
comparison errors with the number of pairs at each level included.

The resulting Theorem~\ref{thm:fixed-cutoff-maximal-tails} bounds the
expected squared length maxima for both the remainder and its exponent
derivative, integrated over \(K\). The largest error over all exponents
still remains to be controlled.

\boldparagraph{\textit{Uniformity in the exponent.}}
An integral over \(p\) can miss a large error concentrated in a short
exponent interval. To control the supremum, we also estimate the
exponent derivative: a one-dimensional Sobolev inequality then bounds
the maximum error using integrated bounds for the function and its
derivative. This requires the additional family
\(\partial_ph_p(x)=\tfrac12x^{p/2}\log x\), with value zero at
\(x=0\). Its estimates do not follow from the ordinary remainder
bounds. At \(p=2\), for example, \(R_{2,d}\) vanishes, whereas
its exponent derivative has a nonzero Laguerre coefficient in every
degree above \(d\). The logarithmic-power functions therefore require
their own estimates, including near zero; see
Section~\ref{sec:coefficient-calculus} and
Appendix~\ref{app:kummer-critical-profile}.

The interpolation weights are independent of \(p\), so differentiation
commutes with polygonal interpolation. We can thus apply the Sobolev
inequality pathwise to the remainder as a \(C(I_T)\)-valued function
of \(p\). From the two integrated length-maximal estimates,
Proposition~\ref{prop:sobolev-prelimit-approximation} yields
\begin{equation}
\limsup_{N\to\infty}\mathbb E
\left\|
\mathcal Z_N^{\mathrm{ec},\mathrm{poly}}
-\mathcal Z_{N,d}^{\mathrm{poly}}
\right\|_{C(K\times I_T)}^2
\le C_{K,T}d^{-p_-}\log^3(d+2).
\label{eq:intro-proof-uniform-approximation}
\end{equation}
Since \(p_->0\), the right-hand side tends to zero as \(d\to\infty\).
Thus, after taking the large-\(N\) upper limit at fixed \(d\), the
finite fields approximate the original field uniformly over the entire
exponent--length rectangle. It remains to join their finite Gaussian
limits to the full limiting field.

\boldparagraph{\textit{Passage to the joint limit.}}
We now have convergence of each finite approximation and a uniform bound
for its difference from the original field. The remaining ingredient is
convergence of the finite Gaussian sums. Coefficient bounds and Gaussian
maximal estimates give \(G_d\to G\) in
\(L^2(\Omega;C^1(K;C(I_T)))\), as proved in
Theorem~\ref{thm:full-gaussian-field}.
Together with \eqref{eq:intro-proof-uniform-approximation}, these facts
give the joint limit in \(C(K\times I_T)\): first take
\(N\to\infty\) with \(d\) fixed, then \(d\to\infty\).
The internal cutoff \(\Lambda_N\) is used to prove the remainder
estimates and does not change this order.

Finally, vanishing adjacent jumps make the floor paths uniformly close
to the continuous interpolation. The bias estimate
\eqref{eq:intro-centering-envelope} makes replacement of exact
expectations by \(m_p\) negligible at the fluctuation scale.
For normalized norms, the power statistics stay uniformly close to
\(m_p\) with high probability, and \(m_p\) is bounded away from
zero on \(K\). A common Taylor neighborhood therefore gives an
expansion of \(x^{1/p}\) uniform in \(p\). We apply the
transformation at each integer length before interpolation, as in
\eqref{eq:poly-floor-W-Cnorm-equi}, obtaining the normalized-norm
limits for both path constructions.

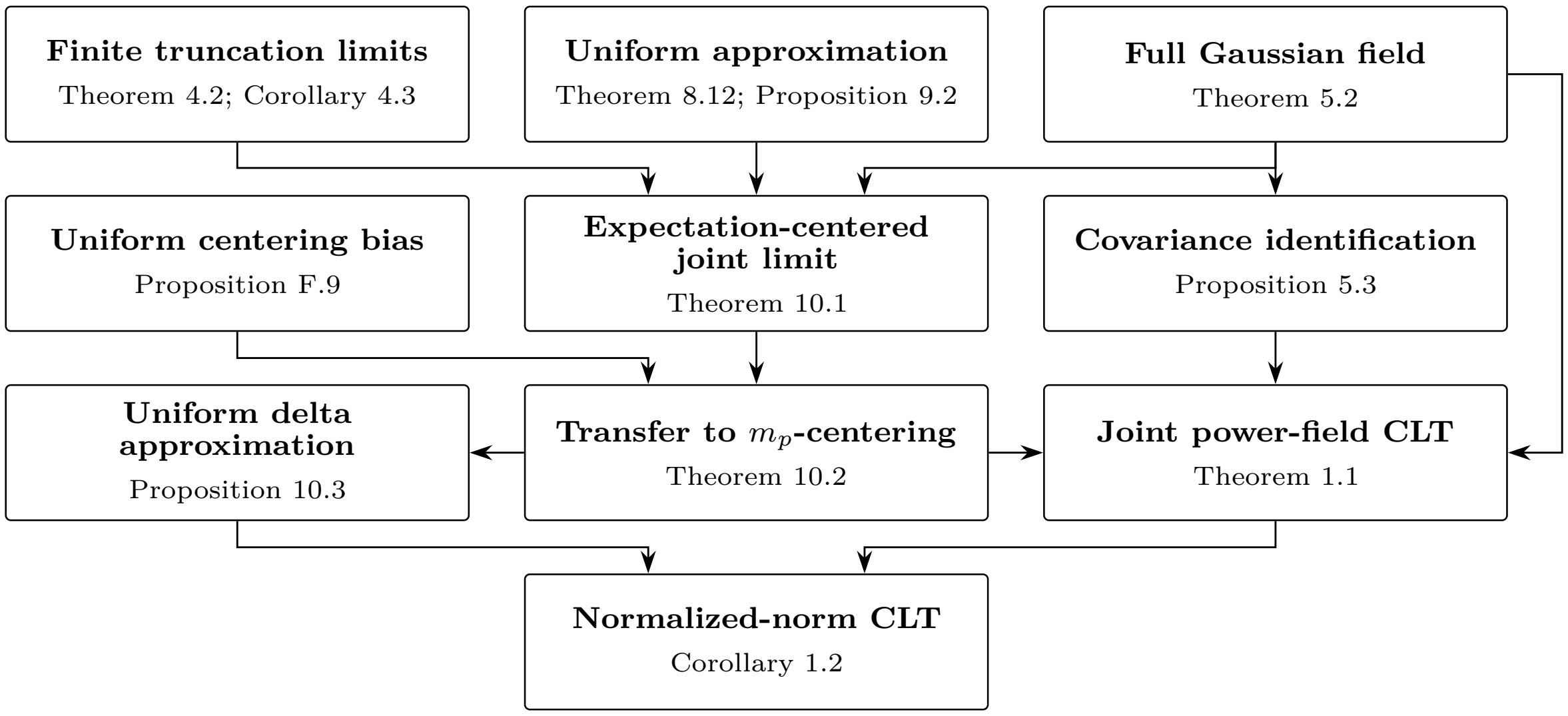
\begin{figure}[!tp]
  \centering
  \resizebox{\textwidth}{!}{
\begin{tikzpicture}[
  x=1cm,y=1cm,
  cltbox/.style={draw=black,line width=.45pt,rounded corners=1.3pt,
    align=flush center,text=black,text width=4cm,minimum height=1.25cm,inner sep=4pt,
    font=\fontsize{7}{8.4}\selectfont,
    execute at begin node={\hyphenpenalty=10000\exhyphenpenalty=10000}},
  cltinput/.style={cltbox},
  cltlimit/.style={cltbox},
  clttransfer/.style={cltbox},
  cltarrow/.style={-{Stealth[length=2.1mm,width=1.4mm]},
    draw=black,line width=.5pt}
]

\node[cltinput,minimum height=1.25cm] (finite) at (0,0) {
  {\fontsize{8}{9.6}\selectfont\bfseries Finite truncation limits}\\[3pt]
  {\fontsize{7}{8.4}\selectfont Theorem~\ref{thm:finite-laguerre-vector};
  Corollary~\ref{cor:finite-power-truncation-joint}}
};

\node[cltinput,minimum height=1.25cm] (tails) at (4.8,0) {
  {\fontsize{8}{9.6}\selectfont\bfseries Uniform approximation}\\[3pt]
  {\fontsize{7}{8.4}\selectfont Theorem~\ref{thm:fixed-cutoff-maximal-tails};
  Proposition~\ref{prop:sobolev-prelimit-approximation}}
};

\node[cltinput,minimum height=1.25cm] (gaussian) at (9.6,0) {
  {\fontsize{8}{9.6}\selectfont\bfseries Full Gaussian field}\\[3pt]
  {\fontsize{7}{8.4}\selectfont Theorem~\ref{thm:full-gaussian-field}}
};

\node[clttransfer,minimum height=1.25cm] (ec) at (4.8,-1.75) {
  {\fontsize{8}{9.6}\selectfont\bfseries Expectation-centered\\joint limit}\\[3pt]
  {\fontsize{7}{8.4}\selectfont Theorem~\ref{thm:expectation-centered-joint-field}}
};

\node[cltinput] (bias) at (0,-1.75) {
  {\fontsize{8}{9.6}\selectfont\bfseries Uniform centering bias}\\[3pt]
  {\fontsize{7}{8.4}\selectfont Proposition~\ref{prop:uniform-deterministic-expectation}}
};

\node[clttransfer,minimum height=1.25cm] (dc) at (4.8,-3.5) {
  {\fontsize{8}{9.6}\selectfont\bfseries Transfer to $m_p$-centering}\\[3pt]
  {\fontsize{7}{8.4}\selectfont Theorem~\ref{thm:deterministic-centering}}
};

\node[cltinput] (covariance) at (9.6,-1.75) {
  {\fontsize{8}{9.6}\selectfont\bfseries Covariance identification}\\[3pt]
  {\fontsize{7}{8.4}\selectfont Proposition~\ref{prop:joint-covariance}}
};

\node[cltlimit,text width=4cm,minimum height=1.25cm] (main) at (9.6,-3.5) {
  {\fontsize{8}{9.6}\selectfont\bfseries Joint power-field CLT}\\[3pt]
  {\fontsize{7}{8.4}\selectfont Theorem~\ref{thm:intro-joint-field}}
};

\node[clttransfer] (delta) at (0,-3.5) {
  {\fontsize{8}{9.6}\selectfont\bfseries Uniform delta approximation}\\[3pt]
  {\fontsize{7}{8.4}\selectfont Proposition~\ref{prop:assembly-uniform-delta-approximation}}
};

\node[cltlimit,minimum height=1.25cm] (norm) at (4.8,-5.25) {
  {\fontsize{8}{9.6}\selectfont\bfseries Normalized-norm CLT}\\[3pt]
  {\fontsize{7}{8.4}\selectfont Corollary~\ref{cor:intro-normalized-norm}}
};

\draw[cltarrow] (finite.south)
  |- ([xshift=-10mm,yshift=2.5mm]ec.north)
  -- ([xshift=-10mm]ec.north);
\draw[cltarrow] (tails.south) -- (ec.north);
\draw[cltarrow] (gaussian.south)
  |- ([xshift=10mm,yshift=2.5mm]ec.north)
  -- ([xshift=10mm]ec.north);

\draw[cltarrow] (ec.south) -- (dc.north);
\draw[cltarrow] (bias.south) -- (0,-2.625) -- (3.8,-2.625)
  -- ([xshift=-1cm]dc.north);
\draw[cltarrow] (gaussian.south) -- (covariance.north);

\draw[cltarrow] (dc.east) -- (main.west);
\draw[cltarrow] (covariance.south) -- (main.north);
\draw[cltarrow] (gaussian.east) -- (12.25,0)
  -- (12.25,-3.5) -- (main.east);
\draw[cltarrow] (dc.west) -- (delta.east);

\draw[cltarrow] (main.south) -- (9.6,-4.375) -- (5.8,-4.375)
  -- ([xshift=1cm]norm.north);
\draw[cltarrow] (delta.south) -- (0,-4.375) -- (3.8,-4.375)
  -- ([xshift=-1cm]norm.north);
\end{tikzpicture}
}
  \caption{Proof structure for Theorem~\ref{thm:intro-joint-field}
  and Corollary~\ref{cor:intro-normalized-norm}.
}
  \label{fig:main-proof-roadmap}
\end{figure}

\subsection{Scope and organization}
\label{sec:intro-scope-organization}

The exponent interval \(K\Subset(0,\infty)\) and the logarithmic-time
window \(I_T=[-T,T]\) remain fixed as \(N\to\infty\).  In particular,
we do not treat the endpoint regimes \(p\downarrow0\) or \(p\uparrow
\infty\), exponent intervals depending on \(N\), or time windows whose
length grows with \(N\).  The range \(0<p<1\) is included; no use is
made of the \(L^p\) triangle inequality.  The limiting fields \(G\)
and \(H\) may be chosen \(C^1\) in the exponent, whereas the prelimit
convergence is asserted in \(C(K\times I_T)\) and
\(D(I_T;C(K))\).  Convergence in a \(C^1\)-topology is not claimed.

The finite-mode theorem also gives, by linearity, expectation-centered
logarithmic-time limits for arbitrary fixed polynomial tests, whereas the
uniform estimates used to remove the spectral truncation are proved for
the power family and its exponent derivative.  All model-specific
arguments needed for the joint \(L^p\)-profile theorem are proved in
the present paper.  In particular, the expectation-comparison
arguments shared with the companion cited in
Section~\ref{sec:intro-first-order-to-fluctuations} are reproduced with
explicit acknowledgement and the estimates uniform over compact
exponent intervals needed here.  This gives a self-contained proof with respect to that
manuscript while recording the common ancestry of these arguments.

Section~\ref{sec:model-main-results} fixes the model and path spaces and
gives the full theorem statements.  Sections~\ref{sec:spectral-framework}--%
\ref{sec:full-gaussian-field} develop the spectral representation,
the finite-mode process limit, and the limiting Gaussian field.
Sections~\ref{sec:coefficient-calculus}--%
\ref{sec:uniform-prelimit-approximation} establish the coefficient and maximal
estimates needed to remove the spectral truncation.  Finally,
Section~\ref{sec:joint-assembly-transfers} proves the joint limit
and handles the two path constructions and changes of centering, followed by
the passage to normalized norms.

To keep the main line of the argument clear, we place several essential estimates of the proof in the appendices, which therefore form an integral part of the paper rather than merely supplementary material.
Appendix~\ref{app:homogeneous-chaos}
develops the homogeneous-chaos counts and finite-mode estimates.
Appendix~\ref{app:finite-replica-occupation} provides the comparison
and estimates at a single length, while
Appendix~\ref{app:block-addition-refinement} proves increment bounds
between two lengths and the temporal refinement used for maxima over
nested lengths.  Appendix~\ref{app:kummer-critical-profile} supplies the
special-function analysis for the logarithmic-power family,
Appendix~\ref{app:moving-labels} treats the moving truncation
boundaries, and Appendix~\ref{app:expectation-envelope} proves the
expectation estimate uniform over compact exponent intervals used for
deterministic centering.

\section{Model, path spaces, and main results}
\label{sec:model-main-results}

This section fixes the model, the interpolation convention, and the path spaces. We then define the prelimit fields and state the full convergence theorem for the power statistics.

\subsection{Random Littlewood polynomials and power statistics}
\label{sec:model-power-statistics}

Let \((\varepsilon_j)_{j\geq0}\) be independent Rademacher random
variables:
\[
\Prob(\varepsilon_j=1)
=
\Prob(\varepsilon_j=-1)
=
\frac12.
\]
For every integer \(n\geq1\), define the random Littlewood polynomial
\begin{equation}
q_n(z)
=
\sum_{j=0}^{n-1}\varepsilon_jz^j.
\label{eq:model-qn}
\end{equation}
All polynomial lengths are coupled through the same infinite sign sequence.
In particular, \(q_m\) is the initial coefficient block of \(q_n\)
whenever \(m<n\).

We write \(\Torus=\mathbb R/ (2\pi\mathbb Z)\) and let \(\mu\) denote the normalized Haar measure on \(\Torus\).
Define the normalized random intensity by
\begin{equation}
S_n(\theta)
=
\frac{|q_n(e^{i\theta})|^2}{n}.
\label{eq:model-Sn}
\end{equation}
For any measurable function
\(f:[0,\infty)\to\R\) for which the integral below is finite, set
\begin{equation}
I_n(f)
=
\int_{\Torus}
f(S_n(\theta))\,\dd\mu(\theta),
\label{eq:model-In}
\end{equation}
and define the corresponding centered statistic on the central-limit
scale by
\begin{equation}
X_n(f)
=
\sqrt n
\left[
I_n(f)-\E I_n(f)
\right].
\label{eq:model-Xn}
\end{equation}

For \(p>0\), let
\begin{equation}
h_p(x)=x^{p/2},
\qquad
x\geq0.
\label{eq:model-hp}
\end{equation}
The normalized \(p\)-th power statistic is
\begin{equation}
T_n^{(p)}
= I_n(h_p)
=\frac{\|q_n\|_p^p}{n^{p/2}}.
\label{eq:model-Tnp}
\end{equation}
Here
\[
\|q_n\|_p
=
\left(
\int_{\Torus}
|q_n(e^{i\theta})|^p\,\dd\mu(\theta)
\right)^{1/p}.
\]
When \(0<p<1\), this notation denotes the usual \(L^p\)-quasinorm.
We fix a compact exponent interval,
\begin{equation}
K=[p_-,p_+]\Subset(0,\infty),
\qquad
p_-<p_+.
\label{eq:model-K}
\end{equation}

\subsection{Canonical logarithmic-time interpolation and path spaces}
\label{sec:model-log-time-interpolation}
Now we turn statistics at integer polynomial lengths into continuous and step paths in logarithmic time.
Fix \(T>0\) and write
\[
I_T=[-T,T].
\]
For \(N\geq1\) and \(t\in I_T\), set
\begin{equation}
n_N(t)
=
\lfloor Ne^t\rfloor.
\label{eq:model-floor-length}
\end{equation}
For all sufficiently large \(N\), define
\begin{equation}
n_N(-T)
=
\lfloor Ne^{-T}\rfloor,
\qquad
n_N(T)
=
\lfloor Ne^T\rfloor.
\label{eq:model-window-endpoints}
\end{equation}
We abbreviate \(n_N(-T)\) as \(n_-\) and \(n_N(T)\) as \(n_+\).
The interpolation knots are defined by
\begin{equation}
\vartheta_{N,n_-}
= -T,\qquad
\vartheta_{N,n}
= \log\left(\frac nN\right),
\qquad
n_-<n\leq n_+.
\label{eq:model-effective-knots}
\end{equation}
The choice
\[
\vartheta_{N,n_-}=-T
\]
rather than \(\log(n_-/N)\) is part of the interpolation convention used
throughout the paper.

\begin{definition}[Canonical \(E\)-valued polygonal interpolation]
\label{def:canonical-polygonal-interpolation}
Let \(E\) be a normed vector space and let
\[
(F_n)_{n_-\leq n\leq n_+}
\subset E.
\]
The canonical polygonal interpolation of \((F_n)\) is the continuous map
\[
t
\longmapsto
\operatorname{Poly}_{N,T}[F_\bullet](t),
\qquad
t\in I_T,
\]
which satisfies
\[
\operatorname{Poly}_{N,T}[F_\bullet](\vartheta_{N,n})
=
F_n,
\qquad
n_-\leq n\leq n_+,
\]
is affine in \(t\) on each interval
\[
[\vartheta_{N,n},\vartheta_{N,n+1}],
\qquad
n_-\leq n<n_+,
\]
and is constant with value \(F_{n_+}\) on
\[
[\vartheta_{N,n_+},T].
\]
\end{definition}

For every deterministic test function \(f\), define its polygonal process by
\begin{equation}
X_{N,f}^{\poly}
=
\operatorname{Poly}_{N,T}
\left[
\bigl(X_n(f)\bigr)_{n_-\leq n\leq n_+}
\right],
\label{eq:paired-generic-polygonal-process}
\end{equation}
and its floor process by
\begin{equation}
X_{N,f}^{\fl}(t)
=
X_{n_N(t)}(f),
\qquad
t\in I_T.
\label{eq:paired-generic-floor-process}
\end{equation}
Since the polygonal process is affine between consecutive grid values and the floor process takes values on the same grid,
\begin{equation}
\left\|X_{N,f}^{\poly}\right\|_{C(I_T)}
=
\left\|X_{N,f}^{\fl}\right\|_{\infty,I_T}
=
\max_{n_-\leq n\leq n_+}
|X_n(f)|.
\label{eq:paired-floor-polygonal-supremum}
\end{equation}

Define the maximum width of the interpolation intervals by
\begin{equation}
\Delta_{N,T}
= \max_{n_-\leq n<n_+}
\left(\vartheta_{N,n+1}
-
\vartheta_{N,n}
\right)\vee
\left(T-\vartheta_{N,n_+}
\right).
\label{eq:model-effective-mesh}
\end{equation}
For every fixed \(T>0\),
\begin{equation}
\Delta_{N,T}
\leq
\frac{C_T}{N}
\label{eq:model-mesh-bound}
\end{equation}
for all sufficiently large \(N\), with \(C_T<\infty\) depending only
on \(T\). 
Indeed, for \(n_-<n<n_+\),
\[
\vartheta_{N,n+1}-\vartheta_{N,n}
=
\log\left(1+\frac1n\right)
\leq
\frac1n,
\]
while the initial and terminal gaps are \(O_T(N^{-1})\) by the
definitions of \(n_-\) and \(n_+\).

When the interpolated grid values take values in \(C(K)\), we view the
polygonal paths as elements of
\[
C(I_T;C(K)),
\]
equipped with the uniform norm
\[
\|F\|_{\infty,I_T;C(K)}
:=
\sup_{t\in I_T}
\|F(t)\|_{C(K)}.
\]
This space is canonically and isometrically identified with
\(C(K\times I_T)\) through
\[
F
\longmapsto
\bigl((p,t)\longmapsto F(t)(p)\bigr).
\]
The corresponding floor paths are viewed as elements of
\[
D(I_T;C(K)),
\]
the space of right-continuous \(C(K)\)-valued paths with left limits.
We equip this space with the Skorokhod \(J_1\)-topology, which permits small continuous, strictly increasing changes of time that fix the endpoints when comparing paths.
Since \(C(K)\) is a separable Banach space,
\(D(I_T;C(K))\) is a Polish space.
Moreover, for \(F_1,F_2\in D(I_T;C(K))\), the \(J_1\)-metric satisfies
\begin{equation}
d_{J_1}^{C(K)}(F_1,F_2)
\leq
\sup_{t\in I_T}
\|F_1(t)-F_2(t)\|_{C(K)}.
\label{eq:model-J1-uniform-domination}
\end{equation}

\subsection{Prelimit joint fields}
\label{sec:model-prelimit-fields}
Now the two centerings and two path constructions give the four fields that appear in our convergence theorem.
For \(p\in K\), set
\[
m_p
=
\Gamma\left(1+\frac p2\right).
\]
For every integer \(n\geq1\), define the expectation-centered and
deterministically centered \(C(K)\)-valued fields at that length by
\begin{equation}
Y_n^{\ec}(p)
=
\sqrt n
\left[
T_n^{(p)}-\E T_n^{(p)}
\right],
\qquad
Y_n^{\dc}(p)
=
\sqrt n
\left[
T_n^{(p)}-m_p
\right].
\label{eq:model-power-grid-fields}
\end{equation}
Their polygonal interpolations under the preceding convention are
\begin{equation}
\begin{aligned}
\mathcal Z_N^{\ec,\poly}
=
\operatorname{Poly}_{N,T}[Y_\bullet^{\ec}],
\qquad
\mathcal Z_N^{\dc,\poly}
=
\operatorname{Poly}_{N,T}[Y_\bullet^{\dc}].
\end{aligned}
\label{eq:model-power-polygonal-fields}
\end{equation}
By Definition~\ref{def:canonical-polygonal-interpolation},
\[
\mathcal Z_N^{\ec,\poly},
\,
\mathcal Z_N^{\dc,\poly}
\in
C(I_T;C(K)).
\]
The corresponding floor processes are
\begin{equation}
\begin{aligned}
\mathcal Z_N^{\ec,\fl}(t)
=
Y_{n_N(t)}^{\ec},
\qquad
\mathcal Z_N^{\dc,\fl}(t)
=
Y_{n_N(t)}^{\dc},
\qquad
t\in I_T.
\end{aligned}
\label{eq:model-power-floor-fields}
\end{equation}
Since \(t\mapsto n_N(t)\) is a finite-step càdlàg function on \(I_T\),
\[
\mathcal Z_N^{\ec,\fl},
\,
\mathcal Z_N^{\dc,\fl}
\in
D(I_T;C(K)).
\]

For fixed \(N\), all the fields above depend only on the finitely many
Rademacher variables
\[
\varepsilon_0,\ldots,\varepsilon_{n_+-1},
\]
and hence are measurable random elements of the indicated path spaces.

\subsection{The joint power-field theorem}
\label{sec:joint-power-field-theorem}

The four fields obtained from the two centerings and the two path
constructions have the same limiting field.  The full convergence theorem
is as follows.

\begin{theorem}[Joint power-field functional central limit theorem]
\label{thm:main-joint-field}
Let 
\(
K=[p_-,p_+]\Subset(0,\infty)
\)
and \(T>0\), and let \(G\) be the centered Gaussian field defined in
\eqref{eq:intro-G-mode}.

\begin{enumerate}[label=\textup{(\roman*)},leftmargin=2.4em]

\item
The expectation-centered polygonal fields satisfy
\begin{equation}
\mathcal Z_N^{\ec,\poly}
\Longrightarrow
\bigl(t\longmapsto G(\cdot,t)\bigr)
\qquad
\text{in }C(I_T;C(K)),
\label{eq:main-ec-poly-convergence}
\end{equation}
and the expectation-centered floor fields satisfy
\begin{equation}
\mathcal Z_N^{\ec,\fl}
\Longrightarrow
\bigl(t\longmapsto G(\cdot,t)\bigr)
\qquad
\text{in }D(I_T;C(K)).
\label{eq:main-ec-floor-convergence}
\end{equation}

\item 
The deterministically centered fields satisfy
\begin{equation}
\mathcal Z_N^{\dc,\poly}
\Longrightarrow
\bigl(t\longmapsto G(\cdot,t)\bigr)
\qquad
\text{in }C(I_T;C(K)),
\label{eq:main-dc-poly-convergence}
\end{equation}
and the deterministically centered floor fields satisfy
\begin{equation}
\mathcal Z_N^{\dc,\fl}
\Longrightarrow
\bigl(t\longmapsto G(\cdot,t)\bigr)
\qquad
\text{in }D(I_T;C(K)).
\label{eq:main-dc-floor-convergence}
\end{equation}
The convergences in \(D(I_T;C(K))\) are with respect to the
Skorokhod \(J_1\)-topology.

\item
For \(p,q\in K\) and \(s,t\in I_T\),
\begin{equation}
\Cov
\left(
G(p,s),
G(q,t)
\right)
=
2
\sum_{a=2}^{\infty}
V_a
c_a(p)c_a(q)
e^{-(a-\frac12)|t-s|}.
\label{eq:main-joint-covariance}
\end{equation}
The covariance series converges absolutely and uniformly on
\[
K\times K\times I_T\times I_T.
\]

\end{enumerate}
\end{theorem}

Two even-exponent specializations are worth recording.
The quadratic case is degenerate, whereas the quartic case is the first
nontrivial specialization and recovers the known quartic fluctuation law
recalled in Section~\ref{sec:intro-first-order-to-fluctuations}.

\begin{remark}
At \(p=2\), Parseval's identity gives
\[
T_n^{(2)}=1,
\]
so both centered fields vanish identically and
\[
G(2,t)=0.
\]

At \(p=4\), one has
\[
m_4=2,
\qquad
\E T_n^{(4)}=2-\frac1n,
\]
and only the second fluctuating Gaussian mode contributes:
\[
G(4,t)=2G^{(2)}(t).
\]
Since \(V_2=2/3\), the covariance formula \eqref{eq:main-joint-covariance} in
Theorem~\ref{thm:main-joint-field} reduces to
\[
\Cov\left(G(4,s),G(4,t)\right)
=
\frac{16}{3}e^{-\frac32|t-s|}.
\]
In particular,
\[
\Var(G(4,t))=\frac{16}{3},
\]
in agreement with the known quartic central limit theorem
\eqref{eq:intro-known-quartic-clt}.
\end{remark}

\section{Laguerre expansions and the semigroup}
\label{sec:spectral-framework}

This section expands the power functions in Laguerre modes and introduces the semigroup used later to control the terms left out of a finite approximation, allowing us to pass from finitely many modes to the full fluctuation field.
The Laguerre expansion separates the power statistics into components of different degrees. The coefficient identities and the Laguerre semigroup are classical; we record their normalization and the estimates that hold uniformly on compact exponent intervals.
We first compute and estimate the Laguerre coefficients of the power family \(h_\alpha(x)=x^\alpha\) and its exponent derivative. We then define truncations in degree and the energies and semigroup used to estimate and smooth their remainders.

\subsection{Laguerre basis and power-function coefficients}
\label{subsec:spectral-power-coefficients}

Now we expand the power functions in the Laguerre basis and estimate how their coefficients vary with the exponent.
Let
\[
\dd\nu(x)=e^{-x}\,\dd x,
\qquad x\geq0.
\]
For every integer \(a\geq0\), define the ordinary Laguerre polynomial by
\begin{equation}
L_a(x)
=
\sum_{j=0}^{a}
(-1)^j
\binom aj
\frac{x^j}{j!},
\qquad a\geq0,
\label{eq:spectral-Laguerre-definition}
\end{equation}
so that
\[
L_0(x)=1,
\qquad
L_1(x)=1-x.
\]
The family \((L_a)_{a\geq0}\) is an orthonormal basis of
\(L^2(\nu)\); see, for example,
\cite{Szego1975,NISTDLMF}.
For \(f\in L^2(\nu)\), we write
\[
\widehat f_a
=
\int_0^\infty
f(x)L_a(x)e^{-x}\,\dd x.
\]
Then,
\[
f
=
\sum_{a=0}^{\infty}
\widehat f_aL_a
\qquad
\text{in } L^2(\nu).
\]

The relation with Hermite expansion is explicit in the Gaussian model.
Let \(X,Y\) be independent standard real Gaussians, put
\(Z=(X+iY)/\sqrt2\), and let \(\operatorname{He}_m\) denote the
probabilists' Hermite polynomials, normalized by
\(e^{tx-t^2/2}=\sum_{m\geq0}\operatorname{He}_m(x)t^m/m!\).
The classical identity \cite[Eq.~(18.18.40)]{NISTDLMF}, after conversion
of the Hermite normalization, gives
\[
L_a(|Z|^2)
=\frac{(-1)^a}{2^a a!}
\sum_{j=0}^a\binom aj
\operatorname{He}_{2j}(X)\operatorname{He}_{2a-2j}(Y).
\]
Thus a radial Laguerre mode belongs to Gaussian chaos of total degree
\(2a\).  The corresponding extraction from the Rademacher field is
proved in Section~\ref{subs:app-A-Walsh-top-chaos}.

Parseval's identity yields
\[
\int_{\Torus}
S_n(\theta)\,\dd\mu(\theta)
=
1.
\]
Since \(L_0(x)=1\) and \(L_1(x)=1-x\), it follows that
\begin{equation}
I_n(L_0)=1,
\qquad
I_n(L_1)=0.
\label{eq:spectral-invisible-mode-values}
\end{equation}
Consequently,
\begin{equation}
X_n(L_0)=X_n(L_1)=0.
\label{eq:spectral-invisible-mode-fluctuations}
\end{equation}
Thus, the degree-zero and degree-one Laguerre terms do not contribute to the centered fluctuation, and all nontrivial fluctuating Laguerre expansions below start at degree \(a=2\).

In the coefficient calculations, it is convenient to reparametrize the exponent by
\[
\alpha=\frac p2,
\]
so that 
\(p\in K=[p_-,p_+]\) 
corresponds to
\(\alpha\in \mathfrak A:=[\alpha_-,\alpha_+]\), where
\begin{equation}
\alpha_-=p_-/2,
\qquad
\alpha_+=p_+/2.
\label{eq:def-alpha-PN}
\end{equation}
In the spectral calculations and the estimates that follow, we use
\(\alpha\) as the exponent parameter and write
\begin{equation}
h_\alpha(x)=x^\alpha,
\qquad x\geq0.
\label{eq:def-halpha}
\end{equation}
For every integer \(a\geq0\), define
\begin{equation}
\gamma_a(\alpha)
:=
\int_0^\infty
h_\alpha(x)L_a(x)e^{-x}\,\dd x,
\label{eq:spectral-gamma-alpha}
\end{equation}
the Laguerre coefficients of \(x^\alpha\).
For \(a\geq2\), these are precisely the coefficients denoted by
\(c_a(p)\) in \eqref{eq:intro-ca}, expressed in the
\(\alpha\)-parametrization:
\begin{equation}
\gamma_a(\alpha)=c_a(2\alpha),
\qquad
c_a(p)=\gamma_a\!\left(\frac p2\right).
\label{eq:equi-p-alpha}
\end{equation}

For \(z\in\C\) and \(a\geq0\), we use the generalized binomial coefficient
\[
\binom za
:=
\frac{z(z-1)\cdots(z-a+1)}{a!},
\qquad
\binom z0=1.
\]

\begin{proposition}
\label{prop:laguerre-coefficients}
For every \(\alpha>0\) and \(a\geq0\),
\begin{equation}
\gamma_a(\alpha)
=
\Gamma(1+\alpha)(-1)^a
\binom{\alpha}{a}.
\label{eq:spectral-gamma-formula}
\end{equation}
Consequently,
\begin{equation}
x^\alpha
=
\sum_{a=0}^{\infty}
\gamma_a(\alpha)L_a
\qquad
\text{in }L^2(\nu).
\label{eq:spectral-power-expansion-alpha}
\end{equation}
\end{proposition}

\begin{proof}
We use the classical generating function for the Laguerre polynomials;
see, for example, \cite[Eq.~(18.12.13)]{NISTDLMF},
\[
\sum_{a=0}^{\infty}
L_a(x)z^a
=
\frac1{1-z}
\exp\left(
-\frac{xz}{1-z}
\right),
\qquad |z|<1.
\]
For \(|z|<1\), integration against \(x^\alpha e^{-x}\,\dd x\) yields
\begin{equation}
\sum_{a=0}^{\infty}
\gamma_a(\alpha)z^a
=
\Gamma(1+\alpha)(1-z)^\alpha.
\label{eq:spectral-power-coefficient-generating-function}
\end{equation}
Comparing coefficients yields
\eqref{eq:spectral-gamma-formula}.
Since \(x^\alpha\in L^2(\nu)\), completeness of the Laguerre basis yields
\eqref{eq:spectral-power-expansion-alpha}.
\end{proof}

In particular, if \(k\in\mathbb Z_{\geq0}\), then
\[
\gamma_a(k)=0,
\qquad
a>k,
\]
so the Laguerre expansion of \(x^k\) terminates at degree \(k\).

We next differentiate the power family with respect to the exponent.
The map
\[
\alpha\longmapsto h_\alpha
\]
is \(C^\infty\) from \((0,\infty)\) into \(L^2(\nu)\).
In particular, set
\begin{equation}
\widetilde h_\alpha(x)
:=
\partial_\alpha h_\alpha(x)
=
x^\alpha\log x,
\qquad
x>0,
\label{eq:def-widetilde-halpha}
\end{equation}
with \(\widetilde h_\alpha(0):=0\).
Moreover, for every compact interval
\(\mathfrak A\Subset(0,\infty)\) and every integer \(r\geq0\),
\begin{equation}
\sup_{\alpha\in\mathfrak A}
\left\|
\partial_\alpha^r  h_\alpha
\right\|_{L^2(\nu)}
<\infty.
\label{eq:spectral-uniform-L2-derivatives}
\end{equation}
These assertions follow directly from dominated convergence.

Define
\begin{equation}
\widetilde\gamma_a(\alpha)
:=
\int_0^\infty
\widetilde h_\alpha(x)L_a(x)e^{-x}\,\dd x,
\qquad a\geq0.
\label{eq:spectral-tilde-gamma-definition}
\end{equation}
Differentiation under the integral sign gives
\begin{equation}
\widetilde\gamma_a(\alpha)
=
\partial_\alpha\gamma_a(\alpha).
\label{eq:spectral-tilde-gamma-derivative}
\end{equation}

\begin{proposition}
\label{prop:compact-coefficient-bounds}
Let
\[
\mathfrak A=[\alpha_-,\alpha_+]\Subset(0,\infty).
\]
There exists \(C_{\mathfrak A}<\infty\) such that, for every
\(\alpha\in\mathfrak A\) and every \(a\geq0\),
\begin{equation}
|\gamma_a(\alpha)|
\leq
C_{\mathfrak A}
(a+1)^{-1-\alpha},
\label{eq:spectral-compact-power-coefficient-bound}
\end{equation}
and
\begin{equation}
|\widetilde\gamma_a(\alpha)|
\leq
C_{\mathfrak A}
(a+1)^{-1-\alpha}
\log(a+2).
\label{eq:spectral-compact-derivative-coefficient-bound}
\end{equation}
Consequently,
\begin{equation}
\sup_{\alpha\in\mathfrak A}
\sum_{a\geq2}
a
\left[
|\gamma_a(\alpha)|^2
+
|\widetilde\gamma_a(\alpha)|^2
\right]
<\infty.
\label{eq:spectral-uniform-weighted-square-summability}
\end{equation}
\end{proposition}

\begin{proof}
It suffices first to consider all sufficiently large \(a\), say
\(a>\alpha_++1\).  From \eqref{eq:spectral-gamma-formula} and
\[
(-1)^a\binom{\alpha}{a}
=
\frac{\Gamma(a-\alpha)}
{\Gamma(-\alpha)\Gamma(a+1)},
\]
we obtain
\begin{equation}
\gamma_a(\alpha)
=
\frac{\Gamma(1+\alpha)}{\Gamma(-\alpha)}
\frac{\Gamma(a-\alpha)}{\Gamma(a+1)}.
\label{eq:gamma-equi-entire}
\end{equation}
Here and below \(1/\Gamma(-\alpha)\) is understood through the reciprocal
Gamma function.  Since \(1/\Gamma\) is entire, the functions
\[
\alpha\longmapsto
\frac{\Gamma(1+\alpha)}{\Gamma(-\alpha)}
\]
and its \(\alpha\)-derivative are bounded on the compact interval
\(\mathfrak A\), including when \(\mathfrak A\) contains positive integers.
Differentiating \eqref{eq:gamma-equi-entire} and using \eqref{eq:spectral-tilde-gamma-derivative} gives
\begin{equation}
\widetilde\gamma_a(\alpha)
=
\partial_\alpha
\left[
\frac{\Gamma(1+\alpha)}{\Gamma(-\alpha)}
\right]
\frac{\Gamma(a-\alpha)}{\Gamma(a+1)}
-
\frac{\Gamma(1+\alpha)}{\Gamma(-\alpha)}
\psi(a-\alpha)
\frac{\Gamma(a-\alpha)}{\Gamma(a+1)},
\label{eq:spectral-tilde-gamma-ratio}
\end{equation}
where \(\psi=\Gamma'/\Gamma\).
The classical gamma-ratio and digamma asymptotics
\cite[Eqs.~(5.11.12)--(5.11.13) and (5.11.2)]{NISTDLMF},
obtained from Stirling's formula with its remainder on a fixed compact set of shifts, yield
\[
\frac{\Gamma(a-\alpha)}{\Gamma(a+1)}
\leq
C_{\mathfrak A}
(a+1)^{-1-\alpha},
\qquad
\alpha\in\mathfrak A,
\]
for all sufficiently large \(a\), while
\[
|\psi(a-\alpha)|
\leq
C_{\mathfrak A}\log(a+2).
\]
Consequently, \eqref{eq:spectral-compact-power-coefficient-bound} and \eqref{eq:spectral-compact-derivative-coefficient-bound} hold uniformly for \(\alpha\in\mathfrak A\) and all sufficiently large \(a\).

For each fixed \(a\), both
\(\gamma_a(\alpha)\) and \(\widetilde\gamma_a(\alpha)\) are continuous in
\(\alpha\).  Since only finitely many smaller values of \(a\) remain,
enlarging \(C_{\mathfrak A}\) if necessary proves
\eqref{eq:spectral-compact-power-coefficient-bound} and
\eqref{eq:spectral-compact-derivative-coefficient-bound} for every
\(a\geq0\).

Finally, since \(\alpha\geq\alpha_->0\), for \(a\geq2\),
\[
a
\left[
|\gamma_a(\alpha)|^2
+
|\widetilde\gamma_a(\alpha)|^2
\right]
\leq
C_{\mathfrak A}
a^{-1-2\alpha_-}\log^2(a+2),
\]
uniformly for \(\alpha\in\mathfrak A\).  The series on the right is
summable, and hence
\eqref{eq:spectral-uniform-weighted-square-summability} follows.
\end{proof}

\subsection{Sharp Laguerre truncations}
\label{subsec:spectral-sharp-truncations}

Now we approximate the power functions by finite Laguerre sums and isolate
the error that must be controlled.
A sharp truncation keeps all terms up to a specified degree, with their
coefficients unchanged, and discards the rest.
The retained terms give the finite approximations studied in
Section~\ref{sec:finite-mode-process}.
For the discarded terms, we must control their centered angular
integrals uniformly over the nested polynomial lengths, as in
Section~\ref{sec:boundaries-fixed-tails}.
The cutoff remains fixed while \(N\to\infty\), and is then sent to infinity.

For every integer \(d\geq0\), define the sharp Laguerre projection by
\begin{equation}
\LagProj_{\leq d}f
=
\sum_{a=0}^{d}\widehat f_aL_a.
\label{eq:spectral-sharp-projection}
\end{equation}
For \(\alpha>0\), set
\begin{equation}
h_{\alpha,\leq d}
:=
\LagProj_{\leq d}h_\alpha
=
\sum_{a=0}^{d}
\gamma_a(\alpha)L_a.
\label{eq:spectral-power-truncation}
\end{equation}
Define the corresponding sharp tail by
\begin{equation}
R_{\alpha,d}
:=
\left(
I-\LagProj_{\leq d}
\right)h_\alpha
=
\sum_{a>d}
\gamma_a(\alpha)L_a.
\label{eq:spectral-power-tail}
\end{equation}
Under the reparametrization \(p=2\alpha\), this is precisely the
remainder denoted by \(R_{p,d}\) in \eqref{eq:intro-proof-remainder}.
Likewise, define
\begin{equation}
\widetilde R_{\alpha,d}
:=
\left(
I-\LagProj_{\leq d}
\right)\widetilde h_\alpha
=
\sum_{a>d}
\widetilde\gamma_a(\alpha)L_a.
\label{eq:spectral-logarithmic-tail}
\end{equation}

Since \(\LagProj_{\leq d}\) is a bounded finite-rank operator on
\(L^2(\nu)\), exponent differentiation commutes with
\(\LagProj_{\leq d}\).
Hence,
\begin{equation}
\partial_\alpha h_{\alpha,\leq d}
=
\LagProj_{\leq d}\widetilde h_\alpha
=
\sum_{a=0}^{d}
\widetilde\gamma_a(\alpha)L_a,
\qquad
\partial_\alpha R_{\alpha,d}
=
\widetilde R_{\alpha,d}.
\label{eq:spectral-sharp-truncation-derivatives}
\end{equation}

\subsection{Laguerre energies, generator, and semigroup}
\label{subsec:spectral-laguerre-semigroup}

We shall use Laguerre degree both to measure the coefficients and to smooth
the power functions.
The energies below weight the squared Laguerre coefficients by powers
of the degree, making their dependence on higher modes explicit.
The generator is diagonal in the Laguerre basis on \(L^2(\nu)\), and
its semigroup smooths the power family by damping the higher degrees.
We use this semigroup to form the Abel blocks in
Section~\ref{subsec:abel-blocks-dyadic-telescoping}.
The corresponding energy bounds enter the estimates over nested
polynomial lengths in Appendix~\ref{app:block-addition-refinement}.

For
\[
f=\sum_{a=0}^{\infty}\widehat f_aL_a
\]
and \(s\geq0\), define the Laguerre energy
\begin{equation}
\mathfrak E_s(f)
=
\sum_{a\geq1}
a^s|\widehat f_a|^2,
\label{eq:spectral-Ds}
\end{equation}
whenever the series converges.
The constant term is omitted.

The ordinary Laguerre differential operator is
\begin{equation}
\LagGen f
:=
xf''+(1-x)f'.
\label{eq:spectral-Laguerre-generator}
\end{equation}
The Laguerre differential equation gives
\begin{equation}
\LagGen L_a=-aL_a.
\label{eq:spectral-Laguerre-eigenvalue}
\end{equation}
For smooth compactly supported \(f,g\), integration by parts gives
\begin{equation}
-\int_0^\infty
f(x)\LagGen g(x)e^{-x}\,\dd x
=
\int_0^\infty
x f'(x)g'(x)e^{-x}\,\dd x.
\label{eq:spectral-Dirichlet-integration-by-parts}
\end{equation}
Consequently, by \eqref{eq:spectral-Laguerre-eigenvalue} and \eqref{eq:spectral-Dirichlet-integration-by-parts},
the corresponding closed quadratic form satisfies
\begin{equation}
\mathfrak E_1(f)
=
\int_0^\infty
x|f'(x)|^2e^{-x}\,\dd x,
\label{eq:spectral-D1-Dirichlet-form}
\end{equation}
whenever \(\mathfrak E_1(f)<\infty\), where the derivative is understood
in the form-domain sense.
These are standard facts from Laguerre polynomial theory;
see, for example, \cite{Szego1975,NISTDLMF}.
It is convenient to introduce the positive number operator
\begin{equation}
\LagNum
:=
-\LagGen,
\qquad
\LagNum L_a=aL_a.
\label{eq:spectral-number-operator}
\end{equation}

The same spectral eigenvalues define the associated Laguerre semigroup.
For \(\tau\geq0\), define the Laguerre semigroup spectrally by
\begin{equation}
P_\tau f
=
\sum_{a=0}^{\infty}
e^{-a\tau}\widehat f_aL_a,
\qquad
P_\tau L_a
=
e^{-a\tau}L_a.
\label{eq:spectral-analytic-semigroup}
\end{equation}
This is the radial form of the classical Gaussian Mehler semigroup.
More explicitly, if \(Z\) is standard circular Gaussian, with density
\(\pi^{-1}e^{-|z|^2}\) on \(\C\), then, with \(r=e^{-\tau}\),
\begin{equation}
P_\tau f(x)
=\E f\!\left(\left|\sqrt r\,\sqrt x+\sqrt{1-r}\,Z\right|^2\right).
\label{eq:spectral-radial-Mehler}
\end{equation}
The associated radial kernel is the Hille--Hardy Poisson kernel
\cite[Eq.~(18.18.27)]{NISTDLMF}; we display it in
Lemma~\ref{lem:uniform-local-semigroup-continuity}.
For completeness, the Gaussian Laplace transform and the Laguerre generating function give, for \(|z|\) sufficiently small,
\[
\E\left[\frac1{1-z}
 \exp\!\left(-\frac{z}{1-z}
 \left|\sqrt r\,\sqrt x+\sqrt{1-r}\,Z\right|^2\right)\right]
=\frac1{1-rz}\exp\!\left(-\frac{rzx}{1-rz}\right).
\]
Comparison of coefficients proves \eqref{eq:spectral-radial-Mehler}
on every \(L_a\).  Gaussian invariance shows that its invariant law is
\(\nu\), since \(|Z|^2\) has density \(e^{-x}\).  Jensen's inequality
therefore extends the identity to \(L^2(\nu)\).  In particular, the
normalization is \(P_\tau L_a=e^{-a\tau}L_a\), with generator
\(xf''+(1-x)f'\). This fixes the normalization of the semigroup parameter.

The spectral representation immediately gives the semigroup law
\[
P_{\tau+\sigma}=P_\tau P_\sigma.
\]
By Parseval's identity,
\begin{equation}
\|P_\tau f\|_{L^2(\nu)}
\leq
\|f\|_{L^2(\nu)},
\qquad
\tau\geq0.
\label{eq:spectral-semigroup-properties}
\end{equation}
Moreover, the semigroup is strongly continuous, that is,
\begin{equation}
P_\tau f
\longrightarrow
f
\qquad
\text{in }L^2(\nu),
\quad
\text{as }\tau\downarrow0.
\label{eq:Ptau-semigroup-strong-continuous}
\end{equation}
Thus, \((P_\tau)_{\tau\geq0}\) is a strongly continuous contraction
semigroup on \(L^2(\nu)\).
Since both \(P_\tau\) and \(\LagProj_{\leq d}\) are diagonal in the
Laguerre basis,
\begin{equation}
P_\tau\LagProj_{\leq d}
=
\LagProj_{\leq d}P_\tau,
\qquad
d\geq0,\quad \tau\geq0.
\label{eq:spectral-semigroup-projection-commutation}
\end{equation}
Finally, \(P_\tau\) is bounded on \(L^2(\nu)\), while
\(\alpha\mapsto h_\alpha\) is \(C^1\) as an \(L^2(\nu)\)-valued map.
Hence exponent differentiation commutes with \(P_\tau\):
\begin{equation}
\partial_\alpha P_\tau h_\alpha
=
P_\tau\widetilde h_\alpha.
\label{eq:spectral-semigroup-p-derivative}
\end{equation}

We use the energies to estimate truncation errors and the semigroup to smooth the power functions. The commutation identities above allow the same operations to be applied to exponent derivatives.

\section{Limits for finite Laguerre truncations}
\label{sec:finite-mode-process}
This section proves joint convergence for finitely many Laguerre modes and combines their limits into the Gaussian approximations used later to recover the full fluctuation field.
We first prove the finite Laguerre vector functional CLT, Theorem~\ref{thm:finite-laguerre-vector}. Each Laguerre statistic is approximated by a polynomial in the signs that keeps only its highest homogeneous degree. This highest-degree part is called the top chaos. The functional de Jong theorem then gives the joint process limit. Linearity yields the fixed-truncation convergence for the power functions in Corollary~\ref{cor:finite-power-truncation-joint}.

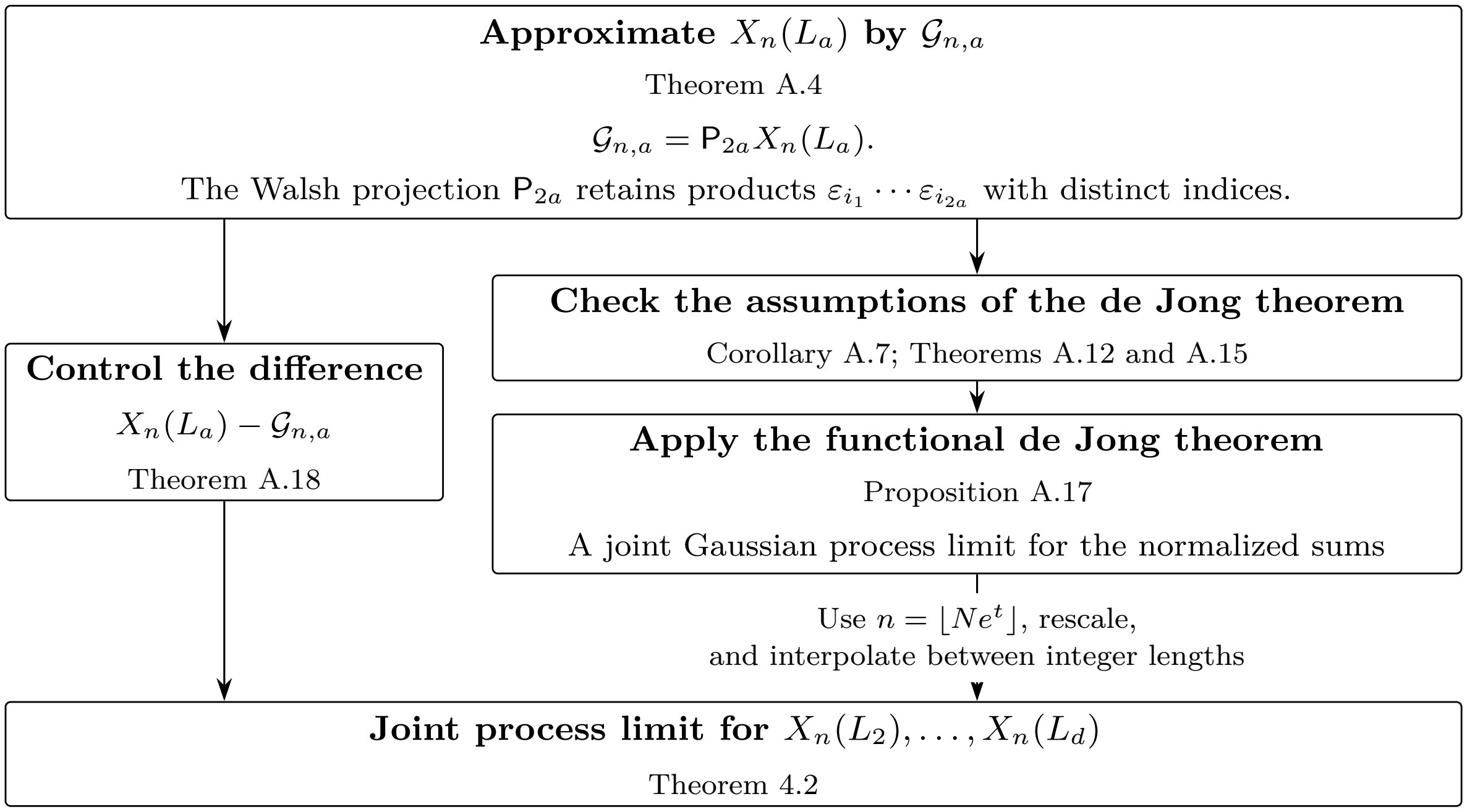
\begin{figure}[htbp]
  \centering
  \resizebox{\textwidth}{!}{%
\begin{tikzpicture}[
  x=1cm,y=1cm,
  foBox/.style={draw=black,line width=.45pt,rounded corners=1.3pt,
    align=flush center,inner xsep=5pt,inner ysep=3pt,
    font=\fontsize{8.2}{10}\selectfont,
    execute at begin node={\hyphenpenalty=10000\exhyphenpenalty=10000}},
  foInput/.style={foBox},
  foStep/.style={foBox},
  foConclusion/.style={foBox},
  foArrow/.style={-{Stealth[length=2.1mm,width=1.4mm]},
    draw=black,line width=.5pt}
]

\node[foStep,text width=12.8cm,minimum height=1.55cm]
(extract) at (5.2,0) {
  {\fontsize{9}{11}\selectfont\bfseries
  Approximate $X_n(L_a)$ by $\mathcal G_{n,a}$}\\[3pt]
  {\fontsize{7.5}{9.5}\selectfont
  Theorem~\ref{thm:fixed-mode-pure-chaos-reduction}}\\[4pt]
  $\mathcal G_{n,a}=\mathsf P_{2a}X_n(L_a).$\\[3pt]
  The Walsh projection $\mathsf P_{2a}$ retains products
  $\varepsilon_{i_1}\cdots\varepsilon_{i_{2a}}$
  with distinct indices.
};

\node[foInput,text width=8.4cm,minimum height=.95cm]
(checks) at (7.4,-1.95) {
  {\fontsize{9}{11}\selectfont\bfseries
  Check the assumptions of the de Jong theorem}\\[3pt]
  {\fontsize{7.5}{9.5}\selectfont
  Corollary~\ref{cor:top-chaos-second-moment};
  Theorems~\ref{thm:fixed-mode-low-influence}
  and~\ref{thm:mixed-fourth-cumulant-top-chaos}}
};

\node[foInput,text width=3.6cm,minimum height=1.4cm]
(error) at (.6,-2.8) {
  {\fontsize{9}{11}\selectfont\bfseries Control the difference}\\[4pt]
  $X_n(L_a)-\mathcal G_{n,a}$\\[4pt]
  {\fontsize{7.5}{9.5}\selectfont
  Theorem~\ref{thm:uniform-pure-chaos-transfer}}
};

\node[foStep,text width=8.4cm,minimum height=1.15cm]
(dejong) at (7.4,-3.45) {
  {\fontsize{9}{11}\selectfont\bfseries
  Apply the functional de Jong theorem}\\[3pt]
  {\fontsize{7.5}{9.5}\selectfont
  Proposition~\ref{prop:top-chaos-functional-de-jong}}\\[4pt]
  A joint Gaussian process limit for the normalized sums
};

\node[foConclusion,text width=12.8cm,minimum height=.9cm]
(limit) at (5.2,-5.8) {
  {\fontsize{9}{11}\selectfont\bfseries
  Joint process limit for $X_n(L_2),\ldots,X_n(L_d)$}\\[4pt]
  {\fontsize{7.5}{9.5}\selectfont
  Theorem~\ref{thm:finite-laguerre-vector}}
};

\draw[foArrow] ([xshift=2.2cm]extract.south) -- (checks.north);
\draw[foArrow] (checks.south) -- (dejong.north);
\draw[foArrow] (dejong.south) -- ([xshift=2.2cm]limit.north)
  node[midway,fill=white,inner sep=3pt,align=flush center,
    font=\fontsize{7.5}{9.5}\selectfont]
    {Use $n=\lfloor Ne^t\rfloor$, rescale,\\
     and interpolate between integer lengths};
\draw[foArrow] ([xshift=-4.6cm]extract.south) -- (error.north);
\draw[foArrow] (error.south) -- ([xshift=-4.6cm]limit.north);
\end{tikzpicture}

}
  \caption{Proof structure for
  Theorem~\ref{thm:finite-laguerre-vector}.
  Fix $d\ge2$. Prove a joint process limit for
  $\mathcal G_{n,2},\ldots,\mathcal G_{n,d}$ and control their differences
  from $X_n(L_2),\ldots,X_n(L_d)$ over the length window.}
  \label{fig:finite-mode-proof-roadmap}
\end{figure}

\subsection{Finite Laguerre vector convergence}
\label{subsec:finite-laguerre-vector-convergence}

We now prove that any fixed finite collection of Laguerre modes converges jointly to independent Gaussian processes in logarithmic time.
For a fixed integer \(d\geq2\), set
\[
\mathbf G_{\leq d}(t)
=
\left(
G^{(2)}(t),
\ldots,
G^{(d)}(t)
\right),
\qquad
t\in I_T,
\]
where \((G^{(a)})_{a\geq2}\) are the independent centered stationary Gaussian processes introduced in Section \ref{subsec:intro-main-results}.
Equivalently, these Gaussian modes may be realized as stationary Ornstein--Uhlenbeck processes:
\[
G^{(a)}(t)
=
2\sqrt{V_a\Bigl(a-\frac12\Bigr)}
\int_{-\infty}^{t}
e^{-(a-1/2)(t-u)}\,\dd W^{(a)}(u),
\qquad t\in\mathbb R,
\]
where \((W^{(a)})_{a\ge2}\) are independent two-sided standard Brownian motions.
This representation has covariance
\eqref{eq:intro-mode-covariance}.
Since the processes \(G^{(a)}\) have
continuous paths,
\[
\mathbf G_{\leq d}
\in
C(I_T;\mathbb R^{d-1})
\qquad\text{almost surely}.
\]

\begin{remark}
The decay rates \(a\) and \(a-1/2\) have different origins.
The eigenvalue identity
\[
\LagNum L_a=aL_a
\]
gives the spectral decay factor \(e^{-a\tau}\) for the \(a\)-th Laguerre component, whereas the rate \(a-1/2\) in the covariance across logarithmic times arises from the normalized covariance across nested polynomial lengths; see Theorem~\ref{thm:nested-top-chaos-covariance}.
\end{remark}

For \(a\geq2\), define
\begin{equation}
X_{N,L_a}^{\poly}
=
\operatorname{Poly}_{N,T}
\bigl[
(X_n(L_a))_{n_-\leq n\leq n_+}
\bigr],
\label{eq:finite-mode-Xa-polygonal}
\end{equation}
where \(\operatorname{Poly}_{N,T}\) is defined in Definition \ref{def:canonical-polygonal-interpolation}.
Set
\begin{equation}
\mathbf X_{N,\leq d}^{\poly}(t)
=
\left(
X_{N,L_2}^{\poly}(t),
\ldots,
X_{N,L_d}^{\poly}(t)
\right),
\label{eq:finite-mode-X-vector-polygonal}
\end{equation}
and the corresponding floor vector
\begin{equation}
\mathbf X_{N,\leq d}^{\fl}(t)
=
\left(
X_{n_N(t)}(L_2),
\ldots,
X_{n_N(t)}(L_d)
\right).
\label{eq:finite-mode-floor-vector}
\end{equation}

\begin{theorem}[Finite Laguerre vector functional CLT]
\label{thm:finite-laguerre-vector}
For every fixed integer \(d\geq2\) and every \(T>0\),
\begin{equation}
\mathbf X_{N,\leq d}^{\poly}
\Longrightarrow
\mathbf G_{\leq d}
\label{eq:finite-mode-polygonal-convergence}
\end{equation}
in
\(
C(I_T;\R^{d-1}).
\)
Moreover,
\begin{equation}
\mathbf X_{N,\leq d}^{\fl}
\Longrightarrow
\mathbf G_{\leq d}
\label{eq:finite-mode-floor-convergence}
\end{equation}
in
\(
D(I_T;\R^{d-1})
\)
with the Skorokhod \(J_1\)-topology.
\end{theorem}

\begin{proof}
For a fixed \(d\), Proposition~\ref{prop:top-chaos-functional-de-jong}
applies the multivariate functional de Jong theorem
\cite[Theorem~1.4]{DoblerKasprzakPeccati2022} to the prefix processes
\[
W_{M,a}(u)=\frac{Q_{\lfloor Mu\rfloor,a}}{\sigma_{M,a}},
\qquad
Q_{n,a}=n^{a-1/2}\mathcal G_{n,a},
\qquad 2\leq a\leq d,
\]
where \(\mathcal G_{n,a}\) and \(\sigma_{M,a}\) are defined in
\eqref{eq:chaos-top-chaos} and
\eqref{eq:chaos-prefix-Q-normalization}, respectively.
It gives joint convergence to independent time-changed Brownian motions
\(B_a(u^{2a-1})\) in \(D([0,1];\R^{d-1})\) with the \(J_1\)-topology.
All limiting paths are continuous.

Put
\[
M_N=\lceil Ne^T\rceil,
\qquad
u_N(t)=\frac{Ne^t}{M_N},
\qquad
\lambda_a=a-\frac12.
\]
Then \(u_N(t)\leq1\),
\(\lfloor M_Nu_N(t)\rfloor=n_N(t)\) exactly, and
\(u_N(t)\to e^{t-T}\) uniformly on \(I_T\).
Moreover,
\[
\mathcal G_{n_N(t),a}
=\frac{\sigma_{M_N,a}}{n_N(t)^{\lambda_a}}
 W_{M_N,a}(u_N(t)),
\qquad
\frac{\sigma_{M_N,a}}{n_N(t)^{\lambda_a}}
\longrightarrow
\sqrt{2V_a}\,e^{-\lambda_a(t-T)}
\]
uniformly on \(I_T\), by
\eqref{eq:chaos-prefix-variance-asymptotic}.
At a continuous limiting path, convergence in \(J_1\) is uniform. Composition with the converging time changes and multiplication by the converging deterministic factors are therefore continuous at that path.
The continuous-mapping theorem therefore yields
\[
\bigl(\mathcal G_{n_N(t),a}\bigr)_{a=2}^{d}
\Longrightarrow
\bigl(\widetilde G^{(a)}(t)\bigr)_{a=2}^{d}
\quad\text{in }D(I_T;\R^{d-1}),
\]
where
\[
\widetilde G^{(a)}(t)
=\sqrt{2V_a}\,e^{-\lambda_a(t-T)}
B_a\bigl(e^{2\lambda_a(t-T)}\bigr).
\]
This logarithmic change of time is an instance of Lamperti's
self-similar/stationary correspondence \cite[Eq.~(7)]{Lamperti1962}.
For \(s\leq t\),
\[
\E\bigl[\widetilde G^{(a)}(s)\widetilde G^{(b)}(t)\bigr]
=\ind_{\{a=b\}}\,2V_a e^{-\lambda_a(t-s)}.
\]
Thus this continuous Gaussian vector has the law of
\(\mathbf G_{\leq d}\).

The same prefix convergence implies
\[
\max_{1\leq n\leq M_N}
\left|W_{M_N,a}(n/M_N)-W_{M_N,a}((n-1)/M_N)\right|
\xrightarrow{\Prob}0,
\qquad
\sup_{0\leq u\leq1}|W_{M_N,a}(u)|=O_{\Prob}(1).
\]
On the length window, the factors
\(\sigma_{M_N,a}n^{-\lambda_a}\) are uniformly bounded and their
adjacent differences are \(O_{a,T}(N^{-1})\).  It follows that
\begin{equation}
\max_{n_-\leq n<n_+}
|\mathcal G_{n+1,a}-\mathcal G_{n,a}|
\xrightarrow{\Prob}0.
\label{eq:finite-mode-de-jong-jumps}
\end{equation}
For the polygonal process
\eqref{eq:chaos-top-polygonal-process}, including its initial interval
and constant final segment, we have
\[
\left\|\mathcal G_{N,a}^{\poly}
      -\mathcal G_{n_N(\cdot),a}\right\|_\infty
\leq
\max_{n_-\leq n<n_+}|\mathcal G_{n+1,a}-\mathcal G_{n,a}|.
\]
Consequently the polygonal top-chaos vector converges in
\(C(I_T;\R^{d-1})\) to \(\mathbf G_{\leq d}\).
Finally, Theorem~\ref{thm:uniform-pure-chaos-transfer} shows that the
full vector of Laguerre statistics and its highest-degree approximation
have a uniform difference tending to zero in probability, for both the polygonal and floor processes. Slutsky's theorem proves
\eqref{eq:finite-mode-polygonal-convergence} and
\eqref{eq:finite-mode-floor-convergence}.
\end{proof}

\subsection{Fixed truncations of the power functions}
\label{subsec:fixed-power-truncations}

We now apply the finite Laguerre vector limit to the finite spectral truncations of the power family \(h_\alpha(x)=x^\alpha\).

For \(d\geq2\) and \(\alpha>0\), recall from
\eqref{eq:spectral-power-truncation} that
\[
h_{\alpha,\leq d}
=
\LagProj_{\leq d}h_\alpha
=
\sum_{a=0}^{d}\gamma_a(\alpha)L_a.
\]
The finite polygonal joint field is
\begin{equation}
\mathcal Z_{N,d}^{\poly}(\alpha,t)
= X_{N,h_{\alpha,\leq d}}^{\poly}(t)= \sum_{a=2}^{d}
\gamma_a(\alpha)X_{N,L_a}^{\poly}(t).
\label{eq:finite-mode-joint-polygonal-field}
\end{equation}
Its floor counterpart is
\begin{equation}
\mathcal Z_{N,d}^{\fl}(t)(\alpha)
= X_{n_N(t)}(h_{\alpha,\leq d})= \sum_{a=2}^{d}
\gamma_a(\alpha)X_{n_N(t)}(L_a).
\label{eq:finite-mode-joint-floor-field}
\end{equation}
Define the finite Gaussian joint field by
\begin{equation}
G_d(\alpha,t)
=
\sum_{a=2}^{d}
\gamma_a(\alpha)G^{(a)}(t).
\label{eq:finite-mode-Gd}
\end{equation}

\begin{corollary}
\label{cor:finite-power-truncation-joint}
For every fixed \(d\geq2\) and every compact interval
\(
\mathfrak A:=[\alpha_-,\alpha_+]\Subset(0,\infty),
\)
\begin{equation}
\mathcal Z_{N,d}^{\poly}
\Longrightarrow
G_d
\label{eq:finite-mode-joint-polygonal-convergence}
\end{equation}
in
\(
C^1(\mathfrak A;C(I_T)).
\)
In particular, the convergence holds in
\(
C(\mathfrak A\times I_T).
\)
Moreover,
\begin{equation}
\mathcal Z_{N,d}^{\fl}
\Longrightarrow
\bigl(t\longmapsto G_d(\cdot,t)\bigr)
\label{eq:finite-mode-joint-floor-convergence}
\end{equation}
in
\(
D(I_T;C(\mathfrak A))
\)
with the \(J_1\)-topology.
\end{corollary}

\begin{proof}
Since \(d\) is fixed and
\(\gamma_a\in C^1(\mathfrak A)\) for \(2\leq a\leq d\), the linear map
\[
(f_2,\ldots,f_d)
\longmapsto
\left[
(\alpha,t)
\longmapsto
\sum_{a=2}^{d}\gamma_a(\alpha)f_a(t)
\right]
\]
is bounded from
\[
C(I_T;\R^{d-1})
\quad\text{to}\quad
C^1(\mathfrak A;C(I_T)).
\]
Hence, the continuous mapping theorem applied to
\eqref{eq:finite-mode-polygonal-convergence}
gives
\eqref{eq:finite-mode-joint-polygonal-convergence}.

For the floor process, the bounded linear map
\[
(x_2,\ldots,x_d)
\longmapsto
\left[
\alpha
\longmapsto
\sum_{a=2}^{d}\gamma_a(\alpha)x_a
\right]
\]
from \(\R^{d-1}\) to \(C(\mathfrak A)\) induces a continuous map
\[
D(I_T;\R^{d-1})
\longrightarrow
D(I_T;C(\mathfrak A))
\]
for the \(J_1\)-topology: the same time changes may be used, and the
error in the path values is controlled by the norm of the state-space map.
Applying the continuous mapping theorem to
\eqref{eq:finite-mode-floor-convergence}
gives
\eqref{eq:finite-mode-joint-floor-convergence}.
\end{proof}

\section{Construction of the full Gaussian field}
\label{sec:full-gaussian-field}

We now pass from the finite Gaussian fields \(G_d\), defined in \eqref{eq:finite-mode-Gd}, to the full spectral series.
To control the truncation tails in
\(C^1(\mathfrak A;C(I_T))\), we combine the coefficient decay from Proposition~\ref{prop:compact-coefficient-bounds} with a maximal estimate
for the individual Gaussian modes.
We then identify the covariance kernel
of the resulting full field.

We use the \(\alpha\)-parametrization fixed in Section~\ref{sec:spectral-framework}.
Accordingly, the field written below as \(G(\alpha,t)\) is the same limiting Gaussian field that is written as \(G(p,t)\) in \eqref{eq:intro-G-mode}, expressed in the coordinates \(\alpha=p/2\).
Indeed, the identity \eqref{eq:equi-p-alpha} translates the Gaussian series in \eqref{eq:intro-G-mode} directly into its \(\alpha\)-parametrized form in \eqref{eq:gaussian-full-series}.
We retain the symbol \(G\) under this
reparametrization.

Throughout this section, we fix
\[
\mathfrak A=[\alpha_-,\alpha_+]\Subset(0,\infty)
\qquad\text{and}\qquad
I_T=[-T,T],
\quad T>0.
\]
Recall from \eqref{eq:intro-Va} that
\[
f_a
=
\ind_{[0,1]}^{*a},
\qquad
V_a
=
\int_0^a f_a(u)^2\,\dd u.
\]
The function \(f_a=\ind_{[0,1]}^{*a}\) is a probability density, so
\(\|f_a\|_{L^1(\mathbb R)}=1\).
Young's inequality also gives
\(\|f_a\|_{L^\infty(\mathbb R)}\leq1\).
Since \(f_a\) is nonzero, it follows that
\begin{equation}
0<V_a\leq1.
\label{eq:gaussian-Va-bound}
\end{equation}

\begin{lemma}
\label{lem:modewise-gaussian-maximal}
For every \(T>0\), there exists \(C_T<\infty\) such that, for every
\(a\geq2\),
\begin{equation}
\left(
\E
\|G^{(a)}\|_{C(I_T)}^2
\right)^{1/2}
\leq
C_T
\sqrt{V_a\log(a+2)}.
\label{eq:modewise-gaussian-maximal}
\end{equation}
\end{lemma}

\begin{proof}
This is the standard Gaussian entropy argument, written on dyadic grids to keep the dependence on the degree explicit.
The canonical metric is the square root of the mean squared increment:
\[
d_a(s,t)^2=4V_a(1-e^{-\lambda_a|t-s|})
\leq4V_a\min\{1,\lambda_a|t-s|\},
\qquad \lambda_a=a-\tfrac12.
\]
The correlation length therefore gives a grid on \(I_T\) with \(O_T(a+1)\) intervals. We control the supremum by bounding the process on this grid and then summing its increments over successively finer grids.
For \(a\geq2\), write
\begin{equation}
\lambda_a:=a-\frac12.
\label{eq:gaussian-lambda-a}
\end{equation}
Choose a finite grid
\[
\mathcal D_{a,0}\subseteq I_T
\]
containing both endpoints and having mesh at most \(\lambda_a^{-1}\).
It may be chosen so that
\[
\#\mathcal D_{a,0}
\leq
C_T(a+1).
\]
By \eqref{eq:intro-mode-covariance}  and the standard maximal inequality for finite centered Gaussian families
(see, e.g., \cite{Janson1997}), we obtain
\begin{equation}
\left(
\E
\max_{t\in\mathcal D_{a,0}}
|G^{(a)}(t)|^2
\right)^{1/2}
\leq
C
\sqrt{
V_a\log\!\left(\#\mathcal D_{a,0}+1\right)
}
\leq
C_T
\sqrt{
V_a\log(a+2)
}.
\label{eq:gaussian-coarse-grid-maximum}
\end{equation}

For \(j\geq0\), let
\(\mathcal D_{a,j+1}\) be obtained from \(\mathcal D_{a,j}\) by bisecting every grid interval.
Each new point
\[
u\in
\mathcal D_{a,j+1}\setminus\mathcal D_{a,j}
\]
may be assigned a point on the preceding grid,
\(\pi_j(u)\in\mathcal D_{a,j}\), satisfying
\[
|u-\pi_j(u)|
\leq
\frac{2^{-j}}{\lambda_a}.
\]
By \eqref{eq:intro-mode-covariance},
\begin{equation}
\E
\left|
G^{(a)}(u)-G^{(a)}(\pi_j(u))
\right|^2
\leq
4V_a2^{-j}.
\label{eq:gaussian-refinement-increment-variance}
\end{equation}
The number of such increments is at most
\[
C_T(a+1)2^j.
\]
Applying the same Gaussian maximal inequality to this increment family gives
\begin{equation}
\left(
\E
\max_{
u\in\mathcal D_{a,j+1}\setminus\mathcal D_{a,j}}
\left|
G^{(a)}(u)-G^{(a)}(\pi_j(u))
\right|^2
\right)^{1/2}
\leq
C_T
\sqrt{V_a}\,
2^{-j/2}
\sqrt{
\log(a+2)+j+1
}.
\label{eq:gaussian-one-refinement-level}
\end{equation}

The union \(\bigcup_{j\geq0}\mathcal D_{a,j}\) is dense in \(I_T\).
By telescoping along these maps \(\pi_j\) and using
continuity of the sample paths, we obtain, almost surely,
\begin{align}
\|G^{(a)}\|_{C(I_T)}
\leq
\max_{t\in\mathcal D_{a,0}}
|G^{(a)}(t)|
+
\sum_{j=0}^{\infty}
\max_{
u\in\mathcal D_{a,j+1}\setminus\mathcal D_{a,j}}
\left|
G^{(a)}(u)-G^{(a)}(\pi_j(u))
\right|.
\label{eq:gaussian-pathwise-chaining}
\end{align}
Minkowski's inequality together with
\eqref{eq:gaussian-coarse-grid-maximum}
and
\eqref{eq:gaussian-one-refinement-level}
therefore yields
\[
\left(
\E
\|G^{(a)}\|_{C(I_T)}^2
\right)^{1/2}
\leq
C_T\sqrt{V_a\log(a+2)}
+
C_T\sqrt{V_a}
\sum_{j=0}^{\infty}
2^{-j/2}
\sqrt{
\log(a+2)+j+1
}.
\]
This proves \eqref{eq:modewise-gaussian-maximal}.
\end{proof}

The maximal estimate and the coefficient decay give a summable bound for the Gaussian series and its exponent derivative.

\begin{theorem}[Full Gaussian field and truncation]
\label{thm:full-gaussian-field}
Let
\(
\mathfrak A=[\alpha_-,\alpha_+]\Subset(0,\infty)
\)
and \(T>0\), and let \(G_d\) be the finite Gaussian fields defined in
\eqref{eq:finite-mode-Gd}.
Then, \((G_d)_{d\geq2}\) converges in
\[
L^2\!\left(
\Omega;
C^1(\mathfrak A;C(I_T))
\right)
\]
to a centered Gaussian random field \(G\).
In particular, \(G\) may be chosen with sample paths in \(C^1(\mathfrak A;C(I_T))\), and
\begin{equation}
G(\alpha,t)
=
\sum_{a=2}^{\infty}
\gamma_a(\alpha)G^{(a)}(t)
\label{eq:gaussian-full-series}
\end{equation}
in
\(
L^2\!\left(
\Omega;
C^1(\mathfrak A;C(I_T))
\right).
\)
Moreover,
\begin{equation}
\partial_\alpha G(\alpha,t)
=
\sum_{a=2}^{\infty}
\widetilde\gamma_a(\alpha)G^{(a)}(t)
\label{eq:gaussian-termwise-alpha-derivative}
\end{equation}
in
\(
L^2\!\left(
\Omega;
C(\mathfrak A;C(I_T))
\right),
\)
where \(\widetilde\gamma_a\) is defined in \eqref{eq:spectral-tilde-gamma-definition} and \eqref{eq:spectral-tilde-gamma-derivative}.
\end{theorem}

Since \(p=2\alpha\), this affine reparametrization identifies
\(C^1(\mathfrak A;C(I_T))\) with \(C^1(K;C(I_T))\); in particular,
\[
\partial_p G(p,t)
=
\frac12\,\partial_\alpha G(\alpha,t)\big|_{\alpha=p/2}.
\]
Thus, the \(C^1\)-regularity established above is equivalent to the \(p\)-parametrized regularity stated in Theorem~\ref{thm:intro-joint-field}.

\begin{proof}
For \(m>d\), Minkowski's inequality gives
\begin{align}
\left(
\E
\|G_m-G_d\|_{C^1(\mathfrak A;C(I_T))}^2
\right)^{1/2}
&\leq
\sum_{a=d+1}^{m}
\left[
\|\gamma_a\|_{C(\mathfrak A)}
+
\|\widetilde\gamma_a\|_{C(\mathfrak A)}
\right]
\left(
\E
\|G^{(a)}\|_{C(I_T)}^2
\right)^{1/2}.
\label{eq:gaussian-Minkowski-C1-tail}
\end{align}
Using Proposition~\ref{prop:compact-coefficient-bounds}, Lemma~\ref{lem:modewise-gaussian-maximal},
and \eqref{eq:gaussian-Va-bound}, we obtain
\begin{equation}
\begin{aligned}
\left(
\E
\|G_m-G_d\|_{C^1(\mathfrak A;C(I_T))}^2
\right)^{1/2}
&\leq
C_{\mathfrak A,T}
\sum_{a>d}
(a+1)^{-1-\alpha_-}
\log^{3/2}(a+2)
\\
&\leq
C_{\mathfrak A,T}
d^{-\alpha_-}
\log^{3/2}(d+2).
\end{aligned}
\label{eq:gaussian-C1-tail-root}
\end{equation}
Thus, \((G_d)_{d\geq2}\) is Cauchy in
\[
L^2\!\left(
\Omega;
C^1(\mathfrak A;C(I_T))
\right),
\]
and hence converges to some
\[
G\in
L^2\!\left(
\Omega;
C^1(\mathfrak A;C(I_T))
\right).
\]
Recalling \eqref{eq:finite-mode-Gd}, this is precisely the convergence
\eqref{eq:gaussian-full-series} in
\(L^2(\Omega;C^1(\mathfrak A;C(I_T)))\).
Since each \(G_d\) is a centered Gaussian field, so is \(G\).

Since
\[
\partial_\alpha G_d(\alpha,t)
=
\sum_{a=2}^{d}
\widetilde\gamma_a(\alpha)G^{(a)}(t),
\]
the convergence in \(C^1(\mathfrak A;C(I_T))\) also yields
\eqref{eq:gaussian-termwise-alpha-derivative}.
\end{proof}

We conclude the construction by identifying the covariance structure of the limiting field.

\begin{proposition}
\label{prop:joint-covariance}
For \(\alpha,\beta\in\mathfrak A\) and \(s,t\in I_T\), the covariance
kernel of \(G\) is
\begin{equation}
\mathcal K(\alpha,\beta;s,t)
:=
\Cov\left(
G(\alpha,s),
G(\beta,t)
\right)
=
2\sum_{a=2}^{\infty}
V_a
\gamma_a(\alpha)\gamma_a(\beta)
e^{-(a-\frac12)|t-s|}.
\label{eq:gaussian-full-covariance-formula}
\end{equation}
The series converges absolutely and uniformly on
\[
\mathfrak A\times\mathfrak A\times I_T\times I_T.
\]
Consequently, \(\mathcal K\) is jointly continuous, and \(G\) is stationary with respect to the logarithmic-time coordinate.
\end{proposition}

\begin{proof}
For every \(d\geq2\), independence of the Gaussian modes and
\eqref{eq:intro-mode-covariance} give
\begin{equation}
\begin{aligned}
\Cov\left(
G_d(\alpha,s),
G_d(\beta,t)
\right)
&=
2
\sum_{a=2}^{d}
V_a
\gamma_a(\alpha)\gamma_a(\beta)
e^{-(a-\frac12)|t-s|}.
\end{aligned}
\label{eq:gaussian-finite-covariance-formula}
\end{equation}
By \eqref{eq:gaussian-Va-bound} and
Proposition~\ref{prop:compact-coefficient-bounds},
\[
2V_a
\left|
\gamma_a(\alpha)\gamma_a(\beta)
\right|
e^{-(a-\frac12)|t-s|}
\leq
C_{\mathfrak A}
(a+1)^{-2-2\alpha_-},
\]
uniformly for
\(\alpha,\beta\in\mathfrak A\) and \(s,t\in I_T\).
The majorant is summable, so the covariance series converges absolutely and uniformly on
\[
\mathfrak A\times\mathfrak A\times I_T\times I_T.
\]
Letting \(d\to\infty\) in
\eqref{eq:gaussian-finite-covariance-formula} and using Theorem~\ref{thm:full-gaussian-field} yields
\eqref{eq:gaussian-full-covariance-formula}.
Joint continuity follows from uniform convergence, while stationarity in logarithmic time is immediate from the dependence on \(s,t\) only through
\(|t-s|\).
\end{proof}

\section{Coefficient and Abel calculus}
\label{sec:coefficient-calculus}

This section gives the deterministic estimates for the ordinary and
logarithmic Abel blocks used in the maximal-tail bounds.  We first derive
coefficient and scaled Abel estimates, then
introduce a decomposition in the intensity variable and obtain the corresponding
Laguerre-energy bounds.  We use the spectral notation of
Section~\ref{sec:spectral-framework} and the special-function estimates
proved in Appendix~\ref{app:kummer-critical-profile}.

Throughout this section, fix
\[
\mathfrak A=[\alpha_-,\alpha_+]\Subset(0,\infty).
\]
All constants are uniform in the Abel scale, the intensity shell, and
\(\alpha\in\mathfrak A\), unless otherwise indicated.

\subsection{Sharp coefficient-energy bounds}
\label{subsec:sharp-coefficient-energy}

We first record the elementary dyadic consequences of the coefficient
bounds in Proposition~\ref{prop:compact-coefficient-bounds}, which are
uniform over the compact exponent interval.
The underlying coefficient identity comes from the classical Laguerre generating function; the estimates below keep explicit the uniformity in the exponent and the additional logarithm from exponent differentiation.

\begin{proposition}
\label{prop:sharp-coefficient-energy}
There exists \(C_{\mathfrak A}<\infty\) such that, uniformly for
\(\alpha\in\mathfrak A\) and dyadic \(A\geq1\),
\begin{align}
\sum_{A<a\leq2A}
a|\gamma_a(\alpha)|^2
&\leq
C_{\mathfrak A}A^{-2\alpha},
\label{eq:calculus-sharp-power-energy}
\\
\sum_{A<a\leq2A}
a|\widetilde\gamma_a(\alpha)|^2
&\leq
C_{\mathfrak A}
A^{-2\alpha}\log^2(A+2).
\label{eq:calculus-sharp-log-energy}
\end{align}
Consequently, for every integer \(d\geq2\), uniformly for
\(\alpha\in\mathfrak A\),
\[
\begin{aligned}
\sum_{a>d}a|\gamma_a(\alpha)|^2
&\leq
C_{\mathfrak A}d^{-2\alpha},
\\
\sum_{a>d}a|\widetilde\gamma_a(\alpha)|^2
&\leq
C_{\mathfrak A}
d^{-2\alpha}\log^2(d+2).
\end{aligned}
\]
\end{proposition}

\begin{proof}
By Proposition~\ref{prop:compact-coefficient-bounds},
uniformly for \(\alpha\in\mathfrak A\),
\[
\sum_{A<a\leq2A}a|\gamma_a(\alpha)|^2
\leq
C_{\mathfrak A}
\sum_{A<a\leq2A}a^{-1-2\alpha}
\leq
C_{\mathfrak A}A^{-2\alpha}.
\]
Since \(\log(a+2)\lesssim\log(A+2)\) on \(A<a\leq2A\),
the same argument gives
\[
\sum_{A<a\leq2A}a|\widetilde\gamma_a(\alpha)|^2
\leq
C_{\mathfrak A}
A^{-2\alpha}\log^2(A+2).
\]

The tail bounds follow by summing these estimates over dyadic blocks.
For the logarithmic family, use
\[
\log(2^j d+2)\lesssim (j+1)\log(d+2)
\]
and \(\alpha\geq\alpha_->0\), so that
\[
\sum_{j\geq0}(j+1)^2\,2^{-2j\alpha_-}<\infty.
\]
\end{proof}

\subsection{Abel blocks and dyadic telescoping}
\label{subsec:abel-blocks-dyadic-telescoping}

Now we recover the power functions from an initial smoothing and a telescoping sum of Abel blocks over dyadic scales.
For \(A\geq2\), define the Abel multiplier
\begin{equation}
m_A(a)
=
e^{-a/(2A)}
-
e^{-a/A}.
\label{eq:calculus-Abel-multiplier}
\end{equation}
It satisfies
\begin{equation}
0\leq m_A(a)
\leq
C
\min
\left\{
\frac aA,
e^{-a/(2A)}
\right\}.
\label{eq:calculus-Abel-multiplier-bound}
\end{equation}
Using the Laguerre semigroup \(P_\tau\) defined in
\eqref{eq:spectral-analytic-semigroup}, define the ordinary
and logarithmic Abel blocks by
\begin{equation}
\begin{aligned}
f_{A,\alpha}
&=
\left(P_{1/(2A)}-P_{1/A}\right)h_\alpha
=
\sum_{a=0}^{\infty}
m_A(a)\gamma_a(\alpha)L_a,
\\
\widetilde f_{A,\alpha}
&=
\left(P_{1/(2A)}-P_{1/A}\right)\widetilde h_\alpha
=
\sum_{a=0}^{\infty}
m_A(a)\widetilde\gamma_a(\alpha)L_a.
\end{aligned}
\label{eq:calculus-Abel-blocks}
\end{equation}

\begin{proposition}
\label{prop:exact-Abel-telescope}
Let \(\Lambda\geq2\) be dyadic. Then
\begin{equation}
\begin{aligned}
\left(I-P_{1/\Lambda}\right)h_\alpha
&=
\sum_{j=0}^{\infty}f_{2^j\Lambda,\alpha},
\qquad
\left(I-P_{1/\Lambda}\right)\widetilde h_\alpha
=
\sum_{j=0}^{\infty}\widetilde f_{2^j\Lambda,\alpha}.
\end{aligned}
\label{eq:calculus-infinite-Abel-telescope}
\end{equation}
The convergence holds in \(L^2(e^{-x}\dd x)\), and is locally uniform in \(x\), uniformly for
\(\alpha\in\mathfrak A\).
\end{proposition}

\begin{proof}
For \(J\geq0\), set
\[
u_J=\frac{1}{2^{J+1}\Lambda}.
\]
The finite telescope yields
\[
\sum_{j=0}^{J}f_{2^j\Lambda,\alpha}
=
P_{u_J}h_\alpha-P_{1/\Lambda}h_\alpha,
\qquad
\sum_{j=0}^{J}\widetilde f_{2^j\Lambda,\alpha}
=
P_{u_J}\widetilde h_\alpha-P_{1/\Lambda}\widetilde h_\alpha.
\]
Hence it remains to let \(u_J\downarrow0\).

Since
\[
\alpha\longmapsto
\bigl(h_\alpha,\widetilde h_\alpha\bigr)
\]
is continuous from the compact interval \(\mathfrak A\) into
\[
L^2(e^{-x}\dd x)\times L^2(e^{-x}\dd x),
\]
its image is compact.  The strong continuity and \(L^2\)-contractivity of
\((P_u)_{u\geq0}\) therefore imply
\[
\sup_{\alpha\in\mathfrak A}
\max\left\{
\|P_uh_\alpha-h_\alpha\|_{L^2(e^{-x}\dd x)},
\,
\|P_u\widetilde h_\alpha-\widetilde h_\alpha\|_{L^2(e^{-x}\dd x)}
\right\}
\longrightarrow0
\]
as \(u\downarrow0\).  This proves
\eqref{eq:calculus-infinite-Abel-telescope}
in \(L^2(e^{-x}\dd x)\), uniformly in \(\alpha\).
The locally uniform convergence follows directly from
Lemma~\ref{lem:uniform-local-semigroup-continuity}.
\end{proof}

\subsection{Low-degree corrections}
\label{subsec:process-equivalent-normalization}

Before applying cutoffs in the intensity variable, we remove the
\(L_1\) component and subtract the value at the origin.  The resulting
functions vanish at the origin and give the same centered statistics as
the original Abel blocks.

Since \(m_A(0)=0\), the Abel blocks have no degree-zero component.
Define the degree-one corrections by
\[
\beta_{A,\alpha}
=
m_A(1)\gamma_1(\alpha),
\qquad
\widetilde\beta_{A,\alpha}
=
m_A(1)\widetilde\gamma_1(\alpha).
\]
Since \(m_A(1)\) is independent of \(\alpha\),
\eqref{eq:spectral-tilde-gamma-derivative} gives
\begin{equation}
\widetilde\beta_{A,\alpha}
=
\partial_\alpha\beta_{A,\alpha}.
\label{eq:calculus-beta-derivative}
\end{equation}
Then
\[
\begin{aligned}
f_{A,\alpha}^{\circ}
&=
f_{A,\alpha}-\beta_{A,\alpha}L_1
=
\sum_{a=2}^{\infty}
m_A(a)\gamma_a(\alpha)L_a,
\\
\widetilde f_{A,\alpha}^{\circ}
&=
\widetilde f_{A,\alpha}
-\widetilde\beta_{A,\alpha}L_1
=
\sum_{a=2}^{\infty}
m_A(a)\widetilde\gamma_a(\alpha)L_a.
\end{aligned}
\]
Consequently,
\begin{equation}
\widetilde f_{A,\alpha}^{\circ}
=
\partial_\alpha f_{A,\alpha}^{\circ}.
\label{eq:calculus-visible-Abel-derivative}
\end{equation}
Since \(L_a(0)=1\),
\[
f_{A,\alpha}^{\circ}(0)
=
\sum_{a=2}^{\infty}
m_A(a)\gamma_a(\alpha),
\qquad
\widetilde f_{A,\alpha}^{\circ}(0)
=
\sum_{a=2}^{\infty}
m_A(a)\widetilde\gamma_a(\alpha).
\]
Both series converge absolutely by
Proposition~\ref{prop:compact-coefficient-bounds}
and the exponential decay of \(m_A(a)\).
By \eqref{eq:spectral-tilde-gamma-derivative},
Proposition~\ref{prop:compact-coefficient-bounds}, and the exponential
decay of \(m_A(a)\), the differentiated series above converges uniformly
on \(\mathfrak A\).
Hence termwise differentiation is justified, and
\begin{equation}
\widetilde f_{A,\alpha}^{\circ}(0)
=
\partial_\alpha
\left[f_{A,\alpha}^{\circ}(0)\right].
\label{eq:calculus-origin-correction-derivative}
\end{equation}

Define the corrected blocks by
\begin{equation}
f_{A,\alpha}^{\ren}
=
f_{A,\alpha}^{\circ}
-
f_{A,\alpha}^{\circ}(0)L_0,
\qquad
\widetilde f_{A,\alpha}^{\ren}
=
\widetilde f_{A,\alpha}^{\circ}
-
\widetilde f_{A,\alpha}^{\circ}(0)L_0.
\label{eq:calculus-renormalized-Abel-block}
\end{equation}
By
\eqref{eq:calculus-visible-Abel-derivative} and
\eqref{eq:calculus-origin-correction-derivative},
\begin{equation}
\widetilde f_{A,\alpha}^{\ren}
=
\partial_\alpha f_{A,\alpha}^{\ren}.
\label{eq:calculus-renormalized-Abel-derivative}
\end{equation}
These blocks vanish at the origin.
By
\eqref{eq:spectral-invisible-mode-fluctuations},
\[
X_n\!\left(f_{A,\alpha}^{\ren}\right)
=
X_n(f_{A,\alpha}),
\qquad
X_n\!\left(\widetilde f_{A,\alpha}^{\ren}\right)
=
X_n(\widetilde f_{A,\alpha}).
\]
We make these corrections before localization in the intensity variable,
since \(X_n(\chi L_0)\) and \(X_n(\chi L_1)\) need not vanish after
multiplication by a cutoff \(\chi\).

\begin{proposition}
\label{prop:correction-coefficient-bounds}
Uniformly for
\(\alpha\in\mathfrak A\) and \(A\geq2\),
\begin{equation}
|\beta_{A,\alpha}|
+
|\widetilde\beta_{A,\alpha}|
\leq
C_{\mathfrak A}A^{-1}.
\label{eq:calculus-beta-bounds}
\end{equation}
Moreover,
\begin{align}
|f_{A,\alpha}^{\circ}(0)|
&\leq
C_{\mathfrak A}
A^{-\min\{\alpha,1\}},
\label{eq:calculus-d0-bound}
\\
|\widetilde f_{A,\alpha}^{\circ}(0)|
&\leq
C_{\mathfrak A}
A^{-\min\{\alpha,1\}}
\log(A+2).
\label{eq:calculus-d1-bound}
\end{align}
\end{proposition}

\begin{proof}
Since \(m_A(1)=O(A^{-1})\) and
\(\gamma_1,\widetilde\gamma_1\) are uniformly bounded on
\(\mathfrak A\), \eqref{eq:calculus-beta-bounds} follows immediately.

For \(c\in\{1/2,1\}\), set \(z=e^{-c/A}\).
Using \eqref{eq:spectral-power-coefficient-generating-function}
and its \(\alpha\)-derivative, together with
\[
1-e^{-c/A}\asymp A^{-1},
\qquad
\bigl|\log(1-e^{-c/A})\bigr|
\leq
C\log(A+2),
\]
we obtain, uniformly for \(\alpha\in\mathfrak A\),
\[
\left|
\sum_{a=0}^{\infty}
\gamma_a(\alpha)e^{-ca/A}
\right|
\leq
C_{\mathfrak A}A^{-\alpha},
\]
and
\[
\left|
\sum_{a=0}^{\infty}
\widetilde\gamma_a(\alpha)e^{-ca/A}
\right|
\leq
C_{\mathfrak A}
A^{-\alpha}\log(A+2).
\]
Since
\[
m_A(a)=e^{-a/(2A)}-e^{-a/A},
\]
subtracting the degree-one contribution and using
\eqref{eq:calculus-beta-bounds} yields
\eqref{eq:calculus-d0-bound}
and
\eqref{eq:calculus-d1-bound}.
\end{proof}

\subsection{Scaled Abel functions}
\label{subsec:scaled-Abel-functions}

Next we rescale the Abel blocks into a common form for the derivative and energy estimates that follow.
The change of variables \(y=Ax\) puts intensity values \(x\) of order
\(A^{-1}\) on a fixed scale.  The functions below express the Abel
blocks and their corrected forms in this variable.

Fix once and for all
\[
0<\varepsilon_0<\frac1{16}.
\]
For
\[
c\in\left\{\frac12,1\right\},
\]
let
\(\mathsf S_{A,c,\alpha}(y)\)
be the scaled semigroup function from
\eqref{eq:def-kummer-scaled-one-semigroup-profile}.  By
\eqref{eq:equi-kummer-scaled-one-semigroup-profile},
\[
\mathsf S_{A,c,\alpha}(y)
=
\Gamma(1+\alpha)
\delta_{A,c}^{\,\alpha}
{}_1F_1
\left(
-\alpha;1;-b_{A,c}y
\right),
\]
where
\[
\delta_{A,c}
=
A(1-e^{-c/A}),
\qquad
b_{A,c}
=
\frac{e^{-c/A}}{\delta_{A,c}},
\]
as in \eqref{eq:kummer-rho-b}, and \({}_1F_1\) denotes the confluent
hypergeometric function.
Define the unmodified ordinary Abel function by
\begin{equation}
F_{A,\alpha}(y)
=
\mathsf S_{A,1/2,\alpha}(y)
-
\mathsf S_{A,1,\alpha}(y).
\label{eq:calculus-unmodified-ordinary-profile}
\end{equation}
Put
\begin{equation}
q_{A,\alpha}
=
A^\alpha\beta_{A,\alpha}
\label{eq:calculus-qA}
\end{equation}
and define the scaled Abel function by
\begin{equation}
F_{A,\alpha}^{\circ}(y)
=
F_{A,\alpha}(y)
-
q_{A,\alpha}L_1(y/A).
\label{eq:calculus-visible-ordinary-profile}
\end{equation}
Thus \(F_{A,\alpha}\) is the unmodified scaled ordinary Abel block;
the corrections subtract \(q_{A,\alpha}L_1(y/A)\) and then the value
at the origin.
Finally, subtract the value at the origin and define
\begin{equation}
\Delta F_{A,\alpha}^{\circ}(y)
=
F_{A,\alpha}^{\circ}(y)
-
F_{A,\alpha}^{\circ}(0).
\label{eq:calculus-renormalized-ordinary-profile}
\end{equation}
Then
\begin{equation}
f_{A,\alpha}^{\ren}(y/A)
=
A^{-\alpha}
\Delta F_{A,\alpha}^{\circ}(y).
\label{eq:calculus-ordinary-profile-scaling}
\end{equation}

\begin{proposition}
\label{prop:ordinary-scaled-Abel-estimates}
Fix
\(0<\delta_0<1/4\)
and let
\(J\Subset\mathfrak A\cap(0,1+\delta_0]\).
For every fixed integer \(M_0\geq1\), uniformly for
\(A\geq 2\),
\(\alpha\in J\),
\(0\leq y\leq8\varepsilon_0A\),
one has
\begin{equation}
\begin{aligned}
|\Delta F_{A,\alpha}^{\circ}(y)|
&\leq
C_J
\left[
1+(1+y)^{\alpha-1}
+A^{-1}(1+y)^\alpha
\right],
\end{aligned}
\label{eq:calculus-ordinary-profile-zero}
\end{equation}
and, for
\(1\leq j\leq M_0\),
\begin{equation}
\left|\partial_y^j
\Delta F_{A,\alpha}^{\circ}(y)\right|
\leq
C_{J,M_0}
\left[
(1+y)^{\alpha-1-j}
+A^{-1}(1+y)^{\alpha-j}
\right].
\label{eq:calculus-ordinary-profile-jets}
\end{equation}
At the critical value \(\alpha=1\),
\begin{equation}
F_{A,1}^{\circ}
\equiv0,
\qquad
\Delta F_{A,1}^{\circ}
\equiv0.
\label{eq:calculus-ordinary-critical-zero}
\end{equation}
\end{proposition}

\begin{proof}
For \(0\leq y\leq1\), by
\eqref{eq:kummer-finite-A-parameters},
the parameters \(\delta_{A,c}\) and \(b_{A,c}\),
\(c\in\{1/2,1\}\), remain in fixed compact subsets of
\((0,\infty)\).  Hence joint smoothness of
\[
{}_1F_1(-\alpha;1;-b_{A,c}y)
\]
gives uniform bounds for every fixed \(y\)-derivative of
\(F_{A,\alpha}\) on this range.
Moreover,
\eqref{eq:calculus-beta-bounds} gives
\[
|q_{A,\alpha}|
\leq
C_JA^{\alpha-1}.
\]
Since \(\alpha\leq1+\delta_0<2\), the affine correction
\[
\frac{q_{A,\alpha}}A y
\]
and its first \(y\)-derivative are uniformly bounded for \(0\leq y\leq1\),
while all higher derivatives vanish. Hence
\eqref{eq:calculus-ordinary-profile-zero} and
\eqref{eq:calculus-ordinary-profile-jets} hold on this range.

Assume
\[
1\leq y\leq8\varepsilon_0A.
\]
By
Corollary~\ref{cor:finite-A-one-profile-expansion},
\[
\mathsf S_{A,c,\alpha}(y)
=
e^{-\alpha c/A}y^\alpha
+\alpha^2\delta_{A,c}
e^{-(\alpha-1)c/A}y^{\alpha-1}
+O_J(y^{\alpha-2}),
\]
uniformly with every fixed \(y\)-derivative.

The leading \(y^\alpha\)-coefficient in the Abel difference is
\[
e^{-\alpha/(2A)}
-
e^{-\alpha/A}
=
O_J(A^{-1}),
\]
whereas the coefficient of \(y^{\alpha-1}\) is uniformly bounded.
The same statements hold after every fixed \(y\)-derivative.

Moreover,
\[
|q_{A,\alpha}|
\leq
C_JA^{\alpha-1}.
\]
Since
\[
L_1(y/A)-L_1(0)
=
-\frac yA,
\]
the affine correction after subtraction at the origin is bounded by
\[
C_JA^{\alpha-2}y
\leq
C_J
\left[
(1+y)^{\alpha-1}
+
A^{-1}(1+y)^\alpha
\right].
\]
Its first \(y\)-derivative is \(O_J(A^{\alpha-2})\), which is absorbed by
\((1+y)^{\alpha-2}\), while all higher derivatives vanish. Together with the uniform bound on \(F_{A,\alpha}(0)\), this proves
\eqref{eq:calculus-ordinary-profile-zero} and
\eqref{eq:calculus-ordinary-profile-jets}.

At \(\alpha=1\),
\[
h_1(x)=x=L_0(x)-L_1(x),
\]
so all Laguerre coefficients of degree at least two vanish. Hence,
\eqref{eq:calculus-ordinary-critical-zero} follows.
\end{proof}

To account for the exponent-dependent scaling, for fixed \(A\geq2\)
define the differential operator
\[
\mathscr D_A
:=
\partial_\alpha-\log A.
\]
Define the corresponding scaled logarithmic function by
\begin{equation}
\mathcal P_{A,\alpha}(y)
=
\mathscr D_A
\Delta F_{A,\alpha}^{\circ}(y).
\label{eq:calculus-log-profile-definition}
\end{equation}
By \eqref{eq:calculus-renormalized-Abel-derivative} and
\eqref{eq:calculus-ordinary-profile-scaling},
\begin{equation}
\widetilde f_{A,\alpha}^{\ren}(y/A)
=
A^{-\alpha}
\mathcal P_{A,\alpha}(y).
\label{eq:calculus-log-profile-scaling}
\end{equation}

\begin{proposition}
\label{prop:logarithmic-scaled-Abel-estimates}
Under the hypotheses of
Proposition~\ref{prop:ordinary-scaled-Abel-estimates}, 
uniformly for
\(A\geq 2\),
\(\alpha\in J\),
\(0\leq y\leq8\varepsilon_0A\),
one has
\begin{equation}
\begin{aligned}
|\mathcal P_{A,\alpha}(y)|
&\leq
C_J
\bigl[\log(A+2)+\log(2+y)\bigr]
\left[
1+(1+y)^{\alpha-1}
+A^{-1}(1+y)^\alpha
\right],
\end{aligned}
\label{eq:calculus-log-profile-zero}
\end{equation}
and, for
\(1\leq j\leq M_0\),
\begin{equation}
\begin{aligned}
\left|\partial_y^j
\mathcal P_{A,\alpha}(y)\right|
&\leq
C_{J,M_0}
\bigl[\log(A+2)+\log(2+y)\bigr]
\left[
(1+y)^{\alpha-1-j}
+
A^{-1}(1+y)^{\alpha-j}
\right].
\end{aligned}
\label{eq:calculus-log-profile-jets}
\end{equation}
\end{proposition}

\begin{proof}
For \(0\leq y\leq1\), joint smoothness of the Kummer representation
\eqref{eq:equi-kummer-scaled-one-semigroup-profile},
together with \eqref{eq:calculus-unmodified-ordinary-profile}
and  \eqref{eq:kummer-finite-A-parameters},
gives, uniformly for \(\alpha\in J\),
\[
\max_{0\leq j\leq M_0}
\left|\mathscr D_A\partial_y^jF_{A,\alpha}(y)\right|
\leq C_{J,M_0}\log(A+2).
\]
Furthermore,
\eqref{eq:calculus-qA},
\eqref{eq:calculus-beta-derivative}, and
\eqref{eq:calculus-beta-bounds} give
\[
|\mathscr D_Aq_{A,\alpha}|
\leq
C_JA^{\alpha-1}.
\]
Since \(\alpha\leq1+\delta_0<2\), the corresponding affine correction and
its first derivative are uniformly bounded on \(0\leq y\leq1\). Therefore,
\[
|\mathcal P_{A,\alpha}(y)|
\leq
C_J\bigl[\log(A+2)+\log(2+y)\bigr],
\]
and, for \(1\leq j\leq M_0\),
\[
|\partial_y^j\mathcal P_{A,\alpha}(y)|
\leq
C_{J,M_0}\bigl[\log(A+2)+\log(2+y)\bigr].
\]

For
\(1\leq y\leq8\varepsilon_0A\),
by
Lemma~\ref{lem:uniform-differentiated-Kummer} and
Corollary~\ref{cor:finite-A-one-profile-expansion},
for \(0\leq j\leq M_0\),
\[
\begin{aligned}
\left|
\mathscr D_A\partial_y^jF_{A,\alpha}(y)
\right|
&\leq
C_{J,M_0}
\bigl[\log(A+2)+\log(2+y)\bigr]
\left[
y^{\alpha-1-j}
+
A^{-1}y^{\alpha-j}
\right].
\end{aligned}
\]
The value at the origin satisfies
\[
\left|
\mathscr D_AF_{A,\alpha}(0)
\right|
\leq
C_J\log(A+2).
\]

For the degree-one correction,
\eqref{eq:calculus-qA},
\eqref{eq:calculus-beta-derivative}, and
\eqref{eq:calculus-beta-bounds} imply
\[
\begin{aligned}
\mathscr D_Aq_{A,\alpha}
=
\left(\partial_\alpha-\log A\right)
\left(A^\alpha\beta_{A,\alpha}\right)
=
A^\alpha\widetilde\beta_{A,\alpha},
\end{aligned}
\]
and hence
\[
|\mathscr D_Aq_{A,\alpha}|
\leq
C_JA^{\alpha-1}.
\]
Since
\[
L_1(y/A)-L_1(0)
=
-\frac yA,
\]
the correction contributes
\[
\left|
\frac yA\mathscr D_Aq_{A,\alpha}
\right|
\leq
C_JA^{\alpha-2}y
\leq
C_J(1+y)^{\alpha-1},
\]
where
\[
\alpha\leq1+\delta_0<2,
\qquad
y\leq8\varepsilon_0A.
\]
Its first \(y\)-derivative is bounded by
\[
C_JA^{\alpha-2}
\leq
C_J(1+y)^{\alpha-2},
\]
and all higher derivatives vanish.

Combining these estimates proves
\eqref{eq:calculus-log-profile-zero} and
\eqref{eq:calculus-log-profile-jets}.

The uniformity across \(\alpha=1\) follows from
Lemma~\ref{lem:uniform-differentiated-Kummer}.
\end{proof}

\begin{remark}
At the critical value \(\alpha=1\), the corrected ordinary scaled Abel
function \(\Delta F_{A,1}^{\circ}\) vanishes identically by
\eqref{eq:calculus-ordinary-critical-zero}, whereas the corresponding
logarithmic scaled function \(\mathcal P_{A,1}\) has a nontrivial limit
as \(A\to\infty\).
More precisely, Appendix~\ref{app:kummer-critical-profile} shows that
\(\mathcal P_{A,1}\) converges locally uniformly, together with every fixed
\(y\)-derivative, to a function \(Q\) satisfying
\[
Q(y)
=
-\frac12\log y+O(1)
\qquad
(y\to\infty);
\]
see Propositions~\ref{prop:critical-Riemann-sum}
and~\ref{prop:critical-logarithmic-function}.
\end{remark}

\subsection{Localization in the intensity variable}
\label{subsec:spatial-localization}

Now we decompose the corrected Abel blocks according to small, intermediate,
and large values of the intensity variable \(x\).  Near zero, derivatives
of the power and logarithmic-power functions may be singular.  Away from
zero these functions are smooth, and the weight \(e^{-x}\) in the
Laguerre-energy integrals supplies decay for large \(x\).  We record the
localization estimates for these three regions below.

Fix a sufficiently large dyadic scale
\[
R_0>4.
\]
Let
\[
\lambda_{\mathrm{low}},
\lambda_{\mathrm{mid}},
\lambda_{\mathrm{up}}
\in C^\infty([0,\infty))
\]
be a smooth partition of unity satisfying
\begin{equation}
0\leq
\lambda_{\mathrm{low}},
\lambda_{\mathrm{mid}},
\lambda_{\mathrm{up}}
\leq1,
\qquad
\lambda_{\mathrm{low}}
+
\lambda_{\mathrm{mid}}
+
\lambda_{\mathrm{up}}
=
1,
\label{eq:calculus-spatial-partition}
\end{equation}
with
\[
\supp\lambda_{\mathrm{low}}
\subseteq[0,4\varepsilon_0],
\qquad
\supp\lambda_{\mathrm{mid}}
\subseteq[\varepsilon_0/2,2R_0],
\qquad
\supp\lambda_{\mathrm{up}}
\subseteq[R_0,\infty).
\]

\boldparagraph{Dyadic resolution of the lower region.}
Here \(A\geq2\) is the Abel scale from
\eqref{eq:calculus-Abel-blocks}.  In the lower region we use the rescaled
variable
\[
y=Ax,
\]
and decompose this variable into dyadic scales \(R=2^k\),
\(k\in\mathbb N_0\).
Let
\[
\chi_0\in C_c^\infty([0,\infty))
\]
be supported in \([0,2]\).  
For these dyadic scales, choose
\[
\omega_R\in C_c^\infty((0,\infty))
\]
so that, for every \(y\geq0\),
\begin{equation}
\chi_0(y)
+
\sum_{k=0}^{\infty}\omega_{2^k}(y)
=
1,
\label{eq:calculus-exact-scaled-dyadic-partition}
\end{equation}
and
\begin{equation}
\supp\omega_R
\subseteq[R/2,4R],
\qquad
\sup_{y\geq0}
\left|\omega_R^{(j)}(y)\right|
\leq
C_jR^{-j}.
\label{eq:calculus-scaled-dyadic-cutoffs}
\end{equation}

Define the lowest-scale piece by
\begin{equation}
r_{A,0,\alpha}(x)
=
\chi_0(Ax)
\lambda_{\mathrm{low}}(x)
f_{A,\alpha}^{\ren}(x),
\label{eq:calculus-lower-core}
\end{equation}
and, for dyadic
\(
1\leq R\leq8\varepsilon_0A,
\)
define the lower shell by
\begin{equation}
r_{A,R,\alpha}(x)
=
\omega_R(Ax)
\lambda_{\mathrm{low}}(x)
f_{A,\alpha}^{\ren}(x).
\label{eq:calculus-lower-shell}
\end{equation}

Since
\(
\supp\omega_R(A\cdot)
\subseteq[R/(2A),4R/A]
\)
and
\(
\supp\lambda_{\mathrm{low}}
\subseteq[0,4\varepsilon_0],
\)
the product
\(
\omega_R(A\cdot)\lambda_{\mathrm{low}}
\)
vanishes whenever \(R>8\varepsilon_0A\).
It follows from
\eqref{eq:calculus-exact-scaled-dyadic-partition},
\eqref{eq:calculus-lower-core}, and
\eqref{eq:calculus-lower-shell} that
\begin{equation}
\begin{aligned}
\lambda_{\mathrm{low}}(x)
f_{A,\alpha}^{\ren}(x)
&=
r_{A,0,\alpha}(x)
+
\sum_{\substack{
R=2^k,\ k\in\mathbb N_0\\
R\leq8\varepsilon_0A}}
r_{A,R,\alpha}(x).
\end{aligned}
\label{eq:calculus-exact-lower-cell-decomposition}
\end{equation}

\boldparagraph{Middle region.}
The middle region remains on a fixed compact set. Define
\begin{equation}
f_{A,\alpha}^{\mathrm{mid}}(x)
=
\lambda_{\mathrm{mid}}(x)
f_{A,\alpha}^{\ren}(x).
\label{eq:calculus-middle-piece}
\end{equation}

\boldparagraph{Dyadic resolution of the upper region.}
Here \(R\) denotes a dyadic scale in the variable \(x\).
For each
\[
R=R_0 2^k,
\qquad
k\in\mathbb N_0,
\]
choose
\[
\varpi_R\in C_c^\infty((0,\infty))
\]
so that
\[
0\leq\varpi_R\leq1,
\qquad
\supp\varpi_R
\subseteq[R/2,2R],
\]
and
\begin{equation}
\sum_{k=0}^{\infty}
\varpi_{R_0 2^k}(x)
=
1,
\qquad
x\geq R_0.
\label{eq:calculus-exact-physical-dyadic-partition}
\end{equation}
Moreover, for every fixed integer \(j\geq0\),
\begin{equation}
\sup_{x\geq0}
\left|\partial_x^j\varpi_R(x)\right|
\leq
C_jR^{-j}.
\label{eq:calculus-physical-dyadic-cutoffs}
\end{equation}

Define the upper shell by
\begin{equation}
f_{A,R,\alpha}^{\mathrm{up}}(x)
=
\varpi_R(x)
\lambda_{\mathrm{up}}(x)
f_{A,\alpha}^{\ren}(x).
\label{eq:calculus-upper-shell}
\end{equation}
Then
\begin{equation}
\lambda_{\mathrm{up}}(x)
f_{A,\alpha}^{\ren}(x)
=
\sum_{k=0}^{\infty}
f_{A,R_0 2^k,\alpha}^{\mathrm{up}}(x).
\label{eq:calculus-exact-upper-shell-decomposition}
\end{equation}

Combining
\eqref{eq:calculus-spatial-partition},
\eqref{eq:calculus-exact-lower-cell-decomposition},
\eqref{eq:calculus-middle-piece}, and
\eqref{eq:calculus-exact-upper-shell-decomposition}, one obtains
\[
f_{A,\alpha}^{\ren}(x)
=
r_{A,0,\alpha}(x)
+
\sum_{\substack{
R=2^k,\ k\in\mathbb N_0\\
R\leq8\varepsilon_0A}}
r_{A,R,\alpha}(x)
+
f_{A,\alpha}^{\mathrm{mid}}(x)
+
\sum_{k=0}^{\infty}
f_{A,R_0 2^k,\alpha}^{\mathrm{up}}(x).
\]

The corresponding logarithmic pieces
\[
\widetilde r_{A,0,\alpha},
\qquad
\widetilde r_{A,R,\alpha},
\qquad
\widetilde f_{A,\alpha}^{\mathrm{mid}},
\qquad
\widetilde f_{A,R,\alpha}^{\mathrm{up}}
\]
are defined by the same formulas with
\(f_{A,\alpha}^{\ren}\) replaced throughout by
\(\widetilde f_{A,\alpha}^{\ren}\).
Thus,
\[
\widetilde f_{A,\alpha}^{\ren}(x)
=
\widetilde r_{A,0,\alpha}(x)
+
\sum_{\substack{
R=2^k,\ k\in\mathbb N_0\\
R\leq8\varepsilon_0A}}
\widetilde r_{A,R,\alpha}(x)
+
\widetilde f_{A,\alpha}^{\mathrm{mid}}(x)
+
\sum_{k=0}^{\infty}
\widetilde f_{A,R_0 2^k,\alpha}^{\mathrm{up}}(x).
\]

\begin{lemma}
\label{lem:fixed-cutoff-Laguerre-localization}
Let
\(0<a<b<\infty\),
\(\chi\in C_c^\infty((a,b))\).
For every fixed integer \(r\geq1\), there exist integers
\(M_r<\infty\) and constants \(C_{\chi,r}<\infty\) such that
\begin{equation}
\mathfrak E_r(\chi f)
\leq
C_{\chi,r}
\sum_{j=0}^{M_r}
\sup_{a/2\leq x\leq2b}
|f^{(j)}(x)|^2
\label{eq:calculus-fixed-cutoff-localization}
\end{equation}
for every smooth \(f\).
Moreover, for every integer \(K\geq0\),
\[
\begin{aligned}
\max_{0\leq j\leq K+1}
\|(\chi f)^{(j)}\|_{L^\infty([0,\infty))}
+\max_{|\beta|\leq K}
\|\partial^\beta\Phi_{\chi f}\|_{L^\infty(\mathbb R^2)}
\leq
C_{\chi,K}
\max_{0\leq j\leq K+1}
\|f^{(j)}\|_{L^\infty([a/2,2b])},
\end{aligned}
\]
where \(\Phi_{\chi f}(z)=(\chi f)'(|z|^2)z\),
\(z\in\mathbb C\simeq\mathbb R^2\), and
\(\partial^\beta\) denotes real partial differentiation.
\end{lemma}

\begin{proof}
For even \(r=2m\) with \(m\geq1\),
\[
\mathfrak E_{2m}(g)
=
\|\LagNum^mg\|_{L^2(e^{-x}\dd x)}^2.
\]
Repeatedly use the commutator identity
\begin{equation}
\begin{aligned}
\LagNum(\chi f)
&=
\chi\LagNum f
-
2x\chi'f'
-
\left[x\chi''+(1-x)\chi'\right]f.
\end{aligned}
\label{eq:calculus-Laguerre-cutoff-commutator}
\end{equation}
On the fixed support, all polynomial coefficients are bounded and \(x\) is
bounded away from zero. Iteration expresses
\[
\LagNum^m(\chi f)
\]
as a finite linear combination of derivatives of \(f\), with bounded
coefficients. Also,
\[
\mathfrak E_0(g)
\leq
\|g\|_{L^2(e^{-x}\dd x)}^2.
\]
For odd energies, including \(r=1\), Cauchy--Schwarz in coefficient space gives
\[
\mathfrak E_{2m+1}(g)^2
\leq
\mathfrak E_{2m}(g)
\mathfrak E_{2m+2}(g),
\qquad
m\geq0.
\]
The derivative bounds follow from Leibniz's rule and the chain rule.
\end{proof}

The following lemma bounds the Laguerre energy on each upper shell using local derivatives, with exponential decay in the shell size.
\begin{lemma}
\label{lem:dyadic-upper-Laguerre-localization}
For
\(
R=R_0 2^k,
\)
\(
k\in\mathbb N_0,
\)
put
\[
\chi_R
=
\varpi_R\lambda_{\mathrm{up}}.
\]
For every fixed integer \(s\geq1\), there exist constants \(K_s,B_s,C_s<\infty\),
such that, for every smooth \(f\),
\begin{equation}
\begin{aligned}
\mathfrak E_s(\chi_Rf)
&\leq
C_s(1+R)^{B_s}e^{-R/2}
\sum_{j=0}^{K_s}
\sup_{R/2\leq x\leq2R}
\left|f^{(j)}(x)\right|^2.
\end{aligned}
\label{eq:calculus-dyadic-upper-localization}
\end{equation}
Moreover, every fixed real derivative of
\[
z\longmapsto
(\chi_Rf)'(|z|^2)z
\]
is bounded by
\begin{equation}
C(1+R)^B
\sum_{j=0}^{M}
\sup_{R/2\leq x\leq2R}
\left|f^{(j)}(x)\right|,
\label{eq:calculus-dyadic-upper-field-jets}
\end{equation}
where \(C,B,M\) depend only on the fixed derivative order and on the fixed
cutoff data.
\end{lemma}

\begin{proof}
By \eqref{eq:calculus-physical-dyadic-cutoffs},
\[
\supp\chi_R\subseteq[R/2,2R],
\qquad
\|\partial_x^j\chi_R\|_\infty
\leq
C_jR^{-j},
\qquad j\geq0.
\]
For every fixed \(m\geq1\), iteration of
\[
\LagNum=-x\partial_x^2-(1-x)\partial_x
\]
yields
\[
\LagNum^m
=
\sum_{j=1}^{2m}
P_{m,j}(x)\partial_x^j,
\]
where the \(P_{m,j}\) are polynomials depending only on \(m\).  Hence,
on \(R/2\leq x\leq2R\), Leibniz's rule yields
\[
\left|
\LagNum^m(\chi_Rf)(x)
\right|
\leq
C_m(1+R)^{B_m}
\sum_{j=0}^{2m}
\sup_{R/2\leq u\leq2R}
|f^{(j)}(u)|
\]
for some \(B_m<\infty\).

Since \(e^{-x}\leq e^{-R/2}\) on \(\supp\chi_R\), whose length is
\(O(R)\), it follows that
\[
\mathfrak E_{2m}(\chi_Rf)
=
\|\LagNum^m(\chi_Rf)\|_{L^2(e^{-x}\dd x)}^2
\leq
C_m(1+R)^{B_m}e^{-R/2}
\sum_{j=0}^{2m}
\sup_{R/2\leq x\leq2R}
|f^{(j)}(x)|^2,
\]
after enlarging \(B_m\).  The same type of estimate for
\(\mathfrak E_0\) is immediate.  The interpolation inequality used in
the proof of Lemma~\ref{lem:fixed-cutoff-Laguerre-localization} then gives
the odd-order bounds.  This proves
\eqref{eq:calculus-dyadic-upper-localization}.

Finally, every fixed real derivative of
\[
z\longmapsto
(\chi_Rf)'(|z|^2)z
\]
is a finite sum of products of derivatives of \(\chi_Rf\) evaluated at
\(|z|^2\) and polynomial factors in the real and imaginary parts of
\(z\).  On the relevant support, \(|z|^2\asymp R\); hence the polynomial
factors contribute only a fixed power of \(1+R\).  Combining this with
Leibniz's rule and the derivative bounds for \(\chi_R\) yields
\eqref{eq:calculus-dyadic-upper-field-jets}.
\end{proof}

\begin{lemma}
\label{lem:compact-upper-semigroup-jets}
Let
\(J\Subset(0,\infty)\)
and let \(M\geq0\) be fixed.
Uniformly for
\(
A\geq2,
\)
\(
\alpha\in J,
\)
the following estimates hold.

\begin{enumerate}[label=\textup{(\roman*)},leftmargin=2.4em]

\item
On every fixed compact interval
\([a,b]\Subset(0,\infty)\),
\[
\max_{0\leq j\leq M}
\sup_{a\leq x\leq b}
\left(
\left|\partial_x^j f_{A,\alpha}(x)\right|
+
\left|\partial_x^j\widetilde f_{A,\alpha}(x)\right|
\right)
\leq
C_{J,a,b,M}A^{-1}.
\]

\item
For every \(\varepsilon>0\), there exist
\(C_{J,M,\varepsilon},M'<\infty\) such that, for
\(x\geq R_0/4\) and \(0\leq j\leq M\),
\[
\left|\partial_x^j f_{A,\alpha}(x)\right|
+
\left|\partial_x^j\widetilde f_{A,\alpha}(x)\right|
\leq
C_{J,M,\varepsilon}
A^{-1}(1+x)^{M'}e^{\varepsilon x}.
\]

\end{enumerate}
\end{lemma}

\begin{proof}
Write
\begin{equation}
\begin{aligned}
f_{A,\alpha}
=
-\int_{1/(2A)}^{1/A}
\partial_uP_uh_\alpha\,\dd u,
\qquad
\widetilde f_{A,\alpha}
=
-\int_{1/(2A)}^{1/A}
\partial_uP_u\widetilde h_\alpha\,\dd u.
\end{aligned}
\label{eq:calculus-Abel-integral-representation}
\end{equation}
For every fixed integer \(j\geq0\), the derivatives
\[
\partial_x^j\partial_uP_uh_\alpha(x),
\qquad
\partial_x^j\partial_uP_u\widetilde h_\alpha(x)
\]
are uniformly bounded for
\[
x\in[a,b]\Subset(0,\infty),
\qquad
0\leq u\leq1,
\qquad
\alpha\in J.
\]
The integration interval has length \(1/(2A)\), which proves part
\textup{(i)}.

For \(x\geq R_0/4\),
Proposition~\ref{prop:kummer-semigroup-power} and
Lemma~\ref{lem:uniform-differentiated-Kummer} imply, uniformly for
\[
0<u\leq1,
\qquad
\alpha\in J,
\qquad
0\leq j\leq M,
\]
that
\[
\left|
\partial_x^j\partial_uP_uh_\alpha(x)
\right|
+
\left|
\partial_x^j\partial_uP_u\widetilde h_\alpha(x)
\right|
\leq
C_{J,M,\varepsilon}
(1+x)^{M'}e^{\varepsilon x}.
\]
Combining this estimate with
\eqref{eq:calculus-Abel-integral-representation}
proves part \textup{(ii)}.
\end{proof}

\subsection{Regional Laguerre-energy estimates}
\label{subsec:regional-Laguerre-energy-estimates}

We now convert the scaled Abel estimates and the localization estimates
from Subsections~\ref{subsec:scaled-Abel-functions}
and~\ref{subsec:spatial-localization} into Laguerre-energy bounds for the pieces
in the lower, middle, and upper intensity regions.

For fixed \(A\) and \(\alpha\), let
\(\Psi\in C^\infty([0,\infty);\mathbb R)\), and define
\[
g_{A,R}^{\Psi,\alpha}(x)
=
A^{-\alpha}\lambda_{\mathrm{low}}(x)
\omega_R(Ax)\Psi(Ax).
\]
Thus \(g_{A,R}^{\Psi,\alpha}\) isolates the part of
\(A^{-\alpha}\Psi(Ax)\) where \(x\) is small and \(Ax\)
is of order \(R\).
Below we take \(\Psi=\Delta F_{A,\alpha}^{\circ}\) for ordinary
blocks and \(\Psi=\mathcal P_{A,\alpha}\) for logarithmic blocks.
For a fixed integer \(s\geq1\), put
\[
M_s=2\left\lceil\frac s2\right\rceil
\]
and define the scaled derivative bound on a shell
\begin{equation}
\mathcal B_{A,R,s}(\Psi)
=
\max_{0\leq k\leq M_s}
R^k
\sup_{R/2\leq y\leq4R}
\left|
\partial_y^k
\left[
\lambda_{\mathrm{low}}(y/A)\omega_R(y)\Psi(y)
\right]
\right|.
\label{eq:calculus-lower-shell-jet-envelope}
\end{equation}
With \(y=Ax\),
\begin{equation}
\LagNum
=
-A
\left[
y\partial_y^2
+
\left(1-\frac yA\right)\partial_y
\right].
\label{eq:calculus-lower-rescaled-number-operator}
\end{equation}
On a shell put \(y=Rw\). Since
\[
1\leq R\leq8\varepsilon_0A,
\]
the variable \(w\) ranges in a fixed compact subset of \((0,\infty)\), and
every application of \(\LagNum\) contributes at most a factor \(CA/R\).
Therefore, for every fixed \(m\geq1\),
\begin{equation}
\mathfrak E_{2m}(g_{A,R}^{\Psi,\alpha})
=
\|\LagNum^mg_{A,R}^{\Psi,\alpha}\|_{L^2(e^{-x}\dd x)}^2
\leq
C_mA^{-2\alpha}
\left(\frac AR\right)^{2m-1}
\mathcal B_{A,R,2m}(\Psi)^2.
\label{eq:calculus-lower-even-energy-scaling}
\end{equation}
Moreover,
\[
\mathfrak E_0(g_{A,R}^{\Psi,\alpha})
\leq
\|g_{A,R}^{\Psi,\alpha}\|_{L^2(e^{-x}\dd x)}^2
\leq
CA^{-2\alpha}
\frac RA
\mathcal B_{A,R,1}(\Psi)^2.
\]
For odd energies, use the interpolation inequality from the proof of
Lemma~\ref{lem:fixed-cutoff-Laguerre-localization}.  Since
\(M_{2m+1}=2m+2\), the same derivative bound controls both neighboring
even energies. Consequently,
\begin{equation}
\mathfrak E_s(g_{A,R}^{\Psi,\alpha})
\leq
C_sA^{-2\alpha}
\left(\frac AR\right)^{s-1}
\mathcal B_{A,R,s}(\Psi)^2
\label{eq:calculus-lower-all-order-energy-scaling}
\end{equation}
for every fixed integer \(s\geq1\).

On the fixed \(y\)-support of the lowest-scale piece, put
\[
g_{A,0}^{\Psi,\alpha}(x)
=
A^{-\alpha}
\chi_0(Ax)
\lambda_{\mathrm{low}}(x)
\Psi(Ax)
\]
and define
\[
\mathcal B_{A,0,s}(\Psi)
=
\max_{0\leq k\leq M_s}
\sup_{0\leq y\leq2}
\left|
\partial_y^k
\left[
\chi_0(y)\lambda_{\mathrm{low}}(y/A)\Psi(y)
\right]
\right|.
\]
For \(m\geq1\), the same rescaling gives
\[
\mathfrak E_{2m}(g_{A,0}^{\Psi,\alpha})
\leq
C_mA^{2m-1-2\alpha}
\mathcal B_{A,0,2m}(\Psi)^2,
\]
while
\[
\mathfrak E_0(g_{A,0}^{\Psi,\alpha})
\leq
\|g_{A,0}^{\Psi,\alpha}\|_{L^2(e^{-x}\dd x)}^2
\leq
CA^{-1-2\alpha}
\mathcal B_{A,0,1}(\Psi)^2.
\]
Interpolation of neighboring even energies therefore yields
\begin{equation}
\mathfrak E_s(g_{A,0}^{\Psi,\alpha})
\leq
C_sA^{s-1-2\alpha}
\mathcal B_{A,0,s}(\Psi)^2.
\label{eq:calculus-core-all-order-energy-scaling}
\end{equation}

Next, we obtain \(A^{\alpha-1}\) bounds for the scaled derivatives above \(\alpha=1\), with an extra logarithmic factor in the logarithmic case.

\begin{lemma}
\label{lem:supercritical-scaled-profile-jets}
Fix
\(0<\delta_0<1/4\)
and let
\[
J_{>}
\Subset
\mathfrak A\cap[1+\delta_0,\infty).
\]
For every fixed integer
\(M\geq0\), uniformly for
\(\alpha\in J_{>}\),
\(A\geq2\),
\(1\leq R\leq8\varepsilon_0A\),
one has
\begin{equation}
\max_{0\leq k\leq M}
R^k
\sup_{y\in\supp[\lambda_{\mathrm{low}}(\,\cdot\,/A)\omega_R]}
\left|
\partial_y^k
\Delta F_{A,\alpha}^{\circ}(y)
\right|
\leq
C_{J_{>},M}
A^{\alpha-1},
\label{eq:calculus-supercritical-ordinary-profile-jets}
\end{equation}
and
\begin{equation}
\max_{0\leq k\leq M}
R^k
\sup_{y\in\supp[\lambda_{\mathrm{low}}(\,\cdot\,/A)\omega_R]}
\left|
\partial_y^k
\mathcal P_{A,\alpha}(y)
\right|
\leq
C_{J_{>},M}
\bigl[\log(A+2)+\log(R+2)\bigr]
A^{\alpha-1}.
\label{eq:calculus-supercritical-log-profile-jets}
\end{equation}
On the microscopic range \(0\leq y\leq2\),
\begin{equation}
\max_{0\leq k\leq M}
\sup_{0\leq y\leq2}
\left|
\partial_y^k
\Delta F_{A,\alpha}^{\circ}(y)
\right|
\leq
C_{J_{>},M}
A^{\alpha-1},
\label{eq:calculus-supercritical-ordinary-core-jets}
\end{equation}
and
\begin{equation}
\max_{0\leq k\leq M}
\sup_{0\leq y\leq2}
\left|
\partial_y^k
\mathcal P_{A,\alpha}(y)
\right|
\leq
C_{J_{>},M}
\log(A+2)
A^{\alpha-1}.
\label{eq:calculus-supercritical-log-core-jets}
\end{equation}
\end{lemma}

\begin{proof}
Put
\[
\mathsf a_{A,\alpha}
=
e^{-\alpha/(2A)}-e^{-\alpha/A}
\]
and
\[
\mathsf b_{A,\alpha}
=
\alpha^2
\left[
\delta_{A,1/2}e^{-(\alpha-1)/(2A)}
-
\delta_{A,1}e^{-(\alpha-1)/A}
\right].
\]
By
\eqref{eq:kummer-finite-A-parameters}, uniformly for
\(\alpha\in J_{>}\),
\begin{equation}
|\mathsf a_{A,\alpha}|
+
|\partial_\alpha\mathsf a_{A,\alpha}|
\leq
C_{J_{>}}A^{-1},
\qquad
|\mathsf b_{A,\alpha}|
+
|\partial_\alpha\mathsf b_{A,\alpha}|
\leq
C_{J_{>}}.
\label{eq:calculus-supercritical-leading-coefficients}
\end{equation}
By
Corollary~\ref{cor:finite-A-one-profile-expansion}, for every fixed
\(0\leq k\leq M\) and \(y\geq1\),
\begin{equation}
\begin{aligned}
\partial_y^kF_{A,\alpha}(y)
&=
\partial_y^k
\left[
\mathsf a_{A,\alpha}y^\alpha
+
\mathsf b_{A,\alpha}y^{\alpha-1}
\right]
+
O_{J_{>},M}\left(y^{\alpha-2-k}\right).
\end{aligned}
\label{eq:calculus-supercritical-F-expansion}
\end{equation}
Hence
\begin{equation}
\left|
\partial_y^kF_{A,\alpha}(y)
\right|
\leq
C_{J_{>},M}
\left[
A^{-1}y^{\alpha-k}
+
y^{\alpha-1-k}
+
y^{\alpha-2-k}
\right].
\label{eq:calculus-supercritical-F-jets}
\end{equation}
Applying \(\mathscr D_A\) to the preceding expansion and using the
\(\alpha\)-differentiated remainder in
Corollary~\ref{cor:finite-A-one-profile-expansion} gives
\begin{equation}
\begin{aligned}
\left|
\mathscr D_A
\partial_y^kF_{A,\alpha}(y)
\right|
&\leq
C_{J_{>},M}
\bigl[\log(A+2)+\log(2+y)\bigr]
\left[
A^{-1}y^{\alpha-k}
+
y^{\alpha-1-k}
+
y^{\alpha-2-k}
\right].
\end{aligned}
\label{eq:calculus-supercritical-DA-F-jets}
\end{equation}

On the support of
\(\lambda_{\mathrm{low}}(y/A)\omega_R(y)\), one has
\[
y\asymp R,
\qquad
R\leq8\varepsilon_0A.
\]
For \(y\geq1\), multiplication of
\eqref{eq:calculus-supercritical-F-jets} by \(R^k\) gives
\[
R^k
\left|
\partial_y^kF_{A,\alpha}(y)
\right|
\leq
C_{J_{>},M}
\left[
A^{-1}R^\alpha
+
R^{\alpha-1}
+
R^{\alpha-2}
\right]
\leq
C_{J_{>},M}A^{\alpha-1}.
\]
The last term is bounded by the same quantity both when \(\alpha\geq2\)
and when \(1+\delta_0\leq\alpha<2\), since \(R\geq1\) and
\(A^{\alpha-1}\geq1\). Similarly,
\eqref{eq:calculus-supercritical-DA-F-jets} gives
\[
R^k
\left|
\mathscr D_A
\partial_y^kF_{A,\alpha}(y)
\right|
\leq
C_{J_{>},M}
\bigl[\log(A+2)+\log(R+2)\bigr]
A^{\alpha-1}.
\]

For \(0\leq y\leq1\), the finite-\(A\) Kummer representation,
\eqref{eq:kummer-finite-A-parameters}, and compactness give
\begin{equation}
\max_{0\leq k\leq M}
\sup_{0\leq y\leq1}
\left|\partial_y^kF_{A,\alpha}(y)\right|
\leq C_{J_{>},M}.
\label{eq:calculus-supercritical-F-bounded-range}
\end{equation}
and
\[
\max_{0\leq k\leq M}
\sup_{0\leq y\leq1}
\left|
\mathscr D_A\partial_y^kF_{A,\alpha}(y)
\right|
\leq
C_{J_{>},M}\log(A+2).
\]
These are bounded by the required right-hand sides because
\(A^{\alpha-1}\geq1\).

By
\eqref{eq:calculus-qA},
\eqref{eq:calculus-beta-bounds}, and
\eqref{eq:calculus-beta-derivative},
\begin{equation}
|q_{A,\alpha}|
\leq
C_{J_{>}}A^{\alpha-1},
\qquad
|\mathscr D_Aq_{A,\alpha}|
=
A^\alpha|\widetilde\beta_{A,\alpha}|
\leq
C_{J_{>}}A^{\alpha-1}.
\label{eq:calculus-supercritical-q-bounds}
\end{equation}
Since
\[
L_1(y/A)-L_1(0)
=
-\frac yA,
\]
\begin{equation}
\Delta F_{A,\alpha}^{\circ}(y)
=
F_{A,\alpha}(y)-F_{A,\alpha}(0)
+
\frac{q_{A,\alpha}}A y.
\label{eq:calculus-supercritical-deltaF-decomposition}
\end{equation}
The bounded-range estimates control the value at the origin.
Together with
\eqref{eq:calculus-supercritical-F-jets},
\eqref{eq:calculus-supercritical-q-bounds}, and
\(R\leq8\varepsilon_0A\), this proves
\eqref{eq:calculus-supercritical-ordinary-profile-jets} and
\eqref{eq:calculus-supercritical-ordinary-core-jets}.

Finally,
\[
\mathcal P_{A,\alpha}
=
\mathscr D_A
\Delta F_{A,\alpha}^{\circ}.
\]
Applying \(\mathscr D_A\) to
\eqref{eq:calculus-supercritical-deltaF-decomposition} gives
\[
\mathcal P_{A,\alpha}(y)
=
\mathscr D_A
\left[
F_{A,\alpha}(y)-F_{A,\alpha}(0)
\right]
+
\frac yA\mathscr D_Aq_{A,\alpha}.
\]
Thus
\eqref{eq:calculus-supercritical-DA-F-jets} and
\eqref{eq:calculus-supercritical-q-bounds} prove
\eqref{eq:calculus-supercritical-log-profile-jets} and
\eqref{eq:calculus-supercritical-log-core-jets}.
\end{proof}

Now we combine the regional estimates into uniform Laguerre energy bounds for ordinary Abel blocks, with exponential decay at large intensities.

\begin{theorem}[Uniform ordinary regional energy bounds]
\label{thm:ordinary-spatial-ledger}
Let
\[
\mathfrak A=[\alpha_-,\alpha_+]\Subset(0,\infty),
\qquad \sigma_*:=\min\{\alpha_-,1\}.
\]
For every fixed integer \(s\geq1\), there exist constants
\(C_{\mathfrak A,s}<\infty\), \(B_s\geq0\), and \(c_s>0\)
such that, uniformly for \(A\geq2\),
\begin{align}
\sup_{\alpha\in\mathfrak A}\mathfrak E_s(r_{A,0,\alpha})
&\leq C_{\mathfrak A,s}A^{s-1-2\sigma_*},
\label{eq:calculus-uniform-core-energy}\\
\sup_{\alpha\in\mathfrak A}\mathfrak E_s(r_{A,R,\alpha})
&\leq C_{\mathfrak A,s}
\left(\frac AR\right)^{s-1}A^{-2\sigma_*},
\label{eq:calculus-uniform-lower-energy}\\
\sup_{\alpha\in\mathfrak A}\mathfrak E_s(f_{A,\alpha}^{\mathrm{mid}})
&\leq C_{\mathfrak A,s}A^{-2\sigma_*},
\label{eq:calculus-uniform-middle-energy}\\
\sup_{\alpha\in\mathfrak A}\mathfrak E_s(f_{A,R,\alpha}^{\mathrm{up}})
&\leq C_{\mathfrak A,s}A^{-2\sigma_*}
(1+R)^{B_s}e^{-c_sR}.
\label{eq:calculus-uniform-upper-energy}
\end{align}
In \eqref{eq:calculus-uniform-lower-energy}, \(R\) ranges over
dyadic scales satisfying \(1\leq R\leq8\varepsilon_0A\).
In \eqref{eq:calculus-uniform-upper-energy},
\(R=R_0 2^k\), \(k\in\mathbb N_0\).
All constants may depend on \(\mathfrak A\), \(s\), and the fixed
cutoff data, but are independent of \(A\), \(R\), and \(\alpha\).
\end{theorem}

\begin{proof}
Fix \(\delta_0=1/8\) and put
\[
J_{\leq}=\mathfrak A\cap(0,1+\delta_0],
\qquad
J_{>}=\mathfrak A\cap[1+\delta_0,\infty).
\]
These compact intervals cover \(\mathfrak A\); an empty interval
is omitted from the argument.

For \(\alpha\in J_{\leq}\),
\eqref{eq:calculus-ordinary-profile-scaling} identifies the scaled
function in the lower region as \(\Delta F_{A,\alpha}^{\circ}\).
On the support of \(\lambda_{\mathrm{low}}(\,\cdot\,/A)\omega_R\),
one has \(y\asymp R\) and \(y\leq4\varepsilon_0A\).
Proposition~\ref{prop:ordinary-scaled-Abel-estimates},
the cutoff bounds \eqref{eq:calculus-scaled-dyadic-cutoffs},
and Leibniz's rule imply
\[
\mathcal B_{A,R,s}(\Delta F_{A,\alpha}^{\circ})
\leq C_{\mathfrak A,s}
\left[1+R^{\alpha-1}+A^{-1}R^\alpha\right].
\]
The lowest-scale cutoff is supported in
\[
\supp\bigl[\chi_0(\,\cdot\,)
\lambda_{\mathrm{low}}(\,\cdot\,/A)\bigr]
\subseteq[0,\min\{2,4\varepsilon_0A\}]
\subseteq[0,8\varepsilon_0A].
\]
The same proposition therefore applies throughout this support and
implies \(\mathcal B_{A,0,s}(\Delta F_{A,\alpha}^{\circ})
\leq C_{\mathfrak A,s}\).
The scaling estimates
\eqref{eq:calculus-lower-all-order-energy-scaling}
and \eqref{eq:calculus-core-all-order-energy-scaling} now yield
\begin{align}
\mathfrak E_s(r_{A,R,\alpha})
&\leq C_{\mathfrak A,s}
\left(\frac AR\right)^{s-1}A^{-2\alpha}
\left[1+R^{2\alpha-2}+A^{-2}R^{2\alpha}\right],
\label{eq:calculus-ordinary-lower-all-orders}\\
\mathfrak E_s(r_{A,0,\alpha})
&\leq C_{\mathfrak A,s}A^{s-1-2\alpha}.
\label{eq:calculus-ordinary-core-all-orders}
\end{align}

For \(\alpha\in J_{>}\),
Lemma~\ref{lem:supercritical-scaled-profile-jets}, the cutoff
bounds, and Leibniz's rule imply
\[
\begin{aligned}
\mathcal B_{A,R,s}(\Delta F_{A,\alpha}^{\circ})
\leq C_{\mathfrak A,s}A^{\alpha-1},\qquad
\mathcal B_{A,0,s}(\Delta F_{A,\alpha}^{\circ})
\leq C_{\mathfrak A,s}A^{\alpha-1}.
\end{aligned}
\]
Substitution into
\eqref{eq:calculus-lower-all-order-energy-scaling}
and \eqref{eq:calculus-core-all-order-energy-scaling} yields
\begin{equation}
\begin{aligned}
\mathfrak E_s(r_{A,R,\alpha})
\leq C_{\mathfrak A,s}\left(\frac AR\right)^{s-1}A^{-2},\qquad
\mathfrak E_s(r_{A,0,\alpha})
\leq C_{\mathfrak A,s}A^{s-3}.
\end{aligned}
\label{eq:calculus-ordinary-supercritical-ledger}
\end{equation}

For \(1\leq R\leq8\varepsilon_0A<A\) and every \(\alpha>0\),
\begin{equation}
A^{-2\alpha}
\left[1+R^{2\alpha-2}+A^{-2}R^{2\alpha}\right]
\leq3A^{-2\min\{\alpha,1\}}.
\label{eq:calculus-regional-uniform-reduction}
\end{equation}
Indeed, if \(\alpha\leq1\), each term in the square brackets
is at most one; if \(\alpha\geq1\), each term after multiplication
by \(A^{-2\alpha}\) is at most \(A^{-2}\).
Since \(\min\{\alpha,1\}\geq\sigma_*\) on \(\mathfrak A\),
\eqref{eq:calculus-ordinary-lower-all-orders}--%
\eqref{eq:calculus-ordinary-supercritical-ledger}, together with
\eqref{eq:calculus-regional-uniform-reduction}, prove
\eqref{eq:calculus-uniform-core-energy}
and \eqref{eq:calculus-uniform-lower-energy}.

For the middle and upper regions, all estimates below are uniform
for \(\alpha\in\mathfrak A\).
For the middle region, choose fixed \(0<a<b<\infty\) such that
\(\supp\lambda_{\mathrm{mid}}\subset(a,b)\).
Lemma~\ref{lem:compact-upper-semigroup-jets} and
Proposition~\ref{prop:correction-coefficient-bounds}
imply, for every fixed nonnegative integer \(M\),
\[
\max_{0\leq m\leq M}\sup_{a/2\leq x\leq2b}
\left|\partial_x^m f_{A,\alpha}^{\ren}(x)\right|
\leq C_{\mathfrak A,M,a,b}A^{-\min\{\alpha,1\}}.
\]
Taking \(M\) at least as large as the derivative order required by
Lemma~\ref{lem:fixed-cutoff-Laguerre-localization} yields
\begin{equation}
\mathfrak E_s(f_{A,\alpha}^{\mathrm{mid}})
\leq C_{\mathfrak A,s}A^{-2\min\{\alpha,1\}}.
\label{eq:calculus-ordinary-middle-ledger}
\end{equation}
This proves \eqref{eq:calculus-uniform-middle-energy}.

For an upper shell, fix \(0<\kappa<1/4\).
Lemma~\ref{lem:compact-upper-semigroup-jets}, applied with
exponential parameter \(\kappa/2\), and
Proposition~\ref{prop:correction-coefficient-bounds}
imply, for every fixed \(M\),
\[
\max_{0\leq m\leq M}\sup_{R/2\leq x\leq2R}
\left|\partial_x^m f_{A,\alpha}^{\ren}(x)\right|
\leq C_{\mathfrak A,M}
A^{-\min\{\alpha,1\}}(1+R)^{C_M}e^{\kappa R}.
\]
Taking \(M\geq K_s\), where \(K_s\) is the derivative order in
Lemma~\ref{lem:dyadic-upper-Laguerre-localization}, and applying
that lemma yields, with \(c_s=1/2-2\kappa>0\) and a suitable
\(B_s\geq0\),
\begin{equation}
\mathfrak E_s(f_{A,R,\alpha}^{\mathrm{up}})
\leq C_{\mathfrak A,s}A^{-2\min\{\alpha,1\}}
(1+R)^{B_s}e^{-c_sR}.
\label{eq:calculus-ordinary-upper-ledger}
\end{equation}
Taking the supremum over \(\alpha\in\mathfrak A\) proves
\eqref{eq:calculus-uniform-upper-energy}.
\end{proof}

Next, we obtain the same regional energy bounds for logarithmic Abel blocks, with an extra factor of \(\log^2(A+2)\).

\begin{theorem}[Uniform logarithmic regional energy bounds]
\label{thm:piecewise-log-shell-ledger}
Let \(\mathfrak A\) and \(\sigma_*\) be as in
Theorem~\ref{thm:ordinary-spatial-ledger}.
For every fixed integer \(s\geq1\), there exist constants
\(C_{\mathfrak A,s}<\infty\), \(B_s\geq0\), and \(c_s>0\)
such that, uniformly for \(A\geq2\),
\begin{align}
\sup_{\alpha\in\mathfrak A}\mathfrak E_s(\widetilde r_{A,0,\alpha})
&\leq C_{\mathfrak A,s}A^{s-1-2\sigma_*}\log^2(A+2),
\label{eq:calculus-uniform-log-core-energy}\\
\sup_{\alpha\in\mathfrak A}\mathfrak E_s(\widetilde r_{A,R,\alpha})
&\leq C_{\mathfrak A,s}
\left(\frac AR\right)^{s-1}A^{-2\sigma_*}\log^2(A+2),
\label{eq:calculus-uniform-log-lower-energy}\\
\sup_{\alpha\in\mathfrak A}
\mathfrak E_s(\widetilde f_{A,\alpha}^{\mathrm{mid}})
&\leq C_{\mathfrak A,s}A^{-2\sigma_*}\log^2(A+2),
\label{eq:calculus-uniform-log-middle-energy}\\
\sup_{\alpha\in\mathfrak A}
\mathfrak E_s(\widetilde f_{A,R,\alpha}^{\mathrm{up}})
&\leq C_{\mathfrak A,s}A^{-2\sigma_*}\log^2(A+2)
(1+R)^{B_s}e^{-c_sR}.
\label{eq:calculus-uniform-log-upper-energy}
\end{align}
In \eqref{eq:calculus-uniform-log-lower-energy}, \(R\) ranges over
dyadic scales satisfying \(1\leq R\leq8\varepsilon_0A\).
In \eqref{eq:calculus-uniform-log-upper-energy},
\(R=R_0 2^k\), \(k\in\mathbb N_0\).
All constants may depend on \(\mathfrak A\), \(s\), and the fixed
cutoff data, but are independent of \(A\), \(R\), and \(\alpha\).
The constants in the two theorems may be chosen in common.
\end{theorem}

\begin{proof}
Use the compact intervals \(J_{\leq}\) and \(J_{>}\) from the
proof of Theorem~\ref{thm:ordinary-spatial-ledger}.
By \eqref{eq:calculus-log-profile-scaling}, the scaled function
for the logarithmic lower pieces is \(\mathcal P_{A,\alpha}\).
For \(\alpha\in J_{\leq}\),
Proposition~\ref{prop:logarithmic-scaled-Abel-estimates},
the cutoff bounds \eqref{eq:calculus-scaled-dyadic-cutoffs},
and Leibniz's rule imply
\[
\begin{aligned}
\mathcal B_{A,R,s}(\mathcal P_{A,\alpha})
&\leq C_{\mathfrak A,s}
\bigl[\log(A+2)+\log(R+2)\bigr]
\left[1+R^{\alpha-1}+A^{-1}R^\alpha\right],\\
\mathcal B_{A,0,s}(\mathcal P_{A,\alpha})
&\leq C_{\mathfrak A,s}\log(A+2).
\end{aligned}
\]
The proposition applies because the shell cutoff is supported where
\(y\asymp R\) and \(y\leq4\varepsilon_0A\), while the lowest-scale
cutoff is supported in \([0,\min\{2,4\varepsilon_0A\}]\).
The scaling estimates
\eqref{eq:calculus-lower-all-order-energy-scaling}
and \eqref{eq:calculus-core-all-order-energy-scaling} yield
\begin{equation}
\begin{aligned}
\mathfrak E_s(\widetilde r_{A,R,\alpha})
\leq C_{\mathfrak A,s}
\left(\frac AR\right)^{s-1}A^{-2\alpha}
\bigl[\log(A+2)+\log(R+2)\bigr]^2
\left[1+R^{2\alpha-2}+A^{-2}R^{2\alpha}\right],
\end{aligned}
\label{eq:calculus-log-lower-all-orders}
\end{equation}
and
\begin{equation}
\mathfrak E_s(\widetilde r_{A,0,\alpha})
\leq C_{\mathfrak A,s}A^{s-1-2\alpha}\log^2(A+2).
\label{eq:calculus-log-core-all-orders}
\end{equation}

For \(\alpha\in J_{>}\),
Lemma~\ref{lem:supercritical-scaled-profile-jets}, the cutoff
bounds, and Leibniz's rule imply
\[
\begin{aligned}
\mathcal B_{A,R,s}(\mathcal P_{A,\alpha})
&\leq C_{\mathfrak A,s}
\bigl[\log(A+2)+\log(R+2)\bigr]A^{\alpha-1},\\
\mathcal B_{A,0,s}(\mathcal P_{A,\alpha})
&\leq C_{\mathfrak A,s}\log(A+2)A^{\alpha-1}.
\end{aligned}
\]
The lower-shell and lowest-scale energy estimates
\eqref{eq:calculus-lower-all-order-energy-scaling}
and \eqref{eq:calculus-core-all-order-energy-scaling} therefore imply
\begin{equation}
\begin{aligned}
\mathfrak E_s(\widetilde r_{A,R,\alpha})
&\leq C_{\mathfrak A,s}\left(\frac AR\right)^{s-1}A^{-2}
\bigl[\log(A+2)+\log(R+2)\bigr]^2,\\
\mathfrak E_s(\widetilde r_{A,0,\alpha})
&\leq C_{\mathfrak A,s}A^{s-3}\log^2(A+2).
\end{aligned}
\label{eq:calculus-log-supercritical-ledger}
\end{equation}
Since \(R\leq A\),
\[
\log(A+2)+\log(R+2)\leq2\log(A+2).
\]
Combining this inequality and
\eqref{eq:calculus-regional-uniform-reduction} with
\eqref{eq:calculus-log-lower-all-orders}--%
\eqref{eq:calculus-log-supercritical-ledger}, and using
\(\min\{\alpha,1\}\geq\sigma_*\), proves
\eqref{eq:calculus-uniform-log-core-energy}
and \eqref{eq:calculus-uniform-log-lower-energy}.

For the middle and upper regions, all estimates below are uniform
for \(\alpha\in\mathfrak A\).
Choose fixed \(0<a<b<\infty\) with
\(\supp\lambda_{\mathrm{mid}}\subset(a,b)\) and fix
\(0<\kappa<1/4\).
Proposition~\ref{prop:correction-coefficient-bounds} and
Lemma~\ref{lem:compact-upper-semigroup-jets} imply, for every fixed
nonnegative integer \(M\),
\[
\max_{0\leq m\leq M}\sup_{a/2\leq x\leq2b}
\left|\partial_x^m\widetilde f_{A,\alpha}^{\ren}(x)\right|
\leq C_{\mathfrak A,M,a,b}
\log(A+2)A^{-\min\{\alpha,1\}}.
\]
and
\[
\begin{aligned}
\max_{0\leq m\leq M}\sup_{R/2\leq x\leq2R}
\left|\partial_x^m\widetilde f_{A,\alpha}^{\ren}(x)\right|\leq C_{\mathfrak A,M}
\log(A+2)A^{-\min\{\alpha,1\}}
(1+R)^{C_M}e^{\kappa R}.
\end{aligned}
\]
The upper-shell estimate uses part \textup{(ii)} of
Lemma~\ref{lem:compact-upper-semigroup-jets} with exponential
parameter \(\kappa/2\).
Choose \(M\) at least as large as the derivative order required by
Lemma~\ref{lem:fixed-cutoff-Laguerre-localization} and with
\(M\geq K_s\), where \(K_s\) is the derivative order in
Lemma~\ref{lem:dyadic-upper-Laguerre-localization}.
Applying these localization lemmas to the middle and upper pieces,
respectively, yields
\begin{align}
\mathfrak E_s(\widetilde f_{A,\alpha}^{\mathrm{mid}})
&\leq C_{\mathfrak A,s}A^{-2\min\{\alpha,1\}}\log^2(A+2),
\label{eq:calculus-log-middle-ledger}\\
\mathfrak E_s(\widetilde f_{A,R,\alpha}^{\mathrm{up}})
&\leq C_{\mathfrak A,s}A^{-2\min\{\alpha,1\}}\log^2(A+2)
(1+R)^{B_s}e^{-c_sR}.
\label{eq:calculus-log-upper-ledger}
\end{align}
Here \(c_s=\frac12-2\kappa>0\).
For the upper piece, the polynomial exponent in
\eqref{eq:calculus-dyadic-upper-localization} is increased
by \(2C_M\), enlarged further if necessary to be nonnegative,
and still denoted by \(B_s\).
Since
\[
A^{-2\min\{\alpha,1\}}\leq A^{-2\sigma_*},
\qquad \alpha\in\mathfrak A,\quad A\geq2,
\]
taking suprema over \(\alpha\in\mathfrak A\) proves
\eqref{eq:calculus-uniform-log-middle-energy} and
\eqref{eq:calculus-uniform-log-upper-energy}.

The constants in
Theorems~\ref{thm:ordinary-spatial-ledger} and
\ref{thm:piecewise-log-shell-ledger}
may be chosen in common by increasing
\(C_{\mathfrak A,s}\) and \(B_s\), and decreasing \(c_s>0\).
For every fixed integer \(S\geq1\), they may also be chosen so that
\[
C_{\mathfrak A,s}\leq C_{\mathfrak A,S},
\qquad B_s\leq B_S,
\qquad 0<c_S\leq c_s,
\qquad 1\leq s\leq S.
\]
Thus \(C_{\mathfrak A,S}\), \(B_S\), and \(c_S\) can be used
in all bounds of both theorems for \(1\leq s\leq S\),
with dependence only on \(\mathfrak A\), \(S\),
and the fixed cutoff data.
\end{proof}

\section{Maximal estimates over intermediate Abel scales}
\label{sec:active-abel-range}

This section proves maximal estimates for the Abel blocks on the
intermediate dyadic range
\(\Lambda_N\leq A\leq N^\beta\); see Theorem~\ref{thm:active-abel-range}.
The proof treats the lower, middle, and upper intensity regions separately,
using the same refinement of the logarithmic-time interval.
The refinement is the usual finite dyadic chaining construction, recorded
with the integer mesh and endpoint conventions in
Appendix~\ref{app:block-addition-refinement}.
Its application here requires uniform regional Laguerre-energy bounds and
summable projected and unprojected comparison errors over all
refinement levels and Abel scales.

Fix
\[
\mathfrak A=[\alpha_-,\alpha_+]\Subset(0,\infty),
\qquad T>0,
\qquad n_-=\lfloor Ne^{-T}\rfloor,
\qquad n_+=\lfloor Ne^T\rfloor.
\]
Set
\begin{equation}
\sigma_*:=\min\{\alpha_-,1\}>0.
\label{eq:active-critical-exponents}
\end{equation}
For a measurable deterministic family
\(f_\bullet=(f_\alpha)_{\alpha\in\mathfrak A}\)
for which the centered statistics are defined, put
\begin{equation}
\mathfrak M_{N,T;\mathfrak A}(f_\bullet)
:=
\left(
\int_{\mathfrak A}
\E\max_{n_-\leq n\leq n_+}|X_n(f_\alpha)|^2
\,\dd\alpha
\right)^{1/2}.
\label{eq:active-integrated-maximal-seminorm}
\end{equation}
Linearity of \(X_n\) and Minkowski's inequality imply that this
seminorm is subadditive under finite sums of such families.

Fix \(\gamma_0>0\) and define
\begin{equation}
\Lambda_N
:=
\min\left\{
A\geq2:
A\text{ is dyadic and }
A>\gamma_0\frac{\log(N+2)}{\log\log(N+3)}
\right\}.
\label{eq:moving-LambdaN}
\end{equation}
Then \(\Lambda_N\to\infty\) and \(\Lambda_N=N^{o(1)}\).
We consider dyadic scales
\begin{equation}
\Lambda_N\leq A\leq N^\beta,
\label{eq:active-Abel-range}
\end{equation}
where \(\beta>0\) will be chosen sufficiently small.

Fix an even integer \(q=2r\) and an integer \(J\) such that
\[
r\geq6,\qquad \frac1q<\sigma_*,
\qquad J\geq3.
\]
These orders are independent of \(N\) and \(A\) and are fixed
before \(\beta\) is chosen.  The condition \(1/q<\sigma_*\) makes
the principal factors \(A^{1/q-\sigma_*}\), including their logarithmic
multipliers, summable over dyadic scales.  Once the moment, expansion,
and derivative orders are fixed, choosing \(\beta\) sufficiently small
allows the negative powers of \(N\) in the comparison errors to dominate
their polynomial growth in \(A\).
Throughout this section, the fixed-time and increment estimates
are used with moment order \(q\), and the expansion remainders
are evaluated at order \(J\).

We use the corrected blocks and fixed intensity cutoffs from
Section~\ref{sec:coefficient-calculus}.
Since \(X_n(L_0)=X_n(L_1)=0\),
\begin{equation}
X_n(f_{A,\alpha})=X_n(f_{A,\alpha}^{\ren}),
\qquad
X_n(\widetilde f_{A,\alpha})
=X_n(\widetilde f_{A,\alpha}^{\ren}).
\label{eq:active-renormalized-process-equivalence}
\end{equation}
Define the lower and upper regional pieces by
\begin{equation}
\begin{aligned}
f_{A,\alpha}^{\mathrm{low}}
&:=\lambda_{\mathrm{low}}f_{A,\alpha}^{\ren},
&\quad
\widetilde f_{A,\alpha}^{\mathrm{low}}
&:=\lambda_{\mathrm{low}}\widetilde f_{A,\alpha}^{\ren},
\\
f_{A,\alpha}^{\mathrm{up}}
&:=\lambda_{\mathrm{up}}f_{A,\alpha}^{\ren},
&
\widetilde f_{A,\alpha}^{\mathrm{up}}
&:=\lambda_{\mathrm{up}}\widetilde f_{A,\alpha}^{\ren}.
\end{aligned}
\label{eq:active-spatial-sector-pieces}
\end{equation}
The middle pieces and the individual lower and upper shells retain
their definitions from Section~\ref{sec:coefficient-calculus}.
Thus
\begin{equation}
\begin{aligned}
f_{A,\alpha}^{\ren}
=f_{A,\alpha}^{\mathrm{low}}
+f_{A,\alpha}^{\mathrm{mid}}
+f_{A,\alpha}^{\mathrm{up}},
\quad
\widetilde f_{A,\alpha}^{\ren}
=\widetilde f_{A,\alpha}^{\mathrm{low}}
+\widetilde f_{A,\alpha}^{\mathrm{mid}}
+\widetilde f_{A,\alpha}^{\mathrm{up}}.
\end{aligned}
\label{eq:active-exact-spatial-partition}
\end{equation}

For \(n_-\leq m<n\leq n_+\), put
\[
\rho_{m,n}:=\frac{n-m}{n}.
\]

For all sufficiently large \(N\), depending only on \(T\),
and each dyadic \(A\geq2\), use the refinement constructed in
Lemma~\ref{lem:active-scale-A-refinement-tree}
with \(\Gamma=A\), and denote its coarse set and edge sets by
\(\mathcal C_{N,A}\) and \(\mathcal E_{N,A,j}\), respectively.

\subsection{Lower region}
\label{subsec:lower-region-maximal}

Next we use the lower-region energy bounds to control the small-intensity contributions uniformly over nested polynomial lengths.
For a deterministic family
\(u_\bullet=(u_\alpha)_{\alpha\in\mathfrak A}\)
of smooth compactly supported functions, we collect the errors in the following expression.
The first term bounds the error on the coarse grid.
The sum bounds the errors in increments along successive refinements;
the factor \((A2^j)^{1/q}\) accounts for the number of pairs at level \(j\).
We define
\begin{equation}
\begin{aligned}
\operatorname{Err}_{N,A}(u_\bullet)
:={}&
C A^{1/q}
\left[
\sup_{\alpha\in\mathfrak A}
\mathcal R_{N,q,J}^{\mathrm{end}}
(u_\alpha)
\right]^{1/2}
\\
&+
C
\sum_{\substack{j\geq0\\
\mathcal E_{N,A,j}\neq\varnothing}}
(A2^j)^{1/q}
\max_{(m,n)\in\mathcal E_{N,A,j}}
\sup_{\alpha\in\mathfrak A}
\left[
\operatorname{Err}_{m,n}^{\mathrm{lead,proj}}(u_\alpha)
+
\operatorname{Err}_{m,n}^{\mathrm{nonlin,dir}}(u_\alpha)
\right].
\end{aligned}
\label{eq:active-refinement-remainder-functional}
\end{equation}
The fixed-length and local remainders are defined in
\eqref{eq:direct-endpoint-remainder},
\eqref{eq:direct-leading-projected-error}, and
\eqref{eq:direct-nonlinear-error}, respectively.
The constant \(C\) is fixed sufficiently large and depends only on
\(\mathfrak A,T\) and the fixed moment and expansion orders.
The term
\(\operatorname{Err}_{m,n}^{\mathrm{diff,proj}}(u_\alpha)\)
is already included in
\(\operatorname{Err}_{m,n}^{\mathrm{nonlin,dir}}(u_\alpha)\);
see \eqref{eq:direct-nonlinear-error}.

Define the total lower-region remainder by
\begin{equation}
\begin{aligned}
\operatorname{Err}^{\mathrm{low}}_{N,A}
:=
\operatorname{Err}_{N,A}(r_{A,0,\bullet})
+
\operatorname{Err}_{N,A}(\widetilde r_{A,0,\bullet})
+
\sum_{\substack{
R=2^k,\ k\in\mathbb N_0\\
R\leq8\varepsilon_0A}}
\left[
\operatorname{Err}_{N,A}(r_{A,R,\bullet})
+
\operatorname{Err}_{N,A}(\widetilde r_{A,R,\bullet})
\right].
\end{aligned}
\label{eq:lower-total-refinement-remainder}
\end{equation}

\begin{lemma}
\label{lem:lower-region-error-summability}
There exists \(\beta_{\mathrm{low}}^{\mathrm{err}}>0\) such that,
for every \(0<\beta\leq\beta_{\mathrm{low}}^{\mathrm{err}}\),
the total refinement remainder in the lower region satisfies
\begin{equation}
\sum_{\substack{
\Lambda_N\leq A\leq N^\beta\\
A\ {\rm dyadic}}}
\operatorname{Err}^{\mathrm{low}}_{N,A}
\longrightarrow0
\qquad
\text{as }N\to\infty.
\label{eq:lower-region-error-summability}
\end{equation}
\end{lemma}
\begin{proof}

Throughout the proof, \(u_\alpha\) denotes any one of
\[
r_{A,0,\alpha},
\qquad
\widetilde r_{A,0,\alpha},
\qquad
r_{A,R,\alpha},
\qquad
\widetilde r_{A,R,\alpha},
\]
where, in the shell cases, \(R\) is dyadic and
\(1\leq R\leq8\varepsilon_0A\).
All estimates below are uniform over these choices.

Every \(u_\alpha\) is supported in the fixed interval
\([0,4\varepsilon_0]\).
Propositions~\ref{prop:ordinary-scaled-Abel-estimates}
and~\ref{prop:logarithmic-scaled-Abel-estimates},
together with
Lemma~\ref{lem:supercritical-scaled-profile-jets}
and the cutoff bounds
\eqref{eq:calculus-scaled-dyadic-cutoffs},
imply that, for every fixed integer \(M\geq0\), there exists
\(d_M\geq0\) such that
\begin{equation}
\sup_{\alpha\in\mathfrak A}
\max_{0\leq k\leq M}
\|\partial_x^k u_\alpha\|_\infty
\leq
C_M A^{d_M}\log(A+2).
\label{eq:lower-remainder-jet-growth}
\end{equation}
Indeed, differentiation under the substitutions \(y=Ax\)
and \(y=Rw\), followed by the Leibniz rule, introduces only
fixed powers of \(A\) and \(R^{-1}\).
Since \(1\leq R\leq8\varepsilon_0A\), these factors are bounded
by fixed powers of \(A\), while
\(\log(R+2)\leq C\log(A+2)\).

Moreover,
Theorems~\ref{thm:ordinary-spatial-ledger}
and~\ref{thm:piecewise-log-shell-ledger}
imply that, for every fixed integer \(s\geq1\),
\begin{equation}
\sup_{\alpha\in\mathfrak A}
\mathfrak E_s(u_\alpha)
\leq
C_s A^{s-1}\log^2(A+2).
\label{eq:lower-remainder-energy-growth}
\end{equation}
Here we have discarded the factor \(A^{-2\sigma_*}\)
and used \(R\geq1\) in the lower-shell estimates.

Since \(\supp u_\alpha\subset[0,4\varepsilon_0]\) and
\(\mathscr Ru_\alpha=u_\alpha'+xu_\alpha''\),
\eqref{eq:replica-compact-jet-seminorm},
\eqref{eq:direct-homogeneous-jet}, and
\eqref{eq:direct-generated-OU-ordinary-jet} give,
for every fixed integer \(K\geq0\),
\[
\begin{aligned}
\sup_{\alpha\in\mathfrak A}
\left[
\mathfrak J_K(u_\alpha)
+\mathfrak H_K(u_\alpha)
+\dot{\mathfrak J}_K(\mathscr Ru_\alpha)
\right]
&\leq
C_K\left(
1+\sup_{\alpha\in\mathfrak A}
\max_{0\leq j\leq2K+7}
\|\partial_x^j u_\alpha\|_\infty
\right)\\
&\leq
C_K A^{d_{2K+7}}\log(A+2),
\end{aligned}
\]
where the last inequality follows from
\eqref{eq:lower-remainder-jet-growth}.
This bounds the derivative quantities in the fixed-length and local
remainders appearing in
\eqref{eq:active-refinement-remainder-functional}.
By the Laguerre lowering identity
\eqref{eq:endpoint-resonant-Laguerre},
\[
\mathfrak E_s(\mathscr Ru_\alpha)
=
\sum_{a\geq2}(a-1)^s a^2
|\widehat{u_\alpha}_a|^2
\leq
\mathfrak E_{s+2}(u_\alpha),
\qquad s\geq1.
\]
Only finitely many energy and derivative orders occur,
because the moment and expansion orders are fixed.
The preceding bounds therefore control all energy and derivative factors
in these remainders by
\(CA^D[\log(A+2)]^M\), for fixed \(D,M\geq0\),
uniformly over the lower-region pieces.

Applying Lemma~\ref{lem:direct-exact-temporal-power-sum}
with \(\Gamma=A\), we obtain, for every fixed \(b\geq0\),
\begin{equation}
\begin{aligned}
&
\sum_{\substack{j\geq0\\
\mathcal E_{N,A,j}\neq\varnothing}}
(A2^j)^{1/q}
\max_{(m,n)\in\mathcal E_{N,A,j}}
\rho_{m,n}^{\,b}
\leq
C_{b,q,T}\log(N+2)
N^{(1/q-b)_+}
A^{\min\{1/q-b,0\}}.
\end{aligned}
\label{eq:lower-remainder-temporal-sum}
\end{equation}

Lemma~\ref{lem:direct-homogeneous-endpoint}
applies to every \(u_\alpha\), with common constants.
Set
\[
a_0
:=
\min\left\{
\frac12-\frac1q,\,
\frac{c_q}{2},\,
\frac{c_{q,J}}{2},\,
\frac{J-1}{2}
\right\}>0,
\]
where the comparison exponents are those in
\eqref{eq:direct-endpoint-remainder}
and \eqref{eq:direct-nonlinear-error}.
Using \eqref{eq:lower-remainder-jet-growth} and
\eqref{eq:lower-remainder-energy-growth} to estimate the terms in
\eqref{eq:direct-endpoint-remainder} and
\eqref{eq:direct-nonlinear-error}, and applying
\eqref{eq:lower-remainder-temporal-sum}
with \(b=0,\frac12,1\), we obtain
\begin{equation}
\begin{aligned}
A^{1/q}&
\left[
\sup_{\alpha\in\mathfrak A}
\mathcal R_{N,q,J}^{\mathrm{end}}
(u_\alpha)
\right]^{1/2}
+
\sum_{\substack{j\geq0\\
\mathcal E_{N,A,j}\neq\varnothing}}
(A2^j)^{1/q}
\max_{(m,n)\in\mathcal E_{N,A,j}}
\sup_{\alpha\in\mathfrak A}
\left[
\operatorname{Err}_{m,n}^{\mathrm{nonlin,dir}}(u_\alpha)
-
\operatorname{Err}_{m,n}^{\mathrm{diff,proj}}(u_\alpha)
\right]\\
&\leq
C N^{-a_0}A^{D_0}
[\log(N+2)]^{M_0},
\qquad 2\leq A\leq N,
\end{aligned}
\label{eq:lower-nonprojected-remainder-bound}
\end{equation}
for fixed \(D_0,M_0\geq0\).

Use the finite composition family \(\mathsf C_r\)
from Appendix~\ref{app:block-addition-refinement}, and write
\[
L=L(\mathbf k),
\qquad
K=K(\mathbf k).
\]
Using the explicit logarithmic factors in
\eqref{eq:direct-projected-B-P} and
\eqref{eq:direct-difference-B-P},
and proceeding as in the derivation of
\eqref{eq:direct-projected-errors-before-growth},
we substitute
\eqref{eq:lower-remainder-jet-growth} and
\eqref{eq:lower-remainder-energy-growth}
and sum by \eqref{eq:lower-remainder-temporal-sum}. Then, uniformly in \(\alpha\), 
\begin{equation}
\begin{aligned}
\sum_{\substack{j\geq0\\
\mathcal E_{N,A,j}\neq\varnothing}}&
(A2^j)^{1/q}
\max_{(m,n)\in\mathcal E_{N,A,j}}
\sup_{\alpha\in\mathfrak A}
\left[
\operatorname{Err}_{m,n}^{\mathrm{lead,proj}}(u_\alpha)
+
\operatorname{Err}_{m,n}^{\mathrm{diff,proj}}(u_\alpha)
\right]
\\
&\quad\leq
C A^{D_1}[\log(N+2)]^{M_1}
\left[
N^{-1/2+1/q}
+
\sum_{\substack{
\mathbf k\in\mathsf C_r\\
L(\mathbf k)\geq1}}
\left(
N^{-a_r(L,K)/q}
+
N^{-b_r(L,K)/q}
\right)
\right]
\\
&\quad\leq
C N^{-1/q}A^{D_1}[\log(N+2)]^{M_1},
\qquad 2\leq A\leq N,
\end{aligned}
\label{eq:lower-projected-remainder-bound}
\end{equation}
for fixed \(D_1,M_1\geq0\).
Indeed, the exponents defined in
\eqref{eq:direct-composition-a}
and \eqref{eq:direct-composition-b}
satisfy
\[
a_r(L,K)\geq1,
\qquad
b_r(L,K)\geq1
\]
for every nonempty composition, while
\[
\frac12-\frac1q\geq\frac1q.
\]
The finiteness of \(\mathsf C_r\) permits common choices
of \(D_1\), \(M_1\), and the multiplicative constant.

Combining
\eqref{eq:lower-nonprojected-remainder-bound}
and \eqref{eq:lower-projected-remainder-bound}
with \eqref{eq:active-refinement-remainder-functional},
we obtain fixed constants
\[
a:=\min\{a_0,1/q\}>0,
\qquad
D\geq1,
\qquad
M\geq0,
\]
such that
\begin{equation}
\operatorname{Err}_{N,A}(u_\bullet)
\leq
C N^{-a}A^D[\log(N+2)]^M,
\qquad 2\leq A\leq N.
\label{eq:lower-uniform-remainder-bound}
\end{equation}
These constants are uniform over all lower-region pieces
and over the ordinary and logarithmic families.

Choose
\[
\beta_{\mathrm{low}}^{\mathrm{err}}
:=
\min\left\{1,\frac{a}{2D}\right\}.
\]
For every
\(0<\beta\leq\beta_{\mathrm{low}}^{\mathrm{err}}\)
and \(A\leq N^\beta\),
\[
N^{-a}A^D\leq N^{-a/2}.
\]
For each \(A\), the number of lower-region pieces,
counting both families, is at most \(C\log(A+2)\).
The number of dyadic scales in
\([\Lambda_N,N^\beta]\) is at most \(C\log(N+2)\).
Therefore, by
\eqref{eq:lower-total-refinement-remainder}
and \eqref{eq:lower-uniform-remainder-bound},
\[
\begin{aligned}
\sum_{\substack{
\Lambda_N\leq A\leq N^\beta\\
A\ {\rm dyadic}}}
\operatorname{Err}^{\mathrm{low}}_{N,A}
\leq
C N^{-a/2}[\log(N+2)]^{M+2}
\longrightarrow0.
\end{aligned}
\]
This proves \eqref{eq:lower-region-error-summability}.
\end{proof}

\begin{theorem}
\label{thm:active-lower-sectors}
There exists \(\beta_{\mathrm{low}}>0\) such that, for every
\(0<\beta\leq\beta_{\mathrm{low}}\),
\begin{equation}
\sum_{\substack{
\Lambda_N\leq A\leq N^\beta\\
A\ {\rm dyadic}}}
\mathfrak M_{N,T;\mathfrak A}
\left(
f_{A,\bullet}^{\mathrm{low}}
\right)
\longrightarrow0
\qquad
\text{as }N\to\infty,
\label{eq:active-ordinary-lower-sector}
\end{equation}
and
\begin{equation}
\sum_{\substack{
\Lambda_N\leq A\leq N^\beta\\
A\ {\rm dyadic}}}
\mathfrak M_{N,T;\mathfrak A}
\left(
\widetilde f_{A,\bullet}^{\mathrm{low}}
\right)
\longrightarrow0
\qquad
\text{as }N\to\infty.
\label{eq:active-log-lower-sector}
\end{equation}
\end{theorem}

\begin{proof}
Set
\(\beta_{\mathrm{low}}:=\beta_{\mathrm{low}}^{\mathrm{err}}\)
as in Lemma~\ref{lem:lower-region-error-summability}, and fix
\(0<\beta\leq\beta_{\mathrm{low}}\).

For each dyadic \(A\in[\Lambda_N,N^\beta]\), apply
Theorem~\ref{thm:refinement-safe-maximal-interface}
with temporal scale \(\Gamma=A\) and parameters
\[
(\Lambda,\mathcal A)
=
\begin{cases}
\bigl(A,A^{1/2-\sigma_*}\bigr),
& \text{for }r_{A,0,\alpha},\\[1mm]
\left(A/R,(A/R)^{1/2}A^{-\sigma_*}\right),
& \text{for }r_{A,R,\alpha}.
\end{cases}
\]
Both choices satisfy \(1\leq\Lambda\leq\Gamma\), since
\(1\leq R\leq8\varepsilon_0A<A\) in the shell case.
Theorem~\ref{thm:ordinary-spatial-ledger} verifies
\eqref{eq:direct-exact-energy-contract} uniformly in
\(\alpha\) and \(R\).
Taking the \(L^2(\mathfrak A)\)-norm and using
\eqref{eq:block-direct-local-error},
\eqref{eq:block-refinement-max-error}, and
\eqref{eq:active-refinement-remainder-functional}, we obtain
\begin{equation}
\begin{aligned}
\mathfrak M_{N,T;\mathfrak A}(r_{A,0,\bullet})
&\leq
CA^{1/q-\sigma_*}
+\operatorname{Err}_{N,A}(r_{A,0,\bullet}),\\
\mathfrak M_{N,T;\mathfrak A}(r_{A,R,\bullet})
&\leq
CA^{1/q-\sigma_*}(1+R^{-1/2})
+\operatorname{Err}_{N,A}(r_{A,R,\bullet}).
\end{aligned}
\label{eq:active-lower-piece-maximal}
\end{equation}

For the logarithmic pieces, replace \(\mathcal A\) by
\(\mathcal A\log(A+2)\).
Theorem~\ref{thm:piecewise-log-shell-ledger} verifies
\eqref{eq:direct-exact-energy-contract}, so the preceding
argument gives
\begin{equation}
\begin{aligned}
\mathfrak M_{N,T;\mathfrak A}(\widetilde r_{A,0,\bullet})
&\leq
C\log(A+2)A^{1/q-\sigma_*}
+\operatorname{Err}_{N,A}(\widetilde r_{A,0,\bullet}),\\
\mathfrak M_{N,T;\mathfrak A}(\widetilde r_{A,R,\bullet})
&\leq
C\log(A+2)A^{1/q-\sigma_*}(1+R^{-1/2})
+\operatorname{Err}_{N,A}(\widetilde r_{A,R,\bullet}).
\end{aligned}
\label{eq:active-lower-log-piece-maximal}
\end{equation}

By \eqref{eq:active-spatial-sector-pieces},
\eqref{eq:calculus-exact-lower-cell-decomposition},
and its logarithmic counterpart,
\[
\begin{aligned}
f_{A,\alpha}^{\mathrm{low}}
=r_{A,0,\alpha}
+\sum_{\substack{R=2^k,\ k\in\mathbb N_0\\R\leq8\varepsilon_0A}}
r_{A,R,\alpha},
\qquad
\widetilde f_{A,\alpha}^{\mathrm{low}}
=\widetilde r_{A,0,\alpha}
+\sum_{\substack{R=2^k,\ k\in\mathbb N_0\\R\leq8\varepsilon_0A}}
\widetilde r_{A,R,\alpha}.
\end{aligned}
\]
There are at most \(C\log(A+2)\) lower shells, whereas
\[
\sum_{\substack{R=2^k,\ k\in\mathbb N_0\\R\leq8\varepsilon_0A}}
R^{-1/2}
\leq\sum_{k=0}^{\infty}2^{-k/2}<\infty.
\]
Subadditivity of \(\mathfrak M_{N,T;\mathfrak A}\), together with
\eqref{eq:active-lower-piece-maximal},
\eqref{eq:active-lower-log-piece-maximal}, and
\eqref{eq:lower-total-refinement-remainder}, therefore implies
\begin{equation}
\begin{aligned}
\mathfrak M_{N,T;\mathfrak A}(f_{A,\bullet}^{\mathrm{low}})
+\mathfrak M_{N,T;\mathfrak A}(\widetilde f_{A,\bullet}^{\mathrm{low}})
\leq C\left[
[\log(A+2)]^2A^{1/q-\sigma_*}
+\log(A+2)A^{1/q-\sigma_*}
\right]
+\operatorname{Err}^{\mathrm{low}}_{N,A}.
\end{aligned}
\label{eq:active-lower-combined-maximal}
\end{equation}

Set \(\delta:=\sigma_*-1/q>0\). The powers of
\(A\) in \eqref{eq:active-lower-combined-maximal} are at most
\(A^{-\delta}\).
For all sufficiently large \(A\),
\([\log(A+2)]^2\leq A^{\delta/2}\).
As \(\Lambda_N\to\infty\), it follows that
\[
\begin{aligned}
\sum_{\substack{\Lambda_N\leq A\leq N^\beta\\A\ {\rm dyadic}}}
\Bigl[
\mathfrak M_{N,T;\mathfrak A}(f_{A,\bullet}^{\mathrm{low}})
+\mathfrak M_{N,T;\mathfrak A}(\widetilde f_{A,\bullet}^{\mathrm{low}})
\Bigr]
&\leq
C\sum_{\substack{A\geq\Lambda_N\\A\ {\rm dyadic}}}A^{-\delta/2}
+\sum_{\substack{\Lambda_N\leq A\leq N^\beta\\A\ {\rm dyadic}}}
\operatorname{Err}^{\mathrm{low}}_{N,A}
\\
&\leq C_\delta\Lambda_N^{-\delta/2}
+\sum_{\substack{\Lambda_N\leq A\leq N^\beta\\A\ {\rm dyadic}}}
\operatorname{Err}^{\mathrm{low}}_{N,A}
\longrightarrow0,
\end{aligned}
\]
where the remainder sum tends to zero by
Lemma~\ref{lem:lower-region-error-summability}.
Nonnegativity proves both
\eqref{eq:active-ordinary-lower-sector} and
\eqref{eq:active-log-lower-sector}.
\end{proof}

\subsection{Middle region}
\label{subsec:middle-region-maximal}

Now we control the middle-intensity contribution on a fixed interval away from zero, uniformly over nested polynomial lengths.

\begin{theorem}
\label{thm:active-middle-sectors}
There exists \(\beta_{\mathrm{mid}}>0\) such that, for every
\(0<\beta\leq\beta_{\mathrm{mid}}\),
\begin{equation}
\sum_{\substack{
\Lambda_N\leq A\leq N^\beta\\
A\ {\rm dyadic}}}
\mathfrak M_{N,T;\mathfrak A}
\left(
f_{A,\bullet}^{\mathrm{mid}}
\right)
\longrightarrow0
\qquad
\text{as }N\to\infty,
\label{eq:active-ordinary-middle-sector}
\end{equation}
and
\begin{equation}
\sum_{\substack{
\Lambda_N\leq A\leq N^\beta\\
A\ {\rm dyadic}}}
\mathfrak M_{N,T;\mathfrak A}
\left(
\widetilde f_{A,\bullet}^{\mathrm{mid}}
\right)
\longrightarrow0
\qquad
\text{as }N\to\infty.
\label{eq:active-log-middle-sector}
\end{equation}
\end{theorem}
\begin{proof}
The middle-region families are supported in the fixed
interval \([\varepsilon_0/2,2R_0]\).
Lemma~\ref{lem:compact-upper-semigroup-jets},
Proposition~\ref{prop:correction-coefficient-bounds},
and Leibniz's rule imply that, for every fixed integer
\(K\geq0\),
\[
\sup_{\alpha\in\mathfrak A}
\max_{0\leq k\leq K}
\left(
\|\partial_x^k f_{A,\alpha}^{\mathrm{mid}}\|_\infty
+
\|\partial_x^k\widetilde f_{A,\alpha}^{\mathrm{mid}}\|_\infty
\right)
\leq C_K,
\qquad A\geq2,
\]
where \(C_K\) is independent of \(A\).
Together with the fixed support, these derivative bounds
control the quantities in
\eqref{eq:replica-compact-jet-seminorm},
\eqref{eq:direct-homogeneous-jet}, and
\eqref{eq:direct-generated-OU-ordinary-jet}
uniformly in \(A\) and \(\alpha\).

Combining these bounds with the energy estimates of
Theorems~\ref{thm:ordinary-spatial-ledger}
and~\ref{thm:piecewise-log-shell-ledger}, we apply the
argument proving \eqref{eq:lower-uniform-remainder-bound}.
Thus there exist constants \(a>0\), \(D\geq1\),
\(M\geq0\), and \(C<\infty\), independent of \(N\) and \(A\),
such that, for all sufficiently large \(N\) and every
dyadic \(A\in[2,N]\),
\[
\begin{aligned}
\operatorname{Err}_{N,A}(f_{A,\bullet}^{\mathrm{mid}})
+\operatorname{Err}_{N,A}
(\widetilde f_{A,\bullet}^{\mathrm{mid}})
\leq
C N^{-a}A^D[\log(N+2)]^M.
\end{aligned}
\]
Choose
\[
\beta_{\mathrm{mid}}
:=\min\left\{1,\frac{a}{2D}\right\}.
\]
For every \(0<\beta\leq\beta_{\mathrm{mid}}\), there are
at most \(C\log(N+2)\) dyadic values of \(A\) satisfying
\(\Lambda_N\leq A\leq N^\beta\), so
\begin{equation}
\begin{aligned}
\sum_{\substack{
\Lambda_N\leq A\leq N^\beta\\
A\ {\rm dyadic}}}
\Bigl[
\operatorname{Err}_{N,A}(f_{A,\bullet}^{\mathrm{mid}})
+\operatorname{Err}_{N,A}(\widetilde f_{A,\bullet}^{\mathrm{mid}})
\Bigr]
\leq C N^{-a/2}[\log(N+2)]^{M+1}
\longrightarrow0.
\end{aligned}
\label{eq:active-middle-error-summability}
\end{equation}

Apply Theorem~\ref{thm:refinement-safe-maximal-interface}
with \(\Gamma=A\), \(\Lambda=1\), and
\[
\mathcal A=
\begin{cases}
A^{-\sigma_*},
& \text{for }f_{A,\bullet}^{\mathrm{mid}},\\
A^{-\sigma_*}\log(A+2),
& \text{for }\widetilde f_{A,\bullet}^{\mathrm{mid}}.
\end{cases}
\]
\eqref{eq:direct-exact-energy-contract}
follows from
Theorems~\ref{thm:ordinary-spatial-ledger}
and~\ref{thm:piecewise-log-shell-ledger}.
Taking the \(L^2(\mathfrak A)\)-norm and using
\eqref{eq:active-refinement-remainder-functional},
we obtain
\begin{equation}
\begin{aligned}
\mathfrak M_{N,T;\mathfrak A}(f_{A,\bullet}^{\mathrm{mid}})
+\mathfrak M_{N,T;\mathfrak A}
(\widetilde f_{A,\bullet}^{\mathrm{mid}})
&\leq C\log(A+2)
\left[
A^{1/q-\sigma_*}
+A^{1/q-1/2-\sigma_*}
\right]
\\
&\quad+
\operatorname{Err}_{N,A}(f_{A,\bullet}^{\mathrm{mid}})
+\operatorname{Err}_{N,A}(\widetilde f_{A,\bullet}^{\mathrm{mid}}).
\end{aligned}
\label{eq:active-middle-combined-maximal}
\end{equation}

Since \(1/q<\sigma_*\),
summing \eqref{eq:active-middle-combined-maximal} over dyadic \(A\in[\Lambda_N,N^\beta]\) and using
\eqref{eq:active-middle-error-summability} proves
\eqref{eq:active-ordinary-middle-sector}
and \eqref{eq:active-log-middle-sector}.
\end{proof}

\subsection{Upper region}
\label{subsec:upper-region-maximal}

Next we combine upper-shell energy bounds with the small probability of very large intensities to control the upper-region contribution uniformly over nested lengths.
For a fixed \(\kappa_{\mathrm{up}}>0\), define
\begin{equation}
\begin{aligned}
\operatorname{Err}^{\mathrm{up}}_{N,A}
:=
\sum_{\substack{
R=R_0 2^k,\ k\in\mathbb N_0\\
R\leq\kappa_{\mathrm{up}}\log(N+2)}}
\Bigl[
\operatorname{Err}_{N,A}
(f_{A,R,\bullet}^{\mathrm{up}})
+
\operatorname{Err}_{N,A}
(\widetilde f_{A,R,\bullet}^{\mathrm{up}})
\Bigr].
\end{aligned}
\label{eq:upper-total-refinement-remainder}
\end{equation}
\begin{lemma}
\label{lem:upper-region-error-summability}
For every fixed \(\kappa_{\mathrm{up}}>0\), there exist
\(\varepsilon>0\) and
\(\beta_{\mathrm{up}}^{\mathrm{err}}>0\),
chosen in this order, such that, for every
\(0<\beta\leq\beta_{\mathrm{up}}^{\mathrm{err}}\),
the total refinement remainder over the compact upper shells satisfies
\begin{equation}
\sum_{\substack{
\Lambda_N\leq A\leq N^\beta\\
A\ {\rm dyadic}}}
\operatorname{Err}^{\mathrm{up}}_{N,A}
\longrightarrow0
\qquad
\text{as }N\to\infty.
\label{eq:upper-region-error-summability}
\end{equation}
\end{lemma}

\begin{proof}
Let \(u_\alpha\) denote either of the two upper-shell families
\[
f_{A,R,\alpha}^{\mathrm{up}},
\qquad
\widetilde f_{A,R,\alpha}^{\mathrm{up}},
\qquad R=R_0 2^k,\quad k\in\mathbb N_0.
\]
All bounds are uniform in \(N,A,R,\alpha\) and the choice of family.

By \eqref{eq:calculus-upper-shell},
\[
\operatorname{supp}u_\alpha\subseteq[R/2,2R].
\]
For every fixed integer \(K\geq0\) and every \(\varepsilon>0\),
Lemma~\ref{lem:compact-upper-semigroup-jets}, applied with
exponential parameter \(\varepsilon/2\),
Proposition~\ref{prop:correction-coefficient-bounds}, and
the cutoff bounds \eqref{eq:calculus-physical-dyadic-cutoffs}
imply
\begin{equation}
\sup_{\alpha\in\mathfrak A}
\max_{0\leq j\leq K}\|\partial_x^j u_\alpha\|_\infty
\leq C_{K,\varepsilon}
A^{-\sigma_*}\log(A+2)
(1+R)^{b_{K,\varepsilon}}e^{\varepsilon R},
\label{eq:upper-remainder-derivative-bound}
\end{equation}
where \(b_{K,\varepsilon}\geq0\).

Since \(A^{-\sigma_*}\log(A+2)\) is bounded for \(A\geq2\), \eqref{eq:upper-remainder-derivative-bound},
\eqref{eq:replica-compact-jet-seminorm},
\eqref{eq:direct-homogeneous-jet}, and
\eqref{eq:direct-generated-OU-ordinary-jet} imply that,
for every fixed integer \(K\geq0\) and every
\(\varepsilon>0\),
\[
\sup_{\alpha\in\mathfrak A}
\max_{v\in\{u_\alpha,\mathscr Ru_\alpha\}}
\left[
\mathfrak J_K(v)
+\dot{\mathfrak J}_K(v)
+\mathfrak H_K(v)
\right]
\leq
C_{K,\varepsilon}(1+R)^{b_{K,\varepsilon}}
e^{\varepsilon R}.
\]
Here,
\[
\mathscr Ru_\alpha=u_\alpha'+xu_\alpha''.
\]
The energies of \(u_\alpha\) are controlled by
\eqref{eq:calculus-uniform-upper-energy} and
\eqref{eq:calculus-uniform-log-upper-energy};
the corresponding bounds for \(\mathscr Ru_\alpha\)
follow from \eqref{eq:endpoint-resonant-Laguerre}.

We now use the remainder argument in the proof of
Lemma~\ref{lem:lower-region-error-summability}, retaining the
dependence on \(R\).
We insert the common energy and derivative bounds before applying the temporal sum
\eqref{eq:lower-remainder-temporal-sum}.
The estimates leading to
\eqref{eq:lower-nonprojected-remainder-bound} and
\eqref{eq:lower-projected-remainder-bound} then yield
\begin{equation}
\operatorname{Err}_{N,A}(u_\bullet)
\leq C_\varepsilon N^{-a}A^D
(1+R)^{B_\varepsilon}
[\log(N+2)]^{M_\varepsilon}
e^{c_*\varepsilon R},
\label{eq:upper-uniform-remainder-bound}
\end{equation}
for all sufficiently large \(N\), uniformly for dyadic
\(2\leq A\leq N\) and \(R=R_0 2^k\).
The fixed powers in these remainder formulas allow
\(a,c_*>0\) and \(D\geq1\) to be chosen independently of
\(\varepsilon\).

Choose \(\varepsilon\) and then
\(\beta_{\mathrm{up}}^{\mathrm{err}}\) so that
\[
0<\varepsilon\leq
\min\left\{1,\frac{a}{4c_*\kappa_{\mathrm{up}}}\right\},
\qquad
\beta_{\mathrm{up}}^{\mathrm{err}}
:=\min\left\{1,\frac{a}{4D}\right\}.
\]
For \(0<\beta\leq\beta_{\mathrm{up}}^{\mathrm{err}}\),
\(A\leq N^\beta\), and
\(R\leq\kappa_{\mathrm{up}}\log(N+2)\), we have
\[
A^D\leq N^{a/4},
\qquad
e^{c_*\varepsilon R}
\leq (N+2)^{c_*\varepsilon\kappa_{\mathrm{up}}}
\leq C N^{a/4}.
\]
Consequently, \eqref{eq:upper-uniform-remainder-bound} implies
\[
\operatorname{Err}_{N,A}(u_\bullet)
\leq C N^{-a/2}
[\log(N+2)]^{B_\varepsilon+M_\varepsilon}.
\]
There are at most \(C\log(N+2)\) dyadic values of \(A\)
in \([\Lambda_N,N^\beta]\), and at most
\(C[1+\log\log(N+3)]\) compact upper shells for each \(A\).
Summing over both families therefore yields
\[
\begin{aligned}
\sum_{\substack{
\Lambda_N\leq A\leq N^\beta\\
A\ {\rm dyadic}}}
\operatorname{Err}^{\mathrm{up}}_{N,A}
\leq C N^{-a/2}
[\log(N+2)]^{B_\varepsilon+M_\varepsilon+1}
[1+\log\log(N+3)]
\longrightarrow0.
\end{aligned}
\]
\end{proof}

For the ordinary and logarithmic upper families, define the far-upper
tails by
\begin{equation}
\begin{aligned}
f_{A,\alpha}^{\mathrm{up,far}}
:=
\sum_{\substack{
R=R_0 2^k,\ k\in\mathbb N_0\\
R>\kappa_{\mathrm{up}}\log(N+2)}}
f_{A,R,\alpha}^{\mathrm{up}},
\qquad
\widetilde f_{A,\alpha}^{\mathrm{up,far}}
:=
\sum_{\substack{
R=R_0 2^k,\ k\in\mathbb N_0\\
R>\kappa_{\mathrm{up}}\log(N+2)}}
\widetilde f_{A,R,\alpha}^{\mathrm{up}}.
\end{aligned}
\label{eq:active-far-upper-label}
\end{equation}

\begin{lemma}
\label{lem:active-upper-intensity-tail}
For \(\kappa_{\mathrm{up}}>0\) sufficiently large, there exists
\(\beta_{\mathrm{far}}>0\) such that, for every
\(0<\beta\leq\beta_{\mathrm{far}}\),
\begin{equation}
\sum_{\substack{
\Lambda_N\leq A\leq N^\beta\\
A\ {\rm dyadic}}}
\left[
\mathfrak M_{N,T;\mathfrak A}
\left(
f_{A,\bullet}^{\mathrm{up,far}}
\right)
+
\mathfrak M_{N,T;\mathfrak A}
\left(
\widetilde f_{A,\bullet}^{\mathrm{up,far}}
\right)
\right]
\longrightarrow0
\qquad
\text{as }N\to\infty.
\label{eq:active-far-upper-tail}
\end{equation}
\end{lemma}
\begin{proof}
Hoeffding's inequality applied to the real and imaginary parts
of \(q_n(e^{i\theta})\) yields
\begin{equation}
\sup_{\substack{n\geq1\\\theta\in\mathbb T}}
\mathbb P\bigl(S_n(\theta)\geq x\bigr)
\leq4e^{-x/4},
\qquad x\geq0.
\label{eq:active-intensity-tail}
\end{equation}
Integrating this bound, we obtain, for \(L\geq0\),
\[
\begin{aligned}
\sup_{\substack{n\geq1\\\theta\in\mathbb T}}
\E\!\left[
e^{S_n(\theta)/8}
\mathbf 1_{\{S_n(\theta)\geq L\}}
\right]
\leq4e^{-L/8}
+\frac12\int_L^\infty e^{-x/8}\,\dd x
=8e^{-L/8}.
\end{aligned}
\]

Consider either family
\[
u_\bullet\in
\bigl\{f_{A,\bullet}^{\mathrm{up,far}},
\widetilde f_{A,\bullet}^{\mathrm{up,far}}\bigr\}.
\]
The support properties of the cutoffs and
\eqref{eq:calculus-exact-physical-dyadic-partition} imply
\[
0\leq
\lambda_{\mathrm{up}}(x)
\sum_{\substack{
R=R_0 2^k,\ k\in\mathbb N_0\\
R>\kappa_{\mathrm{up}}\log(N+2)}}
\varpi_R(x)
\leq
\mathbf 1_{\{x\geq\frac12\kappa_{\mathrm{up}}\log(N+2)\}}.
\]
By \eqref{eq:calculus-upper-shell},
Lemma~\ref{lem:compact-upper-semigroup-jets}\textup{(ii)}
with exponential parameter \(1/32\), and
Proposition~\ref{prop:correction-coefficient-bounds}, we obtain
\begin{equation}
\sup_{\alpha\in\mathfrak A}|u_\alpha(x)|
\leq C A^{-\sigma_*}\log(A+2)e^{x/16}
\mathbf 1_{\{x\geq\frac12\kappa_{\mathrm{up}}\log(N+2)\}},
\label{eq:active-far-upper-growth}
\end{equation}
uniformly for \(A\geq2\) and \(x\geq0\), after absorbing the
fixed polynomial factors into \(e^{x/32}\).
Taking \(L=\frac12\kappa_{\mathrm{up}}\log(N+2)\)
in the exponential-moment bound, we obtain
\begin{equation}
\sup_{\substack{
\alpha\in\mathfrak A\\
n\geq1,\ \theta\in\mathbb T}}
\E\bigl|u_\alpha(S_n(\theta))\bigr|^2
\leq C A^{-2\sigma_*}\log^2(A+2)
(N+2)^{-\kappa_{\mathrm{up}}/16}.
\label{eq:active-far-upper-one-angle}
\end{equation}

By \eqref{eq:model-In}, \eqref{eq:model-Xn}, and Jensen's inequality,
\[
\E|X_n(u_\alpha)|^2
\leq n\E|I_n(u_\alpha)|^2
\leq n\int_{\mathbb T}
\E|u_\alpha(S_n(\theta))|^2\,\dd\mu(\theta).
\]
Together with
\[
\sum_{n=n_-}^{n_+}n\leq C_TN^2,
\]
this yields
\begin{equation}
\begin{aligned}
\mathfrak M_{N,T;\mathfrak A}(u_\bullet)
\leq
\left(
\sum_{n=n_-}^{n_+}
\int_{\mathfrak A}\E|X_n(u_\alpha)|^2\,\dd\alpha
\right)^{1/2}
\leq C N^{1-\kappa_{\mathrm{up}}/32}
A^{-\sigma_*}\log(A+2).
\end{aligned}
\label{eq:active-far-upper-maximal}
\end{equation}

Choose \(\kappa_{\mathrm{up}}>32\) and set
\(\beta_{\mathrm{far}}=1\).
Since \(\sigma_*>0\), for every
\(0<\beta\leq\beta_{\mathrm{far}}\),
\[
\begin{aligned}
\sum_{\substack{
\Lambda_N\leq A\leq N^\beta\\
A\ {\rm dyadic}}}
\Bigl[
\mathfrak M_{N,T;\mathfrak A}
(f_{A,\bullet}^{\mathrm{up,far}})
+
\mathfrak M_{N,T;\mathfrak A}
(\widetilde f_{A,\bullet}^{\mathrm{up,far}})
\Bigr]
&\leq
C N^{1-\kappa_{\mathrm{up}}/32}
\sum_{\substack{A\geq2\\A\ {\rm dyadic}}}
A^{-\sigma_*}\log(A+2)
\\
&\leq C N^{1-\kappa_{\mathrm{up}}/32}
\longrightarrow0.
\end{aligned}
\]
\end{proof}

\begin{theorem}
\label{thm:active-upper-sectors}
There exists \(\beta_{\mathrm{up}}>0\) such that, for every
\(0<\beta\leq\beta_{\mathrm{up}}\),
\begin{equation}
\sum_{\substack{
\Lambda_N\leq A\leq N^\beta\\
A\ {\rm dyadic}}}
\mathfrak M_{N,T;\mathfrak A}
\left(
f_{A,\bullet}^{\mathrm{up}}
\right)
\longrightarrow0
\qquad
\text{as }N\to\infty,
\label{eq:active-ordinary-upper-sector}
\end{equation}
and
\begin{equation}
\sum_{\substack{
\Lambda_N\leq A\leq N^\beta\\
A\ {\rm dyadic}}}
\mathfrak M_{N,T;\mathfrak A}
\left(
\widetilde f_{A,\bullet}^{\mathrm{up}}
\right)
\longrightarrow0
\qquad
\text{as }N\to\infty.
\label{eq:active-log-upper-sector}
\end{equation}
\end{theorem}

\begin{proof}
Choose \(\kappa_{\mathrm{up}}\) sufficiently large for
Lemma~\ref{lem:active-upper-intensity-tail}, and apply
Lemma~\ref{lem:upper-region-error-summability} with this choice.
Set
\[
\beta_{\mathrm{up}}
:=\min\{\beta_{\mathrm{up}}^{\mathrm{err}},\beta_{\mathrm{far}}\},
\]
and fix \(0<\beta\leq\beta_{\mathrm{up}}\).

For each compact upper shell
\[
R=R_0 2^k\leq\kappa_{\mathrm{up}}\log(N+2),
\qquad k\in\mathbb N_0,
\]
let
\[
u_\bullet\in
\bigl\{
f_{A,R,\bullet}^{\mathrm{up}},
\widetilde f_{A,R,\bullet}^{\mathrm{up}}
\bigr\},
\]
and set \(\Gamma=A\), \(\Lambda=1\), and
\[
\mathcal A=
\begin{cases}
A^{-\sigma_*}(1+R)^B e^{-cR},
& u_\bullet=f_{A,R,\bullet}^{\mathrm{up}},\\[1mm]
A^{-\sigma_*}\log(A+2)(1+R)^B e^{-cR},
& u_\bullet=\widetilde f_{A,R,\bullet}^{\mathrm{up}}.
\end{cases}
\]
Choose \(B\geq0\) and \(c>0\), independently of \(N,A,R\),
so that the upper-region bounds in
Theorems~\ref{thm:ordinary-spatial-ledger}
and~\ref{thm:piecewise-log-shell-ledger} imply
\[
\sup_{\alpha\in\mathfrak A}
\mathfrak E_s(u_\alpha)^{1/2}
\leq C_s\mathcal A,
\qquad s=1,2,3.
\]
These bounds verify
\eqref{eq:direct-exact-energy-contract} with \(\Lambda=1\).
Applying Theorem~\ref{thm:refinement-safe-maximal-interface}
with these choices and taking the \(L^2(\mathfrak A)\)-norm,
we obtain
\[
\mathfrak M_{N,T;\mathfrak A}(u_\bullet)
\leq
C A^{1/q}(1+A^{-1/2})\mathcal A
+
\operatorname{Err}_{N,A}(u_\bullet),
\]
where the remainder functional is defined in
\eqref{eq:active-refinement-remainder-functional}.

Since
\[
\sum_{k=0}^{\infty}(1+R_0 2^k)^B e^{-cR_0 2^k}<\infty,
\]
the principal terms can be summed uniformly over the compact
upper shells. Together with
\eqref{eq:calculus-exact-upper-shell-decomposition},
\eqref{eq:active-far-upper-label}, and subadditivity, this yields
\begin{equation}
\begin{aligned}
\mathfrak M_{N,T;\mathfrak A}(f_{A,\bullet}^{\mathrm{up}})
+\mathfrak M_{N,T;\mathfrak A}
(\widetilde f_{A,\bullet}^{\mathrm{up}})
&\leq C\log(A+2)
\left[
A^{1/q-\sigma_*}
+A^{1/q-1/2-\sigma_*}
\right]
+\operatorname{Err}^{\mathrm{up}}_{N,A}
\\
&\quad+
\mathfrak M_{N,T;\mathfrak A}(f_{A,\bullet}^{\mathrm{up,far}})
+\mathfrak M_{N,T;\mathfrak A}
(\widetilde f_{A,\bullet}^{\mathrm{up,far}}),
\end{aligned}
\label{eq:active-upper-combined-maximal}
\end{equation}
where the compact-shell remainders are collected in
\eqref{eq:upper-total-refinement-remainder}.

Since \(1/q<\sigma_*\), both exponents of \(A\) in the principal
term are negative. Hence, as \(\Lambda_N\to\infty\),
\[
\sum_{\substack{A\geq\Lambda_N\\A\ {\rm dyadic}}}
\log(A+2)
\left[
A^{1/q-\sigma_*}
+A^{1/q-1/2-\sigma_*}
\right]
\longrightarrow0.
\]
Summing \eqref{eq:active-upper-combined-maximal} over
\(\Lambda_N\leq A\leq N^\beta\) and applying
Lemmas~\ref{lem:upper-region-error-summability}
and~\ref{lem:active-upper-intensity-tail} completes the proof.
\end{proof}

\subsection{Completion of the intermediate Abel range}
\label{subsec:completion-intermediate-abel-range}

We now combine the three regional estimates to show that the total contribution from the intermediate Abel scales tends to zero.
Choose
\[
0<\beta
\leq
\min\{
\beta_{\mathrm{low}},
\beta_{\mathrm{mid}},
\beta_{\mathrm{up}}
\}.
\]

\begin{theorem}
\label{thm:active-abel-range}
The ordinary Abel blocks satisfy
\begin{equation}
\sum_{\substack{
\Lambda_N\leq A\leq N^\beta\\
A\ {\rm dyadic}}}
\mathfrak M_{N,T;\mathfrak A}
\left(
f_{A,\bullet}
\right)
\longrightarrow0
\qquad
\text{as }N\to\infty,
\label{eq:active-ordinary-complete}
\end{equation}
and the logarithmic Abel blocks satisfy
\begin{equation}
\sum_{\substack{
\Lambda_N\leq A\leq N^\beta\\
A\ {\rm dyadic}}}
\mathfrak M_{N,T;\mathfrak A}
\left(
\widetilde f_{A,\bullet}
\right)
\longrightarrow0
\qquad
\text{as }N\to\infty.
\label{eq:active-log-complete}
\end{equation}
\end{theorem}
\begin{proof}
Subadditivity of \(\mathfrak M_{N,T;\mathfrak A}\) and
Theorems~\ref{thm:active-lower-sectors},
\ref{thm:active-middle-sectors}, and~\ref{thm:active-upper-sectors}
give both assertions.
\end{proof}

\section{Maximal estimates for Laguerre tails}
\label{sec:boundaries-fixed-tails}

In this section we show that keeping sufficiently many Laguerre modes leaves only a small error, uniformly over the polynomial lengths in the time window.
To this end, we establish fixed-cutoff maximal-tail estimates for both the ordinary and logarithmic Laguerre tails.
We first obtain bounds for low Laguerre blocks that are uniform in the coefficients.
These yield sharp-band estimates below the moving cutoff \(\Lambda_N\).
We then combine the Abel semigroup remainder estimate with the
Abel-to-sharp comparison to obtain moving sharp-tail bounds, and finally
pass to a fixed cutoff in Theorem~\ref{thm:fixed-cutoff-maximal-tails}.

The fixed-mode functional limit in
Theorem~\ref{thm:finite-laguerre-vector}, obtained from the functional
de Jong theorem of D\"obler--Kasprzak--Peccati
\cite{DoblerKasprzakPeccati2022}, is used with a fixed number of degrees.
It supplies no bounds uniform in a degree cutoff that grows with \(N\).
The estimates below provide that additional control, within the stated
range \(A\leq\Lambda_N\), before treating the ordinary and logarithmic
power tails.  Moreover, the maxima in this section are over nested
polynomial lengths.  The Laguerre semigroup acts on the intensity
variable and is an analytic device in these estimates; its usual
maximal inequalities do not by themselves estimate those length maxima.

Fix \(\mathfrak A=[\alpha_-,\alpha_+]\Subset(0,\infty)\) and \(T>0\),
and write \(\sigma_*=\min\{\alpha_-,1\}\).
Fix \(0<\beta<1\) sufficiently small that the estimates in
Theorem~\ref{thm:active-abel-range} hold.

The constant \(\gamma_0\) in \eqref{eq:moving-LambdaN} is fixed
sufficiently small for the low-degree estimates below and for
Proposition~\ref{cor:moving-upper-boundary-input}, with the auxiliary
parameters for the latter chosen in the order
\eqref{eq:moving-parameter-order}.
These choices are independent of \(N\) and apply to both families.

\begin{figure}[htbp]
  \centering
  \resizebox{.98\textwidth}{!}{%
\begin{tikzpicture}[
  x=1cm,y=1cm,
  toBox/.style={draw=black!65,line width=.45pt,rounded corners=1.3pt,
    align=flush center,inner xsep=5pt,inner ysep=4pt,
    font=\fontsize{8.2}{10}\selectfont,
    execute at begin node={\hyphenpenalty=10000\exhyphenpenalty=10000}},
  toEstimate/.style={toBox},
  toStep/.style={toBox},
  toConclusion/.style={toBox},
  toArrow/.style={-{Stealth[length=2.1mm,width=1.4mm]},
    draw=black!75,line width=.5pt}
]

\node[toStep,text width=14cm,minimum height=1.65cm]
(split) at (4.9,0) {
  {\fontsize{9}{11}\selectfont\bfseries
  Split the remainder at $\Lambda_N$}\\[3pt]
  Fix $d\ge2$; for large $N$, $d<\Lambda_N$.\\[3pt]
  $\displaystyle R_{\alpha,d}
      =\sum_{d<a\le\Lambda_N}\gamma_a(\alpha)L_a
         +R_{\alpha,\Lambda_N}.$
};

\node[toEstimate,text width=4.2cm,minimum height=1.5cm]
(band) at (0,-4.15) {
  {\fontsize{9}{11}\selectfont\bfseries
  Estimate the terms\\between $d$ and $\Lambda_N$}\\[3pt]
  {\fontsize{7.5}{9.5}\selectfont
  Theorem~\ref{thm:low-moving-sharp-band}}
};

\node[toStep,text width=8.8cm,minimum height=1.55cm]
(splitupper) at (7.5,-2.15) {
  {\fontsize{9}{11}\selectfont\bfseries
  Split the tail above $\Lambda_N$}\\[4pt]
  $\displaystyle R_{\alpha,\Lambda_N}
      =(I-P_{1/\Lambda_N})h_\alpha
       +(P_{1/\Lambda_N}-\Pi^{\mathrm{Lag}}_{\le\Lambda_N})h_\alpha.$
};

\node[toEstimate,text width=3.9cm,minimum height=1.75cm]
(smoothing) at (5.1,-4.15) {
  {\fontsize{9}{11}\selectfont\bfseries Estimate for}\\[3pt]
  $(I-P_{1/\Lambda_N})h_\alpha$\\[3pt]
  {\fontsize{7.5}{9.5}\selectfont
  Theorem~\ref{thm:abel-high-pass}}
};

\node[toEstimate,text width=3.9cm,minimum height=1.75cm]
(comparison) at (9.9,-4.15) {
  {\fontsize{9}{11}\selectfont\bfseries Estimate for}\\[3pt]
  $(P_{1/\Lambda_N}-\Pi^{\mathrm{Lag}}_{\le\Lambda_N})h_\alpha$\\[3pt]
  {\fontsize{7.5}{9.5}\selectfont
  Theorem~\ref{thm:abel-sharp-boundary-comparison}}
};

\node[toStep,text width=8.8cm,minimum height=1.05cm]
(upper) at (7.5,-6.1) {
  {\fontsize{9}{11}\selectfont\bfseries
  Tail above $\Lambda_N$}\\[3pt]
  {\fontsize{7.5}{9.5}\selectfont
  Corollary~\ref{cor:moving-sharp-Laguerre-tails}}
};

\node[toConclusion,text width=14cm,minimum height=1.05cm]
(fixed) at (4.9,-7.85) {
  {\fontsize{9}{11}\selectfont\bfseries
  Fixed-cutoff maximal tails}\\[3pt]
  {\fontsize{7.5}{9.5}\selectfont
  Theorem~\ref{thm:fixed-cutoff-maximal-tails}}
};

\draw[toArrow]
  ([xshift=-4.9cm]split.south) -- (band.north);

\draw[toArrow]
  ([xshift=2.6cm]split.south) -- (splitupper.north);

\draw[toArrow]
  ([xshift=-2.4cm]splitupper.south) -- (smoothing.north);

\draw[toArrow]
  ([xshift=2.4cm]splitupper.south) -- (comparison.north);

\draw[toArrow]
  (smoothing.south) -- ([xshift=-2.4cm]upper.north);

\draw[toArrow]
  (comparison.south) -- ([xshift=2.4cm]upper.north);

\draw[toArrow]
  (band.south) -- ([xshift=-4.9cm]fixed.north);

\draw[toArrow]
  (upper.south) -- ([xshift=2.6cm]fixed.north);

\end{tikzpicture}
}
  \caption{Proof structure for
  Theorem~\ref{thm:fixed-cutoff-maximal-tails}.
  Split the remainder, estimate each part, and combine the bounds.
  The argument treats both $x^\alpha$ and $x^\alpha\log x$.
  The bounds concern maxima over the length window, integrated over
  $\alpha$, with $d$ fixed as $N\to\infty$.}
  \label{fig:fixed-cutoff-tail-roadmap}
\end{figure}
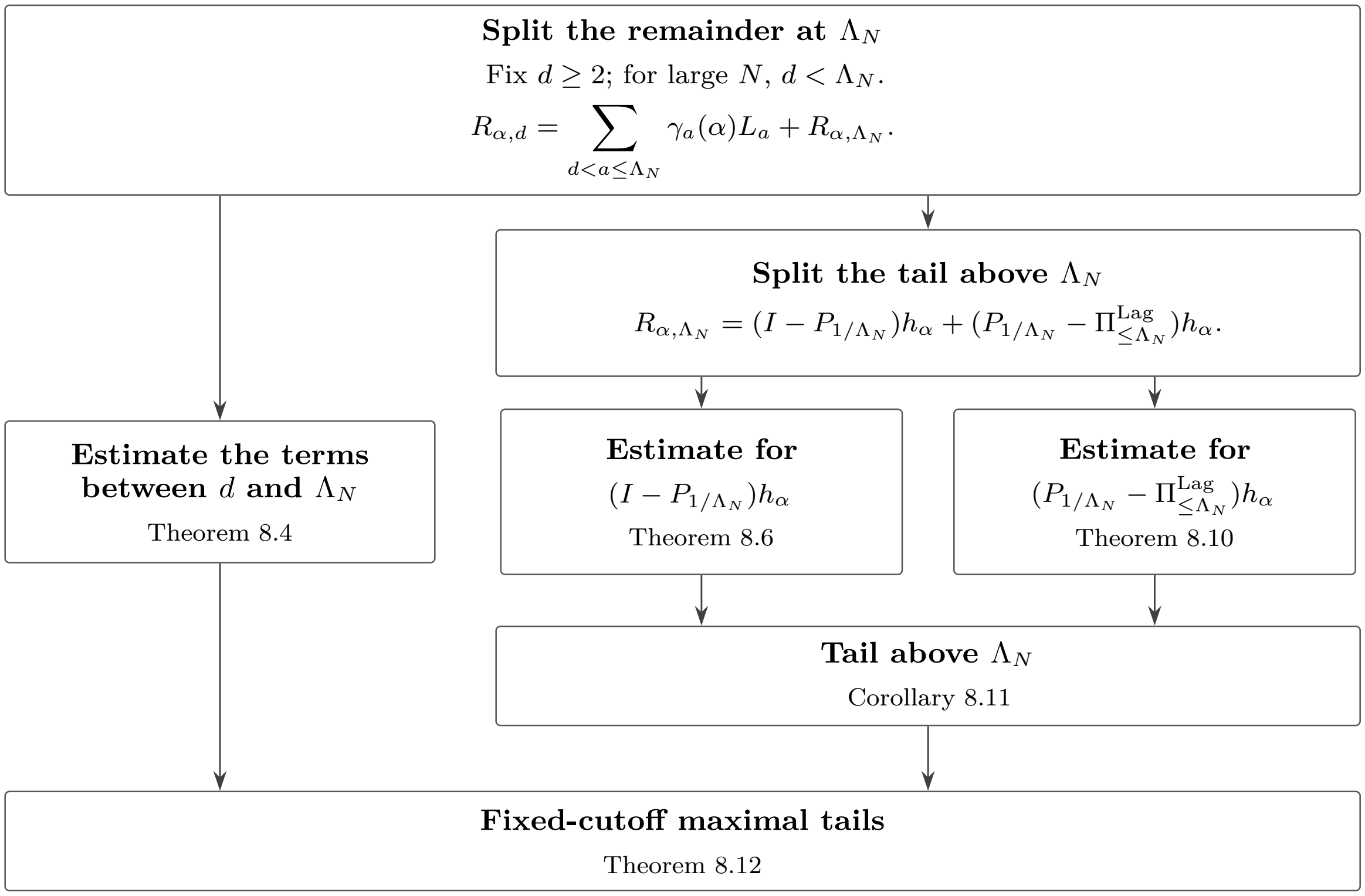

\subsection{Coefficient-uniform dyadic Laguerre blocks}
\label{subsec:coefficient-uniform-low-blocks}

The purpose of this subsection is to obtain a coefficient-uniform maximal estimate for dyadic Laguerre blocks below the moving cutoff \(\Lambda_N\).
Let \(A\geq1\) be dyadic and let
\[
\mathbf u=(u_a)_{A<a\leq2A}
\]
be a real coefficient vector.  Define
\begin{equation}
g_{A,\mathbf u}
=
\sum_{A<a\leq2A}u_aL_a
\label{eq:low-block-gAu}
\end{equation}
and its coefficient energy by
\begin{equation}
\mathcal E_A(\mathbf u)
=
\sum_{A<a\leq2A}a|u_a|^2.
\label{eq:low-block-energy}
\end{equation}
No regularity, monotonicity, sign, or special-function structure is imposed on
\(\mathbf u\).

For \(a\geq2\), recall from \eqref{eq:chaos-top-chaos} that
\(\mathcal G_{n,a}\) is the component of Walsh degree \(2a\) in
\(X_n(L_a)\).
Define
\begin{equation}
\mathcal T_n(g_{A,\mathbf u})
=
\sum_{A<a\leq2A}
u_a\mathcal G_{n,a}
\label{eq:low-block-top-chaos}
\end{equation}
and the extraction remainder
\begin{equation}
\mathcal R_n(g_{A,\mathbf u})
=
X_n(g_{A,\mathbf u})
-
\mathcal T_n(g_{A,\mathbf u}).
\label{eq:low-block-extraction-remainder}
\end{equation}
Thus \(\mathcal T_n(g_{A,\mathbf u})\) keeps the part of highest Walsh degree from each Laguerre mode, with the original weights \(u_a\).

Theorem~\ref{thm:fixed-mode-pure-chaos-reduction} provides an explicit
bound in the Laguerre degree for the single-mode extraction error.
The following proposition extends that bound to arbitrary coefficient
vectors in a dyadic block and then to the maximum over the
logarithmic-time window.

\begin{proposition}
\label{prop:low-block-degree-sensitive-extraction}
There exists \(c_T>0\) such that, for all sufficiently large \(N\),
uniformly for every dyadic
\(
1\leq A\leq\Lambda_N
\)
and every deterministic real coefficient vector
\(\mathbf u=(u_a)_{A<a\leq2A}\),
\begin{equation}
\E
\max_{\lfloor Ne^{-T}\rfloor\leq n\leq\lfloor Ne^T\rfloor}
\left|\mathcal R_n(g_{A,\mathbf u})\right|^2
\leq
N^{-c_T}\mathcal E_A(\mathbf u).
\label{eq:low-block-maximal-extraction}
\end{equation}
\end{proposition}

\begin{proof}
Let \(\mathcal R_{n,a}^{\mathrm{pc}}\) be the single-mode extraction
remainder defined in \eqref{eq:chaos-pure-chaos-remainder}.
By \eqref{eq:chaos-pure-chaos-L2}, uniformly for \(A<a\leq2A\),
\[
\|\mathcal R_{n,a}^{\mathrm{pc}}\|_{L^2}^2
\leq
\frac{\exp\{CA\log(A+2)\}}{n}.
\]
Since
\[
\mathcal R_n(g_{A,\mathbf u})
=
\sum_{A<a\leq2A}u_a\mathcal R_{n,a}^{\mathrm{pc}},
\]
weighted Cauchy--Schwarz gives
\begin{equation}
\|\mathcal R_n(g_{A,\mathbf u})\|_{L^2}^2
\leq
\frac{\exp\{CA\log(A+2)\}}{n}
\mathcal E_A(\mathbf u).
\label{eq:block-extraction-L2-proof}
\end{equation}

Moreover,
\[
\deg \mathcal R_n(g_{A,\mathbf u})\leq4A
\]
by
\eqref{eq:chaos-Laguerre-power-expansion},
\eqref{eq:chaos-even-Walsh-degrees}, and
\eqref{eq:chaos-top-projection}.
Hence Lemma~\ref{lem:mixed-hypercontractivity}, followed by
Lemma~\ref{lem:block-finite-family-maximum}, yields for every fixed
integer \(q\geq2\)
\begin{equation}
\begin{aligned}
\E
\max_{\lfloor Ne^{-T}\rfloor\leq n\leq\lfloor Ne^T\rfloor}
|\mathcal R_n(g_{A,\mathbf u})|^2
\leq
C_T N^{-1+1/q}
\exp\!\left\{
CA\log(A+2)+4A\log(2q-1)
\right\}
\mathcal E_A(\mathbf u).
\end{aligned}
\label{eq:block-extraction-max-precutoff}
\end{equation}

For \(A\leq\Lambda_N\),
\[
A\log(A+2)\leq C\gamma_0\log N,
\qquad
A=o(\log N).
\]
Taking \(q=4\), and using the fixed sufficiently small choice of
\(\gamma_0\), the right-hand side of
\eqref{eq:block-extraction-max-precutoff} is bounded by
\(N^{-c_T}\mathcal E_A(\mathbf u)\), after decreasing \(c_T>0\) if
necessary. This proves the proposition.
\end{proof}

It remains to control the top-chaos term. Only second moments are
needed for this purpose. By
\eqref{eq:chaos-top-second-moment-degree-sensitive} and
\eqref{eq:gaussian-Va-bound},
\[
\sup_{A<a\leq2A}\E\mathcal G_{n,a}^2
\leq 2+\frac{\exp\{CA\log(A+2)\}}{n}\leq C_T
\]
uniformly for \(1\leq A\leq\Lambda_N\) and \(n\asymp_TN\),
for all sufficiently large \(N\).
Indeed, for large \(N\), the dyadic definition
\eqref{eq:moving-LambdaN} gives
\[
\Lambda_N\leq
2\gamma_0\frac{\log(N+2)}{\log\log(N+3)},
\qquad
A\log(A+2)\leq C\gamma_0\log N,
\]
and the fixed choice of \(\gamma_0\) absorbs the absolute constants,
including the factor from \(a\leq2A\).
The variables \(\mathcal G_{n,a}\) for distinct \(a\) lie in distinct
homogeneous Walsh degrees and are orthogonal. Consequently,
\begin{equation}
\|\mathcal T_n(g_{A,\mathbf u})\|_{L^2}^2
=\sum_{A<a\leq2A}|u_a|^2\E\mathcal G_{n,a}^2
\leq C_T\sum_{A<a\leq2A}|u_a|^2.
\label{eq:low-block-unweighted-endpoint}
\end{equation}
We retain the unweighted coefficient norm until the logarithmic-time
window has been covered.

The martingale identity in
Lemma~\ref{lem:low-block-damped-martingale} passes this fixed-length bound
to a maximum on each short logarithmic interval.

\begin{lemma}
\label{lem:low-block-short-interval}
For all sufficiently large \(N\), let
\(1\leq A\leq\Lambda_N\) be dyadic and let
\(\mathbf u=(u_a)_{A<a\leq2A}\) be a deterministic real coefficient
vector. Let
\[
I=[s_0,s_1]\subseteq I_T,
\qquad 0<s_1-s_0\leq A^{-1},
\qquad
\mathcal N_N(I)=\{n\in\mathbb N:\log(n/N)\in I\}.
\]
Then
\begin{equation}
\E\max_{n\in\mathcal N_N(I)}
|\mathcal T_n(g_{A,\mathbf u})|^2
\leq C_T\sum_{A<a\leq2A}|u_a|^2.
\label{eq:low-block-short-interval}
\end{equation}
The maximum over an empty set is zero. The constant and the threshold
for \(N\) are uniform in \(A,I,\mathbf u\); the threshold depends only
on \(T\) and the fixed cutoff constant \(\gamma_0\).
\end{lemma}

\begin{proof}
The empty case is immediate, and the singleton case follows from
\eqref{eq:low-block-unweighted-endpoint}. Otherwise, write
\[
m_I=\min\mathcal N_N(I)<n_I=\max\mathcal N_N(I),
\qquad
\|\mathbf u\|_{\ell^2}^2=\sum_{A<a\leq2A}|u_a|^2.
\]
By Lemma~\ref{lem:low-block-damped-martingale},
\[
M_{n,a}=n^{a-\frac12}\mathcal G_{n,a}
\]
is a martingale. Moreover,
\begin{equation}
\frac{n_I}{m_I}\leq e^{1/A},
\qquad
\sum_{k=m_I+1}^{n_I}\frac Ak
\leq A\log\frac{n_I}{m_I}\leq1.
\label{eq:low-block-short-damping}
\end{equation}
For \(m_I\leq n\leq n_I\), discrete summation by parts gives
\[
\mathcal T_n(g_{A,\mathbf u})
=\mathcal T_{m_I}(g_{A,\mathbf u})+H_n+V_n,
\]
where
\[
\begin{aligned}
H_n&=\sum_{A<a\leq2A}u_a
\sum_{k=m_I+1}^{n}k^{-(a-\frac12)}(M_{k,a}-M_{k-1,a}),\\
V_n&=\sum_{k=m_I+1}^{n}\sum_{A<a\leq2A}u_a
\bigl(k^{-(a-\frac12)}-(k-1)^{-(a-\frac12)}\bigr)M_{k-1,a}.
\end{aligned}
\]
The process \(H_n\) is a martingale. Set
\[
u_a^{(k)}=u_a\left[\left(\frac{k-1}{k}\right)^{a-\frac12}-1\right].
\]
Then
\[
V_n=\sum_{k=m_I+1}^{n}\mathcal T_{k-1}(g_{A,\mathbf u^{(k)}}),
\qquad
\|\mathbf u^{(k)}\|_{\ell^2}
\leq\frac{2A}{k}\|\mathbf u\|_{\ell^2}.
\]
The coefficient vectors are deterministic. Applying Minkowski's
inequality, \eqref{eq:low-block-unweighted-endpoint}, and
\eqref{eq:low-block-short-damping}, we obtain
\[
\left\|\max_{m_I\leq n\leq n_I}|V_n|\right\|_{L^2}
\leq C_T\|\mathbf u\|_{\ell^2}
\sum_{k=m_I+1}^{n_I}\frac{2A}{k}
\leq C_T\|\mathbf u\|_{\ell^2}.
\]
The decomposition at \(n=n_I\), together with the fixed-length bound,
therefore yields
\[
\|H_{n_I}\|_{L^2}
\leq\|\mathcal T_{n_I}(g_{A,\mathbf u})\|_{L^2}
+\|\mathcal T_{m_I}(g_{A,\mathbf u})\|_{L^2}
+\|V_{n_I}\|_{L^2}
\leq C_T\|\mathbf u\|_{\ell^2}.
\]
Doob's \(L^2\) inequality now gives
\[
\left\|\max_{m_I\leq n\leq n_I}|H_n|\right\|_{L^2}
\leq2\|H_{n_I}\|_{L^2}\leq C_T\|\mathbf u\|_{\ell^2}.
\]
Combining these estimates proves the lemma.
\end{proof}

We combine the short-interval bounds by summing their second moments,
and then add the extraction error.

\begin{theorem}[Coefficient-uniform low-block maximum]
\label{thm:coefficient-uniform-low-block}
There exists \(C_T>0\) such that, for all sufficiently large \(N\),
uniformly over every dyadic \(1\leq A\leq\Lambda_N\), every subset
\[
S\subseteq(A,2A]\cap\mathbb N,
\]
and every real coefficient family \((u_a)_{a\in S}\),
\begin{equation}
\begin{aligned}
\left(\E\max_{n_-\leq n\leq n_+}
\left|X_n\!\left(\sum_{a\in S}u_aL_a\right)\right|^2\right)^{1/2}
\leq C_T\left(\sum_{a\in S}a|u_a|^2\right)^{1/2}.
\end{aligned}
\label{eq:coefficient-uniform-low-partial-block}
\end{equation}
\end{theorem}

\begin{proof}
Extend the coefficients by zero to all \(A<a\leq2A\). Partition
\(I_T\) into \(J\leq C_TA\) closed intervals \(I_\ell\) of positive
length at most \(A^{-1}\). Since
\[
\{n_-,\ldots,n_+\}
=\{n_-\}\cup\bigcup_{\ell=1}^{J}\mathcal N_N(I_\ell),
\]
the inequality \(\max_j z_j^2\leq\sum_jz_j^2\),
\eqref{eq:low-block-unweighted-endpoint}, and
Lemma~\ref{lem:low-block-short-interval} give
\begin{equation}
\begin{aligned}
\E\max_{n_-\leq n\leq n_+}|\mathcal T_n(g_{A,\mathbf u})|^2
&\leq\E|\mathcal T_{n_-}(g_{A,\mathbf u})|^2
+\sum_{\ell=1}^{J}\E\max_{n\in\mathcal N_N(I_\ell)}
|\mathcal T_n(g_{A,\mathbf u})|^2\\
&\leq C_T(1+J)\sum_{A<a\leq2A}|u_a|^2
\leq C_T\mathcal E_A(\mathbf u).
\end{aligned}
\label{eq:low-block-top-chaos-prebound}
\end{equation}
The separate endpoint includes the possible value
\(n_-<Ne^{-T}\). The argument also covers \(A=1\) and empty or
singleton sets \(\mathcal N_N(I_\ell)\).
By Proposition~\ref{prop:low-block-degree-sensitive-extraction},
\[
\left\|\max_{n_-\leq n\leq n_+}
|\mathcal R_n(g_{A,\mathbf u})|\right\|_{L^2}
\leq N^{-c_T/2}\mathcal E_A(\mathbf u)^{1/2}.
\]
Finally, \(X_n(g_{A,\mathbf u})=\mathcal T_n(g_{A,\mathbf u})
+\mathcal R_n(g_{A,\mathbf u})\), so the triangle inequality and
\eqref{eq:low-block-top-chaos-prebound} prove the theorem.
\end{proof}

\subsection{Laguerre sums between fixed and moving cutoffs}
\label{subsec:low-moving-sharp-bands}

We now apply the coefficient-uniform estimate from Theorem~\ref{thm:coefficient-uniform-low-block} to the ordinary and logarithmic Laguerre coefficients between a fixed cutoff \(d\) and the moving cutoff \(\Lambda_N\).

Recall from \eqref{eq:def-halpha} and \eqref{eq:def-widetilde-halpha} that
\[
h_\alpha(x)=x^\alpha,
\qquad
\widetilde h_\alpha(x)=x^\alpha\log x.
\]
Here and below, fix \(\mathfrak A=[\alpha_-,\alpha_+]\Subset(0,\infty)\), and 
\[
h_\bullet=(h_\alpha)_{\alpha\in\mathfrak A},
\qquad
\widetilde h_\bullet=(\widetilde h_\alpha)_{\alpha\in\mathfrak A},
\]
with the bullet indicating the \(\alpha\)-indexed family.
We use the maximal seminorm
\(\mathfrak M_{N,T;\mathfrak A}\) from
\eqref{eq:active-integrated-maximal-seminorm}
and the sharp Laguerre projections from
\eqref{eq:spectral-sharp-projection}.

\begin{theorem}
\label{thm:low-moving-sharp-band}
For every fixed integer \(d\geq2\),
\begin{equation}
\begin{aligned}
\limsup_{N\to\infty}
\mathfrak M_{N,T;\mathfrak A}
\left(
(\LagProj_{\leq\Lambda_N}-\LagProj_{\leq d})h_\bullet
\right)
\leq
C_{\mathfrak A,T}
d^{-\alpha_-}
\sqrt{\log(d+2)}.
\end{aligned}
\label{eq:low-band-ordinary-rate}
\end{equation}
Likewise,
\begin{equation}
\begin{aligned}
\limsup_{N\to\infty}
\mathfrak M_{N,T;\mathfrak A}
\left(
(\LagProj_{\leq\Lambda_N}-\LagProj_{\leq d})
\widetilde h_\bullet
\right)
\leq
C_{\mathfrak A,T}
d^{-\alpha_-}
\log^{3/2}(d+2).
\end{aligned}
\label{eq:low-band-logarithmic-rate}
\end{equation}
\end{theorem}

\begin{proof}
We prove \eqref{eq:low-band-ordinary-rate}.
The logarithmic estimate \eqref{eq:low-band-logarithmic-rate}
follows by the same argument, using
\(\widetilde\gamma_a(\alpha)\) and the corresponding logarithmic
coefficient bound from Proposition~\ref{prop:sharp-coefficient-energy}.

Fix \(d\geq2\).
Since \(\Lambda_N\to\infty\), for all sufficiently large \(N\) one has
\(d<\Lambda_N\).
Let \(A_d\) be the unique dyadic integer such that
\[
A_d\leq d<2A_d.
\]
For such \(N\),
\[
(\LagProj_{\leq\Lambda_N}-\LagProj_{\leq d})h_\alpha
=
\sum_{d<a\leq\Lambda_N}\gamma_a(\alpha)L_a.
\]
Partition \((d,\Lambda_N]\cap\mathbb N\) into its intersections with
the dyadic degree blocks \((A,2A]\), where
\(A_d\leq A<\Lambda_N\).
Since Theorem~\ref{thm:coefficient-uniform-low-block} allows arbitrary
subsets of a dyadic block, it applies in particular to the initial
partial block.

By subadditivity of \(\mathfrak M_{N,T;\mathfrak A}\),
Theorem~\ref{thm:coefficient-uniform-low-block}, applied to each dyadic
block, and Proposition~\ref{prop:sharp-coefficient-energy},
\[
\begin{aligned}
\mathfrak M_{N,T;\mathfrak A}
\left(
(\LagProj_{\leq\Lambda_N}-\LagProj_{\leq d})h_\bullet
\right)
&\leq
C_{\mathfrak A,T}
\sum_{\substack{
A_d\leq A<\Lambda_N\\
A\ {\rm dyadic}}}
A^{-\alpha_-}\sqrt{\log(A+2)}
\\
&\leq
C_{\mathfrak A,T}
d^{-\alpha_-}\sqrt{\log(d+2)},
\end{aligned}
\]
where the last estimate follows from standard dyadic summation and
\(A_d\asymp d\).
Taking \(\limsup_{N\to\infty}\) proves
\eqref{eq:low-band-ordinary-rate}.
\end{proof}

\subsection{Global Abel blocks and semigroup remainders}
\label{subsec:global-Abel-high-pass}

This subsection estimates the difference between each ordinary or
logarithmic power function and its Abel smoothing at the moving cutoff
\(\Lambda_N\).  The difference is negligible in
\(\mathfrak M_{N,T;\mathfrak A}\), as stated in
Theorem~\ref{thm:abel-high-pass}.

We now return to the Abel-scale notation of
Sections~\ref{sec:coefficient-calculus} and~\ref{sec:active-abel-range}:
\(A\geq2\) denotes the Abel parameter in
\[
f_{A,\alpha}
=
\left(P_{1/(2A)}-P_{1/A}\right)h_\alpha,
\qquad
\widetilde f_{A,\alpha}
=
\left(P_{1/(2A)}-P_{1/A}\right)\widetilde h_\alpha,
\]
with dyadic restrictions on \(A\) stated explicitly when needed.
We continue to use
\[
\sigma_*=\min\{\alpha_-,1\},
\]
as fixed at the beginning of this section.

\begin{proposition}[Moment bounds for Abel blocks]
\label{prop:global-Abel-endpoint}
Let \(q\geq4\) be fixed and even, and let
\(
\mathfrak A=[\alpha_-,\alpha_+]\Subset(0,\infty).
\)
Uniformly for \(n\geq2\) and \(A\geq2\),
\begin{align}
\sup_{\alpha\in\mathfrak A}
\|X_n(f_{A,\alpha})\|_{L^q}
&\leq
C_{\mathfrak A,q}A^{-\sigma_*},
\label{eq:global-ordinary-Abel-endpoint}
\\
\sup_{\alpha\in\mathfrak A}
\|X_n(\widetilde f_{A,\alpha})\|_{L^q}
&\leq
C_{\mathfrak A,q}
A^{-\sigma_*}\log(A+2).
\label{eq:global-log-Abel-endpoint}
\end{align}
\end{proposition}

\begin{proof}
The derivative estimates must remain stable under the changes of size
\(n^{-1/2}\) in the normalized polynomial field caused by a sign flip.  When \(A>n\),
we therefore separate intensities of order \(n^{-1}\), which exceed
the Abel smoothing scale \(A^{-1}\).  Set
\[
s_{A,n}:=A^{-1}\vee n^{-1}.
\]
Choose
\[
\chi\in C^\infty([0,\infty);[0,1]),
\qquad
\chi=1\ \text{on }[0,1],
\qquad
\chi=0\ \text{on }[2,\infty),
\]
and define
\[
\chi_{A,n}(x)
:=
\begin{cases}
0, & A\leq n,\\[1mm]
\chi(x/s_{A,n}), & A>n.
\end{cases}
\]
For fixed \(\alpha\in\mathfrak A\), set
\begin{align*}
B_0(x)
&=
(x+s_{A,n})^{\alpha-1}
+
|\alpha-1|(x+s_{A,n})^{\alpha-2},
\\
B_1(x)
&=
(x+s_{A,n})^{\alpha-1}
\bigl(1+|\log(x+s_{A,n})|\bigr)
\\
&\quad+
(x+s_{A,n})^{\alpha-2}
\Bigl[
1+
|\alpha-1|
\bigl(1+|\log(x+s_{A,n})|\bigr)
\Bigr],
\\
W_0(x)
&=
(x+s_{A,n})^{2\alpha-1}
+
|\alpha-1|^2(x+s_{A,n})^{2\alpha-3},
\\
W_1(x)
&=
(x+s_{A,n})^{2\alpha-1}
\bigl(1+|\log(x+s_{A,n})|\bigr)^2
\\
&\quad+
(x+s_{A,n})^{2\alpha-3}
\Bigl[
1+
|\alpha-1|
\bigl(1+|\log(x+s_{A,n})|\bigr)
\Bigr]^2.
\end{align*}
The functions \(B_0,B_1\) will bound derivatives of the ordinary and
logarithmic blocks, respectively; \(W_0,W_1\) will bound the corresponding
squared derivatives multiplied by \(x\).
All estimates below are uniform in \(\alpha\in\mathfrak A\).
Lemma~\ref{lem:gae-abel-derivatives} and the product rule yield, for
\(j=1,2\),
\begin{align}
\left|
\partial_x^j
\bigl[(1-\chi_{A,n})f_{A,\alpha}\bigr](x)
\right|
&\leq
CA^{-1}
\left[
(x+s_{A,n})^{1-j}B_0(x)
+
s_{A,n}^{\alpha-1-j}
\mathbf 1_{\{s_{A,n}\leq x\leq2s_{A,n}\}}
\right],
\label{eq:gae-cutoff-ordinary}
\\
\left|
\partial_x^j
\bigl[(1-\chi_{A,n})\widetilde f_{A,\alpha}\bigr](x)
\right|
&\leq
CA^{-1}
\left[
(x+s_{A,n})^{1-j}B_1(x)
\right.
\nonumber\\
&\quad\left.
+
s_{A,n}^{\alpha-1-j}
\bigl(1+\log s_{A,n}^{-1}\bigr)
\mathbf 1_{\{s_{A,n}\leq x\leq2s_{A,n}\}}
\right].
\label{eq:gae-cutoff-log}
\end{align}
When \(\chi_{A,n}\not\equiv0\), one has
\(A^{-1}<n^{-1}=s_{A,n}\), and hence
\[
x+A^{-1}\asymp s_{A,n}
\qquad
\text{for }s_{A,n}\leq x\leq2s_{A,n}.
\]
Consequently,
\begin{equation}
\begin{aligned}
x\left|
\partial_x\bigl[(1-\chi_{A,n})f_{A,\alpha}\bigr](x)
\right|^2
&\leq
CA^{-2}
\left[
W_0(x)
+
s_{A,n}^{2\alpha-3}
\mathbf 1_{\{s_{A,n}\leq x\leq2s_{A,n}\}}
\right],
\\
x\left|
\partial_x\bigl[(1-\chi_{A,n})\widetilde f_{A,\alpha}\bigr](x)
\right|^2
&\leq
CA^{-2}
\left[
W_1(x)
\right.
\\
&\quad\left.
+
s_{A,n}^{2\alpha-3}
\bigl(1+\log s_{A,n}^{-1}\bigr)^2
\mathbf 1_{\{s_{A,n}\leq x\leq2s_{A,n}\}}
\right].
\end{aligned}
\label{eq:gae-square-envelope}
\end{equation}
For \(f\in C^2([0,\infty))\), write
\begin{equation}
H_{n,f}(x)
=
\sup_{\substack{
        r\geq0\\
        |r-\sqrt x|\leq2n^{-1/2}}}
\left(
|f'(r^2)|
+
2r^2|f''(r^2)|
\right).
\label{eq:gae-local-jacobian-envelope}
\end{equation}
If
\(
|r-\sqrt x|\leq2n^{-1/2},
\)
then
\[
r^2+s_{A,n}
\asymp
x+s_{A,n}.
\]
Moreover,
\[
s_{A,n}\leq r^2\leq2s_{A,n}
\quad\Longrightarrow\quad
x\leq16s_{A,n}.
\]
Hence
\begin{equation}
\begin{aligned}
H_{n,(1-\chi_{A,n})f_{A,\alpha}}(x)
&\leq
CA^{-1}
\left[
B_0(x)
+
s_{A,n}^{\alpha-2}
\mathbf 1_{\{x\leq16s_{A,n}\}}
\right],
\\
H_{n,(1-\chi_{A,n})\widetilde f_{A,\alpha}}(x)
&\leq
CA^{-1}
\left[
B_1(x)
+
s_{A,n}^{\alpha-2}
\bigl(1+\log s_{A,n}^{-1}\bigr)
\mathbf 1_{\{x\leq16s_{A,n}\}}
\right].
\end{aligned}
\label{eq:gae-jacobian-envelope}
\end{equation}
Since \(\alpha\) ranges over the compact interval \(\mathfrak A\),
the powers appearing in \(B_i\) and \(W_i\), \(i\in\{0,1\}\), range
over fixed compact intervals.
Proposition~\ref{prop:gae-occupation} and
\eqref{eq:gae-weighted-occupation} therefore imply
\begin{align}
\|I_n(B_i)\|_{L^q}
&\leq C
\left(1+\int_0^1B_i(x)\,\dd x
\right),
\label{eq:gae-B-moment}\\
\|I_n(W_i)\|_{L^{q/2}}
&\leq C
\left(1+\int_0^1W_i(x)\,\dd x
\right).
\label{eq:gae-W-moment}
\end{align}
The same proposition yields
\begin{align}
\left\|I_n\!\left(
s_{A,n}^{2\alpha-3}(1+\log(s_{A,n}^{-1}))^{2i}
\mathbf1_{\{s_{A,n}\leq x\leq2s_{A,n}\}}
\right)
\right\|_{L^{q/2}}^{1/2}
&\leq C s_{A,n}^{\alpha-1}(1+\log(s_{A,n}^{-1}))^i,
\label{eq:gae-cutoff-square-moment}\\
\left\|I_n\!\left(
s_{A,n}^{\alpha-2}(1+\log(s_{A,n}^{-1}))^i
\mathbf1_{\{x\leq16 s_{A,n}\}}
\right)
\right\|_{L^q}
&\leq C s_{A,n}^{\alpha-1}(1+\log(s_{A,n}^{-1}))^i.
\label{eq:gae-cutoff-jacobian-moment}
\end{align}
If \(\alpha_-<1\), direct integration yields
\begin{equation}
\int_0^1B_i(x)\,\dd x
\leq
Cs_{A,n}^{\alpha_--1}
\bigl(1+\log s_{A,n}^{-1}\bigr)^i,
\qquad
\int_0^1W_i(x)\,\dd x
\leq
Cs_{A,n}^{2\alpha_--2}
\bigl(1+\log s_{A,n}^{-1}\bigr)^{2i}.
\label{eq:gae-subcritical-integrals}
\end{equation}
Here the singular terms are controlled by
\[
\int_{s_{A,n}}^1
u^{-1-c}(1+|\log u|)^k\,\dd u
\leq
C_c s_{A,n}^{-c}
\bigl(1+\log s_{A,n}^{-1}\bigr)^k,
\qquad
c>0,
\quad
k\leq2.
\]
Suppose now that \(\alpha_-\geq1\).
Then \(\alpha-1\geq0\), and
\[
(\alpha-1)
\int_{s_{A,n}}^1
u^{\alpha-2}\,\dd u
\leq1,
\qquad
(\alpha-1)^2
\int_{s_{A,n}}^1
u^{2\alpha-3}\,\dd u
\leq
\frac{\alpha-1}{2}.
\]
Thus
\begin{equation}
\int_0^1B_0(x)\,\dd x
+
\int_0^1W_0(x)\,\dd x
\leq C.
\label{eq:gae-critical-ordinary-integrals}
\end{equation}
For \(0<u\leq1\),
\begin{equation}
u^{\alpha-1}
\Bigl[
1+
(\alpha-1)(1+|\log u|)
\Bigr]
\leq
C_{\mathfrak A}.
\label{eq:gae-critical-cancellation}
\end{equation}
Thus the singular integrands in the logarithmic case satisfy
\[
u^{\alpha-2}
\Bigl[
1+
(\alpha-1)(1+|\log u|)
\Bigr]
+
u^{2\alpha-3}
\Bigl[
1+
(\alpha-1)(1+|\log u|)
\Bigr]^2
\leq
Cu^{-1}.
\]
The remaining terms are uniformly integrable, so
\begin{equation}
\int_0^1B_1(x)\,\dd x
+
\int_0^1W_1(x)\,\dd x
\leq
C\bigl(1+\log s_{A,n}^{-1}\bigr).
\label{eq:gae-critical-log-integrals}
\end{equation}
Lemma~\ref{lem:gae-sign-flip}, together with
\eqref{eq:gae-square-envelope}--\eqref{eq:gae-critical-log-integrals},
now yields
\begin{equation}
\begin{aligned}
\sup_{\alpha\in\mathfrak A}
\left\|
X_n\!\left((1-\chi_{A,n})f_{A,\alpha}\right)
\right\|_{L^q}
&\leq
CA^{-1}s_{A,n}^{\sigma_*-1},
\\
\sup_{\alpha\in\mathfrak A}
\left\|
X_n\!\left((1-\chi_{A,n})\widetilde f_{A,\alpha}\right)
\right\|_{L^q}
&\leq
CA^{-1}s_{A,n}^{\sigma_*-1}
\bigl(1+\log s_{A,n}^{-1}\bigr).
\end{aligned}
\label{eq:gae-regular-part}
\end{equation}

It remains to estimate the small-\(x\) parts.
If \(A\leq n\), then \(\chi_{A,n}\equiv0\).
If \(A>n\), then
\[
s_{A,n}=n^{-1},
\qquad
\operatorname{supp}\chi_{A,n}
\subseteq[0,2s_{A,n}].
\]
By
\eqref{eq:gae-value-ordinary}--\eqref{eq:gae-value-log},
\[
\begin{aligned}
\sup_{\alpha\in\mathfrak A}
\|\chi_{A,n}f_{A,\alpha}\|_\infty
&\leq
CA^{-\sigma_*},
\\
\sup_{\alpha\in\mathfrak A}
\|\chi_{A,n}\widetilde f_{A,\alpha}\|_\infty
&\leq
CA^{-\sigma_*}(1+\log A).
\end{aligned}
\]
Since
\[
\|X_n(h)\|_{L^q}
\leq
2\sqrt n\,\|I_n(|h|)\|_{L^q},
\]
Proposition~\ref{prop:gae-occupation}, applied with \(u=2n^{-1}\), yields
\begin{equation}
\begin{aligned}
\sup_{\alpha\in\mathfrak A}
\left\|
X_n(\chi_{A,n}f_{A,\alpha})
\right\|_{L^q}
&\leq
Cn^{-1/2}A^{-\sigma_*},
\\
\sup_{\alpha\in\mathfrak A}
\left\|
X_n(\chi_{A,n}\widetilde f_{A,\alpha})
\right\|_{L^q}
&\leq
Cn^{-1/2}A^{-\sigma_*}(1+\log A).
\end{aligned}
\label{eq:gae-microscopic-part}
\end{equation}

Since
\[
A^{-1}s_{A,n}^{\sigma_*-1}
\leq
A^{-\sigma_*},
\qquad
1+\log s_{A,n}^{-1}
\leq
1+\log A,
\]
and
\[
f_{A,\alpha}
=
(1-\chi_{A,n})f_{A,\alpha}
+
\chi_{A,n}f_{A,\alpha},
\]
with the analogous decomposition for
\(\widetilde f_{A,\alpha}\),
\eqref{eq:gae-regular-part} and
\eqref{eq:gae-microscopic-part} imply
\eqref{eq:global-ordinary-Abel-endpoint} and
\eqref{eq:global-log-Abel-endpoint}.
\end{proof}

\begin{theorem}
\label{thm:abel-high-pass}
For the ordinary family,
\begin{equation}
\mathfrak M_{N,T;\mathfrak A}
\left(
(I-P_{1/\Lambda_N})h_\bullet
\right)
\longrightarrow0.
\label{eq:ordinary-Abel-high-pass}
\end{equation}
For the logarithmic family,
\begin{equation}
\mathfrak M_{N,T;\mathfrak A}
\left(
(I-P_{1/\Lambda_N})\widetilde h_\bullet
\right)
\longrightarrow0.
\label{eq:log-Abel-high-pass}
\end{equation}
\end{theorem}

\begin{proof}
We prove \eqref{eq:ordinary-Abel-high-pass}.
The proof of \eqref{eq:log-Abel-high-pass} is identical, using
\(\widetilde f_{A,\alpha}\) and \(\widetilde h_\alpha\) together with
the corresponding logarithmic estimates.

We first control the dyadic Abel scales \(A>N^\beta\).
Choose a fixed even integer \(q\geq4\) such that
\(
1/q<\beta\sigma_*.
\)
Since
\[
\#\left\{
n\in\mathbb N:
n_-(N,T)\leq n\leq n_+(N,T)
\right\}
\leq C_TN,
\]
Lemma~\ref{lem:block-finite-family-maximum} and
Proposition~\ref{prop:global-Abel-endpoint} give, uniformly for
dyadic \(A>N^\beta\),
\[
\mathfrak M_{N,T;\mathfrak A}(f_{A,\bullet})
\leq
C_{\mathfrak A,T,q}
N^{1/q}A^{-\sigma_*}.
\]
Hence, by dyadic summation,
\begin{equation}
\sum_{\substack{
A>N^\beta\\
A\ {\rm dyadic}}}
\mathfrak M_{N,T;\mathfrak A}(f_{A,\bullet})
\leq
C_{\mathfrak A,T,q}
N^{1/q-\beta\sigma_*}
\longrightarrow0.
\label{eq:high-Abel-scale-disposal}
\end{equation}

We now combine this high-scale estimate with the intermediate-scale
bound from Theorem~\ref{thm:active-abel-range}.
For each fixed \(N\), Proposition~\ref{prop:exact-Abel-telescope}
implies that
\[
\sup_{\alpha\in\mathfrak A}
\sup_{0\leq x\leq n_+(N,T)}
\left|
P_{1/(2^{J+1}\Lambda_N)}h_\alpha(x)-h_\alpha(x)
\right|
\longrightarrow0.
\]
Since \(0\leq S_n(\theta)\leq n_+(N,T)\), the definition of \(X_n\)
therefore yields
\[
\begin{aligned}
&
\mathfrak M_{N,T;\mathfrak A}
\left(
\sum_{j=0}^{J}f_{2^j\Lambda_N,\bullet}
-
(I-P_{1/\Lambda_N})h_\bullet
\right)
\\
&\quad\leq
2\sqrt{|\mathfrak A|\,n_+(N,T)}
\sup_{\alpha\in\mathfrak A}
\sup_{0\leq x\leq n_+(N,T)}
\left|
P_{1/(2^{J+1}\Lambda_N)}h_\alpha(x)-h_\alpha(x)
\right|
\longrightarrow0
\end{aligned}
\]
as \(J\to\infty\).

Applying subadditivity to the finite Abel sum and then letting
\(J\to\infty\), we obtain, for all sufficiently large \(N\),
\[
\mathfrak M_{N,T;\mathfrak A}
\left((I-P_{1/\Lambda_N})h_\bullet\right)
\leq
\sum_{\substack{
\Lambda_N\leq A\leq N^\beta\\
A\ {\rm dyadic}}}
\mathfrak M_{N,T;\mathfrak A}(f_{A,\bullet})
+
\sum_{\substack{
A>N^\beta\\
A\ {\rm dyadic}}}
\mathfrak M_{N,T;\mathfrak A}(f_{A,\bullet}).
\]
The first sum tends to zero by
Theorem~\ref{thm:active-abel-range}, while the second tends to zero by \eqref{eq:high-Abel-scale-disposal}.
This proves \eqref{eq:ordinary-Abel-high-pass}.
\end{proof}

\subsection{Comparison of Abel and sharp Laguerre cutoffs}
\label{subsec:Abel-sharp-boundary-comparison}

The purpose of this subsection is to compare the Abel cutoff
\(P_{1/\Lambda_N}\) with the sharp Laguerre projection
\(\LagProj_{\leq\Lambda_N}\).
We split their difference into terms of Laguerre degree \(a\leq\Lambda_N\)
and \(a>\Lambda_N\).
The terms above the cutoff are controlled by
Proposition~\ref{cor:moving-upper-boundary-input}, while those at or below
the cutoff are estimated dyadically using
Theorem~\ref{thm:coefficient-uniform-low-block}.

Recall that \(\Lambda_N\) denotes the moving dyadic cutoff defined in
\eqref{eq:moving-LambdaN}, and continue to write
\[
\sigma_*=\min\{\alpha_-,1\}.
\]
\begin{definition}
For dyadic \(\Lambda\geq2\), define
\begin{equation}
\begin{aligned}
B_{\Lambda,\alpha}^{-,0}
=
\LagProj_{\leq\Lambda}
(P_{1/\Lambda}-I)h_\alpha,
\qquad
B_{\Lambda,\alpha}^{-,1}
=
\LagProj_{\leq\Lambda}
(P_{1/\Lambda}-I)\widetilde h_\alpha.
\end{aligned}
\label{eq:boundary-lower-definition}
\end{equation}
Define the corresponding terms above the cutoff by
\begin{equation}
\begin{aligned}
B_{\Lambda,\alpha}^{+,0}
=
(I-\LagProj_{\leq\Lambda})
P_{1/\Lambda}h_\alpha,
\qquad
B_{\Lambda,\alpha}^{+,1}
=
(I-\LagProj_{\leq\Lambda})
P_{1/\Lambda}\widetilde h_\alpha.
\end{aligned}
\label{eq:boundary-upper-definition}
\end{equation}
\end{definition}
Thus,
\begin{equation}
\begin{aligned}
(P_{1/\Lambda}-\LagProj_{\leq\Lambda})h_\alpha
&=
B_{\Lambda,\alpha}^{-,0}
+
B_{\Lambda,\alpha}^{+,0},
\\
(P_{1/\Lambda}-\LagProj_{\leq\Lambda})\widetilde h_\alpha
&=
B_{\Lambda,\alpha}^{-,1}
+
B_{\Lambda,\alpha}^{+,1}.
\end{aligned}
\label{eq:boundary-complete-decomposition}
\end{equation}

The terms above the cutoff have already been controlled in
Proposition~\ref{cor:moving-upper-boundary-input}.
It therefore remains to estimate the terms at or below the cutoff.
We decompose them into dyadic Laguerre blocks.
For dyadic \(C\leq\Lambda/2\), define
\begin{equation}
\begin{aligned}
B_{\Lambda,C,\alpha}^{-,0}
&=
\sum_{C<a\leq2C}
\left(e^{-a/\Lambda}-1\right)
\gamma_a(\alpha)L_a,
\\
B_{\Lambda,C,\alpha}^{-,1}
&=
\sum_{C<a\leq2C}
\left(e^{-a/\Lambda}-1\right)
\widetilde\gamma_a(\alpha)L_a.
\end{aligned}
\label{eq:boundary-lower-dyadic-block}
\end{equation}
At the process level, this yields the dyadic decompositions
\begin{equation}
\begin{aligned}
X_n(B_{\Lambda,\alpha}^{-,0})
&=
\sum_{\substack{
C\leq\Lambda/2\\
C\ \mathrm{dyadic}}}
X_n(B_{\Lambda,C,\alpha}^{-,0}),\\
X_n(B_{\Lambda,\alpha}^{-,1})
&=
\sum_{\substack{
C\leq\Lambda/2\\
C\ \mathrm{dyadic}}}
X_n(B_{\Lambda,C,\alpha}^{-,1}).
\end{aligned}
\label{eq:boundary-lower-dyadic-decomposition}
\end{equation}

\begin{lemma}
\label{lem:lower-boundary-coefficient-energies}
Let \(\Lambda\geq2\) and \(C\leq\Lambda/2\) be dyadic.
Uniformly for \(\alpha\in\mathfrak A\),
\begin{equation}
\sum_{C<a\leq2C}
a
\left|
\left(e^{-a/\Lambda}-1\right)
\gamma_a(\alpha)
\right|^2
\leq
C_{\mathfrak A}
\left(\frac C\Lambda\right)^2
C^{-2\alpha_-}.
\label{eq:ordinary-lower-boundary-energy}
\end{equation}
Moreover,
\begin{equation}
\sum_{C<a\leq2C}
a
\left|
\left(e^{-a/\Lambda}-1\right)
\widetilde\gamma_a(\alpha)
\right|^2
\leq
C_{\mathfrak A}
\left(\frac C\Lambda\right)^2
C^{-2\alpha_-}
\log^2(C+2).
\label{eq:log-lower-boundary-energy}
\end{equation}
\end{lemma}

\begin{proof}
For \(C<a\leq2C\),
\[
|1-e^{-a/\Lambda}|
\leq
\frac a\Lambda
\leq
\frac{2C}{\Lambda}.
\]
The two estimates now follow from
\eqref{eq:calculus-sharp-power-energy} and
\eqref{eq:calculus-sharp-log-energy}, respectively.
\end{proof}

We now specialize to the moving cutoff
\(\Lambda=\Lambda_N\) and sum these blockwise estimates to control
the complete contribution at or below the cutoff.

\begin{proposition}
\label{prop:lower-abel-sharp-boundary}
For all sufficiently large \(N\), the ordinary family satisfies
\begin{equation}
\mathfrak M_{N,T;\mathfrak A}
\left(
B_{\Lambda_N,\bullet}^{-,0}
\right)
\leq
C_{\mathfrak A,T}
\Lambda_N^{-\sigma_*}
\log^{3/2}(\Lambda_N+2).
\label{eq:ordinary-lower-boundary}
\end{equation}
For the logarithmic family,
\begin{equation}
\mathfrak M_{N,T;\mathfrak A}
\left(
B_{\Lambda_N,\bullet}^{-,1}
\right)
\leq
C_{\mathfrak A,T}
\Lambda_N^{-\sigma_*}
\log^{5/2}(\Lambda_N+2).
\label{eq:log-lower-boundary}
\end{equation}
\end{proposition}

\begin{proof}
For the ordinary family, by
\eqref{eq:boundary-lower-dyadic-decomposition},
Theorem~\ref{thm:coefficient-uniform-low-block}, and
\eqref{eq:ordinary-lower-boundary-energy},
\[
\mathfrak M_{N,T;\mathfrak A}
\left(
B_{\Lambda_N,\bullet}^{-,0}
\right)
\leq
\frac{C_{\mathfrak A,T}}{\Lambda_N}
\sum_{\substack{
C\leq\Lambda_N/2\\
C\ \mathrm{dyadic}}}
C^{1-\alpha_-}\sqrt{\log(C+2)}.
\]
Write \(C=2^j\) and
\[
J=
\left\lfloor
\log_2\frac{\Lambda_N}{2}
\right\rfloor.
\]
Since
\[
1+j\asymp\log(2^j+2),
\qquad
1+J\asymp\log(\Lambda_N+2),
\]
we have
\[
\sum_{j=0}^{J}
2^{j(1-\alpha_-)}(1+j)^{1/2}
\leq
C_{\mathfrak A}
\begin{cases}
2^{J(1-\alpha_-)}(1+J)^{1/2},
& 0<\alpha_-<1,\\[1mm]
(1+J)^{3/2},
& \alpha_-=1,\\[1mm]
1,
& \alpha_->1.
\end{cases}
\]
Therefore,
\[
\frac1{\Lambda_N}
\sum_{\substack{
C\leq\Lambda_N/2\\
C\ \mathrm{dyadic}}}
C^{1-\alpha_-}\sqrt{\log(C+2)}
\leq
C_{\mathfrak A}
\Lambda_N^{-\sigma_*}
\log^{3/2}(\Lambda_N+2).
\]
This proves \eqref{eq:ordinary-lower-boundary}.

For the logarithmic family,
\eqref{eq:log-lower-boundary-energy} gives
\[
\mathfrak M_{N,T;\mathfrak A}
\left(
B_{\Lambda_N,\bullet}^{-,1}
\right)
\leq
\frac{C_{\mathfrak A,T}}{\Lambda_N}
\sum_{\substack{
C\leq\Lambda_N/2\\
C\ \mathrm{dyadic}}}
C^{1-\alpha_-}\log^{3/2}(C+2).
\]
The same dyadic computation, with the exponent
\(1/2\) replaced by \(3/2\), yields
\[
\frac1{\Lambda_N}
\sum_{\substack{
C\leq\Lambda_N/2\\
C\ \mathrm{dyadic}}}
C^{1-\alpha_-}\log^{3/2}(C+2)
\leq
C_{\mathfrak A}
\Lambda_N^{-\sigma_*}
\log^{5/2}(\Lambda_N+2).
\]
Hence \eqref{eq:log-lower-boundary} follows.
\end{proof}

\begin{theorem}[Comparison of Abel and sharp Laguerre cutoffs]
\label{thm:abel-sharp-boundary-comparison}
For the ordinary family,
\begin{equation}
\mathfrak M_{N,T;\mathfrak A}
\left(
(P_{1/\Lambda_N}
-
\LagProj_{\leq\Lambda_N})h_\bullet
\right)
\longrightarrow0.
\label{eq:ordinary-complete-boundary}
\end{equation}
For the logarithmic family,
\begin{equation}
\mathfrak M_{N,T;\mathfrak A}
\left(
(P_{1/\Lambda_N}
-
\LagProj_{\leq\Lambda_N})
\widetilde h_\bullet
\right)
\longrightarrow0.
\label{eq:log-complete-boundary}
\end{equation}
\end{theorem}

\begin{proof}
By \eqref{eq:boundary-complete-decomposition} and subadditivity,
the terms at or below the cutoff vanish by
Proposition~\ref{prop:lower-abel-sharp-boundary}, while the terms above
the cutoff vanish by
Proposition~\ref{cor:moving-upper-boundary-input}.
This proves both conclusions.
\end{proof}

Combining the Abel--sharp comparison with the Abel semigroup remainder estimate from
Theorem~\ref{thm:abel-high-pass} now gives the corresponding sharp
spectral tail estimate.

\begin{corollary}
\label{cor:moving-sharp-Laguerre-tails}
For the ordinary family,
\begin{equation}
\mathfrak M_{N,T;\mathfrak A}
\left(
(I-\LagProj_{\leq\Lambda_N})h_\bullet
\right)
\longrightarrow0.
\label{eq:ordinary-moving-sharp-tail}
\end{equation}
For the logarithmic family,
\begin{equation}
\mathfrak M_{N,T;\mathfrak A}
\left(
(I-\LagProj_{\leq\Lambda_N})
\widetilde h_\bullet
\right)
\longrightarrow0.
\label{eq:log-moving-sharp-tail}
\end{equation}
\end{corollary}

\begin{proof}
The operator identity
\[
I-\LagProj_{\leq\Lambda_N}
=
(I-P_{1/\Lambda_N})
+
(P_{1/\Lambda_N}-\LagProj_{\leq\Lambda_N})
\]
and Theorems~\ref{thm:abel-high-pass}
and~\ref{thm:abel-sharp-boundary-comparison}
imply both conclusions.
\end{proof}

\subsection{Fixed-cutoff maximal tails}
\label{subsec:fixed-cutoff-maximal-tails}

The purpose of this subsection is to pass from the moving cutoff
\(\Lambda_N\) to a fixed Laguerre cutoff \(d\), see Theorem \ref{thm:fixed-cutoff-maximal-tails}.
We combine the low moving-cutoff sharp-band estimate from
Theorem~\ref{thm:low-moving-sharp-band} with the moving sharp-tail
estimate from Corollary~\ref{cor:moving-sharp-Laguerre-tails}.

For \(d\geq0\), recall from
\eqref{eq:spectral-power-tail} and
\eqref{eq:spectral-logarithmic-tail} that
\[
R_{\alpha,d}
=
(I-\LagProj_{\leq d})h_\alpha,
\qquad
\widetilde R_{\alpha,d}
=
(I-\LagProj_{\leq d})\widetilde h_\alpha.
\]

\begin{theorem}[Fixed-cutoff maximal tails]
\label{thm:fixed-cutoff-maximal-tails}
For every integer \(d\geq2\),
\begin{align}
\limsup_{N\to\infty}
\mathfrak M_{N,T;\mathfrak A}
\left(R_{\bullet,d}\right)
&\leq
C_{\mathfrak A,T}
d^{-\alpha_-}\sqrt{\log(d+2)},
\label{eq:ordinary-maximal-tail-rate}
\\
\limsup_{N\to\infty}
\mathfrak M_{N,T;\mathfrak A}
\left(\widetilde R_{\bullet,d}\right)
&\leq
C_{\mathfrak A,T}
d^{-\alpha_-}\log^{3/2}(d+2).
\label{eq:logarithmic-maximal-tail-alpha-rate}
\end{align}
\end{theorem}

\begin{proof}
We prove \eqref{eq:ordinary-maximal-tail-rate}.
The logarithmic estimate
\eqref{eq:logarithmic-maximal-tail-alpha-rate}
follows by the same argument, using the corresponding logarithmic
bounds in Theorem~\ref{thm:low-moving-sharp-band} and
Corollary~\ref{cor:moving-sharp-Laguerre-tails}.

Fix \(d\geq2\).
Since \(\Lambda_N\to\infty\), one has \(d\leq\Lambda_N\) for all
sufficiently large \(N\), and hence
\[
R_{\alpha,d}
=
\left(
\LagProj_{\leq\Lambda_N}
-
\LagProj_{\leq d}
\right)h_\alpha
+
\left(
I-\LagProj_{\leq\Lambda_N}
\right)h_\alpha.
\]
By subadditivity,
\[
\mathfrak M_{N,T;\mathfrak A}(R_{\bullet,d})
\leq
\mathfrak M_{N,T;\mathfrak A}
\left(
(\LagProj_{\leq\Lambda_N}-\LagProj_{\leq d})h_\bullet
\right)
+
\mathfrak M_{N,T;\mathfrak A}
\left(
(I-\LagProj_{\leq\Lambda_N})h_\bullet
\right).
\]
Taking \(\limsup_{N\to\infty}\), the first term is bounded by
Theorem~\ref{thm:low-moving-sharp-band}, while the second tends to zero
by Corollary~\ref{cor:moving-sharp-Laguerre-tails}.
This proves \eqref{eq:ordinary-maximal-tail-rate}.
\end{proof}

\section{Uniform prelimit approximation}
\label{sec:uniform-prelimit-approximation}

In this section we combine the fixed-cutoff maximal-tail
bounds from Theorem~\ref{thm:fixed-cutoff-maximal-tails} with Sobolev
control in the exponent parameter, and thereby obtain a uniform
approximation of the full prelimit field; see
Proposition~\ref{prop:sobolev-prelimit-approximation}.
Throughout this section, \(N\) is an integer satisfying
\(\lfloor Ne^{-T}\rfloor\geq1\).

Both estimates in Theorem~\ref{thm:fixed-cutoff-maximal-tails}
are needed for uniform control in the exponent parameter.
The ordinary tail controls the field itself, whereas the logarithmic
tail controls its \(\alpha\)-derivative, and the latter cannot in
general be recovered from the former. Indeed, at \(\alpha=1\),
\[
\gamma_a(1)=0,
\qquad
\widetilde\gamma_a(1)=\frac1{a(a-1)},
\qquad
a\geq2.
\]

\begin{lemma}
\label{lem:prelimit-tail-exponent-derivative}
Let
\[
\mathfrak A=[\alpha_-,\alpha_+]\Subset(0,\infty),
\qquad
|\mathfrak A|=\alpha_+-\alpha_->0,
\]
and let \(T>0\).  Fix \(N\) as above and \(d\geq0\).  For every realization, the map
\begin{equation}
\alpha\longmapsto
X_{N,R_{\alpha,d}}^{\poly}
\label{eq:paired-prelimit-tail-map}
\end{equation}
belongs to
\(
C^1(\mathfrak A;C(I_T)),
\)
and
\begin{equation}
\partial_\alpha
X_{N,R_{\alpha,d}}^{\poly}
=
X_{N,\widetilde R_{\alpha,d}}^{\poly}.
\label{eq:paired-prelimit-tail-derivative}
\end{equation}
\end{lemma}

\begin{proof}
For every
\(n_-(N,T)\leq n\leq n_+(N,T)\),
one has
\[
0\leq S_n(\theta)\leq n.
\]
Hence, on the relevant compact range of \((\alpha,S_n(\theta))\),
the functions \(x^\alpha\) and \(x^\alpha\log x\), with the latter
defined to be zero at \(x=0\), are continuous.  Together with
\eqref{eq:spectral-sharp-truncation-derivatives}, this justifies
differentiation under the angular integral and through the finite
expectation, giving
\[
\partial_\alpha X_n(R_{\alpha,d})
=
X_n(\widetilde R_{\alpha,d}).
\]
Since the polygonal interpolation coefficients are independent
of \(\alpha\), the same identity holds after interpolation, proving
\eqref{eq:paired-prelimit-tail-derivative}.
\end{proof}

We now use this derivative identity together with Theorem~\ref{thm:fixed-cutoff-maximal-tails} to obtain uniform control in the exponent parameter.
For this purpose we return to the \(p\)-parametrization.
Recall from
\eqref{eq:def-alpha-PN} that \(p=2\alpha\) identifies \(\mathfrak A\) with \(K\).
The full expectation-centered polygonal field
\(\mathcal Z_N^{\ec,\poly}\) was defined in
\eqref{eq:model-power-polygonal-fields}.  For \(d\geq2\), write the
fixed Laguerre truncation in the \(p\)-parametrization as
\[
\mathcal Z_{N,d}^{\poly}(p,t)
=
\sum_{a=2}^{d}
c_a(p)X_{N,L_a}^{\poly}(t),
\qquad
(p,t)\in K\times I_T.
\]
By \eqref{eq:equi-p-alpha}, this is precisely the fixed truncation
from Section~\ref{sec:finite-mode-process} under the change of variables
\(p=2\alpha\).
\begin{equation}
\mathcal Z_N^{\ec,\poly}(p,t)
-
\mathcal Z_{N,d}^{\poly}(p,t)
=
X_{N,R_{p/2,d}}^{\poly}(t),
\qquad
(p,t)\in K\times I_T.
\label{eq:paired-exact-prelimit-difference}
\end{equation}

\begin{proposition}[Uniform prelimit approximation]
\label{prop:sobolev-prelimit-approximation}
Let
\(
K=[p_-,p_+]\Subset(0,\infty)
\)
and \(T>0\).  For every integer \(d\geq2\),
\begin{equation}
\limsup_{N\to\infty}
\E
\left\|
\mathcal Z_N^{\ec,\poly}
-
\mathcal Z_{N,d}^{\poly}
\right\|_{C(I_T;C(K))}^2
\leq
C_{K,T}
d^{-p_-}
\log^3(d+2).
\label{eq:paired-quantitative-prelimit-approximation}
\end{equation}
In particular,
\begin{equation}
\lim_{d\to\infty}
\limsup_{N\to\infty}
\E
\left\|
\mathcal Z_N^{\ec,\poly}
-
\mathcal Z_{N,d}^{\poly}
\right\|_{C(I_T;C(K))}^2
=
0.
\label{eq:paired-prelimit-approximation-zero}
\end{equation}
\end{proposition}

\begin{proof}
If \(K=\{p_0\}\), first apply the argument below on the interval
\(K'=[p_0,p_0+\eta]\) for a fixed \(\eta>0\), and then restrict
the supremum to \(K\).  The lower endpoint remains \(p_0=p_-\),
so the stated power \(d^{-p_-}\) is preserved.
We may therefore assume \(|\mathfrak A|>0\).
By Lemma~\ref{lem:prelimit-tail-exponent-derivative}, the standard
one-dimensional Sobolev estimate for the \(C(I_T)\)-valued tail gives
\[
\sup_{\alpha\in\mathfrak A}
\left\|
X_{N,R_{\alpha,d}}^{\poly}
\right\|_{C(I_T)}^2
\leq
\frac{2}{|\mathfrak A|}
\int_{\mathfrak A}
\left\|
X_{N,R_{\alpha,d}}^{\poly}
\right\|_{C(I_T)}^2
\,\dd\alpha
+
2|\mathfrak A|
\int_{\mathfrak A}
\left\|
X_{N,\widetilde R_{\alpha,d}}^{\poly}
\right\|_{C(I_T)}^2
\,\dd\alpha .
\]
By \eqref{eq:paired-exact-prelimit-difference} and the correspondence
\(p=2\alpha\),
\[
\sup_{\alpha\in\mathfrak A}
\left\|
X_{N,R_{\alpha,d}}^{\poly}
\right\|_{C(I_T)}
=
\left\|
\mathcal Z_N^{\ec,\poly}
-
\mathcal Z_{N,d}^{\poly}
\right\|_{C(I_T;C(K))}.
\]
Taking expectations and using the definition of
\(\mathfrak M_{N,T;\mathfrak A}\) from
\eqref{eq:active-integrated-maximal-seminorm},
together with
\eqref{eq:paired-floor-polygonal-supremum}, we obtain
\[
\E
\left\|
\mathcal Z_N^{\ec,\poly}
-
\mathcal Z_{N,d}^{\poly}
\right\|_{C(I_T;C(K))}^2
\leq
\frac{2}{|\mathfrak A|}
\mathfrak M_{N,T;\mathfrak A}(R_{\bullet,d})^2
+
2|\mathfrak A|
\mathfrak M_{N,T;\mathfrak A}
(\widetilde R_{\bullet,d})^2.
\]
Taking \(\limsup_{N\to\infty}\) and applying Theorem~\ref{thm:fixed-cutoff-maximal-tails} yields
\[
\limsup_{N\to\infty}
\E
\left\|
\mathcal Z_N^{\ec,\poly}
-
\mathcal Z_{N,d}^{\poly}
\right\|_{C(I_T;C(K))}^2
\leq
C_{K,T}
d^{-p_-}
\log^3(d+2),
\]
because \(p_-=2\alpha_-\) by \eqref{eq:def-alpha-PN}.
This proves
\eqref{eq:paired-quantitative-prelimit-approximation}.
Since \(p_->0\), letting \(d\to\infty\) gives
\eqref{eq:paired-prelimit-approximation-zero}.
\end{proof}

\section{Joint convergence and final transfers}
\label{sec:joint-assembly-transfers}

 In this section we first combine the finite-mode limits with the uniform approximation
estimates to prove the joint functional limit theorem.
We then pass to deterministic centering and normalized \(L^p\)-norms
to obtain the results stated in Section~\ref{sec:model-main-results}.
For the joint limit, we apply the standard converging-together principle
\cite[Theorem~3.2]{Billingsley-1999}, using the quantitative approximation
in Proposition~\ref{prop:sobolev-prelimit-approximation}.
The centering estimate in
Proposition~\ref{prop:uniform-deterministic-expectation}, uniform for
\(p\in K\), justifies the change of centering.
Finally, we obtain the limits for the normalized norms through a Taylor
expansion that is uniform for \(p\in K\), as in the usual delta method.

Throughout this section we return to the \(p\)-parametrization used in
Sections~\ref{sec:introduction} and~\ref{sec:model-main-results}.  Thus,
\[
G(p,t)
=
\sum_{a=2}^{\infty}
c_a(p)G^{(a)}(t),
\qquad
(p,t)\in K\times I_T.
\]
By \eqref{eq:equi-p-alpha}, this is exactly the Gaussian field constructed
in Section~\ref{sec:full-gaussian-field}, written there in the
\(\alpha\)-parametrization with \(\alpha=p/2\).

\subsection{Expectation-centered joint field}
\label{subsec:assembly-expectation-centered}

We first prove the full joint limit with each statistic centered by its exact expectation.
For \(d\geq2\), we recall the fixed Laguerre truncation in the
\(p\)-parametrization as
\[
G_d(p,t)
=
\sum_{a=2}^{d}
c_a(p)G^{(a)}(t).
\]
These are the fields of
Sections~\ref{sec:finite-mode-process} and~\ref{sec:full-gaussian-field}
under the affine change of variables
\(
\alpha=p/2,
\)
by \eqref{eq:equi-p-alpha}.

\begin{theorem}
\label{thm:expectation-centered-joint-field}
Let
\[
K=[p_-,p_+]\Subset(0,\infty)
\]
and \(T>0\).
Then, as \(N\to\infty\),
\begin{equation}
\mathcal Z_N^{\ec,\poly}
\Longrightarrow
\bigl(t\longmapsto G(\cdot,t)\bigr)
\qquad
\text{in }C(I_T;C(K)),
\label{eq:assembly-expectation-centered-polygonal}
\end{equation}
and
\begin{equation}
\mathcal Z_N^{\ec,\fl}
\Longrightarrow
\bigl(t\longmapsto G(\cdot,t)\bigr)
\qquad
\text{in }D(I_T;C(K)),
\label{eq:assembly-expectation-centered-floor}
\end{equation}
where the second convergence is with respect to the Skorokhod
\(J_1\)-topology.
\end{theorem}

\begin{proof}
For every fixed \(d\geq2\), Corollary~\ref{cor:finite-power-truncation-joint},
applied on the exponent interval
\(
\left[p_-/2,p_+/2\right]
\)
and followed by the change of variables \(p=2\alpha\), gives
\[
\mathcal Z_{N,d}^{\poly}
\Longrightarrow
\bigl(t\longmapsto G_d(\cdot,t)\bigr)
\qquad
\text{in }C(I_T;C(K)).
\]

On the prelimit side,
Proposition~\ref{prop:sobolev-prelimit-approximation} yields
\[
\lim_{d\to\infty}
\limsup_{N\to\infty}
\E
\left\|
\mathcal Z_N^{\ec,\poly}
-
\mathcal Z_{N,d}^{\poly}
\right\|_{C(I_T;C(K))}^2
=
0.
\]
Thus the corresponding difference tends to zero in probability in the same iterated limit.
On the limiting side,
Theorem~\ref{thm:full-gaussian-field}, under the same reparametrization,
gives
\[
G_d(p,t)
\longrightarrow
G(p,t)
\qquad
\text{in }
L^2\!\left(\Omega;C(I_T;C(K))\right).
\]
The converging-together theorem therefore yields
\eqref{eq:assembly-expectation-centered-polygonal}.
The order of limits is \(N\to\infty\) at each fixed \(d\), followed
by \(d\to\infty\).  The maximal estimates with a growing cutoff
are used here to justify approximation by fixed truncations.

It remains to pass from polygonal to floor paths.  The convergence
\eqref{eq:assembly-expectation-centered-polygonal} implies tightness in
\(C(I_T;C(K))\).  Since compact subsets of this space are uniformly
equicontinuous and
\(\Delta_{N,T}\to0\) by \eqref{eq:model-mesh-bound}, we have
\[
\sup_{\substack{
s,t\in I_T\\
|t-s|\leq\Delta_{N,T}
}}
\left\|
\mathcal Z_N^{\ec,\poly}(\cdot,t)
-
\mathcal Z_N^{\ec,\poly}(\cdot,s)
\right\|_{C(K)}
\longrightarrow0
\]
in probability.
At the interpolation knots,
\[
\mathcal Z_N^{\ec,\poly}(\cdot,\vartheta_{N,n})
=
Y_n^{\ec}.
\]
Hence the definition of the canonical polygonal interpolation gives
\begin{equation}
\begin{aligned}
\sup_{t\in I_T}
\left\|
\mathcal Z_N^{\ec,\fl}(t)
-
\mathcal Z_N^{\ec,\poly}(\cdot,t)
\right\|_{C(K)}
&\leq
\max_{n_-\leq n<n_+}
\left\|
Y_{n+1}^{\ec}-Y_n^{\ec}
\right\|_{C(K)}
\\
&\leq
\sup_{\substack{
s,t\in I_T\\
|t-s|\leq\Delta_{N,T}
}}
\left\|
\mathcal Z_N^{\ec,\poly}(\cdot,t)
-
\mathcal Z_N^{\ec,\poly}(\cdot,s)
\right\|_{C(K)}
\longrightarrow0
\end{aligned}
\label{eq:assembly-floor-polygonal-zero}
\end{equation}
in probability.
By \eqref{eq:model-J1-uniform-domination},
\eqref{eq:assembly-floor-polygonal-zero} also implies
\[
d_{J_1}^{C(K)}
\left(
\mathcal Z_N^{\ec,\fl},
\mathcal Z_N^{\ec,\poly}
\right)
\longrightarrow0
\]
in probability.
Combining this with
\eqref{eq:assembly-expectation-centered-polygonal} and applying the
converging-together theorem proves
\eqref{eq:assembly-expectation-centered-floor}.
\end{proof}

\subsection{Uniform deterministic centering}
\label{subsec:assembly-deterministic-centering}

Next we replace the exact expectations by \(m_p\), with an error that vanishes uniformly over the exponent and logarithmic time.
For \(n\geq2\), set
\[
b_n(p)
=
\sqrt n
\left[
\E T_n^{(p)}-m_p
\right],
\qquad
p\in K.
\]
Recall that in \eqref{eq:model-power-grid-fields}
\[
Y_n^{\ec}(p)
=
\sqrt n
\left[
T_n^{(p)}-\E T_n^{(p)}
\right],
\qquad
Y_n^{\dc}(p)
=
\sqrt n
\left[
T_n^{(p)}-m_p
\right].
\]
Then,
\begin{equation}
Y_n^{\dc}
=
Y_n^{\ec}+b_n.
\label{eq:Yn-dc-ec-bn}
\end{equation}
Thus the two centering conventions differ only by the deterministic
bias \(b_n\).

\begin{theorem}
\label{thm:deterministic-centering}
Let
\[
K=[p_-,p_+]\Subset(0,\infty)
\]
and \(T>0\).
Then, as \(N\to\infty\),
\begin{equation}
\mathcal Z_N^{\dc,\poly}
\Longrightarrow
\bigl(t\longmapsto G(\cdot,t)\bigr)
\qquad
\text{in }C(I_T;C(K)),
\label{eq:assembly-deterministic-polygonal}
\end{equation}
and
\begin{equation}
\mathcal Z_N^{\dc,\fl}
\Longrightarrow
\bigl(t\longmapsto G(\cdot,t)\bigr)
\qquad
\text{in }D(I_T;C(K)),
\label{eq:assembly-deterministic-floor}
\end{equation}
where the second convergence is with respect to the Skorokhod
\(J_1\)-topology.
\end{theorem}

\begin{proof}
By linearity of the polygonal interpolation and
\eqref{eq:Yn-dc-ec-bn}, both path constructions satisfy the common bound
\[
\begin{aligned}
\max\Bigl\{
&\left\|
\mathcal Z_N^{\dc,\poly}
-
\mathcal Z_N^{\ec,\poly}
\right\|_{C(I_T;C(K))},\\
&\sup_{t\in I_T}
\left\|
\mathcal Z_N^{\dc,\fl}(t)
-
\mathcal Z_N^{\ec,\fl}(t)
\right\|_{C(K)}
\Bigr\}
\leq
\max_{n_-\leq n\leq n_+}\|b_n\|_{C(K)}
\longrightarrow0,
\end{aligned}
\]
where the convergence follows from
Proposition~\ref{prop:uniform-deterministic-expectation}.
The floor-process bound also controls the corresponding
\(J_1\)-distance by \eqref{eq:model-J1-uniform-domination}.
Theorem~\ref{thm:expectation-centered-joint-field} and Slutsky's theorem
in the respective path spaces therefore give
\eqref{eq:assembly-deterministic-polygonal} and
\eqref{eq:assembly-deterministic-floor}.
\end{proof}

\subsection{Transfer to the normalized \texorpdfstring{\(L^p\)}{Lp}-norm}
\label{subsec:assembly-normalized-transfer}

Next we transfer the joint limit from power statistics to normalized \(L^p\)-norms by a uniform Taylor expansion.
Recall from \eqref{eq:intro-first-order-mean} that
\[
m_p
=
\Gamma\left(1+\frac p2\right),
\]
and from \eqref{eq:Hpt-def} that
\[
\mathfrak a_p
=
\frac1p\,m_p^{1/p-1}.
\]
We also recall the normalized-norm fields from
Section~\ref{subsec:intro-main-results}.
For \(n\geq1\),
\[
W_n(p)
=
\sqrt n
\left[
\left(T_n^{(p)}\right)^{1/p}
-
m_p^{1/p}
\right]
=
\sqrt n
\left[
\frac{\|q_n\|_p}{\sqrt n}
-
m_p^{1/p}
\right],
\qquad p\in K,
\]
so that \(W_n\in C(K)\).  Its polygonal and floor realizations are
\[
\mathcal W_N^{\poly}(t)
=
\operatorname{Poly}_{N,T}[W_\bullet](t),
\]
and
\[
\mathcal W_N^{\fl}(t)
=
W_{n_N(t)},
\qquad
t\in I_T.
\]
Equivalently, pointwise in \(p\),
\[
\mathcal W_N^{\fl}(t)(p)
=
W_{n_N(t)}(p).
\]
Thus
\[
\mathcal W_N^{\poly}\in C(I_T;C(K)),
\qquad
\mathcal W_N^{\fl}\in D(I_T;C(K)).
\]

\begin{proposition}
\label{prop:assembly-uniform-delta-approximation}
Let
\[
K=[p_-,p_+]\Subset(0,\infty)
\]
and \(T>0\).
Then, as \(N\to\infty\),
\begin{equation}
\left\|
\mathcal W_N^{\poly}
-
\mathfrak a_\bullet
\mathcal Z_N^{\dc,\poly}
\right\|_{C(I_T;C(K))}
=
o_{\Prob}(1),
\label{eq:assembly-delta-polygonal}
\end{equation}
and
\begin{equation}
\sup_{t\in I_T}
\left\|
\mathcal W_N^{\fl}(t)
-
\mathfrak a_\bullet
\mathcal Z_N^{\dc,\fl}(t)
\right\|_{C(K)}
=
o_{\Prob}(1).
\label{eq:assembly-delta-floor}
\end{equation}
\end{proposition}

\begin{proof}
At the grid points, recall from \eqref{eq:model-power-grid-fields} that
\[
Y_n^{\dc}(p)
=
\sqrt n\left(T_n^{(p)}-m_p\right).
\]
Hence
\[
\begin{aligned}
W_n(p)-\mathfrak a_pY_n^{\dc}(p)
&=
\sqrt n
\left[
\left(T_n^{(p)}\right)^{1/p}
-
m_p^{1/p}
-
\mathfrak a_p
\left(T_n^{(p)}-m_p\right)
\right].
\end{aligned}
\]
By linearity of the canonical polygonal interpolation,
\(\mathcal W_N^{\poly}
-\mathfrak a_\bullet\mathcal Z_N^{\dc,\poly}\)
is the polygonal interpolation of these grid-level remainders.
Since the interpolation attains every grid value and is affine between
adjacent grid points, while the floor process takes precisely the same
grid values,
\begin{equation}
\begin{aligned}
\left\|
\mathcal W_N^{\poly}
-
\mathfrak a_\bullet
\mathcal Z_N^{\dc,\poly}
\right\|_{C(I_T;C(K))}
&=
\sup_{t\in I_T}
\left\|
\mathcal W_N^{\fl}(t)
-
\mathfrak a_\bullet
\mathcal Z_N^{\dc,\fl}(t)
\right\|_{C(K)}
\\
&=
\max_{n_-\leq n\leq n_+}
\sup_{p\in K}
\sqrt n
\left|
\left(T_n^{(p)}\right)^{1/p}
-
m_p^{1/p}
-
\mathfrak a_p
\left(T_n^{(p)}-m_p\right)
\right|.
\end{aligned}
\label{eq:poly-floor-W-Cnorm-equi}
\end{equation}
Thus both \eqref{eq:assembly-delta-polygonal} and
\eqref{eq:assembly-delta-floor} reduce to showing that the common
quantity on the right-hand side is \(o_{\Prob}(1)\).

Since the polygonal interpolation attains every grid value,
Theorem~\ref{thm:deterministic-centering} gives
\[
D_N:=
\max_{n_-\leq n\leq n_+}\sup_{p\in K}
\sqrt n\,|T_n^{(p)}-m_p|
=
\left\|\mathcal Z_N^{\dc,\poly}\right\|_{C(I_T;C(K))}
=O_{\Prob}(1).
\]
Thus the unscaled deviations are at most
\(D_N/\sqrt{n_-}=o_{\Prob}(1)\).
Since \(p\mapsto m_p\) is continuous and strictly positive on \(K\),
we may choose \(0<\delta<\frac12\min_{p\in K}m_p\) and
\(C_K<\infty\) such that, uniformly for \(p\in K\) and
\(|x-m_p|\leq\delta\),
\begin{equation}
\left|
x^{1/p}
-
m_p^{1/p}
-
\mathfrak a_p(x-m_p)
\right|
\leq
C_K|x-m_p|^2.
\label{eq:assembly-uniform-Taylor}
\end{equation}
The event \(D_N/\sqrt{n_-}\leq\delta\) has probability tending to one.
On this event, the common remainder in
\eqref{eq:poly-floor-W-Cnorm-equi} is bounded by
\[
\max_{n_-\leq n\leq n_+}\sup_{p\in K}
C_K\sqrt n\,|T_n^{(p)}-m_p|^2
\leq\frac{C_KD_N^2}{\sqrt{n_-}}
=o_{\Prob}(1).
\]
Together with \eqref{eq:poly-floor-W-Cnorm-equi}, this proves
\eqref{eq:assembly-delta-polygonal} and
\eqref{eq:assembly-delta-floor}.
\end{proof}

\subsection{Proof of the principal results}
\label{subsec:assembly-principal-results}

\begin{proof}[Proof of Theorem~\ref{thm:main-joint-field}]
The expectation-centered convergences
\eqref{eq:main-ec-poly-convergence}
and
\eqref{eq:main-ec-floor-convergence}
are precisely the two conclusions of
Theorem~\ref{thm:expectation-centered-joint-field}.

The deterministically centered convergences
\eqref{eq:main-dc-poly-convergence}
and
\eqref{eq:main-dc-floor-convergence}
follow from
Theorem~\ref{thm:deterministic-centering}.

Finally, Proposition~\ref{prop:joint-covariance}, together with
\[
c_a(p)=\gamma_a\!\left(\frac p2\right),
\]
yields \eqref{eq:main-joint-covariance} after the reparametrization
\(
\alpha=p/2,
\)
\(
\beta=q/2.
\)
The absolute and uniform convergence of the covariance series follows
from the same proposition.
\end{proof}

\begin{proof}[Proof of Corollary~\ref{cor:intro-normalized-norm}]
Recall from \eqref{eq:Hpt-def} that
\[
H(p,t)
=
\mathfrak a_pG(p,t),
\qquad
\mathfrak a_p
=
\frac1p\,m_p^{1/p-1}.
\]
Since \(\mathfrak a_\bullet\in C(K)\), multiplication by
\(\mathfrak a_\bullet\) is a bounded linear map on \(C(K)\), and hence
induces a continuous map on \(C(I_T;C(K))\).
Therefore, by
\eqref{eq:assembly-deterministic-polygonal}
and the continuous mapping theorem,
\[
\mathfrak a_\bullet
\mathcal Z_N^{\dc,\poly}
\Longrightarrow
\bigl(t\longmapsto H(\cdot,t)\bigr)
\qquad
\text{in }C(I_T;C(K)).
\]
Together with
\eqref{eq:assembly-delta-polygonal},
Slutsky's theorem yields
\eqref{eq:intro-norm-polygonal}.

Likewise, multiplication by \(\mathfrak a_\bullet\) induces a continuous
map on \(D(I_T;C(K))\) equipped with the Skorokhod \(J_1\)-topology.
Hence, \eqref{eq:assembly-deterministic-floor} and the continuous mapping
theorem give
\[
\mathfrak a_\bullet
\mathcal Z_N^{\dc,\fl}
\Longrightarrow
\bigl(t\longmapsto H(\cdot,t)\bigr)
\qquad
\text{in }D(I_T;C(K)).
\]
Moreover, by \eqref{eq:model-J1-uniform-domination} and
\eqref{eq:assembly-delta-floor},
\[
d_{J_1}^{C(K)}
\left(
\mathcal W_N^{\fl},
\mathfrak a_\bullet
\mathcal Z_N^{\dc,\fl}
\right)
\longrightarrow0
\]
in probability.
The converging-together theorem therefore yields
\eqref{eq:intro-norm-floor}.

Since \(G\) is centered Gaussian and
\(H(p,t)=\mathfrak a_pG(p,t)\), the field \(H\) is also centered
Gaussian.  Moreover, \(\mathfrak a_\bullet\) is smooth on the compact
interval \(K\), so Theorem~\ref{thm:full-gaussian-field} implies that
\(H\) may be chosen in \(C^1(K;C(I_T))\).  Finally,
\[
\Cov\!\left(H(p,s),H(q,t)\right)
=
\mathfrak a_p\mathfrak a_q
\Cov\!\left(G(p,s),G(q,t)\right),
\]
and the asserted covariance formula follows from
\eqref{eq:main-joint-covariance}.
\end{proof}

\appendix
\section{Homogeneous chaos, balanced counts, and finite-mode convergence}
\label{app:homogeneous-chaos}

This appendix gives the estimates for balanced sums used in the finite-mode
application of the functional de Jong theorem
(Section~\ref{subs:app-A-increments-tightness-transfer}), together with explicit
dependence on degree for the later estimates with growing spectral cutoffs
(Section~\ref{sec:boundaries-fixed-tails}).
Section~\ref{subs:app-A-Walsh-top-chaos} isolates the highest homogeneous
Walsh component and controls the lower-degree remainder.
Section~\ref{subs:app-A-balanced-covariance} establishes the balanced-support
variance count and the exact nested-length martingale identity.
Section~\ref{subs:app-A-influences-diagrams} proves the low-influence and
fourth-cumulant bounds needed for the finite-mode limit.
Section~\ref{subs:app-A-increments-tightness-transfer} recalls the precise
functional limit theorem, verifies its hypotheses for our prefix family,
and proves uniform transfer to the Laguerre statistics.
Finally, Section~\ref{subs:app-A-finite-dimensional-homogeneous-vector}
records the standard mixed hypercontractive inequality used in the
degree-sensitive extraction estimate of Section~\ref{sec:boundaries-fixed-tails}.

\subsection{Walsh decomposition and top-chaos extraction}\label{subs:app-A-Walsh-top-chaos}

We now isolate the highest-degree Walsh component of each Laguerre statistic and show that the remaining terms are negligible at fixed degree.
The Laguerre combination cancels the leading contributions from lower
Walsh degrees. Theorem~\ref{thm:fixed-mode-pure-chaos-reduction} bounds the
remaining terms to order \(n^{-1/2}\) in \(L^2\), with explicit dependence
on the Laguerre degree \(a\).

This cancellation plays the role of the Hermite-to-homogeneous-sum
reduction in Angst and Poly \cite{AngstPoly2023}, where repeated
indices are organized through Newton--Girard identities.
The classical Hermite--Laguerre identity in
Section~\ref{subsec:spectral-power-coefficients} identifies the
corresponding Gaussian chaos degree as \(2a\).
For the common-sign field, the balanced-index counts below give the
Walsh approximation and its degree-dependent remainder.

For an integer \(j\geq0\), define
\begin{equation}
J_{n,j}
=
\int_{\Torus}S_n(\theta)^j\,\dd\mu(\theta).
\label{eq:chaos-Jnj}
\end{equation}
By convention,
\[
J_{n,0}=1,
\]
and Parseval's identity implies
\[
J_{n,1}=1.
\]
For \(j\geq1\), Fourier integration yields
\begin{equation}
\begin{aligned}
J_{n,j}
&= n^{-j}
\sum_{\substack{
0\leq u_1,\ldots,u_j<n\\
0\leq v_1,\ldots,v_j<n\\
u_1+\cdots+u_j=v_1+\cdots+v_j
}}
\varepsilon_{u_1}\cdots\varepsilon_{u_j}
\varepsilon_{v_1}\cdots\varepsilon_{v_j}.
\end{aligned}
\label{eq:chaos-power-statistic-expansion}
\end{equation}
Since
\[
L_a(x)
=
\sum_{j=0}^{a}
(-1)^j\binom aj\frac{x^j}{j!},
\]
\eqref{eq:spectral-invisible-mode-values}
yield, for \(a\geq2\),
\begin{equation}
X_n(L_a)
=
\sqrt n
\sum_{j=2}^{a}
(-1)^j
\binom aj
\frac{J_{n,j}-\E J_{n,j}}{j!}.
\label{eq:chaos-Laguerre-power-expansion}
\end{equation}

For a finite set
\[
E\subseteq\{0,\ldots,n-1\},
\]
write
\[
\varepsilon_E
=
\prod_{r\in E}\varepsilon_r,
\qquad
\varepsilon_{\varnothing}=1.
\]
The Walsh monomials
\[
\bigl\{
\varepsilon_E:
E\subseteq\{0,\ldots,n-1\}
\bigr\}
\]
form an orthonormal basis of
\[
L^2
\bigl(
\sigma(\varepsilon_0,\ldots,\varepsilon_{n-1})
\bigr).
\]
\begin{definition}
For a square-integrable random variable \(F\) measurable with respect to
these signs, define
\begin{equation}
\widehat F(E)
=
\E[F\varepsilon_E]
\label{eq:chaos-Walsh-coefficient}
\end{equation}
and its homogeneous Walsh projection of degree \(r\) by
\begin{equation}
\mathsf P_r^{(n)}F
=
\sum_{\substack{
        E\subseteq\{0,\ldots,n-1\}\\
        |E|=r}}
\widehat F(E)\varepsilon_E.
\label{eq:chaos-Walsh-projection}
\end{equation}
When the ambient length is clear, we write \(\mathsf P_rF\).
\end{definition}
Every monomial in
\eqref{eq:chaos-power-statistic-expansion}
contains \(2j\) Rademacher factors before cancellations.
After using
\(\varepsilon_r^2=1\), its Walsh support consists of the labels occurring with odd multiplicity and therefore has even cardinality.
Hence,
\begin{equation}
J_{n,j}
=
\sum_{s=0}^{j}\mathsf P_{2s}J_{n,j},
\label{eq:chaos-even-Walsh-degrees}
\end{equation}
so only the Walsh degrees
\[
0,2,\ldots,2j
\]
can occur.

We next identify the highest homogeneous Walsh component of \(X_n(L_a)\).

For an integer \(a\geq2\), let \(\mathcal B_{n,a}\) be the family of ordered
tuples
\[
(\mathbf u,\mathbf v)
=
(u_1,\ldots,u_a,v_1,\ldots,v_a)
\in\{0,\ldots,n-1\}^{2a}
\]
such that
\begin{equation}
u_1+\cdots+u_a
=
v_1+\cdots+v_a
\label{eq:chaos-balance-equation}
\end{equation}
and all \(2a\) labels are distinct.  Define
\begin{equation}
B_{n,a}
=
\sum_{(\mathbf u,\mathbf v)\in\mathcal B_{n,a}}
\prod_{\alpha=1}^{a}\varepsilon_{u_\alpha}
\prod_{\beta=1}^{a}\varepsilon_{v_\beta}
\label{eq:chaos-Bna}
\end{equation}
and
\begin{equation}
\mathcal G_{n,a}
=
\frac{(-1)^a}{a!} n^{1/2-a}B_{n,a}.
\label{eq:chaos-top-chaos}
\end{equation}

Since \(J_{n,j}\) has Walsh degree at most \(2j\), all terms with
\(j<a\) vanish under \(\mathsf P_{2a}\).
Hence,
\(\mathsf P_{2a}X_n(L_a)\) is determined entirely by the term \(j=a\)
in \eqref{eq:chaos-Laguerre-power-expansion}.
In that term, degree \(2a\) occurs precisely when all \(2a\) labels are distinct, while Fourier integration imposes the balance relation
\eqref{eq:chaos-balance-equation}.
Since the coefficient of \(J_{n,a}\)
in \(X_n(L_a)\) is
\[
\sqrt n\,\frac{(-1)^a}{a!},
\]
and \(J_{n,a}\) carries the normalization \(n^{-a}\), we obtain
\begin{equation}
\mathcal G_{n,a}
=
\mathsf P_{2a}X_n(L_a).
\label{eq:chaos-top-projection}
\end{equation}

It remains to control the lower homogeneous components.
The leading degree-\(2s\) contribution of an individual power statistic
\(J_{n,j}\) is not small by itself; the required gain appears only after
the full Laguerre combination in
\eqref{eq:chaos-Laguerre-power-expansion}
has been formed.  The following diagonal-insertion estimate isolates the
leading contribution and its remainder.

For a set
\[
E\subseteq\{0,\ldots,n-1\},
\]
put
\begin{equation}
\widehat J_{n,j}(E)
=
\E[J_{n,j}\varepsilon_E].
\label{eq:chaos-Jnj-fixed-support-coefficient}
\end{equation}
 \(\widehat J_{n,j}(E)\) is the weight of the sign product \(\varepsilon_E\) in the Walsh expansion of \(J_{n,j}\).

We insert equal labels in pairs, placing one member of each pair on
each side of the balance equation. These pairs preserve both the balance
relation and the Walsh support. We use the following terminology in this
insertion argument.
The \(2j\) ordered positions
\[
u_1,\ldots,u_j,v_1,\ldots,v_j
\]
are called \emph{slots};
the \(u\)-slots are positive and the \(v\)-slots are negative.
For a fixed assignment, its \emph{equality partition} is the partition of the slots into blocks carrying the same integer label.
For a fixed Walsh support \(E\), a block is called \emph{visible} if its label belongs to \(E\), and \emph{neutral} otherwise.
Thus, visible blocks have odd cardinality and neutral blocks have even cardinality.
A \emph{neutral pair} is a neutral block consisting of one positive and one negative slot.
An assignment is called a \emph{distinct neutral-pair assignment}
if every visible block is a singleton, every neutral block is a
neutral pair, the neutral labels are mutually distinct, and none
of them belongs to \(E\).

Next, we obtain the leading term on each fixed Walsh support by inserting pairs of equal labels, leaving an error of order \(n^{-s-1}\).
\begin{lemma}
\label{lem:fixed-support-diagonal-insertion}
Fix integers
\(
1\leq s\leq j.
\)
For every set
\[
E\subseteq\{0,\ldots,n-1\},
\qquad |E|=2s,
\]
one has
\begin{equation}
\widehat J_{n,j}(E)
=
\binom js^{\!2}(j-s)!
\widehat J_{n,s}(E)
+
r_{n,j,s}(E),
\label{eq:chaos-fixed-support-diagonal-insertion}
\end{equation}
where
\begin{equation}
|r_{n,j,s}(E)|
\leq
(2j)^{4j}n^{-s-1}.
\label{eq:chaos-fixed-support-remainder}
\end{equation}
In particular,
\(
r_{n,s,s}(E)=0.
\)
\end{lemma}

\begin{proof}
By
\eqref{eq:chaos-power-statistic-expansion}
and
\eqref{eq:chaos-Jnj-fixed-support-coefficient},
\(n^j\widehat J_{n,j}(E)\) counts the ordered assignments
\[
(u_1,\ldots,u_j,v_1,\ldots,v_j)
\in\{0,\ldots,n-1\}^{2j}
\]
which satisfy the balance relation
\[
u_1+\cdots+u_j
=
v_1+\cdots+v_j
\]
and whose Walsh support is exactly \(E\).

We only need to prove the case
\(1\le s<j\).
Put
\[
m:=j-s.
\]
Consider first the insertion family obtained from an assignment contributing to \(\widehat J_{n,s}(E)\).
Choose \(s\) positive slots and \(s\) negative slots, and then
pair the remaining \(m\) positive slots bijectively with the
remaining \(m\) negative slots.
Hence, the number of slot choices is
\[
\binom js^{\!2}m!.
\]
Allowing the common label on every inserted positive--negative pair to range over \(\{0,\ldots,n-1\}\) gives \(n^m\) choices.
Therefore, after the normalization \(n^{-j}\), the unrestricted insertion contribution is exactly
\begin{equation}
\binom js^{\!2}m!\,
\widehat J_{n,s}(E).
\label{eq:unrestricted-insertion}
\end{equation}

It remains to bound the difference between the number of
unrestricted insertion choices and the number of assignments
contributing to \(\widehat J_{n,j}(E)\).
Let \(\pi\) be the equality partition of the \(2j\) slots.  The number of
possible equality partitions satisfies
\begin{equation}
\#\{\pi\}
\leq
(2j)^{2j}.
\label{eq:chaos-equality-partition-count}
\end{equation}
Since the Walsh support is \(E\), the labels in \(E\) are attached to the odd blocks of \(\pi\).
For a fixed \(\pi\), the number of assignments of the \(2s\) prescribed labels in \(E\) to these odd blocks is \((2s)!\).
Thus, the number of possible choices of the equality partition and the
labels of its odd blocks is at most
\begin{equation}
(2j)^{2j}(2s)!.
\label{eq:chaos-pattern-visible-count}
\end{equation}

An assignment which is not a distinct neutral-pair assignment falls into
one of the following two cases.
Either some block has larger than its minimal possible size
(that is, some visible block has size at least \(3\), or some neutral block
has size at least \(4\)),
or all blocks have their minimal possible sizes but at least one neutral block consists of two slots on the same side.

In the first case, once the equality partition and the labels of the visible blocks have been fixed, the number \(r\) of neutral blocks whose labels remain free satisfies
\[
r\leq m-1.
\]
Indeed, the \(2s\) visible blocks occupy at least \(2s\) slots and the
\(r\) neutral blocks occupy at least \(2r\) slots;
if one block has larger than its minimal possible size, at least two additional slots are required.
Hence,
\[
2j\geq 2s+2r+2,
\]
and this gives \(r\leq m-1\).
Therefore,
\begin{equation}
\#\{\text{neutral-label assignments}\}
\leq
n^{m-1}.
\label{eq:chaos-bad-pattern-free-label-count}
\end{equation}

It remains to consider the second case, in which every block has its minimal possible size but at least one neutral block is supported entirely on one side.
Let
\[
x_1,\ldots,x_r,
\qquad r\leq m,
\]
be the labels of the neutral blocks.
After the visible labels have been fixed, the balance relation has the form
\begin{equation}
c_1x_1+\cdots+c_rx_r=b_E,
\label{eq:chaos-same-side-balance-relation}
\end{equation}
where
\begin{equation}
c_\nu=
\#\{\text{positive slots in class }\nu\}
-
\#\{\text{negative slots in class }\nu\}.
\label{eq:coefficient-slot}
\end{equation}
Since at least one neutral block is supported entirely on one side,
\[
(c_1,\ldots,c_r)\neq0.
\]
Choose \(\nu_0\) with \(c_{\nu_0}\neq0\).
Once the remaining \(r-1\) labels are fixed,
\eqref{eq:chaos-same-side-balance-relation} determines \(x_{\nu_0}\) up to at most one admissible integer value.
Thus,
\begin{equation}
\#\bigl\{
(x_1,\ldots,x_r)\in\{0,\ldots,n-1\}^r:
\eqref{eq:chaos-same-side-balance-relation}\ \text{holds}
\bigr\}
\leq
n^{m-1}.
\label{eq:chaos-same-side-free-index-loss}
\end{equation}
Combining
\eqref{eq:chaos-pattern-visible-count},
\eqref{eq:chaos-bad-pattern-free-label-count}, and
\eqref{eq:chaos-same-side-free-index-loss}, the number of assignments contributing to \(J_{n,j}\) which are not distinct neutral-pair assignments is at most
\begin{equation}
(2j)^{2j}(2s)!n^{m-1}.
\label{eq:chaos-actual-bad-assignment-count}
\end{equation}

We also have to bound the insertion choices counted by
\eqref{eq:unrestricted-insertion} which fail to produce distinct neutral-pair assignments.
Since an assignment contributing to \(\widehat J_{n,s}(E)\) uses the \(2s\) labels of \(E\) exactly once,
\begin{equation}
\#\{\text{degree-}2s\text{ ordered assignments with support }E\}
\leq
(2s)!.
\label{eq:chaos-base-assignment-count}
\end{equation}
The number of \(m\)-tuples of inserted labels for which either two inserted labels coincide or one inserted label belongs to \(E\) is bounded by
\begin{equation}
\binom m2 n^{m-1}
+
2sm\,n^{m-1}
\leq
2j^2 n^{m-1}.
\label{eq:chaos-insertion-collision-count}
\end{equation}
Hence, by
\eqref{eq:chaos-base-assignment-count}--%
\eqref{eq:chaos-insertion-collision-count},
the contribution of these colliding insertion choices is bounded by
\begin{equation}
(2s)!\binom{j}{s}^2m!j^2n^{m-1}.
\label{eq:chaos-insertion-bad-assignment-count}
\end{equation}

Subtracting the unrestricted insertion contribution \eqref{eq:unrestricted-insertion} from \(\widehat J_{n,j}(E)\), and using
\eqref{eq:chaos-actual-bad-assignment-count} and
\eqref{eq:chaos-insertion-bad-assignment-count}, gives
\begin{align}
|r_{n,j,s}(E)|\leq
n^{-j}(2s)!
\left[(2j)^{2j}
+
\binom{j}{s}^{\!2}m!j^2\right]
n^{m-1}
\leq (2j)^{4j}n^{-s-1}.
\label{eq:chaos-fixed-support-remainder-proof}
\end{align}
This proves
\eqref{eq:chaos-fixed-support-diagonal-insertion} and
\eqref{eq:chaos-fixed-support-remainder}.
\end{proof}

\begin{corollary}
\label{cor:projection-diagonal-insertion}
For \(1\leq s\leq j\), there exists a Walsh polynomial \(R_{n,j,s}\) of degree \(2s\) such that
\begin{equation}
\mathsf P_{2s}J_{n,j}
=
\binom js^{\!2}(j-s)!
\mathsf P_{2s}J_{n,s}
+
R_{n,j,s},
\label{eq:chaos-projection-diagonal-insertion}
\end{equation}
and
\begin{equation}
\E|R_{n,j,s}|^2
\leq
\frac{(2j)^{8j}}{n^2}.
\label{eq:chaos-projection-remainder-L2}
\end{equation}
In particular,
\(
R_{n,s,s}=0.
\)
\end{corollary}

\begin{proof}
For \(j>s\), sum
\eqref{eq:chaos-fixed-support-diagonal-insertion}
over \(|E|=2s\).  By Walsh orthogonality,
\[
\E|R_{n,j,s}|^2
=
\sum_{|E|=2s}|r_{n,j,s}(E)|^2.
\]
Using
\eqref{eq:chaos-fixed-support-remainder},
we obtain
\[
\E|R_{n,j,s}|^2
\leq
\frac{(2j)^{8j}}{n^2}.
\]
When \(j=s\), the stated identity holds with \(R_{n,s,s}=0\).
\end{proof}

\begin{theorem}[Reduction to the highest Walsh degree]
\label{thm:fixed-mode-pure-chaos-reduction}
There exists an absolute constant \(C>0\) such that, for every integer
\(a\geq2\),
\begin{equation}
\E
\left|
X_n(L_a)-\mathcal G_{n,a}
\right|^2
\leq
\frac{\exp\{C a\log(a+2)\}}{n}.
\label{eq:chaos-pure-chaos-L2}
\end{equation}
\end{theorem}

\begin{proof}
Identities
\eqref{eq:chaos-even-Walsh-degrees}
and
\eqref{eq:chaos-top-projection} imply that
\begin{equation}
X_n(L_a)-\mathcal G_{n,a}
=
\sum_{s=1}^{a-1}
\mathsf P_{2s}X_n(L_a).
\label{eq:chaos-pure-chaos-lower-degree-decomposition}
\end{equation}
For \(1\leq s<a\), we consider \(\mathsf P_{2s}X_n(L_a)\).
By
\eqref{eq:chaos-Laguerre-power-expansion}
and \eqref{eq:chaos-even-Walsh-degrees},
\begin{equation}
\mathsf P_{2s}X_n(L_a)
=\sqrt n
\sum_{j=s}^{a}(-1)^j\binom aj
\frac{\mathsf P_{2s}J_{n,j}}{j!}.
\end{equation}
Thus \eqref{eq:chaos-projection-diagonal-insertion} gives
\[
\mathsf P_{2s}X_n(L_a)
=
\sqrt n
\sum_{j=s}^{a}
(-1)^j
\binom aj
\frac1{j!}
\binom js^{\!2}(j-s)!
\mathsf P_{2s}J_{n,s}                       +
\sqrt n
\sum_{j=s}^{a}
(-1)^j
\binom aj
\frac{R_{n,j,s}}{j!}.
\]
The coefficient of the first term vanishes, since
\[
\sum_{j=s}^{a}
(-1)^j
\binom aj
\frac1{j!}
\binom js^{\!2}(j-s)!
=
\frac{a!}{s!^2}
\frac{(-1)^s}{(a-s)!}
\sum_{\ell=0}^{a-s}
(-1)^\ell\binom{a-s}{\ell}
=0.
\]
Since \(R_{n,s,s}=0\),
\begin{equation}
\mathsf P_{2s}X_n(L_a)
=
\sqrt n
\sum_{j=s+1}^{a}
(-1)^j
\binom aj
\frac{R_{n,j,s}}{j!}.
\label{eq:chaos-lower-degree-representation}
\end{equation}
By Minkowski's inequality and
\eqref{eq:chaos-projection-remainder-L2}, \[
\left\|
\mathsf P_{2s}X_n(L_a)
\right\|_{L^2}
\leq\frac{2^a(2a)^{4a}}{\sqrt n},
\]
where we used
\[
\sum_{j=s+1}^{a}\binom aj\frac1{j!}
\leq
\sum_{j=0}^{a}\binom aj
=
2^a.
\]
Therefore, uniformly for \(1\leq s<a\),
\begin{equation}
\E
\left|
\mathsf P_{2s}X_n(L_a)
\right|^2
\leq
\frac{\exp\{C a\log(a+2)\}}{n}.
\label{eq:chaos-lower-degree-L2}
\end{equation}
Finally, summing over \(1\leq s<a\) and using the orthogonality of distinct homogeneous Walsh spaces yields \eqref{eq:chaos-pure-chaos-L2}, after enlarging \(C\) if necessary.
\end{proof}

\subsection{Balanced-support counts and nested covariance}
\label{subs:app-A-balanced-covariance}

Next we count balanced supports to compute variances and use a martingale identity to determine covariances across nested polynomial lengths.
Indeed, we determine the variance and nested-length covariance of the top chaos
\(\mathcal G_{n,a}\). Counting balanced supports gives the second-moment
estimate with explicit dependence on degree. The common Rademacher
sequence then gives a martingale identity and the covariance across lengths.

For \(a\ge 2\) and \(k\in\Z\), define
\begin{equation}
r_{n,a}(k)
=
\#\left\{
(u_1,\ldots,u_a)\in\{0,\ldots,n-1\}^a:
u_1+\cdots+u_a=k
\right\}.
\label{eq:chaos-rnak}
\end{equation}
Recall from \eqref{eq:intro-Va} that
\[
f_a
=
\ind_{[0,1]}^{*a},
\qquad
V_a
=
\int_0^a f_a(u)^2\,\dd u.
\]

For a set
\[
E\subseteq\{0,\ldots,n-1\},
\qquad
|E|=2a,
\]
let \(\nu_{n,a}(E)\) be the number of ordered tuples
\[
(\mathbf u,\mathbf v)\in\mathcal B_{n,a}
\]
whose set of labels is \(E\).  Then
\begin{equation}
B_{n,a}
=
\sum_{\substack{
        E\subseteq\{0,\ldots,n-1\}\\
        |E|=2a}}
\nu_{n,a}(E)\varepsilon_E.
\label{eq:chaos-balanced-support-expansion}
\end{equation}
For a \(2a\)-element support \(E\), a \emph{balanced bipartition} is an
unordered decomposition
\[
E=U\sqcup V,
\qquad
|U|=|V|=a,
\qquad
\sum_{u\in U}u=\sum_{v\in V}v.
\]
Writing
\[
E=\{x_1<\cdots<x_{2a}\},
\]
such a bipartition can be represented by a sign vector
\[
\sigma=(\sigma_1,\ldots,\sigma_{2a})
\in\{-1,1\}^{2a},
\]
where \(\sigma_i=1\) when \(x_i\in U\) and \(\sigma_i=-1\) when \(x_i\in V\).
The cardinality and balance conditions become
\begin{equation}
\#\{i:\sigma_i=1\}
=\#\{i:\sigma_i=-1\}
=a,
\qquad
\sum_{i=1}^{2a}\sigma_i x_i=0.
\label{eq:chaos-balanced-bipartition-sign-encoding}
\end{equation}
Interchanging \(U\) and \(V\) replaces \(\sigma\) by \(-\sigma\).
Thus, an unordered balanced bipartition is represented by a sign vector satisfying \eqref{eq:chaos-balanced-bipartition-sign-encoding},
uniquely up to an overall change of sign.

Each balanced bipartition gives \(2(a!)^2\) ordered tuples in \(\mathcal B_{n,a}\): either part may occupy the \(u\)-side, and the labels on each side may be ordered in \(a!\) ways.

\begin{lemma}
\label{lem:exceptional-balanced-supports}
There exists an absolute constant \(C>0\) such that, for every integer
\(a\geq2\), each of the following quantities is bounded by
\[
\exp\{Ca\log(a+2)\}n^{2a-2}:
\]
\begin{enumerate}[label=\textup{(\roman*)},leftmargin=2.4em]
\item the number of ordered balanced pairs of \(a\)-tuples containing a
repeated label;

\item the number of \(2a\)-element supports
\(E\) that admit two distinct balanced bipartitions.
\end{enumerate}
\end{lemma}

\begin{proof}
We use the elementary rank-counting principle that a fixed affine linear
system with coefficient matrix \(B\) in \(q\) label variables has at most
\[
n^{q-\operatorname{rank}B}
\]
solutions in \(\{0,\ldots,n-1\}^q\).

For \textup{(i)}, fix the partition
\[
\pi=\{B_1,\ldots,B_b\}
\]
of the \(2a\) slots induced by equality of their labels.
Since at least one label is repeated,
\[
b\leq2a-1.
\]
Let \(x_\nu\) be the label carried by \(B_\nu\), and
\[
c_\nu
=
\#\{\text{positive slots in }B_\nu\}
-
\#\{\text{negative slots in }B_\nu\}.
\]
The balance relation becomes
\[
\sum_{\nu=1}^{b}c_\nu x_\nu=0.
\]
If
\(
(c_1,\ldots,c_b)\neq0,
\)
the preceding counting principle gives at most
\[
n^{b-1}\leq n^{2a-2}
\]
label assignments for the fixed partition \(\pi\).
If instead
\[
c_1=\cdots=c_b=0,
\]
then every block contains equally many positive and negative slots.
Since every block is nonempty, each block contains at least two slots, and hence
\[
b\leq a.
\]
Thus, in this case there are at most
\[
n^b\leq n^{2a-2}
\]
label assignments.
By \eqref{eq:chaos-equality-partition-count}, the number
of possible partitions of the \(2a\) slots is at most
\(
(2a)^{2a}.
\)
Therefore, the number of ordered balanced pairs containing a repeated
label is at most
\[
(2a)^{2a}n^{2a-2}
\leq
\exp\{Ca\log(a+2)\}n^{2a-2}.
\]

For \textup{(ii)}, write
\[
E=\{x_1<\cdots<x_{2a}\}.
\]
By the sign-vector encoding
\eqref{eq:chaos-balanced-bipartition-sign-encoding},
two distinct balanced bipartitions are represented by sign vectors
\(
\sigma,\tau\in\{-1,1\}^{2a}
\)
such that
\[
\tau\neq\pm\sigma
\]
and
\[
\sum_{i=1}^{2a}\sigma_i x_i=0,
\qquad
\sum_{i=1}^{2a}\tau_i x_i=0.
\]
The coefficient vectors \(\sigma\) and \(\tau\) are linearly independent:
if \(\tau=c\sigma\), then their \(\{-1,1\}\)-valued entries force
\(c=\pm1\).
Hence, the two displayed equations have rank two.
For each fixed pair \((\sigma,\tau)\), disregarding the ordering and distinctness conditions on \(x_1,\ldots,x_{2a}\), the rank-counting principle gives at most
\[
n^{2a-2}
\]
compatible label assignments.

There are at most
\[
\binom{2a}{a}^{2}
\]
ordered pairs of sign vectors having \(a\) positive and \(a\) negative
entries.
Consequently, the number of supports in \textup{(ii)} is bounded by
\[
\binom{2a}{a}^2n^{2a-2}
\leq
\exp\{Ca\log(a+2)\}n^{2a-2},
\]
after enlarging the absolute constant \(C\) if necessary.
\end{proof}

Now we estimate how many ordered pairs of balanced assignments use the same set of labels.
\begin{proposition}
\label{prop:balanced-support-square-count}
There exists an absolute constant \(C>0\) such that, for every integer
\(a\geq2\),
\begin{equation}
\left|
\sum_{\substack{
E\subseteq\{0,\ldots,n-1\}\\
|E|=2a}}
\nu_{n,a}(E)^2
-
2(a!)^2V_an^{2a-1}
\right|
\leq
\exp\{Ca\log(a+2)\}n^{2a-2}.
\label{eq:chaos-balanced-support-square-count-degree-sensitive}
\end{equation}
\end{proposition}

\begin{proof}
Throughout the proof, \(C>0\) denotes an absolute constant whose value may
increase from line to line.

The leading unrestricted balanced count is a bounded-composition
asymptotic; see G\"unther--Schmidt
\cite[Lemma~4.5]{GuntherSchmidt2017}.
Indeed, replacing each \(v_j\) by \(n-1-v_j\) turns
\(\sum_j u_j=\sum_j v_j\) into a sum of \(2a\) entries equal to
\(a(n-1)\). We give the calculation below to retain the explicit
dependence on \(a\); the repeated-label and multiple-balance corrections
are then treated separately.

We first count all ordered balanced pairs of \(a\)-tuples.  By definition,
their number is
\[
\sum_{k\in\Z}r_{n,a}(k)^2.
\]
Inclusion--exclusion gives
\[
r_{n,a}(k)
=
\sum_{\ell=0}^{a}
(-1)^\ell
\binom a\ell
\binom{k-\ell n+a-1}{a-1},
\]
where the binomial coefficient is interpreted as zero when its upper entry
is smaller than \(a-1\).  On the other hand, the box-spline formula is
\[
f_a(x)
=
\frac1{(a-1)!}
\sum_{\ell=0}^{a}
(-1)^\ell
\binom a\ell
(x-\ell)_+^{a-1},
\]
where \(y_+=\max\{y,0\}\).
Therefore,
\[
\left|
r_{n,a}(k)
-
n^{a-1}f_a(k/n)
\right|
\leq
\exp\{Ca\log(a+2)\}n^{a-2},
\]
uniformly in \(k\).
Moreover,
\[
0\leq f_a\leq1,
\]
and there are at most \(an+1\) relevant values of \(k\).  After squaring
and summing the preceding approximation, we obtain
\begin{equation}
\left|
\sum_{k\in\Z}r_{n,a}(k)^2
-
n^{2a-2}\sum_{k\in\Z}f_a(k/n)^2
\right|
\leq
\exp\{Ca\log(a+2)\}n^{2a-2}.
\label{eq:diff-rnsum-fa}
\end{equation}

We next compare the last sum with the defining integral for \(V_a\).
Since
\[
f_a
=
\ind_{[0,1]}*f_{a-1},
\]
one has
\[
\|f_a\|_\infty\leq1
\]
and, almost everywhere,
\[
f_a'(x)
=
f_{a-1}(x)-f_{a-1}(x-1).
\]
The bounded-variation Riemann-sum estimate therefore gives
\[
\left|
\sum_{k\in\Z}f_a(k/n)^2
-
n\int_0^a f_a(x)^2\,\dd x
\right|
\leq4.
\]
Combining with \eqref{eq:diff-rnsum-fa} yields
\begin{equation}
\left|
\sum_{k\in\Z}r_{n,a}(k)^2
-
V_an^{2a-1}
\right|
\leq
\exp\{Ca\log(a+2)\}n^{2a-2}.
\label{eq:chaos-all-balanced-pair-count}
\end{equation}

The quantity
\[
\sum_{\substack{
        E\subseteq\{0,\ldots,n-1\}\\
        |E|=2a}}
\nu_{n,a}(E)
\]
counts the ordered balanced pairs whose \(2a\) labels are all distinct.
Thus, by part \textup{(i)} of Lemma \ref{lem:exceptional-balanced-supports}, removing the repeated-label
configurations from
\eqref{eq:chaos-all-balanced-pair-count}
gives
\begin{equation}
\left|
\sum_{|E|=2a}\nu_{n,a}(E)
-
V_an^{2a-1}
\right|
\leq
\exp\{Ca\log(a+2)\}n^{2a-2}.
\label{eq:chaos-distinct-balanced-pair-count}
\end{equation}
It remains to pass from
\(
\sum_{|E|=2a}\nu_{n,a}(E)
\)
to
\(
\sum_{|E|=2a}\nu_{n,a}(E)^2.
\)
If \(E\) admits no balanced bipartition, then
\(\nu_{n,a}(E)=0\).  If it admits exactly one balanced bipartition, then
\begin{equation*}
\nu_{n,a}(E)=2(a!)^2.
\end{equation*}
Therefore,
\[
\nu_{n,a}(E)^2-2(a!)^2\nu_{n,a}(E)
\]
vanishes whenever \(E\) admits at most one balanced bipartition.
On the other hand, every support admits at most \(\frac12\binom{2a}{a}\) balanced bipartitions, so
\begin{equation}
\nu_{n,a}(E)
\leq\binom{2a}{a}(a!)^2
=(2a)!.
\label{eq:support-estimate-nu}
\end{equation}
Consequently, whenever \(E\) admits at least two distinct balanced bipartitions,
\[
\left|\nu_{n,a}(E)^2
-
2(a!)^2\nu_{n,a}(E)
\right|
\leq
\exp\{Ca\log(a+2)\}.
\]
Part \textup{(ii)} of Lemma \ref{lem:exceptional-balanced-supports}
now implies
\[
\left|
\sum_{|E|=2a}\nu_{n,a}(E)^2
-
2(a!)^2
\sum_{|E|=2a}\nu_{n,a}(E)
\right|
\leq
\exp\{Ca\log(a+2)\}n^{2a-2}.
\]
Combining this with
\eqref{eq:chaos-distinct-balanced-pair-count}, and enlarging \(C\) if necessary,  proves
\eqref{eq:chaos-balanced-support-square-count-degree-sensitive}.
\end{proof}

\begin{corollary}
\label{cor:top-chaos-second-moment}
There exists an absolute constant \(C>0\) such that, for every integer
\(a\geq2\),
\begin{equation}
\left|
\E\mathcal G_{n,a}^2
-
2V_a
\right|
\leq
\frac{\exp\{Ca\log(a+2)\}}{n}.
\label{eq:chaos-top-second-moment-degree-sensitive}
\end{equation}
\end{corollary}

\begin{proof}
By
\eqref{eq:chaos-top-chaos},
\eqref{eq:chaos-balanced-support-expansion},
and Walsh orthogonality,
\[
\E\mathcal G_{n,a}^2
=
\frac{n^{1-2a}}{(a!)^2}
\sum_{|E|=2a}\nu_{n,a}(E)^2.
\]
The conclusion follows directly from
Proposition~\ref{prop:balanced-support-square-count}.
\end{proof}

For \(n\geq1\), let
\[
\mathcal F_n
=
\sigma(\varepsilon_0,\ldots,\varepsilon_{n-1}).
\]

\begin{lemma}
\label{lem:low-block-damped-martingale}
For \(a\geq2\) and \(m\leq n\),
\begin{equation}
\E
\left[
\mathcal G_{n,a}
\,\middle|\,
\mathcal F_m
\right]
=
\left(\frac mn\right)^{a-\frac12}
\mathcal G_{m,a}.
\label{eq:low-block-damped-martingale}
\end{equation}
Consequently,
\begin{equation}
M_{n,a}
:=
n^{a-\frac12}\mathcal G_{n,a}
\label{eq:low-block-mode-martingale}
\end{equation}
is a martingale with respect to \((\mathcal F_n)_{n\geq1}\).
\end{lemma}

\begin{proof}
It is sufficient to prove \eqref{eq:low-block-damped-martingale}.
By
\eqref{eq:chaos-top-chaos}
and
\eqref{eq:chaos-balanced-support-expansion},
\[
\mathcal G_{n,a}
=
\frac{(-1)^a}{a!}
n^{\frac12-a}
\sum_{\substack{
E\subseteq\{0,\ldots,n-1\}\\
|E|=2a}}
\nu_{n,a}(E)\varepsilon_E.
\]
For every support \(E\),
\[
\E[\varepsilon_E\mid\mathcal F_m]
=
\begin{cases}
\varepsilon_E, & E\subseteq\{0,\ldots,m-1\},\\
0, & \text{otherwise}.
\end{cases}
\]
Moreover, if \(E\subseteq\{0,\ldots,m-1\}\), then
\[
\nu_{n,a}(E)=\nu_{m,a}(E).
\]
This completes the proof of \eqref{eq:low-block-damped-martingale}.
\end{proof}

The preceding conditional-expectation identity immediately yields the exact covariance structure across nested truncation lengths.

\begin{theorem}[Nested covariance of the top chaos]
\label{thm:nested-top-chaos-covariance}
Let \(m\leq n\) and \(a,b\geq2\).  Then,
\begin{equation}
\E
\left[
\mathcal G_{m,a}\mathcal G_{n,b}
\right]
=
\ind_{\{a=b\}}
\left(\frac mn\right)^{a-\frac12}
\E\mathcal G_{m,a}^2.
\label{eq:chaos-exact-nested-covariance}
\end{equation}
Consequently, for fixed \(a,b\geq2\) and \(T>0\), uniformly for
\[
m\leq n,
\qquad
m,n\asymp_TN,
\]
one has
\begin{equation}
\Cov
\left(\mathcal G_{m,a},
\mathcal G_{n,b}
\right)
=
2V_a
\left(\frac mn\right)^{a-\frac12}
\ind_{\{a=b\}}
+
O_{a,b,T}(N^{-1}).
\label{eq:chaos-quantitative-covariance}
\end{equation}
\end{theorem}

\begin{proof}
If \(a\neq b\), then
\(\mathcal G_{m,a}\) and \(\mathcal G_{n,b}\)
belong to different homogeneous Walsh degrees and are therefore orthogonal.
Assume \(a=b\).
Since \(\mathcal G_{m,a}\) is
\(\mathcal F_m\)-measurable, Lemma \ref{lem:low-block-damped-martingale} yields \eqref{eq:chaos-exact-nested-covariance}.

The top chaoses are centered, and Corollary \ref{cor:top-chaos-second-moment}
with \(m\asymp_TN\) gives
\[
\E\mathcal G_{m,a}^2
=
2V_a+O_{a,T}(N^{-1}).
\]
Substituting this into
\eqref{eq:chaos-exact-nested-covariance}
gives
\eqref{eq:chaos-quantitative-covariance}.
\end{proof}

\begin{corollary}
\label{cor:finite-time-covariance-approximation}
Fix \(d\geq2\), \(k\geq1\), \(T>0\), and
\[
t_1,\ldots,t_k\in I_T.
\]
For
\[
F_{N,(a,\ell)}
=
\mathcal G_{\lfloor Ne^{t_\ell}\rfloor,a},
\qquad
2\leq a\leq d,
\quad
1\leq\ell\leq k,
\]
one has
\begin{equation}
\Cov
\left(
F_{N,(a,\ell)},
F_{N,(b,m)}
\right)
=
\Sigma_{(a,\ell),(b,m)}
+
O_{d,k,T}(N^{-1}),
\label{eq:chaos-finite-time-covariance-approximation}
\end{equation}
where
\begin{equation}
\Sigma_{(a,\ell),(b,m)}
=
2V_a
e^{-(a-\frac12)|t_\ell-t_m|}
\ind_{\{a=b\}}.
\label{eq:chaos-target-covariance-matrix}
\end{equation}
\end{corollary}

\begin{proof}
For every \(1\leq\ell,m\leq k\),
\[
\frac{\lfloor Ne^{t_\ell}\rfloor
\wedge
\lfloor Ne^{t_m}\rfloor}
{\lfloor Ne^{t_\ell}\rfloor
\vee
\lfloor Ne^{t_m}\rfloor}
=
e^{-|t_\ell-t_m|}
+
O_T(N^{-1}).
\]
The conclusion follows from
\eqref{eq:chaos-quantitative-covariance}.
\end{proof}

\subsection{Low influences and fourth-cumulant bounds}
\label{subs:app-A-influences-diagrams}

We next show that individual signs have vanishing influence and that the fourth cumulants vanish, as required by the functional de Jong theorem.
Indeed, the scalar case of Theorem~\ref{thm:mixed-fourth-cumulant-top-chaos}
verifies the fourth-moment hypothesis of the functional de Jong theorem.
The accompanying low-influence bound verifies its reinforced Lindeberg
condition.

Write
\begin{equation}
\mathcal G_{n,a}
=
\sum_{\substack{
        E\subseteq\{0,\ldots,n-1\}\\
        |E|=2a}}
b_{n,a}(E)\varepsilon_E,
\label{eq:chaos-top-support-form}
\end{equation}
where
\begin{equation}
b_{n,a}(E)
=
\frac{(-1)^a}{a!}
n^{1/2-a}\nu_{n,a}(E).
\label{eq:chaos-top-coefficient}
\end{equation}
\begin{definition}
For \(0\leq r<n\), define
\begin{equation}
\operatorname{Inf}_r(\mathcal G_{n,a})=
\sum_{\substack{
E\subseteq\{0,\ldots,n-1\}\\
|E|=2a,\ r\in E}}
b_{n,a}(E)^2.
\label{eq:chaos-influence-definition}
\end{equation}
\end{definition}
Thus the influence measures the variance carried by monomials containing
the specified sign.

\begin{theorem}
\label{thm:fixed-mode-low-influence}
For every fixed \(a\geq2\),
\begin{equation}
\max_{0\leq r<n}
\operatorname{Inf}_r(\mathcal G_{n,a})
\leq
\frac{C_a}{n}.
\label{eq:chaos-low-influence}
\end{equation}
\end{theorem}

\begin{proof}
Fix \(0\leq r<n\).
By
\eqref{eq:chaos-influence-definition}
and
\eqref{eq:chaos-top-coefficient},
\[
\operatorname{Inf}_r(\mathcal G_{n,a})
=
\frac{n^{1-2a}}{(a!)^2}
\sum_{\substack{|E|=2a\\r\in E}}
\nu_{n,a}(E)^2.
\]
It is therefore enough to prove that
\begin{equation}
\sum_{\substack{|E|=2a\\r\in E}}
\nu_{n,a}(E)^2
\leq
C_an^{2a-2}.
\label{eq:nu-r-2-up-bound}
\end{equation}

The sum \(\sum_{\substack{|E|=2a\\r\in E}}\nu_{n,a}(E)\) counts the ordered distinct
balanced tuples in \(\mathcal B_{n,a}\) containing \(r\).
Choose the slot occupied by \(r\), then leave the first other slot
unspecified and choose the labels in the remaining \(2a-2\) slots.
The balance equation \eqref{eq:chaos-balance-equation} determines the
unspecified label uniquely.  Ignoring distinctness and the requirement
that the determined label lie in \(\{0,\ldots,n-1\}\) gives
\[
\sum_{\substack{|E|=2a\\r\in E}}\nu_{n,a}(E)
\leq 2a\,n^{2a-2}.
\]
For each support \(E\), there are at most \((2a)!\) orderings of its
labels, so \(\nu_{n,a}(E)\leq(2a)!\).  Consequently,
\[
\sum_{\substack{|E|=2a\\r\in E}}\nu_{n,a}(E)^2
\leq (2a)!\sum_{\substack{|E|=2a\\r\in E}}\nu_{n,a}(E)
\leq 2a(2a)!\,n^{2a-2}.
\]
This proves \eqref{eq:nu-r-2-up-bound}, and the coefficient normalization
gives
\[
\operatorname{Inf}_r(\mathcal G_{n,a})
\leq \frac{2a(2a)!}{(a!)^2}\frac1n,
\]
which proves \eqref{eq:chaos-low-influence}.
\end{proof}

Theorem~\ref{thm:fixed-mode-low-influence}
and a finite summation give the following consequence.
\begin{corollary}
\label{cor:finite-vector-aggregate-influence}
Fix integers \(d\geq2\) and \(k\geq1\), a number \(T>0\), and times
\[
t_1,\ldots,t_k\in I_T.
\]
For \(2\leq a\leq d\) and \(1\leq\ell\leq k\), extend the coefficient
array of
\(
\mathcal G_{\lfloor Ne^{t_\ell}\rfloor,a}
\)
by zero to the common coordinate set
\[
\{0,\ldots,\lfloor Ne^T\rfloor-1\}.
\]
Then, for all sufficiently large \(N\),
\begin{equation}
\max_{0\leq r<\lfloor Ne^T\rfloor}
\sum_{a=2}^{d}
\sum_{\ell=1}^{k}
\operatorname{Inf}_r
\left(
\mathcal G_{\lfloor Ne^{t_\ell}\rfloor,a}
\right)
\leq
\frac{C_{d,k,T}}{N}.
\label{eq:chaos-aggregate-influence}
\end{equation}
\end{corollary}

\begin{proof}
For every \(1\leq\ell\leq k\),
\[
\lfloor Ne^{t_\ell}\rfloor
\geq
\frac12Ne^{-T}
\]
for all sufficiently large \(N\).  Apply
Theorem~\ref{thm:fixed-mode-low-influence}
at each of the finitely many degrees and lengths, and then sum over
\(2\leq a\leq d\) and \(1\leq\ell\leq k\).
\end{proof}

The cancellation of joint cumulants for families that split into independent
groups is classical; see \cite[Chapters~3--4]{PeccatiTaqqu2011}.
We give the support argument in our notation and then prove the
balanced-support rank estimate needed for the present model.

We next turn to fourth-cumulant bounds. For centered random variables
\(Y_1,\ldots,Y_4\), write
\begin{align}
\kappa_4(Y_1,Y_2,Y_3,Y_4)
&=\E[Y_1Y_2Y_3Y_4]-\E[Y_1Y_2]\E[Y_3Y_4]\\
&\quad
-\E[Y_1Y_3]\E[Y_2Y_4]-\E[Y_1Y_4]\E[Y_2Y_3].
\label{eq:chaos-fourth-cumulant}
\end{align}

For nonempty finite supports \(E_1,\ldots,E_s\), define their
\emph{support-intersection graph} to be the simple graph with vertex set
\(\{1,\ldots,s\}\), in which two distinct vertices \(i\) and \(j\) are
joined by an edge if and only if
\[
E_i\cap E_j\neq\varnothing.
\]

We first apply this definition with \(s=4\).  Since every expectation in
\eqref{eq:chaos-fourth-cumulant} has absolute value at most one,
\begin{equation}
\left|
\kappa_4
\left(
\varepsilon_{E_1},
\varepsilon_{E_2},
\varepsilon_{E_3},
\varepsilon_{E_4}
\right)
\right|
\leq4.
\label{eq:cumulate-4}
\end{equation}

If the support-intersection graph is disconnected, its components use
disjoint sets of independent signs. The corresponding Walsh monomials
therefore split into independent groups, so their joint cumulant vanishes.
If a coordinate belongs to exactly one support, averaging over its
centered sign makes every expectation containing that Walsh monomial
vanish. Each term of \eqref{eq:chaos-fourth-cumulant} contains such an
expectation, so the cumulant again vanishes.

We have therefore proved that
\[
\kappa_4
\left(\varepsilon_{E_1},
\varepsilon_{E_2},
\varepsilon_{E_3},
\varepsilon_{E_4}
\right)
\neq0
\]
only if the support-intersection graph is connected and every coordinate
in
\[
E_1\cup E_2\cup E_3\cup E_4
\]
belongs to at least two of the supports \(E_1,\ldots,E_4\).

Fix
\[
a_1,a_2,a_3,a_4\geq2,
\]
and write
\[
\lvert\mathbf a\rvert
:=
a_1+a_2+a_3+a_4.
\]
For \(1\leq i\leq4\), let \(E_i\) be a support with
\[
|E_i|=2a_i.
\]
Choose a balanced bipartition of \(E_i\) and one of its two sign-vector
representatives.  Denote by
\[
\sigma_i(x)\in\{-1,1\},
\qquad x\in E_i,
\]
the sign assigned to \(x\) by this representative.

Write
\[
E_1\cup E_2\cup E_3\cup E_4
=
\{x_1,\ldots,x_m\}.
\]
The four balance equations are represented by the \(4\times m\) matrix
\[
M=(M_{i\nu}),
\]
where
\[
M_{i\nu}
=
\begin{cases}
\sigma_i(x_\nu), & x_\nu\in E_i,\\
0,               & x_\nu\notin E_i.
\end{cases}
\]
Indeed, the \(i\)-th balance equation is
\[
\sum_{\nu=1}^{m}M_{i\nu}x_\nu=0.
\]
The following lemma isolates the counting estimate needed below.
\begin{lemma}
\label{lem:four-copy-rank-deficit}
Assume that the support-intersection graph of
\(E_1,\ldots,E_4\)
is connected and that every coordinate in
\[
E_1\cup E_2\cup E_3\cup E_4
\]
belongs to at least two of the four supports.
Then,
\begin{equation}
m-\operatorname{rank}M
\leq
\lvert\mathbf a\rvert-3.
\label{eq:chaos-rank-deficit}
\end{equation}
\end{lemma}

\begin{proof}
Since
\[
|E_1|+\cdots+|E_4|
=
2\lvert\mathbf a\rvert
\]
and every coordinate in
\(E_1\cup\cdots\cup E_4\)
is counted at least twice in this sum,
\[
m\leq\lvert\mathbf a\rvert.
\]

If \(m\leq \lvert\mathbf a\rvert-2\), then
\(\operatorname{rank}M\geq1\), because every row of \(M\) is nonzero.
Hence
\[
m-\operatorname{rank}M\leq \lvert\mathbf a\rvert-3.
\]

Suppose next that \(m=\lvert\mathbf a\rvert-1\).
We claim that \(\operatorname{rank}M\geq2\).
If \(\operatorname{rank}M=1\), the four nonzero rows of \(M\) are proportional.
Since their components belong to \(\{-1,0,1\}\), their nonzero positions coincide.
Consequently,
\[
E_1=E_2=E_3=E_4,
\]
and all four supports have the same cardinality, say \(2c\).  It follows
that
\[
m=2c,
\qquad
\lvert\mathbf a\rvert=4c.
\]
The identity \(m=\lvert\mathbf a\rvert-1\) would therefore give
\[
2c=4c-1,
\]
which is impossible.

Finally, suppose that
\[
m=\lvert\mathbf a\rvert.
\]
Since the sum of the four support sizes is \(2m\), and every coordinate
occurs in at least two supports, every coordinate occurs in exactly two.
We regard each coordinate as an edge joining those two support indices.
This gives a multigraph on \(\{1,2,3,4\}\) whose underlying simple graph
is the support-intersection graph. In particular, it is connected.
Let \(R_1,\ldots,R_4\) denote the rows of \(M\).
If
\[
c_1R_1+\cdots+c_4R_4=0,
\]
then a coordinate belonging to \(E_i\cap E_j\) gives
\[
c_iM_{i\nu}+c_jM_{j\nu}=0,
\]
and therefore
\[
c_j=\pm c_i.
\]
Connectedness makes the space of row relations at most one-dimensional.
Hence,
\[
\operatorname{rank}M\geq3.
\]
This proves
\eqref{eq:chaos-rank-deficit}
in all three cases.
\end{proof}

\begin{theorem}[Mixed fourth-cumulant bound]
\label{thm:mixed-fourth-cumulant-top-chaos}
Fix \(a_1,a_2,a_3,a_4\geq2\) and \(T>0\).  If
\[
n_i\asymp_TN,
\qquad
1\leq i\leq4,
\]
then
\begin{equation}
\left|
\kappa_4
\left(
\mathcal G_{n_1,a_1},
\mathcal G_{n_2,a_2},
\mathcal G_{n_3,a_3},
\mathcal G_{n_4,a_4}
\right)
\right|
\leq
\frac{C_{a_1,a_2,a_3,a_4,T}}{N}.
\label{eq:chaos-mixed-fourth-cumulant}
\end{equation}
In particular, for every fixed \(a\geq2\),
\begin{equation}
\left|
\E\mathcal G_{n,a}^4
-
3\bigl(\E\mathcal G_{n,a}^2\bigr)^2
\right|
\leq
\frac{C_a}{n}.
\label{eq:chaos-scalar-fourth-cumulant}
\end{equation}
Consequently,
\[
\E\mathcal G_{n,a}^4
\longrightarrow
3(2V_a)^2.
\]
\end{theorem}

\begin{proof}
By multilinearity of the fourth cumulant and
\eqref{eq:chaos-top-support-form},
\begin{align*}
\kappa_4
\left(\mathcal G_{n_1,a_1},
\mathcal G_{n_2,a_2},
\mathcal G_{n_3,a_3},
\mathcal G_{n_4,a_4}
\right)
=
\sum_{\substack{
E_i\subseteq\{0,\ldots,n_i-1\}\\
|E_i|=2a_i,\ 1\leq i\leq4}}
\left(\prod_{i=1}^{4}b_{n_i,a_i}(E_i)
\right)
\kappa_4
\left(\varepsilon_{E_1},
\varepsilon_{E_2},
\varepsilon_{E_3},
\varepsilon_{E_4}
\right).
\end{align*}
By the preceding support observations, a summand can be nonzero only if
the support-intersection graph of \(E_1,\ldots,E_4\) is connected and
every coordinate in
\[
E_1\cup E_2\cup E_3\cup E_4
\]
belongs to at least two of the supports.
Moreover, \eqref{eq:support-estimate-nu},
\eqref{eq:chaos-top-coefficient}, \eqref{eq:cumulate-4},
and \(n_i\asymp_TN\) give
\[
\prod_{i=1}^{4}
|b_{n_i,a_i}(E_i)|
\leq
C_{a_1,a_2,a_3,a_4,T}
N^{2-\lvert\mathbf a\rvert}.
\]
It is therefore enough to show that the number of support quadruples
giving nonzero summands is
\begin{equation}
O_{a_1,a_2,a_3,a_4,T}
\left(
N^{\lvert\mathbf a\rvert-3}
\right).
\label{eq:chaos-contributing-support-quadruple-count}
\end{equation}

For every support occurring in a nonzero summand,
\(b_{n_i,a_i}(E_i)\neq0\), and hence \(E_i\) admits a balanced bipartition.
We overcount the support quadruples by additionally choosing,
for each \(E_i\), a balanced bipartition and one of its two sign-vector
representatives.

Write
\[
E_i
=
\{x_{i,1}<\cdots<x_{i,2a_i}\}.
\]
For fixed \(a_1,\ldots,a_4\), there are only finitely many ways to specify
all equalities
\[
x_{i,p}=x_{j,q},
\qquad
i\neq j,
\]
and to choose the four sign-vector representatives.
Fix one such choice, and let \(m\) be the number of distinct coordinates.
The elementary rank-counting bound and
Lemma~\ref{lem:four-copy-rank-deficit}
give
\[
O_T
\left(
N^{\lvert\mathbf a\rvert-3}
\right)
\]
assignments of values to these coordinates satisfying the four balance equations.
Since the number of choices depends only on \(a_1,\ldots,a_4\), \eqref{eq:chaos-contributing-support-quadruple-count} implies \eqref{eq:chaos-mixed-fourth-cumulant}.

Since \(\mathcal G_{n,a}\) is centered, taking
\[
n_1=\cdots=n_4=n,
\qquad
a_1=\cdots=a_4=a,
\]
gives
\[
\kappa_4
\left(
\mathcal G_{n,a},
\mathcal G_{n,a},
\mathcal G_{n,a},
\mathcal G_{n,a}
\right)
=
\E\mathcal G_{n,a}^4
-
3\bigl(\E\mathcal G_{n,a}^2\bigr)^2.
\]
Thus,
\eqref{eq:chaos-scalar-fourth-cumulant}
follows from
\eqref{eq:chaos-mixed-fourth-cumulant}.
The final convergence follows from
Corollary~\ref{cor:top-chaos-second-moment}.
\end{proof}

\subsection{Functional de Jong theorem and uniform transfer}
\label{subs:app-A-increments-tightness-transfer}

We apply the functional de Jong theorem to the leading Walsh components and transfer their process limit to the Laguerre modes using uniform remainder bounds.
Indeed, we use the following specialization of the multivariate functional de Jong
theorem of D\"obler, Kasprzak and Peccati
\cite[Conditions~1.1--1.3 and Theorem~1.4]{DoblerKasprzakPeccati2022}.
It extends de Jong's fourth-moment and negligibility criterion for
generalized multilinear forms \cite{deJong1990}.
The kernels may depend on the support; no common symmetric kernel across
all supports is required.

\begin{theorem}[Functional de Jong theorem, recalled specialization]
\label{thm:functional-de-jong-recalled}
Fix distinct positive integers \(p_1<\cdots<p_r\).
For all sufficiently large integers \(M\), let
\[
H_{M,j}=\sum_{\substack{E\subseteq\{0,\ldots,M-1\}\\|E|=p_j}}
h_{M,j,E},\qquad \Var(H_{M,j})=1,
\]
be completely degenerate \(U\)-statistics in a common family of independent coordinates,
and let \(H_{M,j}(u)\) denote the sum restricted to
\(E\subseteq\{0,\ldots,\lfloor Mu\rfloor-1\}\), with the empty
prefix interpreted as zero.
Define
\[
\rho_{M,j}^2=\max_{0\leq i<M}\sum_{E\ni i}\E h_{M,j,E}^2,
\qquad
D_{M,j}=\sup_{E:\,\E h_{M,j,E}^2>0}
\frac{\E h_{M,j,E}^4}{(\E h_{M,j,E}^2)^2}.
\]
Suppose that, for each \(j\),
\[
\Var(H_{M,j}(u))\longrightarrow v_j(u)\quad(0\leq u\leq1),
\qquad D_{M,j}\rho_{M,j}^2\longrightarrow0,
\qquad \E H_{M,j}^4\longrightarrow3,
\]
where the continuous functions \(v_j\) describe the limiting prefix variances.
Then \((H_{M,j}(\cdot))_{j=1}^{r}\) converges in
\(D([0,1];\R^r)\) with the \(J_1\)-topology to
\((B_j(v_j(\cdot)))_{j=1}^{r}\), for independent standard Brownian
motions \(B_j\). The limiting paths are continuous.

Here \(h_{M,j,E}\) is a random kernel term depending only on the
coordinates indexed by \(E\); its conditional mean given any proper
subset of these coordinates is zero.
\end{theorem}

The quoted theorem also gives a uniform-topology conclusion. Only its
\(J_1\) conclusion with a continuous limit is used below.

\begin{proposition}[Balanced prefix sums]
\label{prop:top-chaos-functional-de-jong}
For fixed \(d\geq2\), set \(Q_{0,a}=0\) and, for \(n\geq1\),
\begin{equation}
Q_{n,a}=n^{a-1/2}\mathcal G_{n,a}
       =\frac{(-1)^a}{a!}B_{n,a},
\qquad
\sigma_{M,a}^2=\Var(Q_{M,a}).
\label{eq:chaos-prefix-Q-normalization}
\end{equation}
For all sufficiently large \(M\), simultaneously for \(2\leq a\leq d\),
\(\sigma_{M,a}>0\), and
\begin{equation}
\left(\frac{Q_{\lfloor M\cdot\rfloor,a}}{\sigma_{M,a}}\right)_{a=2}^{d}
\Longrightarrow
\left(B_a\bigl((\cdot)^{2a-1}\bigr)\right)_{a=2}^{d}
\quad\text{in }D([0,1];\R^{d-1})
\label{eq:chaos-prefix-functional-limit}
\end{equation}
with the \(J_1\)-topology, where \(B_2,\ldots,B_d\) are independent
standard Brownian motions.
\end{proposition}
 With this normalization, the prefix process \(Q_{\lfloor M\cdot\rfloor,a}/\sigma_{M,a}\) has variance one at time \(1\).
\begin{proof}
By \eqref{eq:chaos-balanced-support-expansion},
\[
Q_{n,a}=\sum_{\substack{E\subseteq\{0,\ldots,n-1\}\\|E|=2a}}
c_a(E)\varepsilon_E,
\qquad c_a(E)=\frac{(-1)^a}{a!}\nu_a(E).
\]
Here \(\nu_a(E)=\nu_{n,a}(E)\) is independent of \(n\) whenever
\(E\subseteq\{0,\ldots,n-1\}\): it counts ordered balanced tuples on
that fixed support. A support with no balance has coefficient zero.
Each nonzero product kernel is completely degenerate, since conditioning
on any proper subset of its independent centered signs leaves at least
one sign of mean zero. Thus terminal normalization gives precisely the
prefix family in Theorem~\ref{thm:functional-de-jong-recalled}, with
one coordinate for each order \(4,6,\ldots,2d\).

For the variance clocks, Corollary~\ref{cor:top-chaos-second-moment} gives
\begin{equation}
\Var(Q_{m,a})=2V_a m^{2a-1}+O_a(m^{2a-2}),
\label{eq:chaos-prefix-variance-asymptotic}
\end{equation}
so \(\sigma_{M,a}>0\) for large \(M\), since \(V_a>0\), and
\[
\frac{\Var(Q_{\lfloor Mu\rfloor,a})}{\sigma_{M,a}^2}
\longrightarrow u^{2a-1}\qquad(0<u\leq1).
\]
At \(u=0\) the variance is exactly zero, giving the continuous clock
\(v_a(u)=u^{2a-1}\).

For the reinforced Lindeberg condition, every nonzero normalized product kernel has
fourth-moment ratio exactly one. Zero kernels are omitted in the
supremum defining \(D_{M,a}\). Hence \(D_{M,a}=1\), while
Theorem~\ref{thm:fixed-mode-low-influence} and
Corollary~\ref{cor:top-chaos-second-moment} imply
\[
\rho_{M,a}^2
=\frac{\max_{0\leq i<M}\operatorname{Inf}_i(\mathcal G_{M,a})}
       {\Var(\mathcal G_{M,a})}
=O_a(M^{-1}).
\]
Finally, the equal-time scalar consequence
\eqref{eq:chaos-scalar-fourth-cumulant} of the connected counting
argument gives
\[
\E\left(\frac{Q_{M,a}}{\sigma_{M,a}}\right)^4
=\frac{\E\mathcal G_{M,a}^4}{(\E\mathcal G_{M,a}^2)^2}
=3+O_a(M^{-1}).
\]
This uses no limit theorem: the fourth-cumulant estimate was proved
directly by support counting.
All three hypotheses of Theorem~\ref{thm:functional-de-jong-recalled}
hold for every member of the fixed finite family, proving
\eqref{eq:chaos-prefix-functional-limit}.
\end{proof}

The preceding application is for fixed \(d\). It supplies no control of
constants as \(d\) grows; the degree-sensitive estimates earlier in this
appendix are still needed in Section~8.

For use with the full Laguerre statistic, define
\begin{equation}
\mathcal G_{N,a}^{\poly}
=\operatorname{Poly}_{N,T}
\bigl[(\mathcal G_{n,a})_{n_-\leq n\leq n_+}\bigr].
\label{eq:chaos-top-polygonal-process}
\end{equation}
The standard Rademacher hypercontractive inequality
\cite[Theorem~9.21]{ODonnell2014} states that, for a Walsh polynomial \(F\) and
\(q\geq2\),
\begin{equation}
\|F\|_{L^q(\Omega)}\leq(q-1)^{\deg(F)/2}\|F\|_{L^2(\Omega)}.
\label{eq:chaos-Rademacher-hypercontractivity}
\end{equation}
We now transfer the top-chaos processes to \(X_n(L_a)\).
For this purpose, set
\begin{equation}
\mathcal R_{n,a}^{\mathrm{pc}}
=X_n(L_a)-\mathcal G_{n,a}.
\label{eq:chaos-pure-chaos-remainder}
\end{equation}
By \eqref{eq:chaos-pure-chaos-lower-degree-decomposition}, \(\mathcal R_{n,a}^{\mathrm{pc}}\) is a Walsh polynomial of degree at most
\(2a-2\).

\begin{theorem}[Uniform approximation on the length window]
\label{thm:uniform-pure-chaos-transfer}
Fix \(a\geq2\) and \(T>0\).  Then
\begin{equation}
\left\|
\max_{n_-\leq n\leq n_+}
\left|
X_n(L_a)-\mathcal G_{n,a}
\right|
\right\|_{L^2}
\longrightarrow0.
\label{eq:chaos-uniform-pure-chaos-transfer}
\end{equation}
Equivalently,
\begin{equation}
\left\|
X_{N,L_a}^{\poly}
-
\mathcal G_{N,a}^{\poly}
\right\|_{C(I_T)}
\longrightarrow0
\label{eq:chaos-polygonal-pure-chaos-transfer}
\end{equation}
in \(L^2(\Omega)\), where
\begin{equation}
X_{N,L_a}^{\poly}
=
\operatorname{Poly}_{N,T}
\bigl[
(X_n(L_a))_{n_-\leq n\leq n_+}
\bigr].
\label{eq:chaos-X-Laguerre-polygonal}
\end{equation}
\end{theorem}

\begin{proof}
Fix an integer \(q>1\).  By Theorem \ref{thm:fixed-mode-pure-chaos-reduction}
and \eqref{eq:chaos-Rademacher-hypercontractivity},
\[
\E\left|
\mathcal R_{n,a}^{\mathrm{pc}}
\right|^{2q}
\leq C_{a,q}n^{-q}.
\]
Therefore,
\[
\begin{aligned}
\E\left[
\max_{n_-\leq n\leq n_+}
\left|\mathcal R_{n,a}^{\mathrm{pc}}
\right|^2\right]
\leq\left(
\sum_{n=n_-}^{n_+}
\E\left|
\mathcal R_{n,a}^{\mathrm{pc}}
\right|^{2q}\right)^{1/q}
\leq C_{a,q,T}N^{1/q-1}.
\end{aligned}
\]
This proves
\eqref{eq:chaos-uniform-pure-chaos-transfer}.
Since the polygonal interpolation of the difference is a convex combination of its endpoint values,
\eqref{eq:chaos-polygonal-pure-chaos-transfer}
follows as well.
\end{proof}

\begin{corollary}
\label{cor:finite-vector-pure-chaos-transfer}
For every fixed \(d\geq2\) and \(T>0\),
\begin{equation}
\max_{2\leq a\leq d}
\left\|
X_{N,L_a}^{\poly}
-
\mathcal G_{N,a}^{\poly}
\right\|_{C(I_T)}
\longrightarrow0
\label{eq:chaos-finite-vector-transfer}
\end{equation}
in probability.
\end{corollary}

\begin{proof}
For each fixed \(2\leq a\leq d\), Theorem \ref{thm:uniform-pure-chaos-transfer}
gives convergence in \(L^2(\Omega)\), hence in probability.
Since \(d\) is fixed, taking the maximum over the finitely many statistics preserves convergence to zero in probability.
\end{proof}

\begin{corollary}
\label{cor:full-Laguerre-maximal-jumps}
Fix \(a\geq2\) and \(T>0\).  Then
\begin{equation}
\max_{n_-\leq n<n_+}
\left|X_{n+1}(L_a)-X_n(L_a)\right|
\longrightarrow0
\label{eq:chaos-full-Laguerre-maximal-jumps}
\end{equation}
in probability.
The same conclusion holds jointly for every fixed finite family \(2\leq a\leq d\).
\end{corollary}

\begin{proof}
By
\eqref{eq:chaos-pure-chaos-remainder},
\[\max_{n_-\leq n<n_+}
\left|X_{n+1}(L_a)-X_n(L_a)\right|
\leq\max_{n_-\leq n<n_+}
\left|\mathcal G_{n+1,a}-\mathcal G_{n,a}
\right|
+2\max_{n_-\leq n\leq n_+}
\left|\mathcal R_{n,a}^{\mathrm{pc}}\right|.
\]
The first term tends to zero in probability by \eqref{eq:finite-mode-de-jong-jumps}, while the second tends to zero in
probability by
\eqref{eq:chaos-uniform-pure-chaos-transfer}.
This proves the scalar assertion.
Since \(d\) is fixed, the corresponding joint assertion follows by a finite union bound.
\end{proof}

\subsection{Mixed hypercontractivity}
\label{subs:app-A-finite-dimensional-homogeneous-vector}

We record the standard hypercontractive estimate for products of
Rademacher and Gaussian measures used in Section~8. It follows from
\cite[Propositions~3.11--3.12 and Theorem~3.13]{MosselODonnellOleszkiewicz2010}.

\begin{lemma}
\label{lem:mixed-hypercontractivity}
Let \(Y_1,\ldots,Y_m\) be independent random variables, each of which is
either a Rademacher variable or a standard Gaussian variable.
If \(Q(Y_1,\ldots,Y_m)\) is multilinear, then, for every \(q\geq2\),
\begin{equation}
\|Q\|_{L^q}
\leq(q-1)^{\deg(Q)/2}\|Q\|_{L^2}.
\label{eq:chaos-mixed-hypercontractivity}
\end{equation}
\end{lemma}

\begin{proof}
Theorem~3.13 of \cite{MosselODonnellOleszkiewicz2010} gives the
\((2,q,(q-1)^{-1/2})\)-hypercontractive inequality for each coordinate.
Independence and Proposition~3.11 of the same reference give this
inequality for their product. Proposition~3.12 then gives the displayed
degree bound.
\end{proof}

\section{Finite-dimensional comparison and endpoint moment estimates}
\label{app:finite-replica-occupation}

To control nested polynomial lengths, we first need a moment bound at
one fixed length.
We begin by comparing finite collections of polynomial values with
Gaussian variables in Section~\ref{subsec:coordinate-path-comparison}.
Weighted derivative bounds, a square-function inequality on the
discrete cube, and an exact coordinate expansion then give the endpoint
moment estimate in Theorem~\ref{thm:compact-endpoint-moment}.
Here an endpoint estimate concerns the statistic at a single polynomial
length.
Appendix~\ref{app:block-addition-refinement} combines these fixed-length
bounds with increment estimates to control maxima over nested lengths.

\subsection{Finite-dimensional comparison along coordinate paths}
\label{subsec:coordinate-path-comparison}

Next we compare finite collections of polynomial values with Gaussian vectors while varying one coefficient and keeping the normalization fixed.
\begin{definition}
Fix \(n\geq2\).
For \(0\leq k<n\) and \(-1\leq t\leq1\), define
\begin{equation}
Z_{n;k,t}(\theta)
=
\frac1{\sqrt n}
\left(
\sum_{\substack{0\leq j<n\\j\neq k}}
\varepsilon_je^{ij\theta}
+
te^{ik\theta}
\right),
\label{eq:occupation-fixed-denominator-path}
\end{equation}
and put
\[
S_{n;k,t}(\theta)=|Z_{n;k,t}(\theta)|^2.
\]
\end{definition}
Thus
\[
Z_{n;k,\varepsilon_k}=Z_n,
\qquad
S_{n;k,\varepsilon_k}=S_n.
\]
Here \(Z_n(\theta)=q_n(e^{i\theta})/\sqrt n\).
The denominator \(\sqrt n\) remains fixed as \(t\) varies.

Let \((G_j)_{j\geq0}\) and \((H_j)_{j\geq0}\)
be independent families of independent standard
real Gaussian variables, independent of the
Rademacher sequence, and put
\begin{equation}
\xi_j=\frac{G_j+iH_j}{\sqrt2}.
\label{eq:replica-circular-coefficients}
\end{equation}
For \(0\leq k<n\), define the fields with
the \(k\)-th coordinate omitted:
\[
\begin{aligned}
Z_{n;k}^{\mathrm R}(\theta)
&=
\frac1{\sqrt n}
\sum_{\substack{0\leq j<n\\j\neq k}}
\varepsilon_je^{ij\theta},\\
Z_{n;k}^{\mathrm M}(\theta)
&=
\frac1{\sqrt n}
\sum_{\substack{0\leq j<n\\j\neq k}}
G_je^{ij\theta},\\
Z_{n;k}^{\circ}(\theta)
&=
\frac1{\sqrt n}
\sum_{\substack{0\leq j<n\\j\neq k}}
\xi_je^{ij\theta}.
\end{aligned}
\]
The superscripts \(\mathrm R,\mathrm M,\circ\) indicate the Rademacher field,
the matched Gaussian field with real coefficients, and the circular complex
Gaussian field, respectively. All three fields are complex-valued. The real
covariance matrices of the \(\mathrm M\) and \(\mathrm R\) models agree.

The covariance kernel of each field with the
\(k\)-th coordinate omitted is
\begin{equation}
\mathscr C_n^{(-k)}(x)
=
\frac1n
\sum_{\substack{0\leq j<n\\j\neq k}}
e^{ijx}.
\label{eq:occupation-deleted-kernel}
\end{equation}
The standard Dirichlet-kernel bound gives,
uniformly in \(0\leq k<n\),
\begin{equation}
\|\mathscr C_n^{(-k)}\|_{L^1(\Torus)}
\leq
C\frac{\log(n+2)}n.
\label{eq:occupation-kernel-L1}
\end{equation}
Moreover,
\begin{equation}
1-\E|Z_{n;k,t}(\theta)|^2
=
\frac{1-t^2}{n},
\label{eq:occupation-variance-defect}
\end{equation}
uniformly for \(0\leq k<n\),
\(-1\leq t\leq1\), and \(\theta\in\Torus\).

For an integer \(r\geq1\) and
\(\boldsymbol\theta=(\theta_1,\ldots,\theta_r)
\in\Torus^r\), write
\[
\mathbf Z_{n;k}^{\diamond}(\boldsymbol\theta)
=
\bigl(
Z_{n;k}^{\diamond}(\theta_1),\ldots,
Z_{n;k}^{\diamond}(\theta_r)
\bigr),
\qquad
\diamond\in\{\mathrm R,\mathrm M,\circ\}.
\]
We identify \(\C^r\) with \(\R^{2r}\).
For \(F\in C_b^4(\R^{2r};\C)\),
let \(\mathrm D^jF\) denote its \(j\)-th
real Fr\'echet differential, and put
\[
\|\mathrm D^jF\|_\infty
=
\sup_{x\in\R^{2r}}
\sup_{\substack{|v_1|\leq1,\ldots,|v_j|\leq1}}
|\mathrm D^jF(x)[v_1,\ldots,v_j]|,
\qquad 1\leq j\leq4.
\]

The next lemma is the smooth covariance-path form of Price's formula;
see \cite[Theorem~1]{VoigtlaenderPrice} for the general formulation on
positive definite covariance matrices.  We include the regularization
argument because singular covariance matrices occur at the endpoints
of our comparisons.

\begin{lemma}
\label{lem:Gaussian-covariance-differentiation}
Let \(\Sigma_0,\Sigma_1\) be symmetric positive semidefinite
\(d\times d\) matrices, and let
\(X_\tau\sim N(0,\Sigma_\tau)\), where
\[
\Sigma_\tau=(1-\tau)\Sigma_0+\tau\Sigma_1,
\qquad
0\leq\tau\leq1.
\]
If \(H\in C^2(\mathbb R^d;\mathbb C)\) has bounded second
derivatives, then
\begin{equation}
\frac{\dd}{\dd\tau}\E[H(X_\tau)]
=
\frac12
\sum_{\alpha,\beta=1}^{d}
(\Sigma_1-\Sigma_0)_{\alpha\beta}
\E[\partial_{\alpha\beta}H(X_\tau)].
\label{eq:replica-Gaussian-covariance-differentiation}
\end{equation}
At \(\tau=0\) and \(\tau=1\), the derivative is understood
as the right and left derivative, respectively.
\end{lemma}

\begin{proof}
Let \(Z\sim N(0,I_d)\) and, for \(\delta>0\), set
\[
X_\tau^\delta
=
(\Sigma_\tau+\delta I_d)^{1/2}Z.
\]
Price's formula \cite[Theorem~1]{VoigtlaenderPrice}, applied to the real and imaginary parts of \(H\),
gives
\[
\begin{aligned}
\E[H(X_t^\delta)]-\E[H(X_0^\delta)]
=
\frac12
\sum_{\alpha,\beta=1}^{d}
(\Sigma_1-\Sigma_0)_{\alpha\beta}
\int_0^t
\E[\partial_{\alpha\beta}H(X_\tau^\delta)]
\,\dd\tau .
\end{aligned}
\]
Since the second derivatives of \(H\) are bounded,
\(|H(x)|\leq C_H(1+|x|^2)\).
The continuity of the matrix square root and dominated convergence
therefore allow \(\delta\downarrow0\) in the preceding identity.

Moreover,
\[
\tau\longmapsto
\E[\partial_{\alpha\beta}H(\Sigma_\tau^{1/2}Z)]
\]
is continuous on \([0,1]\) by bounded convergence.
The fundamental theorem of calculus now yields
\eqref{eq:replica-Gaussian-covariance-differentiation},
with the corresponding one-sided derivatives at the endpoints.
\end{proof}

The following estimate combines the classical Lindeberg coordinate
replacement method (see, for example, \cite{Chatterjee2006Lindeberg})
with Price's formula.  Matching the first three moments gives a
fourth-order Taylor remainder.  We give that calculation explicitly,
together with the deleted-coordinate and angular covariance bounds
that determine the three errors below.

\begin{proposition}
\label{prop:finite-dimensional-comparison}
Fix an integer \(r\geq1\).
There exists \(C_r<\infty\) such that, for every
\(n\geq2\), \(0\leq k<n\), \(-1\leq t\leq1\),
and every deterministic
\(F\in C_b^4(\R^{2r};\C)\),
\begin{equation}
\begin{aligned}
&\left|
\E\int_{\Torus^r}
F\bigl(
Z_{n;k,t}(\theta_1),\ldots,
Z_{n;k,t}(\theta_r)
\bigr)\,\dd\mu^{\otimes r}
-
\E F(\zeta_1,\ldots,\zeta_r)
\right|\\
&\qquad\leq
C_r\left(
n^{-1/2}\|\mathrm DF\|_\infty
+
\frac{\log(n+2)}n\|\mathrm D^2F\|_\infty
+
n^{-1}\|\mathrm D^4F\|_\infty
\right),
\end{aligned}
\label{eq:occupation-comparison-derivative-bound}
\end{equation}
where \(\zeta_1,\ldots,\zeta_r\) are independent
standard circular complex Gaussian variables.
\end{proposition}

\begin{proof}
Fix \(\boldsymbol\theta\in\Torus^r\).
The coordinate path may be written as
\[
\bigl(
Z_{n;k,t}(\theta_1),\ldots,
Z_{n;k,t}(\theta_r)
\bigr)
=
\mathbf Z_{n;k}^{\mathrm R}(\boldsymbol\theta)
+
\frac{t}{\sqrt n}
\bigl(e^{ik\theta_1},\ldots,e^{ik\theta_r}\bigr).
\]
Since
\[
\Big|\frac{t}{\sqrt n}
\bigl(e^{ik\theta_1},\ldots,e^{ik\theta_r}\bigr)\Big|
\leq\sqrt{\frac rn},
\]
the mean-value theorem gives
\[
\left|
\E F\Bigl(
\mathbf Z_{n;k}^{\mathrm R}(\boldsymbol\theta)
+\frac{t}{\sqrt n}
\bigl(e^{ik\theta_1},\ldots,e^{ik\theta_r}\bigr)
\Bigr)
-
\E F\bigl(
\mathbf Z_{n;k}^{\mathrm R}(\boldsymbol\theta)
\bigr)
\right|
\leq
C_rn^{-1/2}\|\mathrm DF\|_\infty.
\]

We first replace the Rademacher coefficients by standard real Gaussian coefficients.
For \(j\neq k\), put
\[
v_j=
\frac1{\sqrt n}
\bigl(
\cos(j\theta_1),\sin(j\theta_1),\ldots,
\cos(j\theta_r),\sin(j\theta_r)
\bigr)\in\R^{2r}.
\]
Then
\[
|v_j|^4\leq\frac{r^2}{n^2}.
\]
Define
\[
W_j
=
\sum_{\substack{0\leq\ell<j\\\ell\neq k}}
G_\ell v_\ell
+
\sum_{\substack{j<\ell<n\\\ell\neq k}}
\varepsilon_\ell v_\ell.
\]
The vector \(W_j\) is independent of
\((\varepsilon_j,G_j)\).
Since their moments agree through order three,
Taylor's formula with fourth-order remainder gives
\[
\begin{aligned}
\left|
\E F(W_j+\varepsilon_jv_j)
-
\E F(W_j+G_jv_j)
\right|
\leq
C\bigl(
\E|\varepsilon_j|^4+\E|G_j|^4
\bigr)
|v_j|^4\|\mathrm D^4F\|_\infty
\leq
C_rn^{-2}\|\mathrm D^4F\|_\infty.
\end{aligned}
\]
Summing over the \(n-1\) coordinates other
than \(k\), we obtain
\[
\left|
\E F(\mathbf Z_{n;k}^{\mathrm R}(\boldsymbol\theta))
-
\E F(\mathbf Z_{n;k}^{\mathrm M}(\boldsymbol\theta))
\right|
\leq
C_rn^{-1}\|\mathrm D^4F\|_\infty.
\]

We next compare the matched Gaussian field \(\mathrm M\), whose coefficients
are real, with the circular complex Gaussian field.
Their covariance kernels agree:
\[
\E\left[
Z_{n;k}^{\mathrm M}(\theta_\nu)
\overline{Z_{n;k}^{\mathrm M}(\theta_\mu)}
\right]
=
\E\left[
Z_{n;k}^{\circ}(\theta_\nu)
\overline{Z_{n;k}^{\circ}(\theta_\mu)}
\right]
=
\mathscr C_n^{(-k)}(\theta_\nu-\theta_\mu).
\]
Their pseudocovariances, however, satisfy
\[
\begin{aligned}
\E\left[
Z_{n;k}^{\mathrm M}(\theta_\nu)
Z_{n;k}^{\mathrm M}(\theta_\mu)
\right]
&=
\mathscr C_n^{(-k)}(\theta_\nu+\theta_\mu),\\
\E\left[
Z_{n;k}^{\circ}(\theta_\nu)
Z_{n;k}^{\circ}(\theta_\mu)
\right]
&=0.
\end{aligned}
\]
Hence Lemma~\ref{lem:Gaussian-covariance-differentiation}
gives
\[
\left|
\E F(\mathbf Z_{n;k}^{\mathrm M}(\boldsymbol\theta))
-
\E F(\mathbf Z_{n;k}^{\circ}(\boldsymbol\theta))
\right|
\leq
C_r\|\mathrm D^2F\|_\infty
\sum_{\nu,\mu=1}^{r}
\left|
\mathscr C_n^{(-k)}(\theta_\nu+\theta_\mu)
\right|.
\]
By Haar invariance, including the case \(\nu=\mu\), and \eqref{eq:occupation-kernel-L1},
\[
\int_{\Torus^r}
\left|
\E F(\mathbf Z_{n;k}^{\mathrm M}(\boldsymbol\theta))
-
\E F(\mathbf Z_{n;k}^{\circ}(\boldsymbol\theta))
\right|
\,\dd\mu^{\otimes r}
\leq
C_r\frac{\log(n+2)}n
\|\mathrm D^2F\|_\infty.
\]

Finally, let \(\zeta_1,\ldots,\zeta_r\) be
independent standard circular complex Gaussian
variables.
The off-diagonal covariance blocks of
\(\mathbf Z_{n;k}^{\circ}(\boldsymbol\theta)\)
are controlled by
\[
\mathscr C_n^{(-k)}(\theta_\nu-\theta_\mu),
\qquad \nu\neq\mu,
\]
whereas the diagonal covariance blocks differ
from those of \((\zeta_1,\ldots,\zeta_r)\)
by \(O(n^{-1})\).
All pseudocovariances vanish.
Thus Lemma~\ref{lem:Gaussian-covariance-differentiation}
yields
\[
\left|
\E F(\mathbf Z_{n;k}^{\circ}(\boldsymbol\theta))
-
\E F(\zeta_1,\ldots,\zeta_r)
\right|
\leq
C_r\|\mathrm D^2F\|_\infty
\left[
\frac1n+
\sum_{\nu\neq\mu}
\left|
\mathscr C_n^{(-k)}(\theta_\nu-\theta_\mu)
\right|
\right].
\]
Integrating over \(\boldsymbol\theta\)
and using \eqref{eq:occupation-kernel-L1}
gives
\[
\int_{\Torus^r}
\left|
\E F(\mathbf Z_{n;k}^{\circ}(\boldsymbol\theta))
-
\E F(\zeta_1,\ldots,\zeta_r)
\right|
\,\dd\mu^{\otimes r}
\leq
C_r\frac{\log(n+2)}n
\|\mathrm D^2F\|_\infty.
\]

Combining the coordinate translation,
the Lindeberg replacement, and the two Gaussian
comparisons, then integrating over
\(\boldsymbol\theta\), proves
\eqref{eq:occupation-comparison-derivative-bound}.
\end{proof}

\subsection{Weighted derivative bounds}
\label{subsec:weighted-derivative-bounds}

We use Gaussian comparison and Laguerre energies to control the weighted derivative terms in the fixed-length moment estimate.
Indeed, we apply the finite-dimensional comparison from \eqref{eq:occupation-comparison-derivative-bound}
to the weighted integrals needed in the endpoint argument.

For \(\psi\in C_c^\infty([0,\infty);\mathbb R)\), set
\begin{equation}
\Phi_\psi(z)
=
\psi'(|z|^2)z,
\qquad
z\in\C\simeq\R^2.
\label{eq:replica-field-multiplier}
\end{equation}
\begin{definition}
For an integer \(K\geq0\), define
\begin{equation}
\begin{aligned}
\mathfrak J_K(\psi)
=
1
+\max_{0\leq j\leq K+1}
\|\psi^{(j)}\|_{L^\infty([0,\infty))}
+\max_{|\beta|\leq K}
\|\partial^\beta\Phi_\psi\|_{L^\infty(\R^2)}.
\end{aligned}
\label{eq:replica-compact-jet-seminorm}
\end{equation}
\end{definition}
Thus \(\mathfrak J_K(\psi)\geq1\). This quantity controls the supremum
norms of the displayed derivatives.
Here \(\beta=(\beta_1,\beta_2)\) is a pair of nonnegative integers,
and \(|\beta|=\beta_1+\beta_2\).  The derivatives
\(\partial^\beta\Phi_\psi\) are taken in the real coordinates of
\(\C\simeq\R^2\), with the Euclidean norm on the values of \(\Phi_\psi\).

Recall the Laguerre energies
\[
\mathfrak E_r(\psi)
=
\sum_{a\geq1}
a^r|\widehat\psi_a|^2.
\]

\begin{proposition}
\label{prop:weighted-derivative-bounds}
Fix an even integer \(q\geq4\).
There exists \(C_q<\infty\) such that, for every
\(n\geq2\) and every deterministic
\(\psi\in C_c^\infty([0,\infty);\mathbb R)\),
\begin{equation}
\left\|
\int_{\Torus}
S_n(\theta)
|\psi'(S_n(\theta))|^2
\,\dd\mu(\theta)
\right\|_{L^{q/2}}
\leq
\mathfrak E_1(\psi)
+
C_qn^{-1/q}\mathfrak J_5(\psi)^2.
\label{eq:occupation-first-weighted-bound}
\end{equation}
Moreover,
\begin{equation}
\left\|
\int_{\Torus}
S_n(\theta)^2
|\psi''(S_n(\theta))|^2
\,\dd\mu(\theta)
\right\|_{L^{q/2}}
\leq
\mathfrak E_2(\psi)
+
C_qn^{-1/q}\mathfrak J_5(\psi)^2.
\label{eq:occupation-second-weighted-bound}
\end{equation}
\end{proposition}

\begin{proof}
Set \(r=q/2\), and define
\[
\Psi_1(x)=x|\psi'(x)|^2,
\qquad
\Psi_2(x)=x^2|\psi''(x)|^2.
\]
By \eqref{eq:spectral-D1-Dirichlet-form},
\begin{equation}
\int_0^\infty
\Psi_1(x)e^{-x}\,\dd x
=
\int_0^\infty
x|\psi'(x)|^2e^{-x}\,\dd x
=
\mathfrak E_1(\psi).
\label{eq:occupation-first-weight-mean}
\end{equation}
We use \(L_j^{(\eta)}\) for the generalized Laguerre polynomial with
parameter \(\eta\), in the normalization of
\cite[Eq.~(18.5.12)]{NISTDLMF}; thus \(L_j^{(0)}=L_j\).
The identity
\[
L_a''=L_{a-2}^{(2)},
\qquad a\geq2,
\]
and the associated Laguerre orthogonality give
\begin{equation}
\int_0^\infty
\Psi_2(x)e^{-x}\,\dd x
=
\int_0^\infty
x^2|\psi''(x)|^2e^{-x}\,\dd x
=
\sum_{a\geq2}
a(a-1)|\widehat\psi_a|^2
\leq
\mathfrak E_2(\psi).
\label{eq:occupation-second-weight-mean}
\end{equation}

For \(i=1,2\), put
\[
g_i(z)=\Psi_i(|z|^2),
\qquad
F_i(z_1,\ldots,z_r)
=
\prod_{\nu=1}^{r}g_i(z_\nu).
\]
Writing \(z=x+iy\) and
\(\partial_{\bar z}=\frac12(\partial_x+i\partial_y)\),
we have
\[
g_1(z)=|\Phi_\psi(z)|^2,
\qquad
\partial_{\bar z}\Phi_\psi(z)
=
\psi''(|z|^2)z^2,
\qquad
g_2(z)=|\partial_{\bar z}\Phi_\psi(z)|^2.
\]
Leibniz's rule and
\eqref{eq:replica-compact-jet-seminorm} give
\[
\max_{\substack{0\leq j\leq4\\i=1,2}}
\|\mathrm D^j g_i\|_\infty
\leq
C\mathfrak J_5(\psi)^2,
\]
and consequently
\begin{equation*}
\max_{\substack{0\leq j\leq4\\i=1,2}}
\|\mathrm D^jF_i\|_\infty
\leq
C_r\mathfrak J_5(\psi)^{2r}.
\end{equation*}

Applying \eqref{eq:occupation-comparison-derivative-bound}
to \(F_i\), for \(i=1,2\), and using the independence
and exponential distributions of
\(|\zeta_1|^2,\ldots,|\zeta_r|^2\), we obtain
\[
\left|
\E\left(
\int_{\Torus}
\Psi_i(S_{n;k,t}(\theta))
\,\dd\mu(\theta)
\right)^r
-
\left(
\int_0^\infty
\Psi_i(x)e^{-x}\,\dd x
\right)^r
\right|
\leq
C_rn^{-1/2}\mathfrak J_5(\psi)^{2r},
\]
uniformly for \(0\leq k<n\) and \(-1\leq t\leq1\).
Here we used
\(\log(n+2)/n\leq Cn^{-1/2}\).

Fix \(0\leq k<n\).
Since
\[
S_n=S_{n;k,\varepsilon_k}
\]
and \(S_{n;k,t}\) is independent of
\(\varepsilon_k\) for each fixed \(t\),
conditioning on \(\varepsilon_k\) gives
\[
\begin{aligned}
\E\left(
\int_{\Torus}
\Psi_i(S_n(\theta))\,\dd\mu(\theta)
\right)^r
&=
\frac12\sum_{t\in\{-1,1\}}
\E\left(
\int_{\Torus}
\Psi_i(S_{n;k,t}(\theta))
\,\dd\mu(\theta)
\right)^r\\
&\leq
\left(
\int_0^\infty
\Psi_i(x)e^{-x}\,\dd x
\right)^r
+
C_rn^{-1/2}\mathfrak J_5(\psi)^{2r}.
\end{aligned}
\]

Since \(\Psi_i\geq0\), taking \(r\)-th roots
and recalling that \(q=2r\), we obtain
\[
\left\|
\int_{\Torus}
\Psi_i(S_n(\theta))\,\dd\mu(\theta)
\right\|_{L^{q/2}}
\leq
\int_0^\infty
\Psi_i(x)e^{-x}\,\dd x
+
C_qn^{-1/q}\mathfrak J_5(\psi)^2,
\qquad i=1,2.
\]
The Laguerre identities
\eqref{eq:occupation-first-weight-mean} and
\eqref{eq:occupation-second-weight-mean}
now give
\eqref{eq:occupation-first-weighted-bound}
and \eqref{eq:occupation-second-weighted-bound}.
\end{proof}

\begin{remark}
The arguments of Propositions~\ref{prop:finite-dimensional-comparison}
and~\ref{prop:weighted-derivative-bounds} also apply uniformly
to an interpolation between two polynomial lengths.

Fix \(T>0\), and let \(m\leq n\) with
\(m,n\asymp_T N\). For \(0\leq s\leq1\), put
\[
\lambda_j(s)=
\begin{cases}
1,&0\leq j<m,\\
s,&m\leq j<n,
\end{cases}
\qquad
\sigma_s^2=m+s^2(n-m).
\]
Define
\[
\mathsf Z_s(\theta)
=
\frac1{\sigma_s}
\sum_{j=0}^{n-1}
\lambda_j(s)\varepsilon_je^{ij\theta},
\qquad
\mathsf Z_{s;k,t}(\theta)
=
\frac1{\sigma_s}
\left[
\sum_{\substack{0\leq j<n\\j\neq k}}
\lambda_j(s)\varepsilon_je^{ij\theta}
+
t\lambda_k(s)e^{ik\theta}
\right],
\]
and write
\[
\mathsf Q_s=|\mathsf Z_s|^2,
\qquad
\mathsf Q_{s;k,t}=|\mathsf Z_{s;k,t}|^2.
\]

For every fixed integer \(r\geq1\), the comparison
estimate
\begin{equation*}
\begin{aligned}
&\left|
\E\int_{\Torus^r}
F\bigl(
\mathsf Z_{s;k,t}(\theta_1),\ldots,
\mathsf Z_{s;k,t}(\theta_r)
\bigr)
\,\dd\mu^{\otimes r}
-
\E F(\zeta_1,\ldots,\zeta_r)
\right|\\
&\qquad\leq
C_{r,T}\left[
N^{-1/2}\|\mathrm DF\|_\infty
+
\frac{\log(N+2)}N\|\mathrm D^2F\|_\infty
+
N^{-1}\|\mathrm D^4F\|_\infty
\right]
\end{aligned}
\end{equation*}
holds for every deterministic
\(F\in C_b^4(\R^{2r};\C)\),
uniformly in \(m,n,s,k,t\), where
\(0\leq k<n\), \(-1\leq t\leq1\), and
\(\zeta_1,\ldots,\zeta_r\) are as in
Proposition~\ref{prop:finite-dimensional-comparison}.

Indeed, uniformly in \(s\) and \(k\),
\[
\sigma_s^2\asymp_TN,
\qquad
\frac{\lambda_k(s)}{\sigma_s}
\leq C_TN^{-1/2},
\qquad
\frac1{\sigma_s^4}
\sum_{j=0}^{n-1}\lambda_j(s)^4
\leq C_TN^{-1}.
\]
The covariance kernel with the \(k\)-th
coordinate omitted is
\[
\mathscr C_s^{(-k)}(x)
=
\frac1{\sigma_s^2}
\sum_{\substack{0\leq j<n\\j\neq k}}
\lambda_j(s)^2e^{ijx}.
\]
It satisfies
\[
\|\mathscr C_s^{(-k)}\|_{L^1(\Torus)}
\leq
C_T\frac{\log(N+2)}N,
\]
by the Dirichlet-kernel bound.
Moreover,
\[
1-\E|\mathsf Z_{s;k,t}(\theta)|^2
=
\frac{(1-t^2)\lambda_k(s)^2}{\sigma_s^2}.
\]
These estimates give the same bounds for the
coordinate translation, the Lindeberg replacement,
and the two Gaussian covariance comparisons
as in Proposition~\ref{prop:finite-dimensional-comparison}.

Consequently, for every even integer \(q\geq4\)
and every deterministic
\(\psi\in C_c^\infty([0,\infty);\R)\),
the weighted bounds also extend uniformly:
\[
\begin{aligned}
\left\|
\int_{\Torus}
\mathsf Q_{s;k,t}
|\psi'(\mathsf Q_{s;k,t})|^2
\,\dd\mu
\right\|_{L^{q/2}}
&\leq
\mathfrak E_1(\psi)
+
C_{q,T}N^{-1/q}\mathfrak J_5(\psi)^2,
\\[4pt]
\left\|
\int_{\Torus}
\mathsf Q_{s;k,t}^2
|\psi''(\mathsf Q_{s;k,t})|^2
\,\dd\mu
\right\|_{L^{q/2}}
&\leq
\mathfrak E_2(\psi)
+
C_{q,T}N^{-1/q}\mathfrak J_5(\psi)^2.
\end{aligned}
\]
The same inequalities hold with
\(\mathsf Q_{s;k,t}\) replaced by
\(\mathsf Q_s\), by conditioning on
the two values of \(\varepsilon_k\).
\end{remark}

\subsection{Coordinate differences and exact sign flips}
\label{subsec:coordinate-differences-sign-flips}

Now we compute exactly how flipping one sign changes the statistic, using derivatives along the corresponding coordinate path.
We use a standard \(L^q\)-Poincar\'e inequality on the discrete cube and an exact Fourier representation of coordinate sign flips.

For a function \(F=F(\varepsilon_0,\ldots,\varepsilon_{n-1})\), let
\(\varepsilon^{(k)}\) be obtained by changing only the sign
\(\varepsilon_k\).  Define the coordinate half-difference
\begin{equation}
\Delta_kF
=
\frac12
\left[
F(\varepsilon)-F(\varepsilon^{(k)})
\right].
\label{eq:endpoint-coordinate-half-difference}
\end{equation}
The dimension-free \(L^q\)-Poincar\'e inequality on the discrete cube
(see \cite[Theorem~1.1]{BenEfraimLustPiquard2008}
or \cite[Eqs.~(2), (3), and (5)]{CorderoErausquinEskenazis2024})
gives, for every \(2\leq q<\infty\),
\begin{equation}
\|F-\E F\|_{L^q}^2
\leq
C_q
\left\|
\sum_{k=0}^{n-1}
|\Delta_kF|^2
\right\|_{L^{q/2}},
\label{eq:endpoint-product-space-square-function}
\end{equation}
where \(C_q\) is independent of \(n\) and \(F\).
The gradient convention in \cite[Eqs.~(2)--(3)]{CorderoErausquinEskenazis2024}
is exactly the half-difference used here; the full-difference convention
in \cite{BenEfraimLustPiquard2008} changes only \(C_q\).
Thus this is an \(L^q\) moment inequality for every finite \(q\geq2\).

We use the coordinate path \(Z_{n;k,t}\) defined
in \eqref{eq:occupation-fixed-denominator-path}.
For \(\psi\in C_c^\infty([0,\infty);\mathbb R)\),
recall from \eqref{eq:intro-proof-Xn} the centered statistic
\[
X_n(\psi)
=
\sqrt n\left[
\int_{\Torus}\psi(S_n(\theta))\,\dd\mu(\theta)
-
\E\int_{\Torus}\psi(S_n(\theta))\,\dd\mu(\theta)
\right].
\]
For \(f\in L^2(\Torus,\dd\mu)\), we use the Fourier convention
\begin{equation}
\widehat f(u)
=
\int_{\Torus}
f(\theta)e^{-iu\theta}\,\dd\mu(\theta),
\qquad
u\in\Z.
\label{eq:replica-Fourier-convention}
\end{equation}

The coordinate half-differences of \(X_n(\psi)\) admit the following exact Fourier representation along this path.
The same coordinate path is used in Lemma~\ref{lem:gae-sign-flip}
to obtain a moment bound through Bessel's inequality.

\begin{lemma}
\label{lem:exact-endpoint-sign-flip}
For every
\(\psi\in C_c^\infty([0,\infty);\mathbb R)\)
and \(0\leq k<n\),
\begin{equation}
\Delta_kX_n(\psi)
=
\varepsilon_k\operatorname{Re}
\int_{-1}^{1}
\widehat{\Phi_\psi(Z_{n;k,t})}(k)
\,\dd t.
\label{eq:endpoint-exact-sign-flip}
\end{equation}
\end{lemma}

\begin{proof}
Since
\[
\frac{\partial}{\partial t}
Z_{n;k,t}(\theta)
=
\frac{e^{ik\theta}}{\sqrt n},
\]
the chain rule gives
\[
\frac{\dd}{\dd t}
\left[
\sqrt n
\int_{\Torus}
\psi(|Z_{n;k,t}(\theta)|^2)
\,\dd\mu(\theta)
\right]
=
2\operatorname{Re}
\widehat{\Phi_\psi(Z_{n;k,t})}(k).
\]
The endpoint values at \(t=\pm1\) correspond to
the two possible signs of the \(k\)-th coordinate.
Their coordinate half-difference is
\(\varepsilon_k/2\) times the difference of these
endpoint values, while the deterministic
centering term cancels.
The fundamental theorem of calculus proves
\eqref{eq:endpoint-exact-sign-flip}.
\end{proof}

\subsection{Expansion about the full field and derivative bounds}
\label{subsec:common-field-expansion}

Next we expand the sign-flip formula around the full field and control the resulting derivative terms for the moment estimate.
The calculations in this subsection use the classical Laguerre
differential and differentiation identities
\cite[\S\S18.8--18.9]{NISTDLMF}, the radial chain rule, and Taylor's
formula.  We retain the explicit identities because their normalization,
weighted derivative orders, and sign-flip factors enter the endpoint
moment bound.  We first record the lowering identities needed for the
zero-frequency term, then expand the coordinate half-differences and control
the second-order remainder.

The scalar zero-frequency contribution is expressed by the lowering
operator
\begin{equation}
\mathscr R\psi
=
\psi'+x\psi''.
\label{eq:endpoint-resonant-lowering}
\end{equation}

\begin{lemma}
\label{lem:resonant-Laguerre-lowering}
For every integer \(a\geq1\),
\begin{equation}
\mathscr RL_a=-aL_{a-1}.
\label{eq:endpoint-resonant-Laguerre}
\end{equation}
Moreover, for every integer \(j\geq0\) and
\(\psi\in C_c^\infty([0,\infty);\mathbb R)\),
\begin{equation}
\mathfrak E_1(\mathscr R^j\psi)
\leq C_j\mathfrak E_{2j+1}(\psi),
\qquad
\mathfrak E_2(\mathscr R^j\psi)
\leq C_j\mathfrak E_{2j+2}(\psi).
\label{eq:endpoint-resonant-energy-transfer}
\end{equation}
Finally,
\begin{equation}
\int_0^\infty
\mathscr R\psi(x)e^{-x}\,\dd x
=
-\widehat\psi_1.
\label{eq:endpoint-resonant-exponential-mean}
\end{equation}
Here \(\widehat\psi_1
=\int_0^\infty\psi(x)(1-x)e^{-x}\,\dd x\)
is the first Laguerre coefficient of \(\psi\).
\end{lemma}

\begin{proof}
The Laguerre equation
\[
xL_a''+(1-x)L_a'+aL_a=0
\]
and the identity
\[
xL_a'=a(L_a-L_{a-1})
\]
give \eqref{eq:endpoint-resonant-Laguerre}.
Iterating this identity yields
\[
\mathscr R^jL_a
=
(-1)^j\frac{a!}{(a-j)!}L_{a-j},
\qquad a\geq j,
\]
with zero when \(a<j\).
Consequently, for \(p=1,2\),
\[
\mathfrak E_p(\mathscr R^j\psi)
=
\sum_{a\geq j+1}
(a-j)^p
\left(\frac{a!}{(a-j)!}\right)^2
|\widehat\psi_a|^2
\leq
\mathfrak E_{p+2j}(\psi),
\]
which proves \eqref{eq:endpoint-resonant-energy-transfer}.

Finally,
\[
\mathscr R\psi
=
(x\psi')'.
\]
Integration by parts yields
\[
\begin{aligned}
\int_0^\infty
\mathscr R\psi(x)e^{-x}\,\dd x
&= \int_0^\infty
x\psi'(x)e^{-x}\,\dd x\\
&= -\int_0^\infty
\psi(x)(1-x)e^{-x}\,\dd x\\
&= -\widehat\psi_1.
\end{aligned}
\]
\end{proof}
\begin{lemma}
\label{lem:resonant-lowering-jet-bounds}
For integers \(K,J\geq0\), set \(K'=2K+2J+1\).
Then
\begin{equation}
\max_{0\leq j\leq J}
\mathfrak J_K(\mathscr R^j\psi)
\leq
C_{K,J}\mathfrak J_{K'}(\psi)
\label{eq:endpoint-resonant-jet-transfer}
\end{equation}
for every \(\psi\in C_c^\infty([0,\infty);\mathbb R)\).
The constant \(C_{K,J}\) is independent of \(\psi\) and its
support.
Moreover, for every integer \(j\geq1\),
\begin{equation}
\|\mathscr R^j\psi\|_\infty
\leq
\mathfrak J_{2j-1}(\psi).
\label{eq:endpoint-resonant-supremum-transfer}
\end{equation}
\end{lemma}

\begin{proof}
For a smooth scalar function \(f\), write
\(g_f(z)=f(|z|^2)\), with \(z\in\mathbb R^2\).
The identities
\[
\nabla_z g_f=2\Phi_f,
\qquad
\Delta_z g_f=4g_{\mathscr R f}
\]
imply
\[
\Phi_{\mathscr R^j\psi}
=
4^{-j}\Delta_z^j\Phi_\psi
\qquad (j\geq0),
\]
and
\[
g_{\mathscr R^j\psi}
=
\frac{1}{2\,4^{j-1}}
\nabla_z\cdot
\bigl(\Delta_z^{j-1}\Phi_\psi\bigr)
\qquad (j\geq1).
\]

Set
\(h_f(t)=g_f(t,0)=f(t^2)\).
For every integer \(a\geq1\), repeated integration by parts
yields
\[
f^{(a)}(t^2)
=
\frac{1}{2^{2a-1}(a-1)!}
\int_0^1
(1-u^2)^{a-1}h_f^{(2a)}(ut)\,du,
\qquad t\geq0.
\]
Consequently,
\[
\|f^{(a)}\|_\infty
\leq
C_a\|\partial_1^{2a}g_f\|_\infty,
\qquad a\geq0,
\]
where the case \(a=0\) is immediate.
Thus, for \(j\geq1\),
\[
\|(\mathscr R^j\psi)^{(a)}\|_\infty
\leq
C_{a,j}
\max_{|\beta|\leq2a+2j-1}
\|\partial^\beta\Phi_\psi\|_\infty,
\]
while
\[
\max_{|\beta|\leq K}
\|\partial^\beta\Phi_{\mathscr R^j\psi}\|_\infty
\leq
C_{K,j}
\max_{|\beta|\leq K+2j}
\|\partial^\beta\Phi_\psi\|_\infty.
\]
For \(a\leq K+1\) and \(j\leq J\), all derivative orders on
the right are at most \(K'=2K+2J+1\).
Together with the case \(j=0\), this proves
\eqref{eq:endpoint-resonant-jet-transfer}.

Finally, expanding \(\Delta_z^{j-1}\) in the divergence
identity yields
\[
\|\mathscr R^j\psi\|_\infty
\leq
2^{1-j}
\max_{|\beta|=2j-1}
\|\partial^\beta\Phi_\psi\|_\infty
\leq
\mathfrak J_{2j-1}(\psi),
\]
which proves
\eqref{eq:endpoint-resonant-supremum-transfer}.
\end{proof}

For \(0\leq k<n\) and \(-1\leq t\leq1\), write
\begin{equation}
Z_{n;k,t}
=
Z_n+\eta_{k,t},
\qquad
\eta_{k,t}(\theta)
=
\frac{t-\varepsilon_k}{\sqrt n}e^{ik\theta}.
\label{eq:endpoint-common-field-perturbation}
\end{equation}

\begin{lemma}
\label{lem:exact-first-order-endpoint-expansion}
For every
\(\psi\in C_c^\infty([0,\infty);\mathbb R)\)
and \(0\leq k<n\),
\begin{equation}
\begin{aligned}
\Delta_kX_n(\psi)
&=
2\varepsilon_k\operatorname{Re}
\widehat{\Phi_\psi(Z_n)}(k)
-
\frac2{\sqrt n}
\int_{\Torus}
\mathscr R\psi(S_n(\theta))\,\dd\mu(\theta)\\
&\quad-
\frac2{\sqrt n}
\operatorname{Re}
\widehat{\psi''(S_n)Z_n^2}(2k)
+
\mathcal E_{n,k}^{(2)}(\psi).
\end{aligned}
\label{eq:endpoint-exact-first-order-expansion}
\end{equation}
Here,
\begin{equation}
\mathcal E_{n,k}^{(2)}(\psi)
=
\varepsilon_k\operatorname{Re}
\int_{-1}^{1}\int_0^1
(1-\tau)
\int_{\Torus}
\mathrm D^2\Phi_\psi
\bigl(Z_n(\theta)+\tau\eta_{k,t}(\theta)\bigr)
[\eta_{k,t}(\theta),\eta_{k,t}(\theta)]
e^{-ik\theta}
\,\dd\mu(\theta)\,\dd\tau\,\dd t.
\label{eq:endpoint-second-order-remainder}
\end{equation}
\end{lemma}

\begin{proof}
By Taylor's formula in \(\R^2\),
\[
\Phi_\psi(Z_{n;k,t})
=
\Phi_\psi(Z_n)
+
\mathrm D\Phi_\psi(Z_n)[\eta_{k,t}]
+
\int_0^1(1-\tau)
\mathrm D^2\Phi_\psi(Z_n+\tau\eta_{k,t})
[\eta_{k,t},\eta_{k,t}]
\,\dd\tau.
\]
The real differential of \(\Phi_\psi\) is
\[
\mathrm D\Phi_\psi(z)[v]
=
\mathscr R\psi(|z|^2)v
+
\psi''(|z|^2)z^2\overline v.
\]
Consequently,
\[
\widehat{
\mathrm D\Phi_\psi(Z_n)[\eta_{k,t}]
}(k)
=
\frac{t-\varepsilon_k}{\sqrt n}
\left[
\int_{\Torus}
\mathscr R\psi(S_n)\,\dd\mu
+
\widehat{\psi''(S_n)Z_n^2}(2k)
\right].
\]
Combining with
Lemma~\ref{lem:exact-endpoint-sign-flip} and
\[
\varepsilon_k
\int_{-1}^{1}(t-\varepsilon_k)\,\dd t=-2
\]
proves \eqref{eq:endpoint-exact-first-order-expansion}.

\end{proof}

\begin{lemma}
\label{lem:endpoint-second-order-ledger}
Fix an even integer \(q\geq4\).
There exists \(C_q<\infty\) such that, for every
\(n\geq2\), every deterministic
\(\psi\in C_c^\infty([0,\infty);\mathbb R)\),
and every deterministic
\(\mathcal I\subseteq\{0,\ldots,n-1\}\),
\begin{equation}
\left\|
\sum_{k\in\mathcal I}
|\mathcal E_{n,k}^{(2)}(\psi)|^2
\right\|_{L^{q/2}}
\leq
\frac{C_q|\mathcal I|}{n^2}
\left[
\mathcal H(\psi)
+
n^{-1/q}\mathfrak J_6(\psi)^2
\right].
\label{eq:endpoint-second-order-ledger}
\end{equation}
Moreover,
\begin{equation}
\mathcal H(\psi)
:=
\int_0^\infty
\left(
x|\psi''(x)|^2+x^3|\psi'''(x)|^2
\right)e^{-x}\,\dd x
\leq
3\mathfrak E_3(\psi).
\label{eq:endpoint-H-energy-bound}
\end{equation}
\end{lemma}

\begin{proof}
Put \(r=q/2\), and write
\[
h_\psi(z)
=
\sum_{a,b=1}^{2}
|\partial_a\partial_b\Phi_\psi(z)|^2,
\qquad
F(z_1,\ldots,z_r)
=
\prod_{\nu=1}^{r}h_\psi(z_\nu).
\]
Direct differentiation and
\eqref{eq:replica-compact-jet-seminorm} give
\[
\begin{aligned}
\int_{\R^2}
h_\psi(z)e^{-|z|^2}\,\frac{\dd z}{\pi}
&\leq
C\mathcal H(\psi),\\
\|F\|_\infty+
\max_{1\leq j\leq4}
\|\mathrm D^jF\|_\infty
&\leq
C_q\mathfrak J_6(\psi)^{2r}.
\end{aligned}
\]
Since \(\zeta_1,\ldots,\zeta_r\) are
independent standard circular complex Gaussian
variables, it follows that
\[
\E F(\zeta_1,\ldots,\zeta_r)
=
\left(
\int_{\R^2}
h_\psi(z)e^{-|z|^2}\,\frac{\dd z}{\pi}
\right)^r
\leq
C_q\mathcal H(\psi)^r.
\]

Applying \eqref{eq:occupation-comparison-derivative-bound}
to \(F\), we obtain
\[
\begin{aligned}
\sup_{\substack{0\leq k<n\\|v|\leq1}}
\E\left[
\left(
\int_{\Torus}
h_\psi(Z_{n;k,v}(\theta))
\,\dd\mu(\theta)
\right)^r
\right]
&\leq
C_q\mathcal H(\psi)^r
+
C_q\left[
n^{-1/2}
+\frac{\log(n+2)}n+n^{-1}
\right]
\mathfrak J_6(\psi)^{2r}\\
&\leq
C_q\left[
\mathcal H(\psi)^r
+n^{-1/2}\mathfrak J_6(\psi)^{2r}
\right].
\end{aligned}
\]
Here we used
\(\log(n+2)/n\leq Cn^{-1/2}\).

Taking \(r\)-th roots and using \(q=2r\)
gives
\begin{equation}
\sup_{\substack{0\leq k<n\\|v|\leq1}}
\left\|
\int_{\Torus}
h_\psi(Z_{n;k,v}(\theta))
\,\dd\mu(\theta)
\right\|_{L^{q/2}}
\leq
C_q\left[
\mathcal H(\psi)
+
n^{-1/q}\mathfrak J_6(\psi)^2
\right].
\label{eq:endpoint-Hessian-occupation}
\end{equation}

We now apply this estimate to the Taylor remainder.
For \(\sigma\in\{-1,1\}\),
\(-1\leq t\leq1\), and \(0\leq\tau\leq1\), put
\[
v_{\sigma,t,\tau}
=
\sigma+\tau(t-\sigma).
\]
Then \(|v_{\sigma,t,\tau}|\leq1\), and
\eqref{eq:endpoint-common-field-perturbation}
gives
\[
Z_n+\tau\eta_{k,t}
=
Z_{n;k,v_{\varepsilon_k,t,\tau}}.
\]

For each deterministic \(v\),
\(Z_{n;k,v}\) is independent of \(\varepsilon_k\).
Thus, for fixed \(t,\tau\), conditioning on
\(\varepsilon_k\) gives
\[
\begin{aligned}
\E\left[
\left(
\int_{\Torus}
h_\psi(Z_{n;k,v_{\varepsilon_k,t,\tau}}(\theta))
\,\dd\mu(\theta)
\right)^r
\right]
&=
\frac12\sum_{\sigma\in\{-1,1\}}
\E\left[
\left(
\int_{\Torus}
h_\psi(Z_{n;k,v_{\sigma,t,\tau}}(\theta))
\,\dd\mu(\theta)
\right)^r
\right]\\
&\leq
C_q\left[
\mathcal H(\psi)^r
+
n^{-1/2}\mathfrak J_6(\psi)^{2r}
\right].
\end{aligned}
\]
Consequently, by
\eqref{eq:endpoint-Hessian-occupation},
\[
\sup_{\substack{0\leq k<n\\
-1\leq t\leq1\\0\leq\tau\leq1}}
\left\|
\int_{\Torus}
h_\psi(Z_n+\tau\eta_{k,t})
\,\dd\mu
\right\|_{L^{q/2}}
\leq
C_q\left[
\mathcal H(\psi)
+
n^{-1/q}\mathfrak J_6(\psi)^2
\right].
\]

For \(z,u\in\R^2\), the Cauchy--Schwarz
inequality gives
\[
|\mathrm D^2\Phi_\psi(z)[u,u]|
\leq
|u|^2h_\psi(z)^{1/2}.
\]
Since
\[
|\eta_{k,t}(\theta)|^2
=
\frac{(t-\varepsilon_k)^2}{n}
\leq\frac4n,
\]
the remainder formula
\eqref{eq:endpoint-second-order-remainder}
and Cauchy--Schwarz on \(\Torus\) imply
\[
\begin{aligned}
|\mathcal E_{n,k}^{(2)}(\psi)|
&\leq
\frac Cn
\int_{-1}^{1}\int_0^1
(1-\tau)
\left(
\int_{\Torus}
h_\psi(Z_n+\tau\eta_{k,t})
\,\dd\mu
\right)^{1/2}
\,\dd\tau\,\dd t.
\end{aligned}
\]
Minkowski's inequality and the preceding
estimate therefore give
\[
\begin{aligned}
\sup_{0\leq k<n}
\|\mathcal E_{n,k}^{(2)}(\psi)\|_{L^q}^2
\leq
\frac{C_q}{n^2}
\left[
\mathcal H(\psi)
+
n^{-1/q}\mathfrak J_6(\psi)^2
\right].
\end{aligned}
\]
Finally, since \(q/2\geq1\),
\[
\left\|
\sum_{k\in\mathcal I}
|\mathcal E_{n,k}^{(2)}(\psi)|^2
\right\|_{L^{q/2}}
\leq
\sum_{k\in\mathcal I}
\|\mathcal E_{n,k}^{(2)}(\psi)\|_{L^q}^2
\leq
\frac{C_q|\mathcal I|}{n^2}
\left[
\mathcal H(\psi)
+
n^{-1/q}\mathfrak J_6(\psi)^2
\right].
\]
This proves \eqref{eq:endpoint-second-order-ledger}.

Integration by parts, Cauchy--Schwarz, and
Laguerre orthogonality yield
\[
\mathcal H(\psi)
\leq
\int_0^\infty
\bigl(x^2|\psi''(x)|^2
+2x^3|\psi'''(x)|^2\bigr)
e^{-x}\,\dd x
\leq
\mathfrak E_2(\psi)+2\mathfrak E_3(\psi)
\leq
3\mathfrak E_3(\psi).
\]
\end{proof}

\subsection{Endpoint moment estimate}
\label{subsec:endpoint-moment-estimate}

Now we combine the coordinate expansion and weighted derivative estimates of the preceding subsections to obtain a uniform moment bound for the centered statistic \(X_n(\psi)\).

\begin{theorem}
\label{thm:compact-endpoint-moment}
Fix an even integer \(q\geq4\) and an integer \(J\geq3\).
There exist an integer \(K_{q,J}\geq1\) and
constants \(C_q,C_{q,J}<\infty\) such that,
for every \(n\geq2\) and every deterministic
\(\psi\in C_c^\infty([0,\infty);\mathbb R)\),
\begin{equation}
\begin{aligned}
\|X_n(\psi)\|_{L^q}^{2}
\leq
C_q\mathfrak E_1(\psi)
+
C_{q,J}\sum_{j=1}^{J}n^{-j}
\left[
\mathfrak E_{2j}(\psi)
+\mathfrak E_{2j+1}(\psi)
\right]
+
C_{q,J}n^{-1/q}
\mathfrak J_{K_{q,J}}(\psi)^2.
\end{aligned}
\label{eq:endpoint-compact-moment}
\end{equation}
The constant \(C_q\) is independent of \(J\).
\end{theorem}

\begin{proof}
By \eqref{eq:endpoint-product-space-square-function},
\[
\|X_n(\psi)\|_{L^q}^{2}
\leq
C_q\left\|
\sum_{k=0}^{n-1}
|\Delta_kX_n(\psi)|^2
\right\|_{L^{q/2}}.
\]
We estimate the four terms in
\eqref{eq:endpoint-exact-first-order-expansion}.
For the terms involving Fourier coefficients at $k$ and $2k$,
Bessel's inequality gives
\[
\begin{aligned}
\sum_{k=0}^{n-1}
|\widehat{\Phi_\psi(Z_n)}(k)|^2
&\leq
\int_{\Torus}
S_n|\psi'(S_n)|^2\,\dd\mu,\\
\sum_{k=0}^{n-1}
|\widehat{\psi''(S_n)Z_n^2}(2k)|^2
&\leq
\int_{\Torus}
S_n^2|\psi''(S_n)|^2\,\dd\mu.
\end{aligned}
\]
By Proposition~\ref{prop:weighted-derivative-bounds},
these two terms contribute at most
\[
C_q\mathfrak E_1(\psi)
+
\frac{C_q}{n}\mathfrak E_2(\psi)
+
C_qn^{-1/q}\mathfrak J_5(\psi)^2
\]
to the square-function estimate.

The zero-frequency term contributes
\[
\left\|
\sum_{k=0}^{n-1}
\left|
\frac2{\sqrt n}
\int_{\Torus}\mathscr R\psi(S_n)\,\dd\mu
\right|^2
\right\|_{L^{q/2}}
=
4\left\|
\int_{\Torus}\mathscr R\psi(S_n)\,\dd\mu
\right\|_{L^q}^2.
\]
For the Taylor remainder,
\eqref{eq:endpoint-second-order-ledger} with
\(\mathcal I=\{0,\ldots,n-1\}\) gives
\[
\left\|
\sum_{k=0}^{n-1}
|\mathcal E_{n,k}^{(2)}(\psi)|^2
\right\|_{L^{q/2}}
\leq
\frac{C_q}{n}
\left[
\mathcal H(\psi)
+
n^{-1/q}\mathfrak J_6(\psi)^2
\right]
\leq
\frac{C_q}{n}\mathfrak J_6(\psi)^2.
\]
For the last inequality, note that
\[
\mathcal H(\psi)
\leq
\|\psi''\|_\infty^2
+
6\|\psi'''\|_\infty^2
\leq
7\mathfrak J_6(\psi)^2.
\]
Combining these estimates, using
\(\mathfrak J_5(\psi)\leq\mathfrak J_6(\psi)\)
and \(n^{-1}\leq n^{-1/q}\), we obtain
\[
\|X_n(\psi)\|_{L^q}^{2}
\leq
C_q\mathfrak E_1(\psi)
+
\frac{C_q}{n}\mathfrak E_2(\psi)
+
C_q\left\|
\int_{\Torus}
\mathscr R\psi(S_n)\,\dd\mu
\right\|_{L^q}^{2}
+
C_qn^{-1/q}\mathfrak J_6(\psi)^2.
\]

We next control the zero-frequency term directly by a product moment.
Put
\[
A_n=\int_{\Torus}\mathscr R\psi(S_n)\,\dd\mu,
\qquad H(z)=\mathscr R\psi(|z|^2),
\]
and
\[
F(z_1,\ldots,z_q)=\prod_{\nu=1}^qH(z_\nu).
\]
Since
\[
H(z)=\frac12\nabla_z\cdot\Phi_\psi(z),
\]
the product rule and \eqref{eq:replica-compact-jet-seminorm} give
\[
\max_{0\leq j\leq4}\|\mathrm D^jF\|_\infty
\leq C_q\mathfrak J_5(\psi)^q.
\]
In particular, \(F\in C_b^4(\R^{2q};\R)\).
Fix one coordinate \(0\leq k<n\).
Use this same deleted coordinate in every factor and apply
\eqref{eq:occupation-comparison-derivative-bound} with \(r=q\).
For each common endpoint \(t\in\{-1,1\}\), this yields
\[
\begin{aligned}
\left|
\E\int_{\Torus^q}
\prod_{\nu=1}^q
H\bigl(Z_{n;k,t}(\theta_\nu)\bigr)
\,\dd\mu^{\otimes q}(\boldsymbol\theta)
-\left(\int_0^\infty\mathscr R\psi(x)e^{-x}\,\dd x\right)^q
\right|\leq C_qn^{-1/2}\mathfrak J_5(\psi)^q,
\end{aligned}
\]
where we used \(\log(n+2)/n\leq Cn^{-1/2}\).
The Gaussian target factors because its \(q\) coordinates are
independent.
Since \(S_n=S_{n;k,\varepsilon_k}\), averaging these two comparisons
over the one common sign \(t\) recovers \(\E A_n^q\) exactly.
By \eqref{eq:endpoint-resonant-exponential-mean}, therefore,
\[
\left|\E A_n^q-(-\widehat\psi_1)^q\right|
\leq C_qn^{-1/2}\mathfrak J_5(\psi)^q.
\]
The integer \(q\) is even, so \(\E A_n^q=\E|A_n|^q\) and
\((-\widehat\psi_1)^q=|\widehat\psi_1|^q\).
Using \((x+y)^{2/q}\leq x^{2/q}+y^{2/q}\) for \(x,y\geq0\),
we obtain
\[
\|A_n\|_{L^q}^2
\leq |\widehat\psi_1|^2
+C_qn^{-1/q}\mathfrak J_5(\psi)^2
\leq\mathfrak E_1(\psi)
+C_qn^{-1/q}\mathfrak J_5(\psi)^2.
\]
Substituting this bound into the preceding square-function estimate
gives the stronger inequality
\[
\|X_n(\psi)\|_{L^q}^2
\leq C_q\mathfrak E_1(\psi)
+\frac{C_q}{n}\mathfrak E_2(\psi)
+C_qn^{-1/q}\mathfrak J_6(\psi)^2.
\]
We may take \(K_{q,J}=2J+11\).
Since \(\mathfrak J_6(\psi)\leq\mathfrak J_{K_{q,J}}(\psi)\)
and the sum in \eqref{eq:endpoint-compact-moment} contains
\(n^{-1}\mathfrak E_2(\psi)\), this proves the stated bound.
The leading constant depends only on \(q\).
\end{proof}

\section{Two-time increment estimates and temporal refinement}
\label{app:block-addition-refinement}

This appendix turns estimates at a single polynomial length into bounds
for maxima over nested lengths.
We combine the fixed-length moment bounds of
Appendix~\ref{app:finite-replica-occupation} with estimates for the changes
caused by adding coefficients.
The Gaussian calculations use the classical complex Wick and Mehler
identities, whose normalization we record below.
We then control the Fourier projections and the errors in Gaussian
replacement, and sum the increment bounds over successively finer grids
of polynomial lengths.
This gives the maximal estimate in
Theorem~\ref{thm:refinement-safe-maximal-interface}, which is applied in
Section~\ref{sec:active-abel-range}.

\begin{figure}[htbp]
  \centering
  \resizebox{.96\textwidth}{!}{%
\begin{tikzpicture}[
  x=1cm,y=1cm,
  cmBox/.style={draw=black!65,line width=.45pt,rounded corners=1.3pt,
    align=flush center,inner xsep=5pt,inner ysep=4pt,
    font=\fontsize{8.2}{10}\selectfont,
    execute at begin node={\hyphenpenalty=10000\exhyphenpenalty=10000}},
  cmInput/.style={cmBox},
  cmStep/.style={cmBox},
  cmConclusion/.style={cmBox},
  cmArrow/.style={-{Stealth[length=2.1mm,width=1.4mm]},
    draw=black!75,line width=.5pt}
]

\node[cmStep,text width=14.4cm,minimum height=1.5cm]
(decomposition) at (5.7,-1.9) {
  {\fontsize{9}{11}\selectfont\bfseries Decompose the increment}\\[4pt]
  $X_n(u)-X_m(u)
      =J_{m,n}^{\mathrm{lin}}(u)+W_{m,n}(u)+V_{m,n}(u).$\\[4pt]
  {\fontsize{7.5}{9.5}\selectfont
  From Proposition~\ref{prop:block-adaptive-residual-hierarchy}
  and~\eqref{eq:direct-conditional-mean-identity}}
};

\node[cmInput,text width=4.2cm,minimum height=1.65cm]
(linear) at (0.8,-4.0) {
  {\fontsize{9}{11}\selectfont\bfseries
  Linear term $\|J_{m,n}^{\mathrm{lin}}(u)\|_{L^q}$}\\[4pt]
  {\fontsize{7.5}{9.5}\selectfont
  Corollary~\ref{cor:direct-leading-leakage}}
};

\node[cmInput,text width=4.2cm,minimum height=1.65cm]
(nonlinear) at (5.7,-4.0) {
  {\fontsize{9}{11}\selectfont\bfseries
  Centered nonlinear \mbox{term $\|W_{m,n}(u)\|_{L^q}$}}\\[4pt]
  {\fontsize{7.5}{9.5}\selectfont
  Proposition~\ref{prop:direct-original-nonlinear-fixed-pair}}
};

\node[cmInput,text width=4.2cm,minimum height=1.65cm]
(mean) at (10.6,-4.0) {
  {\fontsize{9}{11}\selectfont\bfseries
  Change of conditional \mbox{mean $\|V_{m,n}(u)\|_{L^q}$}}\\[4pt]
  {\fontsize{7.5}{9.5}\selectfont
  Proposition~\ref{prop:direct-conditional-mean-drift}}
};

\node[cmStep,text width=9.5cm,minimum height=1.65cm]
(residual) at (8.15,-6.6) {
  {\fontsize{9}{11}\selectfont\bfseries Estimate for $\|X_n(u)-X_m(u)-J_{m,n}^{\mathrm{lin}}(u)\|_{L^q}$}
  \\[4pt]
  {\fontsize{7.5}{9.5}\selectfont
  Proposition~\ref{prop:block-adaptive-residual-hierarchy}}
};

\node[cmInput,text width=3cm,minimum height=2cm]
(endpoint) at (0,-9.3) {
  {\fontsize{9}{11}\selectfont\bfseries Single-length estimate}\\[4pt]
  $\|X_n(u)\|_{L^q}$\\[4pt]
  {\fontsize{7.5}{9.5}\selectfont
  Lemma~\ref{lem:direct-homogeneous-endpoint};\\[2pt]
  error~\eqref{eq:direct-endpoint-remainder}}
};

\node[cmInput,text width=10.6cm,minimum height=2cm]
(local) at (7.6,-9.3) {
  {\fontsize{9}{11}\selectfont\bfseries
  Estimate for $\|X_n(u)-X_m(u)\|_{L^q}$}\\[3pt]
  {\fontsize{7.5}{9.5}\selectfont
  Corollary~\ref{cor:direct-exact-energy-increment};
  error~\eqref{eq:block-direct-local-error}}
};

\node[cmStep,text width=14.4cm,minimum height=1.2cm]
(maxima) at (5.7,-11.9) {
  {\fontsize{9}{11}\selectfont\bfseries
  Apply the maximum bound on the length grids}\\[3pt]
  {\fontsize{7.5}{9.5}\selectfont
  Lemmas~\ref{lem:block-finite-family-maximum}
  and~\ref{lem:active-scale-A-refinement-tree}}
};

\node[cmConclusion,text width=14.4cm,minimum height=1.6cm]
(maximum) at (5.7,-14.0) {
  {\fontsize{9}{11}\selectfont\bfseries Maximum over all lengths}\\[4pt]
  $\bigl\|\max_{n_-\le n\le n_+}|X_n(u)|\bigr\|_{L^2}$\\[4pt]
  {\fontsize{7.5}{9.5}\selectfont
  Theorem~\ref{thm:refinement-safe-maximal-interface}}
};

\draw[-, draw=black!75, line width=.5pt]
  (linear.north     |- decomposition.south) -- (linear.north)
  (nonlinear.north  |- decomposition.south) -- (nonlinear.north)
  (mean.north       |- decomposition.south) -- (mean.north);
\draw[cmArrow] (nonlinear.south) -- ([xshift=-2.45cm]residual.north);
\draw[cmArrow] (mean.south) -- ([xshift=2.45cm]residual.north);
\draw[cmArrow] (linear.south) -- (0.8,-7.8) -- (3.8,-7.8)
  -- ([xshift=-3.8cm]local.north);
\draw[cmArrow] (residual.south) -- ([xshift=.55cm]local.north);
\draw[cmArrow] (endpoint.south) -- ([xshift=-5.7cm]maxima.north);
\draw[cmArrow] (local.south) -- ([xshift=1.9cm]maxima.north);
\draw[cmArrow] (maxima.south) -- (maximum.north);
\end{tikzpicture}

}
  \caption{Proof structure for
  Theorem~\ref{thm:refinement-safe-maximal-interface}.}
  \label{fig:local-to-maximal-roadmap}
\end{figure}
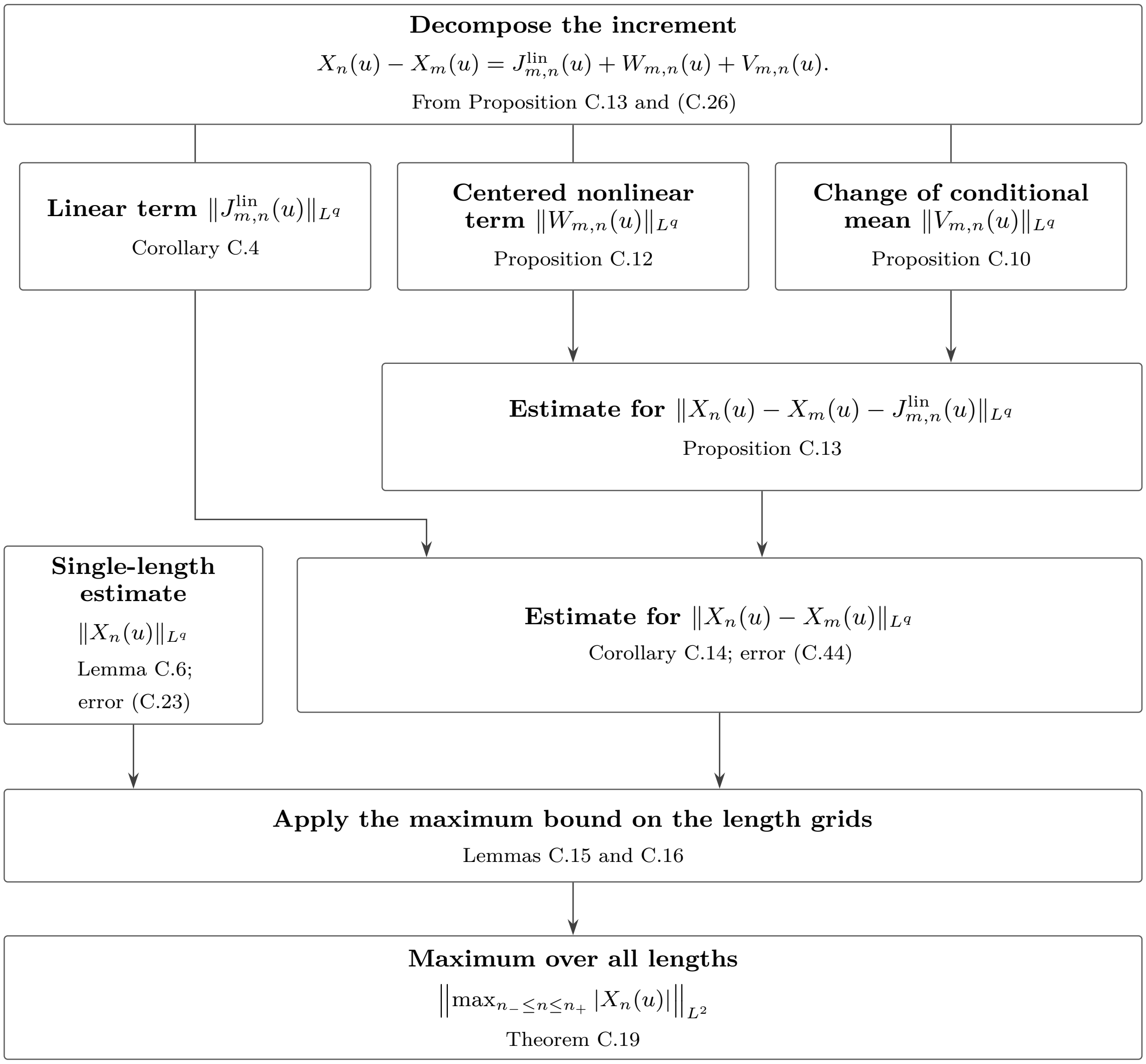

Throughout this appendix, \(\Pi_I\), \(\Phi_u\), \(\mathfrak J_K(u)\),
\(\mathsf C_s\), and \(\mathscr Ru\) have the meanings fixed in
\eqref{eq:replica-Fourier-projection},
\eqref{eq:replica-field-multiplier},
\eqref{eq:replica-compact-jet-seminorm},
\eqref{eq:replica-composition-family}, and
\eqref{eq:endpoint-resonant-lowering}, respectively.

Fix \(T>0\) and integers
\begin{equation}
\lfloor Ne^{-T}\rfloor
\leq m<n
\leq\lfloor Ne^T\rfloor .
\label{eq:block-window-lengths}
\end{equation}
Put
\begin{equation}
h=n-m,
\qquad
\rho_{m,n}=\frac{h}{n},
\qquad
r_{m,n}=\frac{m}{n}=1-\rho_{m,n}.
\label{eq:block-h-rho}
\end{equation}

Let
\begin{equation}
\Delta q_{m,n}(e^{i\theta})
=
\sum_{j=m}^{n-1}
\varepsilon_j e^{ij\theta}.
\label{eq:block-added-polynomial}
\end{equation}

For an integrable random variable \(Y\), write
\begin{equation}
[Y]_0=Y-\E Y.
\label{eq:block-centering-bracket}
\end{equation}

Fix an even integer
\begin{equation}
q=2r,
\qquad
r\geq6.
\label{eq:block-local-moment}
\end{equation}
Whenever it appears below, \(J\) denotes a fixed integer
with \(J\geq3\). 

For each \(N\), fix a deterministic function
\(u\in C_c^\infty([0,\infty);\mathbb R)\).
The function \(u\) is held fixed as \(m,n\) vary.
Recall the natural filtration
\[
\mathcal F_m
=
\sigma(\varepsilon_0,\ldots,\varepsilon_{m-1}).
\]
For \(\ell\geq1\), put
\begin{equation}
F_\ell(u)=\sqrt\ell I_\ell(u)
=
\sqrt\ell
\int_{\Torus}
u(S_\ell(\theta))
\,\dd\mu(\theta),
\qquad
X_\ell(u)=F_\ell(u)-\E F_\ell(u).
\label{eq:direct-F-X}
\end{equation}
Define
\begin{equation}
\begin{aligned}
U_{m,n}(u)
&= F_n(u)
-
\E(F_n(u)\mid\mathcal F_m),\\
V_{m,n}(u)
&= \left[\E(F_n(u)\mid\mathcal F_m)
-
F_m(u)
\right]_0 .
\end{aligned}
\label{eq:direct-U-V}
\end{equation}
Then
\begin{equation}
X_n(u)-X_m(u)
=
U_{m,n}(u)+V_{m,n}(u).
\label{eq:direct-two-time-conditional-split}
\end{equation}

We control the conditionally centered term \(U_{m,n}\)
using Fourier projection estimates for the field multiplier
\(\Phi_u(z)=u'(|z|^2)z\).
For a complex-valued field \(Z\) on \(\Torus\) and a finite
frequency set \(I\subset\Z\), the corresponding projection
energy is
\[
\|\Pi_I\Phi_u(Z)\|_{L^2(\Torus)}^2,
\]
where \(\Pi_I\) is the orthogonal projection onto
\(\operatorname{span}\{e^{ik\theta}:k\in I\}\).
The conditional mean \(V_{m,n}\) is represented through
the Laguerre semigroup.

For later use, put
\[
\mathcal I_{m,n}=\{m,\ldots,n-1\},
\qquad
Z_\ell=\frac{q_\ell}{\sqrt\ell}.
\]
The term linear in the new signs in the first-order Taylor
expansion of \(\sqrt n\,u(S_n)\) about \(S_m\), after integration,
is
\begin{equation}
J_{m,n}^{\mathrm{lin}}(u)
=
\frac2{\sqrt n}
\operatorname{Re}
\int_{\Torus}
u'(S_m)q_m
\overline{\Delta q_{m,n}}
\,\dd\mu .
\label{eq:block-leading-linearized-term}
\end{equation}
Conditionally on \(\mathcal F_m\), this is a centered linear polynomial in
\(\varepsilon_m,\ldots,\varepsilon_{n-1}\).

\subsection{Fourier projection estimates and endpoint bounds}\label{subsec:block-fourier-projection}

Next we combine Gaussian identities with Fourier projection estimates to obtain the fixed-length bounds needed for the later increment analysis.

\medskip

\boldparagraph{Circular Gaussian Wick identities.}
The radial identity below is the Laguerre form of the Wick covariance
formula used by Sodin and Tsirelson \cite[\S2.2, Eq.~(2.6)]{SodinTsirelson2004}.
We also record its charge-one version.  Here charge one means that
rotating the complex argument multiplies the function by the same phase.
Both identities follow from the complex Wick formula; see \cite[Lemma~2.1.7]{HoughKrishnapurPeresVirag2009}.

\begin{lemma}[Radial and charge-one Gaussian covariances]
\label{lem:circular-Laguerre-Wick}
Let \((Z,W)\) be a centered jointly circular complex Gaussian pair with
\(\E|Z|^2=\E|W|^2=1\), and put \(r=\E[Z\overline W]\).
Here joint circularity means
\(\E Z^2=\E W^2=\E[ZW]=0\).
For all integers \(a,b\geq0\),
\begin{align}
\E[L_a(|Z|^2)L_b(|W|^2)]
&=\mathbf1_{\{a=b\}}|r|^{2a},
\label{eq:circular-radial-Wick}\\
\E\!\left[ZL_a^{(1)}(|Z|^2)
\overline{WL_b^{(1)}(|W|^2)}\right]
&=\mathbf1_{\{a=b\}}(a+1)r|r|^{2a}.
\label{eq:circular-charge-one-Wick}
\end{align}
These formulas include \(|r|=1\).
\end{lemma}

\begin{proof}
Define complex Hermite polynomials by
\[
\exp(tz+s\overline z-ts)
=\sum_{p,q\geq0}H_{p,q}(z)
\frac{t^ps^q}{p!q!}.
\]
The Gaussian exponential formula gives, as an identity of power series,
\[
\E\!\left[
 e^{tZ+s\overline Z-ts}
 e^{v\overline W+wW-vw}
\right]
=\exp(tv r+sw\overline r).
\]
Coefficient comparison therefore yields
\[
\E[H_{p,q}(Z)\overline{H_{p',q'}(W)}]
=\mathbf1_{\{p=p',\ q=q'\}}p!q!r^p\overline r^{\,q}.
\]
The finite Laguerre expansions give
\[
H_{a,a}(z)=(-1)^aa!L_a(|z|^2),\qquad
H_{a+1,a}(z)=(-1)^aa!zL_a^{(1)}(|z|^2).
\]
Substitution proves both identities.  This polynomial calculation uses
no inverse covariance matrix and so also covers degenerate pairs.
\end{proof}

\boldparagraph{Fourier projection estimates on integer intervals.}
We use the Fourier convention
\eqref{eq:replica-Fourier-convention}.
For a finite set \(I\subseteq\mathbb Z\), put
\begin{equation}
\Pi_I f
=
\sum_{u\in I}\widehat f(u)e^{iu\theta},
\qquad
D_I(x)
=
\sum_{u\in I}e^{iux}.
\label{eq:replica-Fourier-projection}
\end{equation}
Then
\begin{equation}
\|\Pi_I f\|_{L^2(\Torus)}^2
=
\int_{\Torus^2}
D_I(\varphi-\theta)
f(\theta)\overline{f(\varphi)}
\,\dd\mu(\theta)\,\dd\mu(\varphi).
\label{eq:replica-projection-kernel-identity}
\end{equation}
For \(\ell\geq1\), write
\[
D_\ell(x)
=
\sum_{j=0}^{\ell-1}e^{ijx},
\qquad
\|D_\ell\|_{L^1(\Torus)}
\leq
C\log(\ell+2).
\]

For an integer interval \(I\subset\Z\) and
\(\diamond\in\{\mathrm R,\mathrm M,\circ\}\), set
\begin{equation}
\begin{aligned}
Z_n^{\mathrm R}(\theta)
&= \frac1{\sqrt n}\sum_{j=0}^{n-1}
\varepsilon_j e^{ij\theta},\\
Z_n^{\mathrm M}(\theta)
&= \frac1{\sqrt n}\sum_{j=0}^{n-1}
G_j e^{ij\theta},\\
Z_n^\circ(\theta)
&= \frac1{\sqrt n}\sum_{j=0}^{n-1}
\xi_j e^{ij\theta},
\end{aligned}
\label{eq:direct-three-fields}
\end{equation}
where \(G_j\) are the real Gaussian coefficients of the matched model
and \(\xi_j\) are the circular complex Gaussian coefficients fixed in
\eqref{eq:replica-circular-coefficients}.

Define
\begin{equation}
\mathcal Q_{n,I}^\diamond(u)
=
\left\|
\Pi_I\Phi_u(Z_n^\diamond)
\right\|_{L^2(\Torus)}^2 .
\label{eq:direct-projected-Q}
\end{equation}

For \(q=2r\) and \(|I|=h\), define
\begin{equation}
\begin{aligned}
B_{n,I}(u)
&= \mathfrak J_4(u)^2
\left[\frac{\log(n+2)\log(h+2)}n
+
\frac1n
\right],\\
P_{n,I}(u)
&= C_q
\left[\frac hn\mathfrak E_1(u)
+
B_{n,I}(u)
\right].
\end{aligned}
\label{eq:direct-projected-B-P}
\end{equation}
\begin{definition}
For an integer \(r\geq0\), let \(\mathsf C_r\) consist of the
empty composition and all ordered tuples
\begin{equation}
\mathbf k=(k_1,\ldots,k_L),
\qquad
k_i\in\{1,2,3,4\},
\qquad
K(\mathbf k):=\sum_{i=1}^{L}k_i\leq r.
\label{eq:replica-composition-family}
\end{equation}
For nonempty \(\mathbf k\), let \(L(\mathbf k)\) denote its length,
and set
\[
L(\mathbf k)=K(\mathbf k)=0
\]
for the empty composition. 
\end{definition}
The tuple entries record successive decreases
in moment order when the recurrence
\eqref{eq:direct-Rademacher-Gaussian-recursion} is iterated. The tuple length counts
the number of error terms selected from that recurrence.

\begin{lemma}
\label{lem:direct-arbitrary-interval-leakage}
For every fixed integer \(r\geq1\),
\begin{equation}
\begin{aligned}
\E
\mathcal Q_{n,I}^{\mathrm R}(u)^r
\leq C_r
\sum_{\mathbf k\in\mathsf C_r}
n^{-L(\mathbf k)}
\mathfrak J_4(u)^{2K(\mathbf k)}
P_{n,I}(u)^{r-K(\mathbf k)} .
\end{aligned}
\label{eq:direct-arbitrary-interval-leakage}
\end{equation}
\end{lemma}

\begin{proof}
Write
\[
h=|I|.
\]
For a deterministic coefficient vector
\(\mathbf x=(x_0,\ldots,x_{n-1})\in\C^n\), define
\[
Z_{\mathbf x}(\theta)
=
\frac1{\sqrt n}
\sum_{j=0}^{n-1}x_je^{ij\theta}.
\]
For the circular field, differentiation of the Laguerre expansion in
\(L^2(xe^{-x}\,\dd x)\) gives the charge-one expansion
\[
\Phi_u(Z)
=
\sum_{a\geq0}
b_a(u)ZL_a^{(1)}(|Z|^2),
\qquad
\sum_{a\geq0}
(a+1)|b_a(u)|^2
=
\mathfrak E_1(u),
\qquad b_a(u)=-\widehat u_{a+1}.
\]
By \eqref{eq:circular-charge-one-Wick}, the covariance of the degree-\(a\)
term is \((a+1)r|r|^{2a}\), with
\(r=D_n(\theta-\varphi)/n\).
Set
\[
p_{a,n}(v)
=
\widehat{
\frac{D_n}{n}
\left|\frac{D_n}{n}\right|^{2a}
}(v).
\]
Since
\[
0\leq p_{a,n}(v)\leq\frac1n,
\]
orthogonality yields
\[
\begin{aligned}
 \E
 \mathcal Q_{n,I}^{\circ}(u)
 = \sum_{a\geq0}
 (a+1)|b_a(u)|^2
 \sum_{v\in I}p_{a,n}(v)
\leq \frac hn\mathfrak E_1(u).
\end{aligned}
\]
Moreover,
\[
\mathbf x
\longmapsto
\left\|
\Pi_I\Phi_u(Z_{\mathbf x})
\right\|_2
\]
has Euclidean Lipschitz constant at most
\[
Cn^{-1/2}\mathfrak J_4(u).
\]
Gaussian concentration therefore yields
\[
\E
\mathcal Q_{n,I}^{\circ}(u)^r
\leq
C_r
\left(
\frac hn\mathfrak E_1(u)
\right)^r
+
C_rn^{-r}\mathfrak J_4(u)^{2r}.
\]

Let
\[
a_\tau
=
\left(\frac{1+\tau}{2}\right)^{1/2},
\qquad
b_\tau
=
\left(\frac{1-\tau}{2}\right)^{1/2},
\qquad
0\leq\tau\leq1,
\]
and let \((g_j)\) and \((g'_j)\) be independent standard real Gaussian
families.  Set
\[
Z_\tau(\theta)
=
\frac1{\sqrt n}
\sum_{j=0}^{n-1}
\bigl(
a_\tau g_j+ib_\tau g'_j
\bigr)e^{ij\theta}.
\]
Thus
\[
Z_0=Z^\circ,
\qquad
Z_1=Z^{\mathrm M},
\]
and
\[
\begin{aligned}
\E
\left[Z_\tau(\theta)
\overline{Z_\tau(\varphi)}
\right]
= \frac{D_n(\theta-\varphi)}{n},
\qquad
\E
\left[Z_\tau(\theta)Z_\tau(\varphi)
\right]
= \frac{\tau D_n(\theta+\varphi)}{n}.
\end{aligned}
\]
Put
\[
\mathcal Q_{n,I,\tau}(u)
=
\left\|
\Pi_I\Phi_u(Z_\tau)
\right\|_2^2.
\]
By \eqref{eq:replica-projection-kernel-identity},
\[
\begin{aligned}
\mathcal Q_{n,I,\tau}(u)
&= \int_{\Torus^2}
D_I(\varphi-\theta)
\Phi_u(Z_\tau(\theta))
\overline{\Phi_u(Z_\tau(\varphi))}
\,\dd\mu(\theta)\,\dd\mu(\varphi).
\end{aligned}
\]
Lemma~\ref{lem:Gaussian-covariance-differentiation} and
\[
\max_{|\boldsymbol{\nu}|\leq2}
\|\partial^{\boldsymbol{\nu}}\Phi_u\|_\infty
\leq
\mathfrak J_4(u)
\]
imply
\[
\begin{aligned}
\left|\frac{\dd}{\dd\tau}
\E
\mathcal Q_{n,I,\tau}(u)
\right|
&\leq C\mathfrak J_4(u)^2
\int_{\Torus^2}
|D_I(\varphi-\theta)|
\frac{|D_n(\theta+\varphi)|}{n}
\,\dd\mu(\theta)\,\dd\mu(\varphi)\\
&\leq C
\frac{\log(n+2)\log(h+2)}{n}
\mathfrak J_4(u)^2.
\end{aligned}
\]
Integration over \(0\leq\tau\leq1\), followed by the preceding circular
estimate and Gaussian concentration, yields
\[
\E
\mathcal Q_{n,I}^{\mathrm M}(u)^r
\leq
C_rP_{n,I}(u)^r.
\]

For \(\mathbf x\in\R^n\), set
\[
\mathcal S_I(\mathbf x)
=
\left\|
\Pi_I\Phi_u(Z_{\mathbf x})
\right\|_{L^2(\Torus)}^2.
\]
Let \(g_0,\ldots,g_{n-1}\) be independent standard real Gaussian
variables, independent of
\(\varepsilon_0,\ldots,\varepsilon_{n-1}\), and set
\begin{equation}
\begin{aligned}
\mathbf X^{(j)}
&= (
\varepsilon_0,\ldots,\varepsilon_{j-1},
g_j,\ldots,g_{n-1}
),
\qquad
0\leq j\leq n.
\end{aligned}
\label{eq:direct-Lindeberg-hybrid}
\end{equation}
Thus
\[
\mathcal S_I(\mathbf X^{(0)})
=
\mathcal Q_{n,I}^{\mathrm M}(u),
\qquad
\mathcal S_I(\mathbf X^{(n)})
=
\mathcal Q_{n,I}^{\mathrm R}(u).
\]

Fix
\[
1\leq s\leq r,
\qquad
0\leq j<n.
\]
Conditioning on all coordinates except the \(j\)-th, set
\[
\begin{aligned}
F_{j,s}(t)
&= \mathcal S_I
(
\varepsilon_0,\ldots,\varepsilon_{j-1},
t,
g_{j+1},\ldots,g_{n-1}
)^s.
\end{aligned}
\]
Since
\[
\E\varepsilon_j^\nu
=
\E g_j^\nu,
\qquad
0\leq\nu\leq3,
\]
the third-order Taylor formula with fourth-order remainder yields
\begin{equation}
\begin{aligned}
\left|\E
\mathcal S_I(\mathbf X^{(j+1)})^s
-
\E
\mathcal S_I(\mathbf X^{(j)})^s
\right|
\leq \frac1{3!}
\int_0^1
(1-\lambda)^3
\E
\left[|F_{j,s}^{(4)}(\lambda\varepsilon_j)|
+
|g_j|^4
|F_{j,s}^{(4)}(\lambda g_j)|
\right]
\,\dd\lambda .
\end{aligned}
\label{eq:direct-one-coordinate-replacement}
\end{equation}

For \(1\leq\nu\leq4\),
\[
\left\|
\partial_{x_j}^{\nu}
\Pi_I\Phi_u(Z_{\mathbf x})
\right\|_2
\leq
C_\nu n^{-\nu/2}\mathfrak J_4(u).
\]
Leibniz's rule and Cauchy--Schwarz yield
\begin{equation}
\begin{aligned}
\left|\partial_{x_j}^4
\mathcal S_I(\mathbf x)^s
\right|
&\leq C_s n^{-2}
\sum_{a=1}^{\min\{4,s\}}
\mathfrak J_4(u)^{2a}
\mathcal S_I(\mathbf x)^{s-a}.
\end{aligned}
\label{eq:direct-fourth-derivative-power-bound}
\end{equation}

The Taylor remainder evaluates a coordinate at an intermediate multiple
of its value. Absorbing the fourth-power Gaussian weight will also enlarge
its variance. We use the following family to keep these changes within the
comparison as the moment order decreases.
For \(0\leq s\leq r\), let \(\mathcal Y_s\) be the family of
independent-coordinate random vectors
\[
\mathbf Y=(c_0\zeta_0,\ldots,c_{n-1}\zeta_{n-1}),
\]
where each \(\zeta_j\) is either \(\varepsilon_j\) or \(g_j\).
These choices and the coefficients \(c_j\) are deterministic.
We require that there exist a set
\(D\subseteq\{0,\ldots,n-1\}\) such that
\begin{equation}
\begin{gathered}
|D|\leq r-s,
\qquad
c_j=1 \quad (j\notin D),\\
0\leq c_j\leq 2^{(r-s)/2}
\quad (j\in D).
\end{gathered}
\label{eq:direct-enlarged-hybrid-family}
\end{equation}
The Rademacher and Gaussian coordinates may occur in any order.
Set
\[
\mathfrak P_s
=
\sup_{\mathbf Y\in\mathcal Y_s}
\E\mathcal S_I(\mathbf Y)^s
\quad (1\leq s\leq r),
\qquad
\mathfrak P_0=1.
\]
These quantities are finite, since
\[
\mathcal S_I(\mathbf x)
\leq
\|\Phi_u\|_\infty^2
\leq
\mathfrak J_4(u)^2.
\]
Throughout the following argument, constants may depend on the
fixed integer \(r\), uniformly for \(0\leq s\leq r\).

We first record a Gaussian density bound.
Let \(\varphi_v\) denote the density of \(N(0,v)\).
Since
\[
\frac{x^4\varphi_1(x)}{\varphi_2(x)}
=
\sqrt2\,x^4e^{-x^2/4}
\]
is bounded, every nonnegative Borel function \(F\) satisfies
\begin{equation}
\E\bigl[|g|^4F(g)\bigr]
\leq
C\,\E F(\sqrt2g),
\qquad g\sim N(0,1).
\label{eq:direct-Gaussian-weight-absorption}
\end{equation}
The same inequality holds conditionally on any random variables
independent of \(g\).

Fix \(1\leq s\leq r\), \(\mathbf Y\in\mathcal Y_s\), and
\(0\leq j<n\).  Write \(\mathbf Y^{[j]}(t)\) for the vector obtained
by replacing its \(j\)-th coordinate by \(t\).
For \(0\leq\lambda\leq1\) and
\(1\leq a\leq\min\{4,s\}\), both
\[
\mathbf Y^{[j]}(\lambda c_j\varepsilon_j)
\quad\text{and}\quad
\mathbf Y^{[j]}(\sqrt2\lambda c_jg_j)
\]
belong to \(\mathcal Y_{s-a}\).
Their exceptional coordinates are contained in
\(D\cup\{j\}\), and
\[
|D\cup\{j\}|
\leq r-s+1
\leq r-(s-a),
\]
while
\[
\sqrt2\lambda c_j
\leq 2^{(r-s+1)/2}
\leq 2^{(r-(s-a))/2}.
\]
All unchanged exceptional coefficients satisfy the latter bound
as well.  Consequently,
\begin{equation}
\begin{aligned}
\E\mathcal S_I
\bigl(\mathbf Y^{[j]}(\lambda c_j\varepsilon_j)\bigr)^{s-a}
&\leq \mathfrak P_{s-a},\\
\E\!\left[
|g_j|^4
\mathcal S_I
\bigl(\mathbf Y^{[j]}(\lambda c_jg_j)\bigr)^{s-a}
\right]
&\leq C\mathfrak P_{s-a}.
\end{aligned}
\label{eq:direct-weighted-hybrid-moment-bound}
\end{equation}
For the second inequality, condition on the other coordinates
and apply \eqref{eq:direct-Gaussian-weight-absorption}.
When \(s=a\), both estimates are understood with
\(\mathcal S_I^0=1\).

The variables \(c_j\varepsilon_j\) and \(c_jg_j\) have matching
moments through order three.  The same Taylor expansion as in
\eqref{eq:direct-one-coordinate-replacement}, together with
\eqref{eq:direct-fourth-derivative-power-bound}, therefore yields
\[
\begin{aligned}
&
\left|
\E\mathcal S_I
\bigl(\mathbf Y^{[j]}(c_j\varepsilon_j)\bigr)^s
-
\E\mathcal S_I
\bigl(\mathbf Y^{[j]}(c_jg_j)\bigr)^s
\right|\\
&\quad\leq
C_s n^{-2}c_j^4
\sum_{a=1}^{\min\{4,s\}}
\mathfrak J_4(u)^{2a}
\int_0^1(1-\lambda)^3
\Bigl\{
\E\mathcal S_I
\bigl(\mathbf Y^{[j]}(\lambda c_j\varepsilon_j)\bigr)^{s-a}\\
&\hspace{14em}
+
\E\!\left[
|g_j|^4
\mathcal S_I
\bigl(\mathbf Y^{[j]}(\lambda c_jg_j)\bigr)^{s-a}
\right]
\Bigr\}\,\dd\lambda\\
&\quad\leq
C_r n^{-2}
\sum_{a=1}^{\min\{4,s\}}
\mathfrak J_4(u)^{2a}\mathfrak P_{s-a},
\end{aligned}
\]
where we used \(c_j^4\leq 2^{2r}\).

Replacing a Rademacher coordinate by a Gaussian coordinate with
the same coefficient \(c_j\) leaves the vector in
\(\mathcal Y_s\).
Thus, starting from any \(\mathbf Y\in\mathcal Y_s\), we may
replace all its Rademacher coordinates successively.
There are at most \(n\) replacements.
The resulting Gaussian vector is
\[
\mathbf Y_{\mathrm G}
=
(c_0g_0,\ldots,c_{n-1}g_{n-1}),
\]
and telescoping yields
\begin{equation}
\begin{aligned}
\E\mathcal S_I(\mathbf Y)^s
\leq
\E\mathcal S_I(\mathbf Y_{\mathrm G})^s
+
\frac{C_r}{n}
\sum_{a=1}^{\min\{4,s\}}
\mathfrak J_4(u)^{2a}\mathfrak P_{s-a}.
\end{aligned}
\label{eq:direct-enlarged-hybrid-comparison}
\end{equation}

It remains to bound the Gaussian vectors obtained after these
replacements, uniformly over \(\mathcal Y_s\).
The corresponding covariance kernel is
\[
\mathscr C_{\mathbf c}(x)
=
\frac1n\sum_{j=0}^{n-1}c_j^2e^{ijx}.
\]
Since only the coordinates in \(D\) have nonunit coefficients,
\[
\left\|
\mathscr C_{\mathbf c}-\frac{D_n}{n}
\right\|_\infty
\leq
\frac1n\sum_{j\in D}|c_j^2-1|
\leq
\frac{C_r}{n}.
\]
For the real Gaussian coefficient fields, the covariance kernel and
the pseudocovariance kernel (formed without complex conjugation) are,
respectively,
\(\mathscr C_{\mathbf c}(\theta-\varphi)\) and
\(\mathscr C_{\mathbf c}(\theta+\varphi)\).
Hence all entries of the real covariance matrix of
\[
\bigl(Z_{\mathbf Y_{\mathrm G}}(\theta),
      Z_{\mathbf Y_{\mathrm G}}(\varphi)\bigr)
\]
differ from their matched Gaussian counterparts by at most
\(C_r/n\).
Applying Lemma~\ref{lem:Gaussian-covariance-differentiation} to the
projection-kernel representation yields
\[
\begin{aligned}
\left|
\E\mathcal S_I(\mathbf Y_{\mathrm G})
-
\E\mathcal Q_{n,I}^{\mathrm M}(u)
\right|
&\leq
\frac{C_r}{n}\mathfrak J_4(u)^2
\int_{\Torus^2}|D_I(\theta-\varphi)|
\,\dd\mu(\theta)\,\dd\mu(\varphi)\\
&\leq
C_r\frac{\log(h+2)}{n}\mathfrak J_4(u)^2.
\end{aligned}
\]
Moreover, the map
\[
\mathbf g\longmapsto
\mathcal S_I(c_0g_0,\ldots,c_{n-1}g_{n-1})^{1/2}
\]
has Euclidean Lipschitz constant at most
\[
C_rn^{-1/2}\mathfrak J_4(u),
\]
because \(\max_jc_j\leq2^{r/2}\).
The matched Gaussian mean bound proved above and Gaussian
concentration therefore imply
\[
\begin{aligned}
\E\mathcal S_I(\mathbf Y_{\mathrm G})^s
&\leq
C_r\left[
\bigl(\E\mathcal S_I(\mathbf Y_{\mathrm G})\bigr)^s
+n^{-s}\mathfrak J_4(u)^{2s}
\right]\\
&\leq C_rP_{n,I}(u)^s,
\qquad 1\leq s\leq r.
\end{aligned}
\]

Taking the supremum in
\eqref{eq:direct-enlarged-hybrid-comparison} now yields
\begin{equation}
\begin{aligned}
\mathfrak P_s
&\leq
C_rP_{n,I}(u)^s
+
\frac{C_r}{n}
\sum_{a=1}^{\min\{4,s\}}
\mathfrak J_4(u)^{2a}\mathfrak P_{s-a},
\qquad 1\leq s\leq r,\\
\mathfrak P_0&=1.
\end{aligned}
\label{eq:direct-Rademacher-Gaussian-recursion}
\end{equation}
Each error term in this recurrence lowers the moment order by an integer
from one to four. Its factor \(n^{-1}\) comes from summing at most \(n\)
coordinate errors of order \(n^{-2}\). Induction on \(s\) in
\eqref{eq:direct-Rademacher-Gaussian-recursion} gives
\[
\mathfrak P_s
\leq
C_r
\sum_{\mathbf k\in\mathsf C_s}
n^{-L(\mathbf k)}
\mathfrak J_4(u)^{2K(\mathbf k)}
P_{n, I}(u)^{s-K(\mathbf k)},
\qquad 0\leq s\leq r.
\]
Since the all-Rademacher vector belongs to \(\mathcal Y_r\),
\[
\E\mathcal Q_{n, I}^{\mathrm R}(u)^r
\leq\mathfrak P_r.
\]
Taking \(s=r\) proves
\eqref{eq:direct-arbitrary-interval-leakage}.
\end{proof}

Define
\begin{equation}
\begin{aligned}
\operatorname{Err}_{n,I}^{\mathrm{proj}}(u)
:={}
C_qB_{n,I}(u)^{1/2}                        +
C_q
\sum_{\substack{
\mathbf k\in\mathsf C_r\\
L(\mathbf k)\geq1}}
n^{-L(\mathbf k)/q}
\mathfrak J_4(u)^{K(\mathbf k)/r}
P_{n,I}(u)^{(r-K(\mathbf k))/q}.
\end{aligned}
\label{eq:direct-projected-error}
\end{equation}
Then
\begin{equation}
\|
\mathcal Q_{n,I}^{\mathrm R}(u)
\|_{L^r}^{1/2}
\leq
C_q
\left[
\frac hn\mathfrak E_1(u)
\right]^{1/2}
+
\operatorname{Err}_{n,I}^{\mathrm{proj}}(u).
\label{eq:direct-projected-root}
\end{equation}

\boldparagraph{The leading linearized term.}
With \(\mathcal I_{m,n}\) and \(Z_m\) as above, the leading term may be written as
\begin{equation}
J_{m,n}^{\mathrm{lin}}(u)
=
2\sqrt{\frac mn}
\sum_{k=m}^{n-1}
\varepsilon_k
\operatorname{Re}
\widehat{\Phi_u(Z_m)}(k).
\label{eq:direct-leading-Fourier-form}
\end{equation}

\begin{corollary}
\label{cor:direct-leading-leakage}
For \(q=2r\),
\begin{equation}
\begin{aligned}
\|
J_{m,n}^{\mathrm{lin}}(u)
\|_{L^q}
\leq C_{q,T}
[\rho_{m,n}\mathfrak E_1(u)]^{1/2}
+
\operatorname{Err}_{m,n}^{\mathrm{lead,proj}}(u),
\end{aligned}
\label{eq:direct-leading-leakage}
\end{equation}
where
\begin{equation}
\operatorname{Err}_{m,n}^{\mathrm{lead,proj}}(u)
=
C_T
\operatorname{Err}_{m,\mathcal I_{m,n}}^{\mathrm{proj}}(u).
\label{eq:direct-leading-projected-error}
\end{equation}
\end{corollary}

\begin{proof}
Condition on \(\mathcal F_m\) and use Khintchine's inequality:
\[
\begin{aligned}
\|
J_{m,n}^{\mathrm{lin}}(u)
\|_q
&\leq C_{q,T}
\left\|\Pi_{\mathcal I_{m,n}}
\Phi_u(Z_m)
\right\|_{L^q(\Omega;L^2(\Torus))}.
\end{aligned}
\]
Applying \eqref{eq:direct-projected-root} with base length \(m\) and using
\(h/m\asymp_T\rho_{m,n}\) proves
\eqref{eq:direct-leading-leakage}.
\end{proof}

\boldparagraph{Endpoint estimates with homogeneous seminorms.}

The next estimate does not require a fixed normalization of the test
function. 
\begin{definition}
    Define the homogeneous derivative seminorm
\begin{equation}
\begin{aligned}
\dot{\mathfrak J}_K(f)
:={}
\inf_{c\in\R}
\|f-c\|_\infty
+
\max_{1\leq a\leq K+1}
\|f^{(a)}\|_\infty                  +
\max_{|\boldsymbol{\nu}|\leq K}
\|
\partial^{\boldsymbol{\nu}}
\Phi_f
\|_\infty .
\end{aligned}
\label{eq:direct-homogeneous-jet}
\end{equation}
\end{definition}

For fixed \(q\) and \(J\geq3\), take \(K_{q,J}\)
as in Theorem~\ref{thm:compact-endpoint-moment} and choose
\(C_{q,J}\) sufficiently large. 
Fix \(0<c_{q,J}<1/q\). We keep the higher Laguerre energies and the
derivative seminorms separate because later applications estimate them
differently. Define
\begin{equation}
\begin{aligned}
\mathcal R_{N,q,J}^{\mathrm{end}}(f)
:={}&
C_{q,J}
\sum_{j=1}^{J}
N^{-j}
\left[\mathfrak E_{2j}(f)+\mathfrak E_{2j+1}(f)
\right]\\
&+
C_{q,J}
N^{1-J}
\dot{\mathfrak J}_{K_{q,J}}(f)^2\\
&+
N^{-c_{q,J}}
\dot{\mathfrak J}_{K_{q,J}}(f)^2 .
\end{aligned}
\label{eq:direct-endpoint-remainder}
\end{equation}

\begin{lemma}
\label{lem:direct-homogeneous-endpoint}
There exists \(N_0=N_0(q,J,T,c_{q,J},C_{q,J})\) such that,
for all \(N\geq N_0\) and \(n\asymp_TN\), every deterministic
\(f\in\R+C_c^\infty([0,\infty);\R)\) satisfies
\begin{equation}
\|X_n(f)\|_{L^q}^2
\leq
C_{q,J}\mathfrak E_1(f)
+
\mathcal R_{N,q,J}^{\mathrm{end}}(f).
\label{eq:direct-homogeneous-endpoint}
\end{equation}
\end{lemma}

\begin{proof}
Let \(c_0\in\R\) be such that
\(g:=f-c_0\in C_c^\infty([0,\infty))\).
Since \(f=c_0\) outside a compact set, for every \(c\in\R\),
\[
|c-c_0|\leq\|f-c\|_\infty,
\qquad
\|g\|_\infty
\leq\|f-c\|_\infty+|c-c_0|.
\]
Consequently,
\[
\|g\|_\infty
\leq2\inf_{c\in\R}\|f-c\|_\infty.
\]

Put \(d=\dot{\mathfrak J}_{K_{q,J}}(f)\).
If \(d=0\), then \(f\) is constant and \(X_n(f)=0\).
Otherwise, set \(\widetilde f=g/d\).  Then
\(\widetilde f\in C_c^\infty([0,\infty))\), and
\eqref{eq:direct-homogeneous-jet} yields
\[
\mathfrak J_{K_{q,J}}(\widetilde f)\leq C_{q,J}.
\]
Moreover, for every \(s\geq1\),
\[
X_n(f)=dX_n(\widetilde f),
\qquad
\mathfrak E_s(f)=d^2\mathfrak E_s(\widetilde f).
\]
Applying Theorem~\ref{thm:compact-endpoint-moment} to
\(\widetilde f\) and multiplying by \(d^2\), we obtain
\[
\begin{aligned}
\|X_n(f)\|_q^2
\leq
C_{q,J}\mathfrak E_1(f)
+
C_{q,J}\sum_{j=1}^{J}N^{-j}
\bigl[\mathfrak E_{2j}(f)+\mathfrak E_{2j+1}(f)\bigr]+
C_{q,J}N^{-1/q}d^2.
\end{aligned}
\]
By \eqref{eq:direct-endpoint-remainder}, this proves
\eqref{eq:direct-homogeneous-endpoint}.
\end{proof}

\subsection{Conditional mean via the Laguerre semigroup}\label{subsec:block-conditional-mean}

We next express the main part of the conditional mean through the Laguerre semigroup and bound the Gaussian replacement error.
We use the standard circular Gaussian Mehler representation of the
Laguerre semigroup to identify the conditional-mean term.  The error in
Lemma~\ref{lem:direct-conditional-mean-identity} is obtained by coordinate
replacement for the added Rademacher block and by covariance
interpolation for its noncircular Gaussian counterpart; its displayed
block-length and angular factors are retained in the refinement sum.

Put
\begin{equation}
\tau_{m,n}
=
-\log r_{m,n}.
\label{eq:direct-semigroup-time}
\end{equation}
Let \((P_\tau)_{\tau\geq0}\) be the Laguerre semigroup from
Section~\ref{sec:spectral-framework}.

For independent standard circular Gaussian \(\xi\), the complex Mehler
identity implies
\begin{equation}
\E
u\left(
\left|
\sqrt{r_{m,n}}z
+
\sqrt{\rho_{m,n}}\xi
\right|^2
\right)
=
P_{\tau_{m,n}}u(|z|^2).
\label{eq:direct-complex-Mehler}
\end{equation}
The expectation is taken over \(\xi\), with \(z\) fixed.

Define
\begin{equation}
\psi_{m,n}
=
\sqrt{\frac nm}
P_{\tau_{m,n}}u
-
u.
\label{eq:direct-semigroup-drift-label}
\end{equation}
Here \(\psi_{m,n}\) records the change in \(u\) after smoothing and rescaling.

\begin{lemma}
\label{lem:direct-conditional-mean-identity}
One has
\begin{equation}
V_{m,n}(u)
=
X_m(\psi_{m,n})
+
\mathcal E_{m,n}^{\mathrm{mean}}(u),
\label{eq:direct-conditional-mean-identity}
\end{equation}
where
\begin{equation}
\begin{aligned}
\|
\mathcal E_{m,n}^{\mathrm{mean}}(u)
\|_{L^q}
\leq C_{q,T}N^{-1/2}
\left[\rho_{m,n}\mathfrak J_4(u)
+
\log(h+2)\mathfrak J_2(u)
\right].
\end{aligned}
\label{eq:direct-conditional-mean-error}
\end{equation}
\end{lemma}

\begin{proof}
Recall that
\[
\frac{q_n(e^{i\theta})}{\sqrt n}
=
\sqrt{r_{m,n}}\,Z_m(\theta)
+
\frac1{\sqrt n}
\sum_{j=m}^{n-1}\varepsilon_je^{ij\theta},
\]
where
\[
r_{m,n}
=
\frac mn
=
1-\rho_{m,n},
\qquad
Z_m(\theta)
=
\frac{q_m(e^{i\theta})}{\sqrt m}.
\]

For \(z\in\C\), set \(g_u(z)=u(|z|^2)\), and let
\[
\begin{aligned}
\mathcal K_{m,n}^{\mathrm R}(\theta,z)
&= \E
g_u\left(
\sqrt{r_{m,n}}\,z
+
\frac1{\sqrt n}
\sum_{j=m}^{n-1}\varepsilon_je^{ij\theta}
\right),\\
\mathcal K_{m,n}^{\mathrm M}(\theta,z)
&= \E
g_u\left(
\sqrt{r_{m,n}}\,z
+
\frac1{\sqrt n}
\sum_{j=m}^{n-1}G_je^{ij\theta}
\right),\\
\mathcal K_{m,n}^{\circ}(\theta,z)
&= \E
g_u\left(
\sqrt{r_{m,n}}\,z
+
\frac1{\sqrt n}
\sum_{j=m}^{n-1}\xi_je^{ij\theta}
\right),
\end{aligned}
\]
where \(G_j\) are standard real Gaussians and
\(\xi_j\) are standard circular complex Gaussians. In the above kernel definitions, \(\theta\) and \(z\) are fixed,
and the expectation is taken over the added coefficients indexed
by \(m,\ldots,n-1\).

For \(\mathbf y=(y_m,\ldots,y_{n-1})\in\R^h\), define
\[
f_{\theta,z}(\mathbf y)
=
g_u\left(
\sqrt{r_{m,n}}\,z
+
\frac1{\sqrt n}
\sum_{j=m}^{n-1}y_je^{ij\theta}
\right).
\]
For \(m\leq k\leq n\), let
\[
\mathbf Y^{(k)}
=
\bigl(Y_m^{(k)},\ldots,Y_{n-1}^{(k)}\bigr),
\qquad
Y_j^{(k)}
=
        \begin{cases}
G_j, & m\leq j<k,\\
\varepsilon_j, & k\leq j<n.
        \end{cases}
\]
Thus \(\mathbf Y^{(m)}\) is the Rademacher vector and
\(\mathbf Y^{(n)}\) is the matched Gaussian vector.

The definition of \(\mathfrak J_K\) implies
\[
\max_{|\boldsymbol{\nu}|=4}
\|\partial^{\boldsymbol{\nu}} g_u\|_\infty
\leq
C\mathfrak J_4(u).
\]
For fixed \((\theta,z)\) and \(m\leq k<n\), the one-coordinate replacement
\(\mathbf Y^{(k)}\mapsto\mathbf Y^{(k+1)}\) therefore satisfies, by
third-order Taylor expansion with fourth-order remainder and
\[
\E\varepsilon_k^\nu
=
\E G_k^\nu,
\qquad
0\leq\nu\leq3,
\]
the estimate
\[
\left|
\E f_{\theta,z}(\mathbf Y^{(k)})
-
\E f_{\theta,z}(\mathbf Y^{(k+1)})
\right|
\leq
Cn^{-2}\mathfrak J_4(u).
\]
Summation over the \(h=n-m\) new coordinates yields
\begin{equation}
\left|
\mathcal K_{m,n}^{\mathrm R}(\theta,z)
-
\mathcal K_{m,n}^{\mathrm M}(\theta,z)
\right|
\leq
C\frac{h}{n^2}\mathfrak J_4(u).
\label{eq:direct-mean-Rademacher-matched-bound}
\end{equation}

For the Gaussian comparison, put
\[
\eta_\theta^{\mathrm M}
=
\frac1{\sqrt n}
\sum_{j=m}^{n-1}G_je^{ij\theta},
\qquad
\eta_\theta^\circ
=
\frac1{\sqrt n}
\sum_{j=m}^{n-1}\xi_je^{ij\theta}.
\]
Then
\[
\E|\eta_\theta^{\mathrm M}|^2
=
\E|\eta_\theta^\circ|^2
=
\rho_{m,n},
\]
whereas
\[
\E(\eta_\theta^{\mathrm M})^2
=
\frac{D_n(2\theta)-D_m(2\theta)}{n},
\qquad
\E(\eta_\theta^\circ)^2
=
0,
\]
Since \(h=n-m\),
\[
D_n(2\theta)-D_m(2\theta)
=
e^{2im\theta}D_h(2\theta).
\]
Under \(\C\simeq\R^2\),
\[
\Sigma_\theta^{\mathrm M}
-
\Sigma_\theta^\circ
=
\frac1{2n}
        \begin{pmatrix}
\operatorname{Re}\left(e^{2im\theta}D_h(2\theta)\right)
&
\operatorname{Im}\left(e^{2im\theta}D_h(2\theta)\right)\\
\operatorname{Im}\left(e^{2im\theta}D_h(2\theta)\right)
&
-\operatorname{Re}\left(e^{2im\theta}D_h(2\theta)\right)
        \end{pmatrix},
\]
and hence
\[
\sum_{i_1,i_2=1}^{2}
\left|
(\Sigma_\theta^{\mathrm M}
-
\Sigma_\theta^\circ)_{i_1i_2}
\right|
\leq
\frac{C}{n}|D_h(2\theta)|.
\]
Lemma~\ref{lem:Gaussian-covariance-differentiation}, applied to
\[
(x_1,x_2)
\longmapsto
g_u\left(
\sqrt{r_{m,n}}\,z+x_1+ix_2
\right),
\]
together with
\[
\max_{|\boldsymbol{\nu}|=2}
\|\partial^{\boldsymbol{\nu}} g_u\|_\infty
\leq
C\mathfrak J_2(u),
\]
yields
\begin{equation}
\left|
\mathcal K_{m,n}^{\mathrm M}(\theta,z)
-
\mathcal K_{m,n}^{\circ}(\theta,z)
\right|
\leq
\frac{C}{n}|D_h(2\theta)|\mathfrak J_2(u).
\label{eq:direct-mean-matched-circular-pointwise}
\end{equation}
By Haar invariance and the standard bound,
\begin{equation}
\int_{\Torus}
\frac{|D_h(2\theta)|}{n}
\,\dd\mu(\theta)
\leq
C\frac{\log(h+2)}n.
\label{eq:direct-mean-anomalous-L1-bound}
\end{equation}
Consequently, uniformly over measurable \(z:\Torus\to\C\),
\begin{equation}
\begin{aligned}
\int_{\Torus}
\left|\mathcal K_{m,n}^{\mathrm R}(\theta,z(\theta))
-
\mathcal K_{m,n}^{\circ}(\theta,z(\theta))
\right|
\,\dd\mu(\theta)
\leq C\frac{h}{n^2}\mathfrak J_4(u)
+
C\frac{\log(h+2)}n\mathfrak J_2(u).
\end{aligned}
\label{eq:direct-mean-integrated-kernel-comparison}
\end{equation}

For the circular block,
\[
\E
\left|
\frac1{\sqrt n}
\sum_{j=m}^{n-1}\xi_je^{ij\theta}
\right|^2
=
\rho_{m,n},
\qquad
\E
\left(
\frac1{\sqrt n}
\sum_{j=m}^{n-1}\xi_je^{ij\theta}
\right)^2
=
0.
\]
Since
\[
e^{-\tau_{m,n}}
=
r_{m,n},
\]
the complex Mehler identity
\eqref{eq:direct-complex-Mehler} yields
\[
\mathcal K_{m,n}^{\circ}(\theta,z)
=
P_{\tau_{m,n}}u(|z|^2).
\]

Define
\begin{equation}
\begin{aligned}
\mathcal A_{m,n}^{\mathrm{mean}}(u)
&= \sqrt n
\int_{\Torus}
\left[\mathcal K_{m,n}^{\mathrm R}
(\theta,Z_m(\theta))
-
\mathcal K_{m,n}^{\circ}
(\theta,Z_m(\theta))
\right]
\,\dd\mu(\theta),
\end{aligned}
\label{eq:direct-mean-raw-error-definition}
\end{equation}
and set
\[
\mathcal E_{m,n}^{\mathrm{mean}}(u)
=
\left[
\mathcal A_{m,n}^{\mathrm{mean}}(u)
\right]_0.
\]
Conditional Fubini and the independence of
\(\varepsilon_m,\ldots,\varepsilon_{n-1}\) from \(\mathcal F_m\) imply
\[
\begin{aligned}
\E(F_n(u)\mid\mathcal F_m)
&= \sqrt n
\int_{\Torus}
P_{\tau_{m,n}}u(S_m(\theta))
\,\dd\mu(\theta)
+
\mathcal A_{m,n}^{\mathrm{mean}}(u).
\end{aligned}
\]
Using
\[
\psi_{m,n}
=
\sqrt{\frac nm}\,
P_{\tau_{m,n}}u-u,
\]
one has
\[
\begin{aligned}
 \sqrt n
 \int_{\Torus}
 P_{\tau_{m,n}}u(S_m)
 \,\dd\mu
 -
 F_m(u)
= \sqrt m
 \int_{\Torus}
 \psi_{m,n}(S_m)
 \,\dd\mu.
\end{aligned}
\]
Hence, by
\eqref{eq:direct-U-V} and \eqref{eq:direct-F-X},
\[
\begin{aligned}
V_{m,n}(u)
= \left[\E(F_n(u)\mid\mathcal F_m)-F_m(u)
\right]_0
= X_m(\psi_{m,n})
+
\mathcal E_{m,n}^{\mathrm{mean}}(u).
\end{aligned}
\]

Finally,
\eqref{eq:direct-mean-integrated-kernel-comparison} with
\(z(\theta)=Z_m(\theta)\) yields
\[
\begin{aligned}
|\mathcal A_{m,n}^{\mathrm{mean}}(u)|
&\leq C\sqrt n
\left[\frac{h}{n^2}\mathfrak J_4(u)
+
\frac{\log(h+2)}n\mathfrak J_2(u)
\right]\\
&= \frac C{\sqrt n}
\left[\rho_{m,n}\mathfrak J_4(u)
+
\log(h+2)\mathfrak J_2(u)
\right].
\end{aligned}
\]
Therefore
\[
\|
\mathcal E_{m,n}^{\mathrm{mean}}(u)
\|_{L^q}
\leq
2
\|
\mathcal A_{m,n}^{\mathrm{mean}}(u)
\|_{L^\infty}.
\]
Since \(m,n\asymp_TN\), the last two estimates imply
\eqref{eq:direct-conditional-mean-error}.
\end{proof}

If
\[
u=\sum_{a\geq0}\widehat u_aL_a,
\]
then
\begin{equation}
\psi_{m,n}
=
\sum_{a\geq0}
\left(
r_{m,n}^{\,a-\frac12}-1
\right)
\widehat u_aL_a.
\label{eq:direct-semigroup-Laguerre-expansion}
\end{equation}

Let
\[
\mathscr A_{\mathrm{OU}}
=\frac14\Delta_z-\frac12z\cdot\nabla_z
\]
be the two-dimensional Ornstein--Uhlenbeck generator.
For \(g_u(z)=u(|z|^2)\),
\[
\mathscr A_{\mathrm{OU}}g_u(z)
=(\LagGen u)(|z|^2).
\]

For \(M\geq0\), define
\begin{equation}
\begin{aligned}
\mathfrak H_M(u)
:={}&
\max_{0\leq a\leq M}
\left[\|u^{(a)}\|_\infty
+
\|(\LagGen u)^{(a)}\|_\infty
\right]\\
&+
\max_{|\boldsymbol{\nu}|\leq2M+4}
\left[\|\partial^{\boldsymbol{\nu}} g_u\|_\infty
+
\|
\mathscr A_{\mathrm{OU}}
\partial^{\boldsymbol{\nu}} g_u
\|_\infty
\right].
\end{aligned}
\label{eq:direct-generated-OU-seminorm}
\end{equation}
Here \(\mathfrak H_M(u)\) collects uniform bounds on the derivatives of \(u\) and \(g_u\), including the terms produced by the two generators.

\begin{lemma}
\label{lem:direct-generated-OU-jets}
If
\(\supp u\subset[0,B]\),
then
\begin{equation}
\mathfrak H_M(u)
\leq
C_M
(1+B)^{M+3}
\max_{0\leq a\leq2M+7}
\|u^{(a)}\|_\infty .
\label{eq:direct-generated-OU-ordinary-jet}
\end{equation}
Moreover,
\begin{equation}
\dot{\mathfrak J}_K(\psi_{m,n})
\leq
C_{K,T}
\rho_{m,n}
\mathfrak H_{K+2}(u).
\label{eq:direct-semigroup-jet-gain}
\end{equation}
\end{lemma}

\begin{proof}
Let
\[
x=|z|^2,
\qquad
g_u(z)=u(x).
\]
For every multi-index $\boldsymbol{\nu}$, the chain rule yields
\[
\partial^{\boldsymbol{\nu}} g_u(z)
=
\sum_{\ell=0}^{|\boldsymbol{\nu}|}
P_{\boldsymbol{\nu},\ell}(z)u^{(\ell)}(x),
\]
where
\[
\deg P_{\boldsymbol{\nu},\ell}\le |\boldsymbol{\nu}|.
\]
Since $\supp u\subset[0,B]$, one has
\(|z|^2\le B\)
on $\supp g_u$, and hence
\[
\max_{|\boldsymbol{\nu}|\le2M+4}
|\partial^{\boldsymbol{\nu}} g_u|_\infty
\le
C_M(1+B)^{M+2}
\max_{0\le a\le2M+4}|u^{(a)}|_\infty .
\]

Moreover,
\[
\mathscr A_{\mathrm{OU}}g_u
=
\LagGen u(|z|^2),
\]
and applying $\mathscr A_{\mathrm{OU}}$ to the above radial expansion yields
\[
\mathscr A_{\mathrm{OU}}
\partial^{\boldsymbol{\nu}} g_u(z)
=
\sum_{\ell=0}^{|\boldsymbol{\nu}|+2}
\widetilde P_{\boldsymbol{\nu},\ell}(z)
u^{(\ell)}(x),
\]
with
\[
\deg\widetilde P_{\boldsymbol{\nu},\ell}\le|\boldsymbol{\nu}|+2 .
\]
Consequently,
\[
\max_{|\boldsymbol{\nu}|\le2M+4}
|\mathscr A_{\mathrm{OU}}\partial^{\boldsymbol{\nu}} g_u|_\infty
\le
C_M(1+B)^{M+3}
\max_{0\le a\le2M+7}|u^{(a)}|_\infty .
\]
Together with the scalar part in the definition of $\mathfrak H_M$,
this proves \eqref{eq:direct-generated-OU-ordinary-jet}.

For the second assertion, write
\[
\psi_{m,n}
=
\left(\sqrt{\frac nm}-1\right)u
+
\sqrt{\frac nm}
\int_0^{\tau_{m,n}}
\LagGen P_su\,\dd s .
\]
Using
\[
\tau_{m,n}\le C_T\rho_{m,n},
\qquad
\left|\sqrt{\frac nm}-1\right|
\le C_T\rho_{m,n},
\]
it remains to control the two terms.

The Laguerre semigroup satisfies
\[
\LagGen P_su=P_s\LagGen u,
\qquad
\|P_sf\|_\infty\le \|f\|_\infty .
\]
For the scalar derivatives, let \(P_s^{[a]}\), \(a\geq0\),
denote the Markov Laguerre semigroup on \([0,\infty)\)
with generator
\[
\mathscr L^{[a]}f
=
xf''+(a+1-x)f',
\qquad
P_s^{[0]}=P_s.
\]
For every nonnegative integer \(a\),
\[
\frac{\dd^a}{\dd x^a}\LagGen f
=
\bigl(\mathscr L^{[a]}-a\bigr)f^{(a)}.
\]
Differentiating the semigroup equation therefore yields, for
\(f\in C_c^\infty([0,\infty))\),
\[
(P_sf)^{(a)}
=
e^{-as}P_s^{[a]}f^{(a)}.
\]
Since \(P_s^{[a]}\) is a Markov contraction,
\begin{equation}
\|(P_sf)^{(a)}\|_\infty
\leq
e^{-as}\|f^{(a)}\|_\infty,
\qquad s\geq0.
\label{eq:direct-scalar-derivative-contraction}
\end{equation}
In particular,
\[
\begin{aligned}
\max_{0\leq a\leq K+1}
\left[
\|(P_s\LagGen u)^{(a)}\|_\infty
+
\|u^{(a)}\|_\infty
\right]\leq
\max_{0\leq a\leq K+1}
\left[
\|(\LagGen u)^{(a)}\|_\infty
+
\|u^{(a)}\|_\infty
\right]
\leq
\mathfrak H_{K+2}(u).
\end{aligned}
\]
For \(F\in C_c^\infty(\R^2;\C)\), define
\[
(Q_sF)(z)
=
\E F\bigl(e^{-s/2}z+\sqrt{1-e^{-s}}\,\xi\bigr),
\qquad s\geq0,
\]
where \(\R^2\) is identified with \(\C\) and \(\xi\) is
standard circular complex Gaussian. The expectation is taken over \(\xi\), with \(z\) fixed.
For vector-valued functions, \(Q_s\) acts componentwise.  By the complex Mehler identity,
\(g_{P_sf}=Q_sg_f\).  Since \(\Phi_f=\frac12\nabla g_f\), differentiation
under the Gaussian integral yields, for every real multi-index \(\boldsymbol{\nu}\),
\[
\partial^{\boldsymbol{\nu}}\Phi_{P_sf}
=
e^{-(|\boldsymbol{\nu}|+1)s/2}Q_s(\partial^{\boldsymbol{\nu}}\Phi_f),
\qquad
\|\partial^{\boldsymbol{\nu}}\Phi_{P_sf}\|_\infty
\leq
\|\partial^{\boldsymbol{\nu}}\Phi_f\|_\infty.
\]
Together with \eqref{eq:direct-scalar-derivative-contraction} and
\(P_sc=c\) for constants \(c\), this yields
\[
\dot{\mathfrak J}_K(P_s\LagGen u)
+
\dot{\mathfrak J}_K(u)
\le
C_K\mathfrak H_{K+2}(u).
\]
Hence,
\[
\begin{aligned}
\dot{\mathfrak J}_K(\psi_{m,n})
\le C_T\rho_{m,n}\dot{\mathfrak J}_K(u)
+
C_T
\int_0^{\tau_{m,n}}
\dot{\mathfrak J}_K(\LagGen P_su)\,\dd s\le C_{K,T}
\rho_{m,n}
\mathfrak H_{K+2}(u).
\end{aligned}
\]
\end{proof}

\begin{lemma}
\label{lem:direct-semigroup-endpoint}
For every fixed even \(q\) and fixed \(J\geq3\),
\begin{equation}
\begin{aligned}
\|X_m(\psi_{m,n})\|_q^2
\leq{}&
C_{q,T}\rho_{m,n}^2\mathfrak E_3(u)+
C_{q,T}\rho_{m,n}^2
\sum_{j=1}^{J}
N^{-j}
\left[\mathfrak E_{2j+2}(u)
+
\mathfrak E_{2j+3}(u)
\right]\\
&+
C_{q,T}\rho_{m,n}^2
\left[N^{1-J}
+
N^{-c_{q,J}}
\right]
\mathfrak H_{K_{q,J}+2}(u)^2 .
\end{aligned}
\label{eq:direct-semigroup-endpoint}
\end{equation}
\end{lemma}

\begin{proof}
Fix \(N,m,n,u\), and set
\[
\tau=\tau_{m,n}>0,
\qquad
\psi=\psi_{m,n}
=\sqrt{\frac nm}\,P_\tau u-u.
\]
Choose \(B\geq1\) such that \(\supp u\subset[0,B]\).

We first justify approximation by compactly supported functions.
By the Hille--Hardy formula
\cite[Eq.~(18.18.27)]{NISTDLMF},
\[
P_\tau u(x)
=
\int_0^B K_\tau(x,y)u(y)e^{-y}\,\dd y,
\]
where, with \(d_\tau=1-e^{-\tau}\),
\[
K_\tau(x,y)
=
\frac1{d_\tau}
\exp\left(-\frac{e^{-\tau}(x+y)}{d_\tau}\right)
I_0\left(\frac{2e^{-\tau/2}\sqrt{xy}}{d_\tau}\right).
\]
Here \(I_j\) denotes the modified Bessel function of the first kind of
order \(j\), normalized as in \cite[Eq.~(10.25.2)]{NISTDLMF}.
For integers \(j\geq0\), the identities and bounds
\[
\frac{\dd^j}{\dd x^j}I_0(2\sqrt{cx})
=
\left(\frac cx\right)^{j/2}I_j(2\sqrt{cx}),
\qquad
|I_j(t)|\leq e^t
\quad(t\geq0),
\]
with the first identity extended to \(c=0\) by continuity,
yield
\[
\begin{aligned}
\sup_{0\leq y\leq B}
|\partial_x^jK_\tau(x,y)|
\leq
C_{\tau,B,j}
\exp\left(
-\frac{e^{-\tau}}{d_\tau}x
+\frac{2e^{-\tau/2}\sqrt B}{d_\tau}\sqrt x
\right)\leq
C_{\tau,B,j}e^{-c_\tau x},
\qquad x\geq1,
\end{aligned}
\]
where \(c_\tau=e^{-\tau}/(2d_\tau)>0\), after enlarging the
constant.  Differentiation under the integral is therefore
justified, and
\[
|\partial_x^jP_\tau u(x)|
\leq C_{\tau,u,j}e^{-c_\tau x},
\qquad x\geq1.
\]
Since \(u\) is compactly supported, the same exponential tail
bound holds for every derivative of \(\psi\).
These bounds are used only with \(\tau\) and \(u\) fixed.

Choose a fixed function
\[
\chi\in C_c^\infty([0,\infty)),
\qquad
0\leq\chi\leq1,
\qquad
\chi=1\ \text{on }[0,1],
\qquad
\supp\chi\subset[0,2],
\]
and define
\[
\chi_R(x)=\chi(x/R),
\qquad
\psi_R=\chi_R\psi,
\qquad
f_R=\psi_R-\psi,
\qquad R\geq1.
\]
For each integer \(j\geq0\),
\[
\|\chi_R^{(j)}\|_\infty\leq C_jR^{-j}.
\]
Leibniz's rule and the exponential tail bounds imply that,
for every pair of nonnegative integers \(M,j\),
\begin{equation}
\sup_{x\geq0}
(1+x)^M|f_R^{(j)}(x)|
\longrightarrow0
\qquad(R\to\infty).
\label{eq:direct-semigroup-cutoff-derivatives}
\end{equation}

Each real derivative of
\(\Phi_{f_R}(z)=zf_R'(|z|^2)\)
is a finite sum of polynomial factors in \(z\) multiplied by
scalar derivatives of \(f_R\).
Thus \eqref{eq:direct-semigroup-cutoff-derivatives} yields
\[
\dot{\mathfrak J}_K(f_R)\longrightarrow0
\qquad(R\to\infty)
\]
for every fixed \(K\).
In particular,
\[
\dot{\mathfrak J}_K(\psi_R)
\longrightarrow
\dot{\mathfrak J}_K(\psi),
\]
by the triangle inequality for this seminorm.

For completeness, fix \(s\geq1\) and choose an integer \(b\geq1\)
such that \(2b\geq s\).
The differential operator \((-\LagGen)^b\) has polynomial
coefficients.  Hence
\eqref{eq:direct-semigroup-cutoff-derivatives} implies
\[
\|(-\LagGen)^b f_R\|_{L^2(\nu)}
\longrightarrow0.
\]
Moreover, \(f_R\) vanishes near \(0\), and all its derivatives
decay exponentially at infinity.
Repeated integration by parts against \(L_a\) therefore yields
\[
\widehat{(-\LagGen)^bf_R}_a
=
a^b\widehat{f_R}_a.
\]
Bessel's inequality now yields
\[
\mathfrak E_s(f_R)
\leq
\sum_{a\geq1}a^{2b}|\widehat{f_R}_a|^2
\leq
\|(-\LagGen)^bf_R\|_{L^2(\nu)}^2
\longrightarrow0.
\]
Consequently,
\[
\mathfrak E_s(\psi_R)\longrightarrow\mathfrak E_s(\psi)
\]
for every fixed \(s\geq1\).

Since \(0\leq S_m(\theta)\leq m\), one has
\[
X_m(\psi_R)=X_m(\psi)
\qquad\text{whenever }R\geq m.
\]
Apply Lemma~\ref{lem:direct-homogeneous-endpoint} to
\(\psi_R\in C_c^\infty([0,\infty))\).
Its constants are independent of \(R\).
Letting \(R\to\infty\) and using the preceding convergences yields
\[
\begin{aligned}
\|X_m(\psi)\|_{L^q}^2
\leq{}&
C_{q,J,T}\mathfrak E_1(\psi)+
C_{q,J,T}
\sum_{j=1}^{J}N^{-j}
\left[
\mathfrak E_{2j}(\psi)+\mathfrak E_{2j+1}(\psi)
\right]\\
&+
C_{q,J,T}
\left[N^{1-J}+N^{-c_{q,J}}\right]
\dot{\mathfrak J}_{K_{q,J}}(\psi)^2.
\end{aligned}
\]

Since
\(0\leq1-r_{m,n}^{\,a-\frac12}
\leq1-r_{m,n}^{\,a}
\leq a\rho_{m,n}\) for \(a\geq1\),
\eqref{eq:direct-semigroup-Laguerre-expansion} gives
\[
\begin{aligned}
\mathfrak E_s(\psi)
=
\sum_{a\geq1}
a^s
\left(1-r_{m,n}^{\,a-\frac12}\right)^2
|\widehat u_a|^2\leq
\rho_{m,n}^2\mathfrak E_{s+2}(u),
\qquad s\geq1,
\end{aligned}
\]
while \eqref{eq:direct-semigroup-jet-gain} yields
\[
\dot{\mathfrak J}_{K_{q,J}}(\psi)
\leq
C_{q,J,T}\rho_{m,n}
\mathfrak H_{K_{q,J}+2}(u).
\]
Substituting these bounds proves
\eqref{eq:direct-semigroup-endpoint}, with dependence on the
fixed \(J\) suppressed in the displayed constants.
\end{proof}

\begin{proposition}
\label{prop:direct-conditional-mean-drift}
For every fixed even \(q\) and every fixed integer \(J\geq3\),
\begin{equation}
\begin{aligned}
\|V_{m,n}(u)\|_q
\leq{}&
C_{q,T}
\rho_{m,n}\mathfrak E_3(u)^{1/2}+
C_{q,T}
\rho_{m,n}
\left\{\sum_{j=1}^{J}
N^{-j}
\left[\mathfrak E_{2j+2}(u)
+
\mathfrak E_{2j+3}(u)
\right]
\right\}^{1/2}\\
&+
C_{q,T}
\rho_{m,n}
\left[N^{(1-J)/2}
+
N^{-c_{q,J}/2}
\right]
\mathfrak H_{K_{q,J}+2}(u)\\
&+
C_{q,T}N^{-1/2}
\left[\rho_{m,n}\mathfrak J_4(u)
+
\log(h+2)\mathfrak J_2(u)
\right].
\end{aligned}
\label{eq:direct-conditional-mean-drift}
\end{equation}
\end{proposition}

\begin{proof}
Combine
\eqref{eq:direct-conditional-mean-identity},
\eqref{eq:direct-conditional-mean-error},
and
\eqref{eq:direct-semigroup-endpoint}.
\end{proof}

\subsection{The conditionally centered term and the local increment estimate}

Next we bound the conditionally centered term and combine it with the conditional-mean estimate to control increments between nearby polynomial lengths.
We estimate the nonlinear part of the conditionally centered term
\(U_{m,n}\) by comparing Fourier projections at the two lengths.

\boldparagraph{Two-time projected difference.}
Construct the nested fields
\[
Z_m^\diamond,
\qquad
Z_n^\diamond,
\qquad
\diamond\in\{\mathrm R,\mathrm M,\circ\},
\]
from the same coefficient sequence. We first compare the values at \(m\) and \(n\) after rescaling, then measure the energy of their difference in the chosen frequency range. Put
\begin{equation}
\Theta_{m,n}^\diamond[u]
=
\Phi_u(Z_n^\diamond)
-
\sqrt{r_{m,n}}
\Phi_u(Z_m^\diamond),
\label{eq:direct-two-time-Theta}
\end{equation}
and
\begin{equation}
\mathcal Q_{m,n}^{\mathrm{diff},\diamond}(u)
=
\left\|
\Pi_{\mathcal I_{m,n}}
\Theta_{m,n}^\diamond[u]
\right\|_2^2.
\label{eq:direct-difference-leakage}
\end{equation}

Define
\begin{equation}
\begin{aligned}
B_{m,n}^{\mathrm{diff}}(u)
&= \mathfrak J_4(u)^2
\left[\frac{\log(N+2)\log(h+2)}N
+
\frac1N
\right],\\
P_{m,n}^{\mathrm{diff}}(u)
&= C_{q,T}
\left[\rho_{m,n}\mathfrak E_2(u)
+
B_{m,n}^{\mathrm{diff}}(u)
\right].
\end{aligned}
\label{eq:direct-difference-B-P}
\end{equation}

\begin{lemma}
\label{lem:direct-two-time-projected-leakage}
Let \(q=2r\).  Then
\begin{equation}
\begin{aligned}
\E
\left[\mathcal Q_{m,n}^{\mathrm{diff},\mathrm R}(u)^r
\right]
\leq C_{q,T}
\sum_{\mathbf k\in\mathsf C_r}
N^{-L(\mathbf k)}
\mathfrak J_4(u)^{2K(\mathbf k)}
P_{m,n}^{\mathrm{diff}}(u)^{r-K(\mathbf k)} .
\end{aligned}
\label{eq:direct-two-time-projected-leakage}
\end{equation}
Consequently,
\begin{equation}
\begin{aligned}
\|
\mathcal Q_{m,n}^{\mathrm{diff},\mathrm R}(u)
\|_{L^r}^{1/2}
\leq C_{q,T}
[\rho_{m,n}\mathfrak E_2(u)]^{1/2}
+
\operatorname{Err}_{m,n}^{\mathrm{diff,proj}}(u),
\end{aligned}
\label{eq:direct-two-time-projected-root}
\end{equation}
where
\begin{equation}
\begin{aligned}
\operatorname{Err}_{m,n}^{\mathrm{diff,proj}}(u)
:={}&
C_{q,T}
B_{m,n}^{\mathrm{diff}}(u)^{1/2}\\
&+
C_{q,T}
\sum_{\substack{
\mathbf k\in\mathsf C_r\\
L(\mathbf k)\geq1}}
N^{-L(\mathbf k)/q}
\mathfrak J_4(u)^{K(\mathbf k)/r}
P_{m,n}^{\mathrm{diff}}(u)^{(r-K(\mathbf k))/q}.
\end{aligned}
\label{eq:direct-two-time-projected-error}
\end{equation}
\end{lemma}

\begin{proof}
For the circular pair,
\[
\E
\left[
Z_n^\circ\overline{Z_m^\circ}
\right]
=
\sqrt{r_{m,n}},
\qquad
\E[Z_n^\circ Z_m^\circ]=0.
\]
\[
\Phi_u(z)
=
u'(|z|^2)z
=
-
\sum_{a\geq0}
\widehat u_{a+1}
zL_a^{(1)}(|z|^2),
\]
and \eqref{eq:circular-charge-one-Wick}, applied with
\(r=\sqrt{r_{m,n}}\), yield
\[
\begin{aligned}
\E
\left|\Theta_{m,n}^\circ[u]
\right|^2
&= \sum_{a\geq0}
(a+1)
|\widehat u_{a+1}|^2
\left(1+r_{m,n}-2r_{m,n}^{a+1}
\right).
\end{aligned}
\]
Since
\[
1+r_{m,n}-2r_{m,n}^{a+1}
\leq
C\rho_{m,n}(a+1),
\]
one has
\[
\begin{aligned}
\E
\mathcal Q_{m,n}^{\mathrm{diff},\circ}(u)
\leq C\rho_{m,n}
\sum_{a\geq0}
(a+1)^2|\widehat u_{a+1}|^2
\leq C\rho_{m,n}\mathfrak E_2(u).
\end{aligned}
\]

The covariance comparison between the matched Gaussian model \(\mathrm M\)
and the circular Gaussian model \(\circ\) in
Lemma~\ref{lem:direct-arbitrary-interval-leakage}, applied to
\(\Theta_{m,n}[u]\), yields
\[
\begin{aligned}
\left|\E
\mathcal Q_{m,n}^{\mathrm{diff},\mathrm M}(u)
-
\E
\mathcal Q_{m,n}^{\mathrm{diff},\circ}(u)
\right|
&\leq C_{q,T}
B_{m,n}^{\mathrm{diff}}(u).
\end{aligned}
\]
Moreover,
\[
\mathbf x
\longmapsto
\left\|
\Pi_{\mathcal I_{m,n}}
\Theta_{m,n}^{\mathrm M}[u](\mathbf x)
\right\|_2
\]
has Euclidean Lipschitz constant at most
\[
C_TN^{-1/2}\mathfrak J_4(u).
\]
Hence Gaussian concentration and
\eqref{eq:direct-difference-B-P} imply
\[
\E
\left[
\mathcal Q_{m,n}^{\mathrm{diff},\mathrm M}(u)^s
\right]
\leq
C_{s,T}
P_{m,n}^{\mathrm{diff}}(u)^s,
\qquad
1\leq s\leq r.
\]

Set
\[
\mathcal S_{m,n}^{\mathrm{diff}}(\mathbf x)
=
\left\|
\Pi_{\mathcal I_{m,n}}
\Theta_{m,n}[u](\mathbf x)
\right\|_2^2.
\]
Since
\[
m\asymp_Tn\asymp_TN,
\]
the coordinate derivative bounds used in
Lemma~\ref{lem:direct-arbitrary-interval-leakage} remain valid:
\[
\left\|
\partial_{x_j}^{\nu}
\Pi_{\mathcal I_{m,n}}
\Theta_{m,n}[u](\mathbf x)
\right\|_2
\leq
C_TN^{-\nu/2}\mathfrak J_4(u),
\qquad
1\leq\nu\leq4.
\]
Thus
\[
\begin{aligned}
\left|\partial_{x_j}^4
\mathcal S_{m,n}^{\mathrm{diff}}(\mathbf x)^s
\right|
&\leq C_{s,T}N^{-2}
\sum_{a=1}^{\min\{4,s\}}
\mathfrak J_4(u)^{2a}
\mathcal S_{m,n}^{\mathrm{diff}}(\mathbf x)^{s-a}.
\end{aligned}
\]
Use the families \(\mathcal Y_s\) defined in
\eqref{eq:direct-enlarged-hybrid-family}, now with \(n\) coordinates,
and set
\[
\mathfrak P_s^{\mathrm{diff}}
=
\sup_{\mathbf Y\in\mathcal Y_s}
\E\mathcal S_{m,n}^{\mathrm{diff}}(\mathbf Y)^s,
\qquad 1\leq s\leq r,
\qquad
\mathfrak P_0^{\mathrm{diff}}=1.
\]

We first estimate the Gaussian vectors in these enlarged families.
For
\(\mathbf Y_{\mathrm G}=(c_0g_0,\ldots,c_{n-1}g_{n-1})\), define
\[
Z_{\ell,\mathbf c}(\theta)
=
\frac1{\sqrt\ell}
\sum_{j=0}^{\ell-1}c_jg_je^{ij\theta},
\qquad \ell\in\{m,n\}.
\]
For \(a,b\in\{m,n\}\), the changes in the covariance and
pseudocovariance kernels relative to \(c_j\equiv1\) are
\[
\begin{aligned}
&
\E\!\left[
Z_{a,\mathbf c}(\theta)
\overline{Z_{b,\mathbf c}(\varphi)}
\right]
-
\E\!\left[
Z_a^{\mathrm M}(\theta)
\overline{Z_b^{\mathrm M}(\varphi)}
\right]
=
\frac1{\sqrt{ab}}
\sum_{j<\min\{a,b\}}
(c_j^2-1)e^{ij(\theta-\varphi)},\\
&
\E\!\left[
Z_{a,\mathbf c}(\theta)Z_{b,\mathbf c}(\varphi)
\right]
-
\E\!\left[
Z_a^{\mathrm M}(\theta)Z_b^{\mathrm M}(\varphi)
\right]
=
\frac1{\sqrt{ab}}
\sum_{j<\min\{a,b\}}
(c_j^2-1)e^{ij(\theta+\varphi)}.
\end{aligned}
\]
At most \(r\) summands are nonzero, their coefficients are bounded
in terms of \(r\), and \(a,b\asymp_TN\).
Both kernel differences are therefore bounded by \(C_{r,T}/N\),
uniformly in the angles.

Apply Lemma~\ref{lem:Gaussian-covariance-differentiation} to the
real Gaussian vector consisting of the two nested fields at
\(\theta\) and \(\varphi\).
The projection-kernel representation and the derivative bounds
for \(\Phi_u\) yield
\[
\begin{aligned}
\left|
\E\mathcal S_{m,n}^{\mathrm{diff}}(\mathbf Y_{\mathrm G})
-
\E
\mathcal Q_{m,n}^{\mathrm{diff},\mathrm M}(u)
\right|
\leq
C_{r,T}\frac{\log(h+2)}N\mathfrak J_4(u)^2.
\end{aligned}
\]
The map
\[
\mathbf g\longmapsto
\mathcal S_{m,n}^{\mathrm{diff}}
(c_0g_0,\ldots,c_{n-1}g_{n-1})^{1/2}
\]
has Euclidean Lipschitz constant at most
\(C_{r,T}N^{-1/2}\mathfrak J_4(u)\).
Consequently, the matched Gaussian estimate above and Gaussian
concentration yield
\[
\E\mathcal S_{m,n}^{\mathrm{diff}}(\mathbf Y_{\mathrm G})^s
\leq
C_{r,T}P_{m,n}^{\mathrm{diff}}(u)^s,
\qquad 1\leq s\leq r.
\]

The enlarged families are preserved by the coordinate replacements
and coefficient rescaling used above.  The Gaussian weight bound
\eqref{eq:direct-Gaussian-weight-absorption} also applies.
Using the fourth-coordinate derivative bound above, each
replacement of \(c_j\varepsilon_j\) by \(c_jg_j\) contributes at most
\[
C_{r,T}N^{-2}
\sum_{a=1}^{\min\{4,s\}}
\mathfrak J_4(u)^{2a}\mathfrak P_{s-a}^{\mathrm{diff}}.
\]
There are at most \(n\leq C_TN\) replacements.
Taking the supremum over \(\mathcal Y_s\) therefore yields
\[
\begin{aligned}
\mathfrak P_s^{\mathrm{diff}}
&\leq
C_{r,T}P_{m,n}^{\mathrm{diff}}(u)^s
+
\frac{C_{r,T}}N
\sum_{a=1}^{\min\{4,s\}}
\mathfrak J_4(u)^{2a}
\mathfrak P_{s-a}^{\mathrm{diff}},\\
\mathfrak P_0^{\mathrm{diff}}&=1.
\end{aligned}
\]
Iterating this recursion as in
Lemma~\ref{lem:direct-arbitrary-interval-leakage}, and using the
all-Rademacher vector in \(\mathcal Y_r\), yields
\[
\begin{aligned}
\E\!\left[
\mathcal Q_{m,n}^{\mathrm{diff},\mathrm R}(u)^r
\right]
\leq \mathfrak P_r^{\mathrm{diff}}\leq
C_{r,T}
\sum_{\mathbf k\in\mathsf C_r}
N^{-L(\mathbf k)}
\mathfrak J_4(u)^{2K(\mathbf k)}
P_{m,n}^{\mathrm{diff}}(u)^{r-K(\mathbf k)}.
\end{aligned}
\]
This proves \eqref{eq:direct-two-time-projected-leakage}.

Since \(q=2r\), taking the \(1/q\)-power and separating the empty
composition yields
\[
\begin{aligned}
\|
\mathcal Q_{m,n}^{\mathrm{diff},\mathrm R}(u)
\|_{L^r}^{1/2}
&\leq C_{q,T}
P_{m,n}^{\mathrm{diff}}(u)^{1/2}\\
&\quad+
C_{q,T}
\sum_{\substack{
\mathbf k\in\mathsf C_r\\
L(\mathbf k)\geq1}}
N^{-L(\mathbf k)/q}
\mathfrak J_4(u)^{K(\mathbf k)/r}
P_{m,n}^{\mathrm{diff}}(u)^{(r-K(\mathbf k))/q}.
\end{aligned}
\]
By \eqref{eq:direct-difference-B-P},
\[
P_{m,n}^{\mathrm{diff}}(u)^{1/2}
\leq
C_{q,T}
[\rho_{m,n}\mathfrak E_2(u)]^{1/2}
+
C_{q,T}
B_{m,n}^{\mathrm{diff}}(u)^{1/2}.
\]
Together with
\eqref{eq:direct-two-time-projected-error},
this proves
\eqref{eq:direct-two-time-projected-root}.
\end{proof}

\boldparagraph{The conditionally centered nonlinear part.}

Define
\begin{equation}
W_{m,n}(u)
=
U_{m,n}(u)
-
J_{m,n}^{\mathrm{lin}}(u).
\label{eq:direct-W-definition}
\end{equation}
Since
\[
\E(
J_{m,n}^{\mathrm{lin}}(u)
\mid
\mathcal F_m
)
=
0,
\]
one has
\begin{equation}
\begin{aligned}
W_{m,n}(u)
= [F_n(u)-J_{m,n}^{\mathrm{lin}}(u)]
-
\E[
F_n(u)-J_{m,n}^{\mathrm{lin}}(u)
\mid\mathcal F_m].
\end{aligned}
\label{eq:direct-W-conditional-centering}
\end{equation}

With the lowering operator \(\mathscr R\) recalled above, the exact coordinate
half-difference formula yields, for \(m\leq k<n\),
\begin{equation}
\begin{aligned}
\Delta_k
(F_n-J_{m,n}^{\mathrm{lin}})
={}&
2\varepsilon_k
\operatorname{Re}
\widehat{
\Theta_{m,n}^{\mathrm R}[u]
}(k)\\
&-
\frac2{\sqrt n}
\int_{\Torus}
\mathscr Ru(S_n)
\,\dd\mu\\
&-
\frac2{\sqrt n}
\operatorname{Re}
\widehat{
u''(S_n)Z_n^2
}(2k)
+
\mathcal E_{n,k}^{(2)}(u).
\end{aligned}
\label{eq:direct-W-coordinate-difference}
\end{equation}

\begin{proposition}
\label{prop:direct-original-nonlinear-fixed-pair}
For \(q=2r\geq12\) and every fixed integer \(J\geq3\),
\begin{equation}
\begin{aligned}
\|W_{m,n}(u)\|_q
\leq{}&
C_{q,T}
[\rho_{m,n}\mathfrak E_2(u)]^{1/2}
+
C_{q,T}
[\rho_{m,n}\mathfrak E_1(u)]^{1/2}\\
&+
\operatorname{Err}_{m,n}^{\mathrm{diff,proj}}(u)\\
&+
C_{q,T}
\left(\frac{\rho_{m,n}}n
\right)^{1/2}
\left[\mathfrak E_3(u)
+
\mathcal R_{N,q,J}^{\mathrm{end}}(\mathscr Ru)
\right]^{1/2}\\
&+
C_{q,T}
\rho_{m,n}^{1/2}
N^{-c_q/2}
\mathfrak J_{K_q}(u)^{C_q/2}\\
&+
C_{q,T}
n^{-1/2}
\left[\mathfrak E_2(u)
+
N^{-c_q}
\mathfrak J_{K_q}(u)^{C_q}
\right]^{1/2}\\
&+
C_{q,T}
\left(\frac{\rho_{m,n}}n
\right)^{1/2}
\left[\mathfrak E_3(u)
+
N^{-c_q}
\mathfrak J_{K_q}(u)^{C_q}
\right]^{1/2}.
\end{aligned}
\label{eq:direct-W-estimate}
\end{equation}
\end{proposition}

\begin{proof}
Condition on \(\mathcal F_m\).  \eqref{eq:endpoint-product-space-square-function},
\eqref{eq:direct-W-conditional-centering}, and
\eqref{eq:direct-W-coordinate-difference} imply
\[
\begin{aligned}
\|W_{m,n}(u)\|_q
\leq C_q\Bigg\{&
\|
\mathcal Q_{m,n}^{\mathrm{diff},\mathrm R}(u)
\|_{L^r}^{1/2}
+
\rho_{m,n}^{1/2}\|A_n(u)\|_q+
\left\|\frac1n
\sum_{k=m}^{n-1}
\left|
\widehat{
u''(S_n)Z_n^2
}(2k)
\right|^2
\right\|_{L^{q/2}}^{1/2}\\
&+
\left\|\sum_{k=m}^{n-1}
|\mathcal E_{n,k}^{(2)}(u)|^2
\right\|_{L^{q/2}}^{1/2}
\Bigg\},
\end{aligned}
\]
where
\[
A_n(u)
=
\int_{\Torus}
\mathscr Ru(S_n)\,\dd\mu.
\]
By Lemma~\ref{lem:direct-two-time-projected-leakage},
\[
\begin{aligned}
\|
\mathcal Q_{m,n}^{\mathrm{diff},\mathrm R}(u)
\|_{L^r}^{1/2}
&\leq C_{q,T}
[\rho_{m,n}\mathfrak E_2(u)]^{1/2}
+
\operatorname{Err}_{m,n}^{\mathrm{diff,proj}}(u).
\end{aligned}
\]

For \(A_n(u)\), set \(F(z)=\mathscr Ru(|z|^2)\).
By \eqref{eq:replica-compact-jet-seminorm},
\[
\max_{|\beta|\leq4}
\|\partial^\beta F\|_{L^\infty(\mathbb R^2)}
\leq C\mathfrak J_5(u).
\]
Applying \eqref{eq:occupation-comparison-derivative-bound}
with \(r=1\) and the test function \(F\),
and averaging over \(t=\pm1\), we obtain
\[
\left|
\E A_n(u)
-
\int_0^\infty \mathscr Ru(x)e^{-x}\,\dd x
\right|
\leq C_TN^{-1/2}\mathfrak J_5(u).
\]
By \eqref{eq:endpoint-resonant-exponential-mean},
\[
\left|\int_0^\infty
\mathscr Ru(x)e^{-x}\,\dd x\right|^2\leq\mathfrak E_1(u).
\]
Thus, 
\[
|\E A_n(u)|^2
\leq
2\mathfrak E_1(u)
+C_TN^{-1}\mathfrak J_5(u)^2.
\]
Moreover,
\[
A_n(u)-\E A_n(u)
=\frac1{\sqrt n}X_n(\mathscr Ru).
\]
Therefore,
\eqref{eq:direct-homogeneous-endpoint} and
\eqref{eq:endpoint-resonant-energy-transfer} imply
\[
\|A_n(u)-\E A_n(u)\|_q^2
\leq
\frac{C_{q,T}}n
\left[
\mathfrak E_3(u)
+\mathcal R_{N,q,J}^{\mathrm{end}}(\mathscr Ru)
\right].
\]
Consequently,
\[
\begin{aligned}
\|A_n(u)\|_q^2
\leq{}&
C_{q,T}\mathfrak E_1(u)
+\frac{C_{q,T}}n
\left[
\mathfrak E_3(u)
+\mathcal R_{N,q,J}^{\mathrm{end}}(\mathscr Ru)
\right]\\
&+C_{q,T}N^{-1}\mathfrak J_5(u)^2.
\end{aligned}
\]

Next, Bessel's inequality and
\eqref{eq:occupation-second-weighted-bound} yield
\[
\left\|
\frac1n\sum_{k=m}^{n-1}
\left|
\widehat{u''(S_n)Z_n^2}(2k)
\right|^2
\right\|_{L^{q/2}}
\leq
\frac{C_{q,T}}n
\left[
\mathfrak E_2(u)
+N^{-1/q}\mathfrak J_5(u)^2
\right].
\]
Finally, Lemma~\ref{lem:endpoint-second-order-ledger},
with \(\psi=u\) and \(\mathcal I=\mathcal I_{m,n}\),
together with \eqref{eq:endpoint-H-energy-bound}, gives
\[
\left\|
\sum_{k=m}^{n-1}|\mathcal E_{n,k}^{(2)}(u)|^2
\right\|_{L^{q/2}}
\leq
C_{q,T}\frac{\rho_{m,n}}n
\left[
\mathfrak E_3(u)
+N^{-1/q}\mathfrak J_6(u)^2
\right].
\]

Taking \(K_q\geq6\), \(C_q\geq2\), and
\(0<c_q<1/q\), we have
\[
N^{-1}\mathfrak J_5(u)^2
\leq N^{-c_q}\mathfrak J_{K_q}(u)^{C_q},
\qquad
N^{-1/q}\mathfrak J_6(u)^2
\leq N^{-c_q}\mathfrak J_{K_q}(u)^{C_q}.
\]
Combining these bounds with
\eqref{eq:endpoint-product-space-square-function},
\eqref{eq:direct-W-conditional-centering}, and
\eqref{eq:direct-W-coordinate-difference},
and using \(\mathfrak J_5(u)\leq\mathfrak J_6(u)\),
proves \eqref{eq:direct-W-estimate}.
\end{proof}

\boldparagraph{Combining the increment estimates.}

Combining the conditionally centered term with the conditional mean
gives an increment estimate for a fixed pair of lengths.

We keep the different powers of the relative added length
\(\rho_{m,n}\) separate, since they give different sums over refinement
levels. Define the nonlinear remainder by
\begin{equation}
\begin{aligned}
\operatorname{Err}_{m,n}^{\mathrm{nonlin,dir}}(u)
:={}&
\operatorname{Err}_{m,n}^{\mathrm{diff,proj}}(u)\\
&+
C
\left(\frac{\rho_{m,n}}n
\right)^{1/2}
\mathcal R_{N,q,J}^{\mathrm{end}}(\mathscr Ru)^{1/2}\\
&+
C\rho_{m,n}^{1/2}
N^{-c_q/2}
\mathfrak J_{K_q}(u)^{C_q/2}\\
&+
Cn^{-1/2}
\left[\mathfrak E_2(u)^{1/2}
+
N^{-c_q/2}
\mathfrak J_{K_q}(u)^{C_q/2}
\right]\\
&+
C
\left(\frac{\rho_{m,n}}n
\right)^{1/2}
\left[\mathfrak E_3(u)^{1/2}
+
N^{-c_q/2}
\mathfrak J_{K_q}(u)^{C_q/2}
\right]\\
&+
C\rho_{m,n}
\left\{\sum_{j=1}^{J}
N^{-j}
[
\mathfrak E_{2j+2}(u)+\mathfrak E_{2j+3}(u)
]
\right\}^{1/2}\\
&+
C\rho_{m,n}
[
N^{(1-J)/2}
+
N^{-c_{q,J}/2}
]
\mathfrak H_{K_{q,J}+2}(u)\\
&+
CN^{-1/2}
[
\rho_{m,n}\mathfrak J_4(u)
+
\log(h+2)\mathfrak J_2(u)
].
\end{aligned}
\label{eq:direct-nonlinear-error}
\end{equation}

\begin{proposition}
\label{prop:block-adaptive-residual-hierarchy}
One has the exact identity
\begin{equation}
X_n(u)-X_m(u)-J_{m,n}^{\mathrm{lin}}(u)
=
W_{m,n}(u)
+
X_m(\psi_{m,n})
+
\mathcal E_{m,n}^{\mathrm{mean}}(u).
\label{eq:block-direct-hierarchy-identity}
\end{equation}
Consequently,
\begin{equation}
\begin{aligned}
\|
X_n(u)-X_m(u)-J_{m,n}^{\mathrm{lin}}(u)
\|_q
\leq{}&
C_{q,T}
[\rho_{m,n}\mathfrak E_2(u)]^{1/2}\\
&+
C_{q,T}
[\rho_{m,n}\mathfrak E_1(u)]^{1/2}\\
&+
C_{q,T}
\rho_{m,n}\mathfrak E_3(u)^{1/2}\\
&+
\operatorname{Err}_{m,n}^{\mathrm{nonlin,dir}}(u).
\end{aligned}
\label{eq:block-direct-hierarchy-bound}
\end{equation}
\end{proposition}

\begin{proof}
Use
\eqref{eq:direct-two-time-conditional-split},
the definition of \(W_{m,n}\), and
\eqref{eq:direct-conditional-mean-identity}.
Then apply
Propositions~\ref{prop:direct-original-nonlinear-fixed-pair} and~\ref{prop:direct-conditional-mean-drift}.
\end{proof}

Fix an integer \(J\geq3\) in
\eqref{eq:direct-nonlinear-error}, and define
\begin{equation}
\operatorname{Err}_{m,n}^{\mathrm{loc,dir}}(u)
=
\operatorname{Err}_{m,n}^{\mathrm{lead,proj}}(u)
+
\operatorname{Err}_{m,n}^{\mathrm{nonlin,dir}}(u).
\label{eq:block-direct-local-error}
\end{equation}

\begin{corollary}[Local increment estimate]
\label{cor:direct-exact-energy-increment}
Fix constants \(C_0,C_1,C_2,C_3\).
Suppose that \(\Lambda\geq1\), \(\mathcal A\geq0\), and
\begin{equation}
\mathfrak E_s(u)^{1/2}
\leq C_s\Lambda^{(s-2)/2}\mathcal A,
\qquad s=1,2,3.
\label{eq:direct-exact-energy-contract}
\end{equation}
For all sufficiently large \(N\), every pair \(m<n\) satisfying
\eqref{eq:block-window-lengths} and
\(\rho_{m,n}\Lambda\leq C_0\) obeys
\begin{equation}
\|X_n(u)-X_m(u)\|_{L^q}
\leq
C\rho_{m,n}^{1/2}\mathcal A
+
\operatorname{Err}_{m,n}^{\mathrm{loc,dir}}(u).
\label{eq:direct-exact-energy-increment}
\end{equation}
Here \(C\) and the large-\(N\) threshold depend only on
\(q,J,T,C_0,C_1,C_2,C_3\).
\end{corollary}

\begin{proof}
By \eqref{eq:direct-exact-energy-contract},
\(\Lambda\geq1\), and \(\rho_{m,n}\Lambda\leq C_0\),
\[
\rho_{m,n}^{1/2}
\left[\mathfrak E_1(u)^{1/2}+\mathfrak E_2(u)^{1/2}\right]
+\rho_{m,n}\mathfrak E_3(u)^{1/2}
\leq C\rho_{m,n}^{1/2}\mathcal A.
\]
The conclusion follows from
Corollary~\ref{cor:direct-leading-leakage},
Proposition~\ref{prop:block-adaptive-residual-hierarchy}, and
\eqref{eq:block-direct-local-error}.
\end{proof}

\subsection{Temporal refinement and summation of the projection errors}\label{subsec:block-temporal-refinement}

Now we sum the local projection errors over successively finer time grids to control their accumulation across nested polynomial lengths.
We apply the elementary finite-family \(L^q\) maximum bound and a
dyadic refinement tree to the preceding local estimates, with temporal
scale \(\Gamma\) and spectral scale \(\Lambda\).  The two scales remain
distinct.  Every local error is multiplied by the number of refinement
edges at its level, raised to the reciprocal moment order, before
summation.  Under the growth conditions stated below, the projection
errors tend to zero after this summation.

\boldparagraph{Finite-family maximum.}
\begin{lemma}
\label{lem:block-finite-family-maximum}
Let \(q\geq2\), let \(\mathcal I\) be a nonempty finite index set, and let
\((Y_\iota)_{\iota\in\mathcal I}\)
be random variables in \(L^q\).  Then
\begin{equation}
\left\|
\max_{\iota\in\mathcal I}
|Y_\iota|
\right\|_{L^2}
\leq
|\mathcal I|^{1/q}
\max_{\iota\in\mathcal I}
\|Y_\iota\|_{L^q}.
\label{eq:block-finite-family-maximum}
\end{equation}
\end{lemma}

\begin{proof}
Since \(q\geq2\),
\[
\begin{aligned}
\left\|\max_{\iota\in\mathcal I}|Y_\iota|\right\|_{L^2}
&\leq \left\|\max_{\iota\in\mathcal I}|Y_\iota|\right\|_{L^q}
\leq \left(\sum_{\iota\in\mathcal I}\|Y_\iota\|_{L^q}^q\right)^{1/q}\\
&\leq |\mathcal I|^{1/q}\max_{\iota\in\mathcal I}\|Y_\iota\|_{L^q}.
\end{aligned}
\]
This proves \eqref{eq:block-finite-family-maximum}.
\end{proof}

\boldparagraph{Refinement geometry and exact temporal summation.}

Recall that
\[
n_-=\lfloor Ne^{-T}\rfloor,
\qquad
n_+=\lfloor Ne^T\rfloor.
\]
Now, we reach every polynomial length from a coarse grid, using at most one step at each finer level.
\begin{lemma}[Temporal refinement at scale \(\Gamma\)]
\label{lem:active-scale-A-refinement-tree}
Fix \(T>0\).
There exists \(C_T<\infty\) such that, for all sufficiently large
\(N\), depending only on \(T\), and every real \(\Gamma\geq1\),
there exist a nonempty coarse set
\[
\mathcal C_{N,\Gamma}\subseteq\{n_-,\ldots,n_+\}
\]
and refinement edge sets
\[
\mathcal E_{N,\Gamma,\ell}
\subseteq
\{(m,n)\in\mathbb N^2:n_-\leq m<n\leq n_+\},
\qquad \ell\in\mathbb N_0,
\]
with the following properties.
The coarse set satisfies
\[
\#\mathcal C_{N,\Gamma}\leq C_T\Gamma.
\]
Writing \(M_\ell=\Gamma2^\ell\), every nonempty refinement level
satisfies
\begin{equation}
\begin{gathered}
\#\mathcal E_{N,\Gamma,\ell}\leq C_TM_\ell,
\qquad
\max_{(m,n)\in\mathcal E_{N,\Gamma,\ell}}\rho_{m,n}
\leq\frac{C_T}{M_\ell},
\\
M_\ell\leq C_TN.
\end{gathered}
\label{eq:direct-refinement-level-geometry}
\end{equation}
Only finitely many refinement levels are nonempty.
Every \(n\in\{n_-,\ldots,n_+\}\) can be reached from some point
of \(\mathcal C_{N,\Gamma}\) by a directed chain whose edges use
strictly increasing refinement levels, with at most one edge
from each level.
The chain may have no edges when
\(n\in\mathcal C_{N,\Gamma}\).
\end{lemma}

\begin{proof}
Write \(D=n_+-n_-\) and set \(B=\lceil\Gamma\rceil\).
Then
\[
\Gamma\leq B\leq2\Gamma.
\]
For all sufficiently large \(N\), depending only on \(T\),
\[
n_-\geq c_TN>0,
\qquad
1\leq D\leq C_TN.
\]
For \(\ell\in\mathbb N_0\) and integers \(0\leq k\leq B2^\ell\),
define
\[
v_{\ell,k}
:=
n_-+\left\lfloor\frac{kD}{B2^\ell}\right\rfloor.
\]
These grids are nested, since
\[
v_{\ell+1,2k}=v_{\ell,k}.
\]
Moreover,
\[
v_{\ell,0}=n_-,
\qquad
v_{\ell,B2^\ell}=n_+,
\]
and, for integers \(0\leq k<B2^\ell\),
\[
0\leq v_{\ell,k+1}-v_{\ell,k}
\leq\left\lceil\frac{D}{B2^\ell}\right\rceil.
\]

Set
\[
\mathcal C_{N,\Gamma}
:=
\{v_{0,k}:0\leq k\leq B\},
\qquad
J_*:=\min\{\ell\in\mathbb N_0:B2^\ell\geq D\}.
\]
For \(0\leq\ell<J_*\), let
\[
\mathcal E_{N,\Gamma,\ell}
:=
\bigl\{
(v_{\ell,k},v_{\ell+1,2k+1}):
0\leq k<B2^\ell,\;
v_{\ell,k}<v_{\ell+1,2k+1}
\bigr\},
\]
and set \(\mathcal E_{N,\Gamma,\ell}=\varnothing\) for
\(\ell\geq J_*\).
The coarse set is nonempty, and
\[
\#\mathcal C_{N,\Gamma}\leq B+1\leq3\Gamma,
\qquad
\#\mathcal E_{N,\Gamma,\ell}\leq B2^\ell\leq2M_\ell.
\]

If \((m,n)\in\mathcal E_{N,\Gamma,\ell}\), then
\(\ell<J_*\), so \(B2^\ell<D\).
The grid construction gives
\[
n-m
\leq\frac{D}{B2^{\ell+1}}+1
\leq\frac{3D}{2B2^\ell}.
\]
Consequently,
\[
\rho_{m,n}
=\frac{n-m}{n}
\leq\frac{3D}{2n_-B2^\ell}
\leq\frac{C_T}{M_\ell}.
\]
On every nonempty refinement level, we also have
\[
M_\ell=\Gamma2^\ell\leq B2^\ell<D\leq C_TN.
\]
This proves \eqref{eq:direct-refinement-level-geometry}.

At the terminal level \(J_*\), successive grid values differ
by at most one, because \(D/(B2^{J_*})\leq1\).
Since the grid is nondecreasing and has endpoints \(n_-\) and \(n_+\),
\[
\{v_{J_*,k}:0\leq k\leq B2^{J_*}\}
=\{n_-,\ldots,n_+\}.
\]
Fix \(n\) in this set and choose an integer \(k_{J_*}\) such that
\(v_{J_*,k_{J_*}}=n\).
Define
\[
k_\ell=\left\lfloor\frac{k_{\ell+1}}2\right\rfloor,
\qquad
n^{(\ell)}=v_{\ell,k_\ell},
\qquad \ell=J_*-1,\ldots,0,
\]
with \(n^{(J_*)}=n\).
Then \(n^{(0)}\in\mathcal C_{N,\Gamma}\).
If \(k_{\ell+1}=2k_\ell\), nestedness gives
\(n^{(\ell+1)}=n^{(\ell)}\).
If \(k_{\ell+1}=2k_\ell+1\), either the two values coincide or
\[
(n^{(\ell)},n^{(\ell+1)})
\in\mathcal E_{N,\Gamma,\ell}.
\]
Deleting consecutive repetitions from
\(n^{(0)},\ldots,n^{(J_*)}\) gives the required directed chain,
whose edge levels are strictly increasing.
If \(J_*=0\), the coarse grid coincides with the terminal grid
and therefore contains every integer from \(n_-\) to \(n_+\);
all refinement edge sets are empty. This completes the proof.
\end{proof}

With this construction, we bound the total contribution across refinement levels, accounting for both the number and size of the steps.
\begin{lemma}
\label{lem:direct-exact-temporal-power-sum}
Fix \(T>0\) and \(q\geq2\).
For every \(b\geq0\), there exists \(C_{b,q,T}<\infty\) such that,
for all sufficiently large \(N\), depending only on \(T\),
and every \(\Gamma\geq1\), the refinement edge sets
\(\mathcal E_{N,\Gamma,\ell}\) constructed in
Lemma~\ref{lem:active-scale-A-refinement-tree} satisfy
\begin{equation}
\begin{aligned}
\sum_{\substack{\ell\geq0\\
\mathcal E_{N,\Gamma,\ell}\neq\varnothing}}
(\Gamma2^\ell)^{1/q}
\max_{(m,n)\in\mathcal E_{N,\Gamma,\ell}}
\rho_{m,n}^{\,b}
\leq
C_{b,q,T}\log(N+2)
N^{(1/q-b)_+}
\Gamma^{\min\{1/q-b,0\}}.
\end{aligned}
\label{eq:direct-exact-temporal-power-sum}
\end{equation}
\end{lemma}

\begin{proof}
If all refinement edge sets are empty, the left-hand side is zero.
Otherwise, \eqref{eq:direct-refinement-level-geometry} gives
\[
\begin{aligned}
\sum_{\substack{\ell\geq0\\
\mathcal E_{N,\Gamma,\ell}\neq\varnothing}}
(\Gamma2^\ell)^{1/q}
\max_{(m,n)\in\mathcal E_{N,\Gamma,\ell}}
\rho_{m,n}^{\,b}
\leq
C_{b,T}
\sum_{\substack{\ell\geq0\\
\Gamma2^\ell\leq C_TN}}
(\Gamma2^\ell)^{1/q-b}.
\end{aligned}
\]
Summing the geometric progression and using \(\Gamma\geq1\)
yields
\[
\sum_{\substack{\ell\geq0\\
\Gamma2^\ell\leq C_TN}}
(\Gamma2^\ell)^{1/q-b}
\leq
C_{b,q,T}
\begin{cases}
N^{1/q-b}, & b<1/q,\\[2mm]
\log(N+2), & b=1/q,\\[2mm]
\Gamma^{1/q-b}, & b>1/q.
\end{cases}
\]
Combining these estimates proves
\eqref{eq:direct-exact-temporal-power-sum}.
\end{proof}

The applications use the following scales:
\begin{equation}
\begin{array}{c|c|c}
\text{family}
& \text{temporal scale }\Gamma
& \text{energy scale }\Lambda\\ \hline
\text{moving boundary} & \Lambda_N & \Lambda_N\\
\text{lowest-scale piece} & A & A\\
\text{lower shell} & A & A/R\\
\text{middle piece or upper shell} & A & 1
\end{array}
\label{eq:direct-two-scale-specialization}
\end{equation}
In the lower-shell row, \(R\) is the dyadic scale in \(y=Ax\),
with \(1\leq R\leq A\) in the applications.
All scales in the table satisfy \(1\leq\Lambda\leq\Gamma\).
For these scales, the refinement constructed in
Lemma~\ref{lem:active-scale-A-refinement-tree} satisfies
\[
\rho_{m,n}\Lambda
\leq C_T\frac{\Lambda}{\Gamma}2^{-\ell}
\leq C_T2^{-\ell},
\qquad
(m,n)\in\mathcal E_{N,\Gamma,\ell}.
\]

For the nonprojected refinement estimate below, retain the fixed
integers \(r\geq6\) and \(J\geq3\), with \(q=2r\).
Let \(\Gamma,\Lambda\geq1\) and \(\mathcal A\geq0\), and assume that
\[
\mathfrak E_s(u)^{1/2}
\leq C_s\Lambda^{(s-2)/2}\mathcal A,
\qquad s=1,\ldots,2J+3,
\]
with fixed constants \(C_s\).
Dependence on these constants and on \(J\) is suppressed.
Since \(\Lambda\geq1\), for \(1\leq d\leq J\),
\begin{equation}
\left[
\mathfrak E_{2d+2}(u)+\mathfrak E_{2d+3}(u)
\right]^{1/2}
\leq C_d\Lambda^{d+1/2}\mathcal A.
\label{eq:direct-high-energy-substitution}
\end{equation}

For \(s\geq1\), the Laguerre lowering identity yields
\[
\mathfrak E_s(\mathscr Ru)
=
\sum_{a\geq2}(a-1)^sa^2|\widehat u_a|^2
\leq
\mathfrak E_{s+2}(u).
\]
Using
\[
g_{\mathscr Ru}=\frac14\Delta g_u,
\qquad
\Phi_{\mathscr Ru}=\frac18\nabla\Delta g_u,
\]
the proof of Lemma~\ref{lem:resonant-lowering-jet-bounds} and
\eqref{eq:direct-generated-OU-seminorm} yield
\[
\dot{\mathfrak J}_K(\mathscr Ru)
\leq
C_K
\max_{|\boldsymbol{\nu}|\leq2K+4}
\|\partial^{\boldsymbol{\nu}}g_u\|_\infty
\leq
C_K\mathfrak H_{K+2}(u).
\]

The next expression collects bounds for the nonprojected local errors
after summation over refinement levels. At level \(\ell\), the factor
\((\Gamma2^\ell)^{1/q}\) accounts for the number of pairs before the
errors are summed. Define
\begin{equation}
\begin{aligned}
\mathcal S_N^{\mathrm{dir}}
:={}&
\mathcal A
\left[
N^{-1/2+1/q}
+
N^{-1/2}\Gamma^{1/q-1/2}\Lambda^{1/2}
\right]\\
&+
N^{-c_q/2}\mathfrak J_{K_q}(u)^{C_q/2}
\left[
N^{-1/2+1/q}
+
(1+N^{-1/2})\Gamma^{1/q-1/2}
\right]\\
&+
\left[
N^{-1/2}\Gamma^{1/q-1/2}
+
\Gamma^{1/q-1}
\right]
\mathcal A
\sum_{d=1}^{J}N^{-d/2}\Lambda^{d+1/2}\\
&+
\left[
N^{-1/2}\Gamma^{1/q-1/2}
+
\Gamma^{1/q-1}
\right]
\left[
N^{(1-J)/2}+N^{-c_{q,J}/2}
\right]
\mathfrak H_{K_{q,J}+2}(u)\\
&+
N^{-1/2}\Gamma^{1/q-1}\mathfrak J_4(u)
+
N^{-1/2+1/q}\log(N+2)\mathfrak J_2(u).
\end{aligned}
\label{eq:direct-complete-refinement-budget}
\end{equation}
Using the refinement edge sets constructed in
Lemma~\ref{lem:active-scale-A-refinement-tree},
we apply the preceding bounds and
Lemma~\ref{lem:direct-exact-temporal-power-sum}
to \eqref{eq:direct-nonlinear-error} and obtain
\begin{equation}
\begin{aligned}
\sum_{\substack{\ell\geq0\\
\mathcal E_{N,\Gamma,\ell}\neq\varnothing}}
(\Gamma2^\ell)^{1/q}
\max_{(m,n)\in\mathcal E_{N,\Gamma,\ell}}
\left[
\operatorname{Err}_{m,n}^{\mathrm{nonlin,dir}}(u)
-
\operatorname{Err}_{m,n}^{\mathrm{diff,proj}}(u)
\right]
\leq
C_{q,T}\log(N+2)\mathcal S_N^{\mathrm{dir}}.
\end{aligned}
\label{eq:direct-nonprojected-refinement-bound}
\end{equation}

\boldparagraph{Summation of the projected errors.}

Use the refinement edge sets \(\mathcal E_{N,\Gamma,\ell}\)
constructed in Lemma~\ref{lem:active-scale-A-refinement-tree}.
Assume that, for fixed \(\delta_E,\delta_J\geq0\),
\begin{equation}
\mathfrak E_2(u)\leq N^{\delta_E+o(1)},
\label{eq:direct-compact-projected-input}
\end{equation}
and
\begin{equation}
\mathfrak J_4(u)\leq N^{\delta_J+o(1)},
\label{eq:direct-effective-jet-bound}
\end{equation}
where both \(o(1)\) terms are uniform over the test functions
under consideration.

Fix \(q=2r\) with \(r\geq6\).
For a nonempty
\[
\mathbf k\in\mathsf C_r,
\]
put
\[
L=L(\mathbf k),
\qquad
K=K(\mathbf k),
\]
and define
\begin{equation}
a_r(L,K)
=
L-\max\{K-r+1,0\},
\label{eq:direct-composition-a}
\end{equation}
\begin{equation}
b_r(L,K)
=
L+r-K-1.
\label{eq:direct-composition-b}
\end{equation}
For every nonempty composition,
\[
a_r(L,K)>0,
\qquad
b_r(L,K)>0.
\]

Under these growth conditions, the projected errors vanish even after summing over all refinement levels.
\begin{proposition}
\label{prop:direct-projected-composition-absorption}
Retain \(q=2r\), where \(r\geq6\) is an integer, and assume
\eqref{eq:direct-compact-projected-input} and
\eqref{eq:direct-effective-jet-bound}.
For each \(\Gamma\geq1\), use the refinement edge sets
\(\mathcal E_{N,\Gamma,\ell}\) constructed in
Lemma~\ref{lem:active-scale-A-refinement-tree}.
Suppose that, for every nonempty \(\mathbf k\in\mathsf C_r\),
\begin{equation}
2K\delta_J+(r-K)\delta_E<a_r(L,K),
\label{eq:direct-difference-composition-margin}
\end{equation}
where \(L=L(\mathbf k)\) and \(K=K(\mathbf k)\).
Then there exists \(c>0\), depending only on the fixed parameters,
such that, for all sufficiently large \(N\),
\begin{equation}
\begin{aligned}
\sum_{\substack{\ell\geq0\\
\mathcal E_{N,\Gamma,\ell}\neq\varnothing}}
(\Gamma2^\ell)^{1/q}
\max_{(m,n)\in\mathcal E_{N,\Gamma,\ell}}
\Bigl[
\operatorname{Err}_{m,n}^{\mathrm{lead,proj}}(u)
+
\operatorname{Err}_{m,n}^{\mathrm{diff,proj}}(u)
\Bigr]
\leq N^{-c+o(1)},
\end{aligned}
\label{eq:direct-projected-refinement-decay}
\end{equation}
uniformly for \(\Gamma\geq1\) and over the test functions
under consideration.
\end{proposition}

\begin{proof}
If all refinement edge sets are empty, the left-hand side of
\eqref{eq:direct-projected-refinement-decay} is zero.
Otherwise, since \(m\asymp_Tn\asymp_TN\) and
\(\mathfrak E_1(u)\leq\mathfrak E_2(u)\),
\eqref{eq:direct-projected-B-P} and
\eqref{eq:direct-difference-B-P} give
\[
\begin{aligned}
B_{m,\mathcal I_{m,n}}(u)
+B_{m,n}^{\mathrm{diff}}(u)
&\leq
N^{-1+o(1)}\mathfrak J_4(u)^2,\\
P_{m,\mathcal I_{m,n}}(u)
+P_{m,n}^{\mathrm{diff}}(u)
&\leq
C_{q,T}
\left[
\rho_{m,n}\mathfrak E_2(u)
+
N^{-1+o(1)}\mathfrak J_4(u)^2
\right].
\end{aligned}
\]
For \(K<r\), use
\[
(x+y)^{(r-K)/q}
\leq x^{(r-K)/q}+y^{(r-K)/q},
\qquad x,y\geq0,
\]
since \(0<(r-K)/q\leq1/2\).
For \(K=r\), the corresponding zero power of \(P\) is \(1\).
Substituting the preceding bounds into
\eqref{eq:direct-projected-error},
\eqref{eq:direct-leading-projected-error}, and
\eqref{eq:direct-two-time-projected-error}, and using \(q=2r\),
we obtain
\begin{equation}
\begin{aligned}
\operatorname{Err}_{m,n}^{\mathrm{lead,proj}}(u)
+\operatorname{Err}_{m,n}^{\mathrm{diff,proj}}(u)
&\leq
C_{q,T}N^{-1/2+o(1)}\mathfrak J_4(u)
\\
&\quad+
C_{q,T}
\sum_{\substack{
\mathbf k\in\mathsf C_r\\
L(\mathbf k)\geq1}}
\Bigl[
N^{-L/q}
\rho_{m,n}^{(r-K)/q}
\mathfrak J_4(u)^{K/r}
\mathfrak E_2(u)^{(r-K)/q}
\\
&\quad\quad+
N^{-(L+r-K)/q+o(1)}\mathfrak J_4(u)
\Bigr].
\end{aligned}
\label{eq:direct-projected-errors-before-growth}
\end{equation}
Take the maximum over each nonempty
\(\mathcal E_{N,\Gamma,\ell}\), multiply by
\((\Gamma2^\ell)^{1/q}\), and sum over \(\ell\).
Apply Lemma~\ref{lem:direct-exact-temporal-power-sum}
with \(b=0\) and \(b=(r-K)/q\), together with
\eqref{eq:direct-compact-projected-input} and
\eqref{eq:direct-effective-jet-bound}.
Since
\[
\Gamma^{\min\{1/q-b,0\}}\leq1,
\qquad \Gamma\geq1,
\]
and logarithmic factors can be absorbed into \(N^{o(1)}\),
the left-hand side of
\eqref{eq:direct-projected-refinement-decay} is bounded by
\[
\begin{aligned}
N^{-1/2+1/q+\delta_J+o(1)}
+
\sum_{\substack{
\mathbf k\in\mathsf C_r\\
L(\mathbf k)\geq1}}
\Bigl[
N^{-a_r(L,K)/q+(K/r)\delta_J
+(r-K)\delta_E/q+o(1)}
+
N^{-b_r(L,K)/q+\delta_J+o(1)}
\Bigr].
\end{aligned}
\]

By \eqref{eq:direct-difference-composition-margin},
\[
-\frac{a_r(L,K)}q+\frac Kr\delta_J
+\frac{r-K}q\delta_E<0.
\]
Applying \eqref{eq:direct-difference-composition-margin}
with \(K=r\) and \(L=\lceil r/4\rceil\), and using
\eqref{eq:replica-composition-family} and
\eqref{eq:direct-composition-b}, we obtain
\[
\delta_J<\frac{\lceil r/4\rceil-1}{q}
\leq\frac{b_r(L,K)}q,
\qquad
-\frac12+\frac1q+\delta_J<0.
\]
Thus all limiting exponents are strictly negative.
The finiteness of \(\mathsf C_r\) and the uniformity of the
\(o(1)\) terms yield
\eqref{eq:direct-projected-refinement-decay}
for some \(c>0\), uniformly for \(\Gamma\geq1\) and over the
test functions under consideration. This completes the proof.
\end{proof}

\subsection{A maximal estimate from the increment bounds}

Now we combine the fixed-length and increment bounds to control the maximum over all nested polynomial lengths in the time window.
The endpoint moment estimate controls the coarse grid.  Summing the
increment bounds over successively finer grids then controls every
length in the window, with the error sum kept explicit.

\begin{theorem}[Maximal estimate from Laguerre-energy bounds]
\label{thm:refinement-safe-maximal-interface}
Fix \(T>0\), integers \(r\geq6\) and \(J\geq3\), and constants
\(C_1,C_2,C_3>0\). Set \(q=2r\).
There exist \(C<\infty\) and \(N_0\geq1\) such that, for every
\(N\geq N_0\), \(\Gamma\geq1\), \(1\leq\Lambda\leq\Gamma\),
\(\mathcal A\geq0\), and every deterministic
\(u\in C_c^\infty([0,\infty);\mathbb R)\) satisfying
\[
\mathfrak E_s(u)^{1/2}
\leq C_s\Lambda^{(s-2)/2}\mathcal A,
\qquad s=1,2,3,
\]
we have
\begin{equation}
\begin{aligned}
\left\|
\max_{n_-\leq n\leq n_+}|X_n(u)|
\right\|_{L^2}
\leq
C\Gamma^{1/q}
\bigl(\Lambda^{-1/2}+\Gamma^{-1/2}\bigr)\mathcal A
+\operatorname{Err}^{\mathrm{max}}(u),
\end{aligned}
\label{eq:block-refinement-safe-maximal}
\end{equation}
where
\begin{equation}
\begin{aligned}
\operatorname{Err}^{\mathrm{max}}(u)
:=
C\Gamma^{1/q}
\bigl[\mathcal R_{N,q,J}^{\mathrm{end}}(u)\bigr]^{1/2}
+
C\sum_{\substack{\ell\geq0\\
\mathcal E_{N,\Gamma,\ell}\neq\varnothing}}
(\Gamma2^\ell)^{1/q}
\max_{(m,n)\in\mathcal E_{N,\Gamma,\ell}}
\operatorname{Err}_{m,n}^{\mathrm{loc,dir}}(u).
\end{aligned}
\label{eq:block-refinement-max-error}
\end{equation}
The sets \(\mathcal C_{N,\Gamma}\) and
\(\mathcal E_{N,\Gamma,\ell}\) are those constructed in
Lemma~\ref{lem:active-scale-A-refinement-tree}.
The remainders \(\mathcal R_{N,q,J}^{\mathrm{end}}\) and
\(\operatorname{Err}_{m,n}^{\mathrm{loc,dir}}\) are defined in
\eqref{eq:direct-endpoint-remainder} and
\eqref{eq:block-direct-local-error}, respectively.

\end{theorem}

\begin{proof}
The chain property in
Lemma~\ref{lem:active-scale-A-refinement-tree} and Minkowski's
inequality give
\[
\begin{aligned}
\left\|
\max_{n_-\leq n\leq n_+}|X_n(u)|
\right\|_{L^2}
\leq
\left\|
\max_{c\in\mathcal C_{N,\Gamma}}|X_c(u)|
\right\|_{L^2}
+
\sum_{\substack{\ell\geq0\\
\mathcal E_{N,\Gamma,\ell}\neq\varnothing}}
\left\|
\max_{(m,n)\in\mathcal E_{N,\Gamma,\ell}}
|X_n(u)-X_m(u)|
\right\|_{L^2}.
\end{aligned}
\]
By Lemma~\ref{lem:direct-homogeneous-endpoint},
the \(s=1\) case of \eqref{eq:direct-exact-energy-contract}, and
Lemma~\ref{lem:block-finite-family-maximum},
the coarse-set cardinality bound in
Lemma~\ref{lem:active-scale-A-refinement-tree} yields
\[
\begin{aligned}
\left\|
\max_{c\in\mathcal C_{N,\Gamma}}|X_c(u)|
\right\|_{L^2}
\leq
C\Gamma^{1/q}\Lambda^{-1/2}\mathcal A
+C\Gamma^{1/q}
\bigl[\mathcal R_{N,q,J}^{\mathrm{end}}(u)\bigr]^{1/2}.
\end{aligned}
\]

On every refinement edge, \eqref{eq:direct-refinement-level-geometry}
and \(\Lambda\leq\Gamma\) give
\[
\rho_{m,n}\Lambda
\leq C_T\frac{\Lambda}{\Gamma}2^{-\ell}
\leq C_T.
\]
Corollary~\ref{cor:direct-exact-energy-increment}, with \(C_0=C_T\),
therefore gives
\begin{equation}
\|X_n(u)-X_m(u)\|_{L^q}
\leq
C\rho_{m,n}^{1/2}\mathcal A
+\operatorname{Err}_{m,n}^{\mathrm{loc,dir}}(u).
\label{eq:block-refinement-local-contract}
\end{equation}
For each nonempty level, applying
Lemma~\ref{lem:block-finite-family-maximum} and
\eqref{eq:direct-refinement-level-geometry} to
\eqref{eq:block-refinement-local-contract} yields
\[
\begin{aligned}
\left\|
\max_{(m,n)\in\mathcal E_{N,\Gamma,\ell}}
|X_n(u)-X_m(u)|
\right\|_{L^2}
\leq
C\mathcal A(\Gamma2^\ell)^{1/q-1/2}
+
C(\Gamma2^\ell)^{1/q}
\max_{(m,n)\in\mathcal E_{N,\Gamma,\ell}}
\operatorname{Err}_{m,n}^{\mathrm{loc,dir}}(u).
\end{aligned}
\]
Since \(q>2\), summing the geometric series gives
\[
\sum_{\substack{\ell\geq0\\
\mathcal E_{N,\Gamma,\ell}\neq\varnothing}}
(\Gamma2^\ell)^{1/q-1/2}
\leq C_q\Gamma^{1/q-1/2}.
\]
Combining these estimates and taking \(C\) sufficiently large
proves \eqref{eq:block-refinement-safe-maximal}, with the error
defined in \eqref{eq:block-refinement-max-error}.
\end{proof}

\section{Kummer estimates and the critical logarithmic function}
\label{app:kummer-critical-profile}

This appendix supplies the special-function estimates used in
Section \ref{sec:coefficient-calculus}.
Section \ref{subsec:kummer-semigroup-representation} gives the Kummer
representation of the Laguerre semigroup and continuity uniform on
compact parameter sets; see Proposition \ref{prop:kummer-semigroup-power} and Lemma~\ref{lem:uniform-local-semigroup-continuity}.
In Section \ref{subsec:kummer-uniform-asymptotics}, we establish the uniform
Kummer expansion; see Lemma \ref{lem:uniform-differentiated-Kummer}.
Finally, Section~\ref{subsec:kummer-mehler-heine-critical} identifies the
critical logarithmic function through the classical complementary exponential integral and proves its limit directly from the Kummer formula; see
Proposition~\ref{prop:critical-Riemann-sum} and Proposition~\ref{prop:critical-logarithmic-function}.

\subsection{The Kummer representation of the Laguerre semigroup}\label{subsec:kummer-semigroup-representation}

We now express the smoothed power functions through Kummer's function, providing explicit formulas for the later Abel-block estimates.
Recall the confluent hypergeometric function
\begin{equation}
{}_1F_1(a;b;z)
=
\sum_{m=0}^{\infty}
\frac{(a)_m}{(b)_m}
\frac{z^m}{m!},
\label{eq:kummer-1F1-definition}
\end{equation}
where \(b\notin\{0,-1,-2,\ldots\}\) and
\[
(a)_0=1,
\qquad
(a)_m=a(a+1)\cdots(a+m-1).
\]

For \(\alpha>0\), set
\[
h_\alpha(x)=x^\alpha,
\qquad
\widetilde h_\alpha(x)=x^\alpha\log x,
\]
with \(\widetilde h_\alpha(0)=0\).

\begin{proposition}
\label{prop:kummer-semigroup-power}
Let \(\alpha>0\).
Then, for every \(\tau>0\) and \(x\geq0\),
\begin{equation}
P_\tau h_\alpha(x)
=
\Gamma(1+\alpha)
\left(1-e^{-\tau}\right)^\alpha
{}_1F_1
\left(
-\alpha;1;-\frac{xe^{-\tau}}{1-e^{-\tau}}
\right).
\label{eq:kummer-semigroup-power}
\end{equation}
Moreover,
\begin{equation}
P_\tau\widetilde h_\alpha
=
\partial_\alpha P_\tau h_\alpha.
\label{eq:kummer-semigroup-log-derivative}
\end{equation}
For each fixed \(\tau>0\), the differentiation with respect to
\(\alpha\) is locally uniform for
\(\alpha\in \mathfrak A\Subset(0,\infty)\) and \(x\) in compact intervals.
\end{proposition}

\begin{proof}
We evaluate the radial Mehler representation
\eqref{eq:spectral-radial-Mehler}.  Put \(r=e^{-\tau}\) and
\(q=rx/(1-r)\).  Polar integration of the circular Gaussian density gives
\[
P_\tau h_\alpha(x)
=\frac{e^{-q}}{1-r}\int_0^\infty
 y^\alpha e^{-y/(1-r)}
 I_0\!\left(\frac{2\sqrt{rxy}}{1-r}\right)\,\dd y.
\]
Expanding \(I_0\) in its nonnegative power series and integrating each term gives
\[
P_\tau h_\alpha(x)
=\Gamma(1+\alpha)(1-r)^\alpha e^{-q}
 {}_1F_1(1+\alpha;1;q).
\]
Kummer's classical transformation
\cite[Eq.~(13.2.39)]{NISTDLMF} now proves
\eqref{eq:kummer-semigroup-power}.

The map \(\alpha\mapsto h_\alpha\) is differentiable in
\(L^2(\nu)\), with derivative \(\widetilde h_\alpha\).
The boundedness of \(P_\tau\) proves
\eqref{eq:kummer-semigroup-log-derivative} in that space.
For fixed \(\tau>0\), the radial Gaussian representation also permits differentiation under the integral uniformly for \(\alpha\) and \(x\) in compact sets: near zero, \(y^{\alpha_-}|\log y|\) is bounded, and at infinity the Gaussian kernel dominates every fixed power of \(y\) times \(|\log y|\).  Equivalently, the hypergeometric series and its parameter derivative converge locally uniformly.  This proves the stated pointwise identity and local uniformity.
\end{proof}

We shall also use semigroup continuity uniformly for \(\alpha\) and
\(x\) in compact sets.

\begin{lemma}
\label{lem:uniform-local-semigroup-continuity}
Let
\(\mathfrak A\Subset(0,\infty)\) and
\(0<M<\infty.\)
Then, as \(\tau\downarrow0\),
\begin{equation}
\sup_{\substack{\alpha\in \mathfrak A\\0\leq x\leq M}}
\left|
P_\tau h_\alpha(x)-h_\alpha(x)
\right|
\longrightarrow0,
\label{eq:kummer-power-semigroup-continuity}
\end{equation}
and
\begin{equation}
\sup_{\substack{\alpha\in \mathfrak A\\0\leq x\leq M}}
\left|
P_\tau\widetilde h_\alpha(x)-\widetilde h_\alpha(x)
\right|
\longrightarrow0.
\label{eq:kummer-log-semigroup-continuity}
\end{equation}
\end{lemma}

\begin{proof}
By the Hille--Hardy formula
\cite[Eq.~(6.2.25)]{AndrewsAskeyRoy1999}
(see also \cite[Eq.~(18.18.27)]{NISTDLMF}),
\(P_\tau\) has the positive kernel
\[
K_\tau(x,y)
=
\frac1{1-e^{-\tau}}
\exp\left(
-\frac{e^{-\tau}(x+y)}{1-e^{-\tau}}
\right)
I_0\left(
\frac{2e^{-\tau/2}\sqrt{xy}}
{1-e^{-\tau}}
\right)
\]
with respect to \(e^{-y}\dd y\), where \(I_0\) is the modified Bessel function of the first kind of order zero.
That is,
\[
P_\tau f(x)
=
\int_0^\infty
K_\tau(x,y)f(y)e^{-y}\,\dd y.
\]

The two families
\[
(\alpha,x)\longmapsto h_\alpha(x),
\qquad
(\alpha,x)\longmapsto
\widetilde h_\alpha(x)
\]
are uniformly continuous on
\(\mathfrak A\times[0,4M]\), where
\(\widetilde h_\alpha(0)=0\).
Using \(I_0(r)\leq e^r\) for \(r\geq0\), we obtain
\begin{equation}
K_\tau(x,y)e^{-y}
\leq
\frac1{1-e^{-\tau}}
\exp\left(
-\frac{
\left(
\sqrt y-e^{-\tau/2}\sqrt x
\right)^2
}{1-e^{-\tau}}
\right).
\label{eq:kummer-kernel-estimate}
\end{equation}

For continuity, use the radial Mehler representation
\eqref{eq:spectral-radial-Mehler}.  Put \(r=e^{-\tau}\) and, using the
same standard circular Gaussian \(Z\) for every \(x\), set
\[
Y_r(x)=\left|\sqrt r\,\sqrt x+\sqrt{1-r}\,Z\right|^2.
\]
Then
\[
\sup_{0\leq x\leq M}|Y_r(x)-x|
\leq M(1-r)+2\sqrt M\sqrt{1-r}\,|Z|+(1-r)|Z|^2
\longrightarrow0
\]
almost surely as \(r\uparrow1\).
For each realization of \(Z\), both \(x\) and \(Y_r(x)\) lie in
\([0,2M+2|Z|^2]\).  The two families \(h_\alpha\) and
\(\widetilde h_\alpha\), extended by zero at the origin, are jointly
continuous in \(\alpha\) and \(x\) on compact subsets of
\((0,\infty)\times[0,\infty)\).  Their differences at \(Y_r(x)\) and
\(x\) therefore tend to zero uniformly for
\(\alpha\in\mathfrak A\) and \(0\leq x\leq M\), almost surely.
Writing \(\alpha_+=\sup\mathfrak A\), the sum of their absolute
differences is bounded, uniformly in these parameters and \(0<r<1\), by
\[
C_{\mathfrak A,M}\bigl(1+|Z|^{2(\alpha_++1)}\bigr).
\]
This bound is integrable.  Dominated convergence applied to the
supremum of the differences, followed by
\(\sup|\E F|\leq\E\sup|F|\), proves
\eqref{eq:kummer-power-semigroup-continuity} and
\eqref{eq:kummer-log-semigroup-continuity}.
\end{proof}

\subsection{Uniform Kummer asymptotics}
\label{subsec:kummer-uniform-asymptotics}

Next we obtain Kummer asymptotics with uniform control of the derivatives needed for the ordinary and logarithmic Abel-block bounds.
In particular, we derive the large-\(z\) Kummer expansion needed below, uniformly on compact exponent sets and after finitely many \(z\)-derivatives and one derivative with respect to \(\alpha\).

\begin{lemma}
\label{lem:uniform-differentiated-Kummer}
Let \(\mathfrak A\Subset(0,\infty)\) and let
\(r\in\mathbb N_0\).
There exists \(C_{\mathfrak A,r}<\infty\) such that, for
\[
\alpha\in\mathfrak A,
\qquad
z\geq1,
\qquad
0\leq j\leq r,
\qquad
k\in\{0,1\},
\]
one has
\begin{equation}
\left|
\partial_\alpha^k
\partial_z^j
\left[
\Gamma(1+\alpha)
{}_1F_1(-\alpha;1;-z)
-z^\alpha
-\alpha^2z^{\alpha-1}
\right]
\right|
\leq
C_{\mathfrak A,r}
\left(1+\log(2+z)\right)^k
z^{\alpha-2-j}.
\label{eq:kummer-uniform-two-term}
\end{equation}
\end{lemma}

\begin{proof}
By Kummer's transformation
\cite[Eq.~(13.2.39)]{NISTDLMF},
\[
\Gamma(1+\alpha)
{}_1F_1(-\alpha;1;-z)
=\Gamma(1+\alpha)
e^{-z}
{}_1F_1(1+\alpha;1;z).
\]
Substituting
\[
a=1+\alpha,
\qquad
b=1,
\]
into the Poincar\'e-type expansion
\cite[Eq.~(13.7.1)]{NISTDLMF} gives
\[
\Gamma(1+\alpha)e^{-z}
{}_1F_1(1+\alpha;1;z)
=
z^\alpha
\left[
1+\frac{\alpha^2}{z}
+O_{\mathfrak A}(z^{-2})
\right],
\]
uniformly for \(\alpha\) in a sufficiently small fixed complex
neighborhood of \(\mathfrak A\), as \(z\to\infty\) along
the positive real axis.

We use the \({}_1F_1\) normalization, with the factor
\(\Gamma(b)\) retained.
In the differentiated formulas below, \(b=1+j\),
\(0\le j\le r\), so this factor is \(\Gamma(1+j)=j!\).
The uniform remainder bound follows from the connection formula
\cite[Eq.~(13.2.41)]{NISTDLMF}
and the error bounds
\cite[Eqs.~(13.7.4)--(13.7.9)]{NISTDLMF}.
In the parameter ranges used here, \(a\) stays in a fixed
compact set and \(b\) takes only finitely many values.
The parameter factors in these bounds are therefore uniformly
bounded for sufficiently large positive \(z\).
The reciprocal gamma factors are entire, so integer exponents
cause no singularity.
The additional term in the connection formula is exponentially
small after multiplication by \(e^{-z}\) and is included
in the remainder.
Using the standard differentiation formula for
\(e^{-z}{}_1F_1(a;b;z)\)
\cite[Eq.~(13.3.20)]{NISTDLMF} gives
\[
\partial_z^j
\left[
e^{-z}{}_1F_1(1+\alpha;1;z)
\right]
=
(-1)^j\frac{(-\alpha)_j}{j!}
e^{-z}{}_1F_1(1+\alpha;1+j;z).
\]
Applying the same Poincar\'e-type expansion with
\(a=1+\alpha\) and \(b=1+j\), for \(0\leq j\leq r\), then yields the differentiated leading terms
\(\partial_z^j z^\alpha\) and
\(\alpha^2\partial_z^j z^{\alpha-1}\): the factor
\(\Gamma(1+j)=j!\) from the unregularized Kummer function cancels the denominator above.  Consequently,
\[
\left|
\partial_z^j
\left[
\Gamma(1+\alpha)
{}_1F_1(-\alpha;1;-z)
-z^\alpha
-\alpha^2z^{\alpha-1}
\right]
\right|
\leq
C_{\mathfrak A,r}
z^{\alpha-2-j}
\]
for all sufficiently large \(z\).

For the exponent derivative, use the analyticity in \(\alpha\) and
Cauchy's formula on a circle centered at \(\alpha\), contained in the
same fixed complex neighborhood of \(\mathfrak A\), whose radius is
comparable to
\[
\frac1{1+\log(2+z)}.
\]
On this circle,
\[
z^{\operatorname{Re}\beta-\alpha}
\leq C,
\]
and Cauchy's estimate therefore yields
\[
\left|
\partial_\alpha
\partial_z^j
\left[
\Gamma(1+\alpha)
{}_1F_1(-\alpha;1;-z)
-z^\alpha
-\alpha^2z^{\alpha-1}
\right]
\right|
\leq
C_{\mathfrak A,r}
\left(1+\log(2+z)\right)
z^{\alpha-2-j}.
\]
This proves \eqref{eq:kummer-uniform-two-term} for all sufficiently
large \(z\).
The remaining compact range \(1\leq z\leq z_0\) follows from joint smoothness and compactness.
\end{proof}

For
\[
A\geq2,
\qquad
c\in\left\{\frac12,1\right\},
\]
set
\begin{equation}
\delta_{A,c}
=
A\left(1-e^{-c/A}\right),
\qquad
b_{A,c}
=
\frac{e^{-c/A}}{\delta_{A,c}}.
\label{eq:kummer-rho-b}
\end{equation}
Taylor expansion gives, uniformly in \(c\in\{1/2,1\}\) and \(A\geq2\),
\begin{equation}
\delta_{A,c}
=
c+O(A^{-1}),
\qquad
b_{A,c}
=
\frac1c+O(A^{-1}).
\label{eq:kummer-finite-A-parameters}
\end{equation}
In particular,
\[
\delta_{A,c}\asymp1,
\qquad
b_{A,c}\asymp1.
\]

Define the scaled semigroup function by
\begin{equation}
\mathsf S_{A,c,\alpha}(y)
=
A^\alpha
P_{c/A}h_\alpha(y/A)
\label{eq:def-kummer-scaled-one-semigroup-profile}
\end{equation}
By \eqref{eq:kummer-semigroup-power},
\begin{equation}
\mathsf S_{A,c,\alpha}(y)
=
\Gamma(1+\alpha)
\delta_{A,c}^{\,\alpha}
{}_1F_1
\left(
-\alpha;1;-b_{A,c}y
\right).
\label{eq:equi-kummer-scaled-one-semigroup-profile}
\end{equation}
Applying Lemma~\ref{lem:uniform-differentiated-Kummer} with
\(z=b_{A,c}y\), using \eqref{eq:kummer-finite-A-parameters}, and absorbing the remaining bounded range by joint smoothness and compactness, gives the following consequence.

\begin{corollary}
\label{cor:finite-A-one-profile-expansion}
Let \(\mathfrak A\Subset(0,\infty)\) and \(r\in\mathbb N_0\).
Uniformly for
\[
A\geq2,
\qquad
\alpha\in\mathfrak A,
\qquad
c\in\left\{\frac12,1\right\},
\qquad
y\geq1,
\qquad
0\leq j\leq r,
\]
one has
\begin{equation}
\partial_y^j
\mathsf S_{A,c,\alpha}(y)
=
\partial_y^j
\left[
e^{-\alpha c/A}y^\alpha
+
\alpha^2
\delta_{A,c}
e^{-(\alpha-1)c/A}
y^{\alpha-1}
\right]
+
O_{\mathfrak A,r}
\left(
y^{\alpha-2-j}
\right).
\label{eq:kummer-one-profile-expansion}
\end{equation}
After one derivative with respect to \(\alpha\), the remainder is
\[
O_{\mathfrak A,r}
\left(
\left[1+\log(2+y)\right]
y^{\alpha-2-j}
\right).
\]
\end{corollary}

\subsection{The critical logarithmic limit and the complementary exponential integral}
\label{subsec:kummer-mehler-heine-critical}

We next compute the nonzero rescaled logarithmic limit at \(\alpha=1\), where the corrected ordinary Abel block vanishes.
The special function needed at the critical exponent is the classical
complementary exponential integral
\cite[Eqs.~(6.2.3)--(6.2.4)]{NISTDLMF},
\[
\operatorname{Ein}(z)=\int_0^z\frac{1-e^{-v}}v\,\dd v,
\qquad
\operatorname{Ein}(y)=\gamma_{\!E}+\log y+E_1(y)\quad(y>0),
\]
where \(\gamma_{\!E}\) is Euler's constant and
\(E_1(y)=\int_y^\infty e^{-v}\,\dd v/v\).
The apparent singularity of the integrand at zero is removable, and
\(\operatorname{Ein}\) is entire.  Define
\begin{equation}
Q(y)=-\int_{1/2}^{1}\operatorname{Ein}(y/c)\,\dd c,
\qquad y\geq0.
\label{eq:critical-Ein-representation}
\end{equation}
This expresses the critical limiting function in terms of a classical
special function.  The averaging over \(c\in[1/2,1]\) comes from the
Abel multiplier used here.
In the Bessel normalization, the same function is
\begin{equation}
Q(y)
=\int_0^\infty
\frac{e^{-s/2}-e^{-s}}{s^2}
\left[J_0(2\sqrt{sy})-1\right]\,\dd s,
\label{eq:calculus-Q}
\end{equation}
where \(J_0\) is the Bessel function of order zero; this identity is
verified below.

At \(\alpha=1\), the reciprocal-gamma formula
\eqref{eq:gamma-equi-entire} and its derivative give
\(\widetilde\gamma_a(1)=1/[a(a-1)]\) for \(a\geq2\).
Thus, at finite scale, the spectral expression is
\begin{equation}
\mathcal P_{A,1}(y)
=A\sum_{a=2}^{\infty}\frac{m_A(a)}{a(a-1)}
 \left[L_a(y/A)-1\right].
\label{eq:calculus-critical-finite-A-profile}
\end{equation}
The Kummer representation gives a direct proof of its limit.

We now show that the profile and its derivatives get uniformly close to \(Q\) and its derivatives on each bounded interval as \(A\) grows.
\begin{proposition}
\label{prop:critical-Riemann-sum}
Let \(0<M<\infty\) and \(r\in\mathbb N_0\).
Then, as \(A\to\infty\),
\begin{equation}
\max_{0\leq j\leq r}
\sup_{0\leq y\leq M}
\left|
\partial_y^j\mathcal P_{A,1}(y)
-
\partial_y^jQ(y)
\right|
\longrightarrow0.
\label{eq:kummer-critical-Riemann-sum-convergence}
\end{equation}
Moreover, the integral in \eqref{eq:calculus-Q} may be differentiated up to order \(r\) under the integral sign, uniformly for
\(0\leq y\leq M\).
\end{proposition}

\begin{proof}
Set
\[
H_\alpha(c,y)=\Gamma(1+\alpha)c^\alpha
 {}_1F_1(-\alpha;1;-y/c),\qquad c>0.
\]
For \(\alpha\) in a fixed neighborhood of \(1\), the explicit
parameters in \eqref{eq:kummer-finite-A-parameters} and the entire
Kummer series imply
\[
\mathsf S_{A,c,\alpha}\longrightarrow H_\alpha(c,\cdot),
\qquad c\in\{1/2,1\},
\]
in \(C^r([0,M])\), together with one \(\alpha\)-derivative.
This is a compact-parameter assertion, including \(y=0\).

The affine normalization in
\eqref{eq:calculus-visible-ordinary-profile} gives exactly
\[
\Delta F_{A,\alpha}^{\circ}(y)
=F_{A,\alpha}(y)-F_{A,\alpha}(0)
 +A^{\alpha-1}\beta_{A,\alpha}y.
\]
Since \(\Delta F_{A,1}^{\circ}=0\), its shifted derivative
\(\mathscr D_A\) at \(\alpha=1\) equals its ordinary
\(\alpha\)-derivative.  The derivative of the last term is
\[
 y\bigl[\log A\,\beta_{A,1}
           +\widetilde\beta_{A,1}\bigr],
\]
which tends to zero in \(C^r([0,M])\) by
\eqref{eq:calculus-beta-bounds}.  Consequently,
\(\mathcal P_{A,1}\) converges in that space to
\[
\left.\partial_\alpha
 \bigl[H_\alpha(1/2,y)-H_\alpha(1,y)
       -H_\alpha(1/2,0)+H_\alpha(1,0)\bigr]
\right|_{\alpha=1}.
\]
To identify this limit, differentiation of the Kummer series gives
\[
\partial_c H_\alpha(c,y)=\alpha^2H_{\alpha-1}(c,y).
\]
This identity holds near \(\alpha=1\); the expression on the right
is analytic there.  At exponent zero the same series yields
\[
\left.\partial_\beta H_\beta(c,y)\right|_{\beta=0}
=\psi(1)+\log c+\operatorname{Ein}(y/c).
\]
Indeed, \(\partial_\beta(-\beta)_m|_{\beta=0}=-(m-1)!\)
for \(m\geq1\).  After subtracting the value at \(y=0\),
differentiating the preceding \(c\)-identity at \(\alpha=1\)
therefore gives \(\operatorname{Ein}(y/c)\).
Integration from \(c=1/2\) to \(c=1\) identifies the limit with
\eqref{eq:critical-Ein-representation}, proving
\eqref{eq:kummer-critical-Riemann-sum-convergence}.

For completeness, the Bessel integral and its fixed \(y\)-derivatives
are uniformly dominated on \([0,M]\).  Near \(s=0\),
\(J_0(2\sqrt{sy})-1=O_M(s)\), and its \(j\)-th
\(y\)-derivative is \(O_{M,j}(s^j)\) for \(j\geq1\).
At infinity, the Bessel series bounds these derivatives by a fixed
power of \(1+s\) times \(e^{C_M\sqrt s}\), which is absorbed by
\(e^{-s/2}\).  These bounds justify differentiation under the integral;
the next proof verifies its equality with \(Q\).
\end{proof}

\begin{proposition}
\label{prop:critical-logarithmic-function}
For \(y>0\),
\begin{equation}
Q'(y)
=
-\frac1y
\int_{1/2}^{1}
\left(1-e^{-y/c}\right)
\,\dd c.
\label{eq:calculus-Qprime}
\end{equation}
Consequently,
\begin{equation}
|Q'(y)|
\leq
\frac{C}{1+y},
\qquad
|Q(y)|
\leq
C\log(2+y),
\label{eq:calculus-Q-bounds}
\end{equation}
and
\begin{equation}
Q(y)
=
-\frac12\log y
+
O(1)
\qquad
(y\to\infty).
\label{eq:calculus-Q-asymptotic}
\end{equation}
\end{proposition}

\begin{proof}
We first verify the Bessel representation \eqref{eq:calculus-Q}. The domination established above justifies the following operations.
Use
\begin{equation}
e^{-s/2}-e^{-s}
=
s\int_{1/2}^{1}e^{-cs}\,\dd c.
\label{eq:calculus-Abel-c-integral}
\end{equation}
For \(c>0\), set
\[
\Phi_c(y)
=
\int_0^\infty
e^{-cs}
\frac{J_0(2\sqrt{sy})-1}{s}
\,\dd s.
\]
Then the right-hand side of \eqref{eq:calculus-Q} is
\[
\int_{1/2}^{1}\Phi_c(y)\,\dd c.
\]
Since
\[
y\partial_yJ_0(2\sqrt{sy})
=
s\partial_sJ_0(2\sqrt{sy}),
\]
integration by parts yields
\[
y\Phi_c'(y)
=
\int_0^\infty
e^{-cs}\partial_sJ_0(2\sqrt{sy})\,\dd s
=
-1
+
c\int_0^\infty
e^{-cs}J_0(2\sqrt{sy})\,\dd s.
\]

The Laplace--Bessel identity
\begin{equation}
\int_0^\infty
e^{-cs}J_0(2\sqrt{sy})\,\dd s
=
\frac1c e^{-y/c}
\label{eq:calculus-Laplace-Bessel}
\end{equation}
therefore implies
\[
\Phi_c'(y)
=
-\frac{1-e^{-y/c}}{y}.
\]
Since \(\Phi_c(0)=0\), the definition of \(\operatorname{Ein}\) gives
\(\Phi_c(y)=-\operatorname{Ein}(y/c)\). This proves
\eqref{eq:calculus-Q} and, after integration over \(c\in[1/2,1]\),
\eqref{eq:calculus-Qprime}.

For \(0<y\leq1\),
\[
1-e^{-y/c}\leq Cy,
\]
whereas for \(y\geq1\),
\[
1-e^{-y/c}\leq1.
\]
Hence
\[
|Q'(y)|
\leq
\frac{C}{1+y}.
\]
Since \(Q(0)=0\), integration proves
\eqref{eq:calculus-Q-bounds}.

Finally, the classical identity
\(\operatorname{Ein}(v)=\gamma_{\!E}+\log v+E_1(v)\),
together with \(0\leq E_1(v)\leq e^{-v}/v\), gives the logarithmic asymptotic directly from
\eqref{eq:critical-Ein-representation}. Equivalently, retaining the derivative information,
\[
\begin{aligned}
\int_{1/2}^{1}
\left(1-e^{-y/c}\right)
\,\dd c
&=
\frac12+O(e^{-y}),
\end{aligned}
\]
and therefore
\[
Q'(y)
=
-\frac1{2y}
+
O\left(\frac{e^{-y}}{y}\right).
\]
Integration proves
\eqref{eq:calculus-Q-asymptotic}.
\end{proof}

\section{Maximal estimates above a moving Laguerre cutoff}
\label{app:moving-labels}

This appendix controls the Laguerre modes above a cutoff that grows with
\(N\), uniformly over nested polynomial lengths.
The maximal argument uses the comparison and temporal refinement
estimates of Appendices~\ref{app:finite-replica-occupation}
and~\ref{app:block-addition-refinement}.
We verify the hypotheses needed to apply those estimates to the modes
above the moving cutoff after Abel smoothing.
We record the parameter choices explicitly so that these verifications
remain uniform in the exponent and cover both the ordinary and
logarithmic families as the cutoff varies with \(N\).
The resulting bounds are used in the Abel--sharp comparison of
Section~\ref{subsec:Abel-sharp-boundary-comparison}.

Fix \(\mathfrak A=[\alpha_-,\alpha_+]\Subset(0,\infty)\)
and \(T>0\), and let \(\Lambda_N\) be the dyadic cutoff
defined in \eqref{eq:moving-LambdaN}.
Take \(\chi\in C_c^\infty([0,\infty))\) with
\(0\leq\chi\leq1\), \(\chi=1\) on \([0,1]\), and
\(\chi=0\) on \([2,\infty)\).
For \(L>0\), set
\begin{align}
\ell_N
&=1+\log(\Lambda_N+2),
\label{eq:moving-ellN}\\
\mathsf R_N
&=L\log(N+2),
\qquad
\chi_N(x)=\chi(x/\mathsf R_N).
\label{eq:moving-spatial-cutoff}
\end{align}
Then \(\Lambda_N=N^{o(1)}\) and \(\ell_N=N^{o(1)}\).

For \(\alpha\in\mathfrak A\), the ordinary term above the moving Laguerre cutoff
in \eqref{eq:boundary-upper-definition} and its logarithmic counterpart
have the representations
\begin{equation}
\begin{aligned}
B_{\Lambda_N,\alpha}^{+,0}
&=
\left(I-\LagProj_{\leq\Lambda_N}\right)
P_{1/\Lambda_N}h_\alpha
=
\sum_{a>\Lambda_N}
e^{-a/\Lambda_N}\gamma_a(\alpha)L_a,\\
B_{\Lambda_N,\alpha}^{+,1}
&=
\left(I-\LagProj_{\leq\Lambda_N}\right)
P_{1/\Lambda_N}\widetilde h_\alpha
=
\sum_{a>\Lambda_N}
e^{-a/\Lambda_N}\widetilde\gamma_a(\alpha)L_a,
\end{aligned}
\label{eq:moving-upper-boundaries}
\end{equation}
where \(\widetilde h_\alpha(x)=x^\alpha\log x\) for \(x>0\),
with \(\widetilde h_\alpha(0)=0\).
For \(i\in\{0,1\}\), define
\begin{equation}
u_{N,\alpha}^{(i)}
=
\chi_N B_{\Lambda_N,\alpha}^{+,i},
\qquad
v_{N,\alpha}^{(i)}
=
(1-\chi_N)B_{\Lambda_N,\alpha}^{+,i}.
\label{eq:moving-boundary-split}
\end{equation}
Since \(\Lambda_N\) and \(\chi_N\) are independent of \(\alpha\),
\[
\begin{aligned}
\partial_\alpha B_{\Lambda_N,\alpha}^{+,0}
=B_{\Lambda_N,\alpha}^{+,1},\qquad
\partial_\alpha u_{N,\alpha}^{(0)}
=u_{N,\alpha}^{(1)},
\qquad
\partial_\alpha v_{N,\alpha}^{(0)}
=v_{N,\alpha}^{(1)}.
\end{aligned}
\]

\subsection{Derivative and energy estimates}

Now we bound the derivatives and Laguerre energies of the smoothed terms above the moving cutoff, preparing the maximal estimate.

\begin{proposition}
\label{prop:moving-upper-boundary-spectral-energies}
For every fixed \(s\geq0\),
\begin{equation}
\mathfrak E_s\!\left(B_{\Lambda_N,\alpha}^{+,i}\right)
\leq
C_{\mathfrak A,s}\ell_N^{2i}
\Lambda_N^{s-2\alpha_--1},
\label{eq:moving-upper-energy}
\end{equation}
uniformly for \(\alpha\in\mathfrak A\) and \(i\in\{0,1\}\).
\end{proposition}

\begin{proof}
By Proposition~\ref{prop:compact-coefficient-bounds},
\[
\begin{aligned}
\mathfrak E_s\!\left(B_{\Lambda_N,\alpha}^{+,i}\right)
\leq
C_{\mathfrak A}
\sum_{a>\Lambda_N}
a^{s-2\alpha_--2}
\log^{2i}(a+2)e^{-2a/\Lambda_N}\leq
C_{\mathfrak A,s}\ell_N^{2i}
\Lambda_N^{s-2\alpha_--1}.
\end{aligned}
\]
\end{proof}

\begin{lemma}
\label{lem:moving-semigroup-derivatives}
For every fixed integer \(j\geq0\), there exist \(C,\iota>0\),
depending only on \(\mathfrak A\) and \(j\), such that
\begin{equation}
\sup_{\alpha\in\mathfrak A}
\max_{i\in\{0,1\}}
\left|
\partial_x^j\partial_\alpha^iP_\tau h_\alpha(x)
\right|
\leq C(1+\tau^{-1}+x)^\iota,
\qquad 0<\tau\leq1,\quad x\geq0.
\label{eq:moving-semigroup-derivative}
\end{equation}
\end{lemma}

\begin{proof}
Put
\[
d_\tau=1-e^{-\tau},
\qquad
y_\tau=\frac{e^{-\tau}x}{d_\tau},
\qquad
\pi_j(\alpha)=\prod_{\ell=0}^{j-1}(\alpha-\ell),
\]
with \(\pi_0=1\).
By \eqref{eq:kummer-semigroup-power}, the differentiation formula
for \({}_1F_1\), and Kummer's transformation,
\[
\partial_x^jP_\tau h_\alpha(x)
=
e^{-j\tau}d_\tau^{\alpha-j}\pi_j(\alpha)
\E\!\left[
\frac{\Gamma(Y+\alpha+1)}{\Gamma(Y+j+1)}
\right],
\qquad
Y\sim\operatorname{Poisson}(y_\tau).
\]
Uniformly for \(\alpha\in\mathfrak A\) and \(k\geq0\),
\[
\max_{i\in\{0,1\}}
\left|
\partial_\alpha^i
\frac{\Gamma(k+\alpha+1)}{\Gamma(k+j+1)}
\right|
\leq C(k+1)^{\alpha_++1}.
\]
Since the law of \(Y\) is independent of \(\alpha\), this bound
justifies differentiation under the expectation.
Polynomial moment bounds for \(Y\), together with
\(d_\tau\asymp\tau\) and \(y_\tau\leq Cx/\tau\), give
\[
\begin{aligned}
\left|
\partial_x^j\partial_\alpha^iP_\tau h_\alpha(x)
\right|
\leq
C\tau^{-j}(1+|\log\tau|)
(1+x/\tau)^{\alpha_++1}\leq C(1+\tau^{-1}+x)^\iota
\end{aligned}
\]
for a sufficiently large \(\iota\), uniformly over the indicated
indices and parameters.
The logarithmic case is identified by
\eqref{eq:kummer-semigroup-log-derivative}.
\end{proof}

The finite expansion of the Laguerre polynomials gives, for every
fixed integer \(j\geq0\),
\begin{equation}
|\partial_x^jL_a(x)|
\leq C_j(a+1)^j(2+x)^a,
\qquad a\in\mathbb N_0,\quad x\geq0.
\label{eq:moving-crude-Laguerre-derivative}
\end{equation}

The two terms above the Laguerre cutoff admit the finite-sum representations
\begin{equation}
\begin{aligned}
B_{\Lambda_N,\alpha}^{+,0}
&=
P_{1/\Lambda_N}h_\alpha
-
\sum_{a=0}^{\Lambda_N}
e^{-a/\Lambda_N}\gamma_a(\alpha)L_a,\\
B_{\Lambda_N,\alpha}^{+,1}
&=
P_{1/\Lambda_N}\widetilde h_\alpha
-
\sum_{a=0}^{\Lambda_N}
e^{-a/\Lambda_N}\widetilde\gamma_a(\alpha)L_a.
\end{aligned}
\label{eq:moving-finite-sum}
\end{equation}

\begin{proposition}
\label{prop:moving-upper-boundary-growth}
Fix an integer \(j\geq0\), \(L>0\),
\(0<\varepsilon<1/4\), and \(\delta>0\).
For sufficiently small \(\gamma_0>0\), there exist \(C,\iota>0\)
such that, for all sufficiently large \(N\),
\begin{equation}
\sup_{\alpha\in\mathfrak A}
\max_{i\in\{0,1\}}
\left|\partial_x^jB_{\Lambda_N,\alpha}^{+,i}(x)\right|
\leq
CN^\delta(1+x)^\iota e^{\varepsilon x},
\qquad x\geq\mathsf R_N.
\label{eq:moving-boundary-growth}
\end{equation}
Moreover, uniformly for \(\alpha\in\mathfrak A\) and \(i\in\{0,1\}\),
\begin{equation}
\left|\partial_x^jv_{N,\alpha}^{(i)}(x)\right|
\leq
CN^\delta(1+x)^\iota e^{\varepsilon x}
\ind_{\{x\geq\mathsf R_N\}},
\qquad x\geq0.
\label{eq:moving-remainder-growth}
\end{equation}
\end{proposition}

\begin{proof}
By \eqref{eq:moving-finite-sum},
Lemma~\ref{lem:moving-semigroup-derivatives},
\eqref{eq:moving-crude-Laguerre-derivative}, and
Proposition~\ref{prop:compact-coefficient-bounds},
\[
\begin{aligned}
\sup_{\alpha\in\mathfrak A}
\max_{i\in\{0,1\}}
\left|\partial_x^bB_{\Lambda_N,\alpha}^{+,i}(x)\right|
\leq
C(1+\Lambda_N+x)^\iota+
C(\Lambda_N+1)^{j+1}\ell_N(2+x)^{\Lambda_N},
\end{aligned}
\]
uniformly for \(0\leq b\leq j\) and \(x\geq0\).
The first term is bounded by \(N^{o(1)}(1+x)^\iota\).
For the second, using
\(\Lambda_N\log(\Lambda_N+2)\leq C\gamma_0\log(N+2)\),
we obtain
\[
\begin{aligned}
C(\Lambda_N+1)^{j+1}\ell_N(2+x)^{\Lambda_N}
&\leq
C(\Lambda_N+1)^{j+1}\ell_N
\exp\!\left\{
C_\varepsilon\Lambda_N\log(\Lambda_N+2)
\right\}e^{\varepsilon x}\\
&\leq
N^{C_\varepsilon'\gamma_0+o(1)}e^{\varepsilon x}
\leq N^\delta e^{\varepsilon x},
\end{aligned}
\]
after choosing \(C_\varepsilon'\gamma_0<\delta\).
This proves \eqref{eq:moving-boundary-growth}, with the same
bound at every derivative order \(0,\ldots,j\).
Finally, Leibniz's rule, the support of \(1-\chi_N\), and
\(\|\chi_N^{(a)}\|_\infty\leq C_a\mathsf R_N^{-a}\)
give \eqref{eq:moving-remainder-growth}.
\end{proof}

We use the quantities \(\mathfrak J_K\) and \(\mathfrak H_M\)
defined in \eqref{eq:replica-compact-jet-seminorm} and
\eqref{eq:direct-generated-OU-seminorm}, respectively.

\begin{proposition}
\label{prop:moving-boundary-finite-jets}
Fix integers \(K,M\geq0\) and real numbers \(L,\delta>0\).
For sufficiently small \(\gamma_0>0\),
\begin{equation}
\sup_{\substack{\alpha\in\mathfrak A\\i\in\{0,1\}}}
\max\left\{
\mathfrak J_K(u_{N,\alpha}^{(i)}),
\mathfrak H_M(u_{N,\alpha}^{(i)})
\right\}
\leq N^\delta
\label{eq:moving-finite-jet-bounds}
\end{equation}
for all sufficiently large \(N\).
\end{proposition}

\begin{proof}
By \eqref{eq:moving-LambdaN} and
\eqref{eq:moving-spatial-cutoff},
\[
\Lambda_N\log(2+2\mathsf R_N)
\leq C_L\gamma_0\log(N+2)
\]
for all sufficiently large \(N\).
Using \eqref{eq:moving-finite-sum},
Lemma~\ref{lem:moving-semigroup-derivatives},
\eqref{eq:moving-crude-Laguerre-derivative}, and
Proposition~\ref{prop:compact-coefficient-bounds},
Leibniz's rule gives, for each fixed integer \(j\geq0\),
\begin{equation}
\begin{aligned}
\sup_{\substack{\alpha\in\mathfrak A\\i\in\{0,1\}}}
\left\|(u_{N,\alpha}^{(i)})^{(j)}\right\|_\infty
&\leq
C(1+\Lambda_N+\mathsf R_N)^\iota+
C(\Lambda_N+1)^{j+1}\ell_N
(2+2\mathsf R_N)^{\Lambda_N}\\
&\leq N^{C_L\gamma_0+o(1)}.
\end{aligned}
\label{eq:moving-compactified-scalar-jets}
\end{equation}
Here the constants and the \(o(1)\) term are uniform in
\(\alpha\) and \(i\).

Since
\(\operatorname{supp}u_{N,\alpha}^{(i)}
\subseteq[0,2\mathsf R_N]\),
the chain rule, \eqref{eq:replica-compact-jet-seminorm}, and
\eqref{eq:direct-generated-OU-ordinary-jet} yield
\[
\begin{aligned}
\mathfrak J_K(u_{N,\alpha}^{(i)})
&\leq
1+C_K(1+\mathsf R_N)^{(K+1)/2}
\max_{0\leq j\leq K+1}
\left\|(u_{N,\alpha}^{(i)})^{(j)}\right\|_\infty,\\
\mathfrak H_M(u_{N,\alpha}^{(i)})
&\leq
C_M(1+\mathsf R_N)^{M+3}
\max_{0\leq j\leq2M+7}
\left\|(u_{N,\alpha}^{(i)})^{(j)}\right\|_\infty.
\end{aligned}
\]
By \eqref{eq:moving-compactified-scalar-jets} and
\(\mathsf R_N=N^{o(1)}\), the left-hand side of
\eqref{eq:moving-finite-jet-bounds} is at most
\(N^{C_{K,M,L}\gamma_0+o(1)}\).
Choosing \(C_{K,M,L}\gamma_0<\delta\) proves the assertion.
\end{proof}

\begin{lemma}
\label{lem:moving-family-energy-contract}
For every fixed integer \(s\geq1\), choosing \(L\) sufficiently large
and then \(\gamma_0>0\) sufficiently small gives
\begin{equation}
\mathfrak E_s\!\left(u_{N,\alpha}^{(i)}\right)^{1/2}
\leq C_{\mathfrak A,s}\ell_N^i
\Lambda_N^{(s-1)/2-\alpha_-},
\qquad \alpha\in\mathfrak A,\quad i\in\{0,1\},
\label{eq:moving-family-energy-contract}
\end{equation}
for all sufficiently large \(N\).
\end{lemma}

\begin{proof}
Write \(B=B_{\Lambda_N,\alpha}^{+,i}\), \(u=\chi_NB\), and
\(m=\lceil s/2\rceil\).
Since \(B\) has no constant Laguerre mode and \(0\leq\chi_N\leq1\),
\[
\mathfrak E_0(u)
\leq\|\chi_NB\|_{L^2(e^{-x}\dd x)}^2
\leq\|B\|_{L^2(e^{-x}\dd x)}^2
=\mathfrak E_0(B).
\]
For \(\LagNum=-\LagGen\) and smooth \(f\), the product rule gives
\[
\LagNum(\chi_Nf)-\chi_N\LagNum f
=-2x\chi_N'f'-[x\chi_N''+(1-x)\chi_N']f.
\]
On iterating this identity, every commutator term contains a derivative
of \(\chi_N\) and is supported in \([\mathsf R_N,2\mathsf R_N]\).
Proposition~\ref{prop:moving-upper-boundary-growth}, with
\(\delta=1\), any fixed \(0<\varepsilon<1/4\), and derivative
orders \(0,\ldots,2m\), therefore yields
\[
\left\|
\LagNum^m(\chi_NB)-\chi_N\LagNum^mB
\right\|_{L^2(e^{-x}\dd x)}^2
\leq CN^2\mathsf R_N^C e^{-\mathsf R_N/2}
\leq N^{-2},
\]
after choosing \(L\) sufficiently large and then \(\gamma_0\)
sufficiently small.
Using \(\mathfrak E_{2m}(f)=\|\LagNum^mf\|_{L^2(e^{-x}\dd x)}^2\)
and \eqref{eq:moving-upper-energy}, we obtain
\[
\begin{aligned}
\mathfrak E_{2m}(u)
\leq 2\mathfrak E_{2m}(B)+2N^{-2}\leq C_{\mathfrak A,s}\ell_N^{2i}
\Lambda_N^{2m-2\alpha_--1},
\end{aligned}
\]
since \(N^{-2}=o(\ell_N^{2i}\Lambda_N^{2m-2\alpha_--1})\).
For \(s<2m\), Hölder's inequality applied to the Laguerre
coefficients, together with \eqref{eq:moving-upper-energy}
at order zero, gives
\[
\begin{aligned}
\mathfrak E_s(u)
=\sum_{a\geq1}a^s|\widehat u_a|^2\leq
\mathfrak E_0(u)^{1-s/(2m)}
\mathfrak E_{2m}(u)^{s/(2m)}\leq C_{\mathfrak A,s}\ell_N^{2i}
\Lambda_N^{s-2\alpha_--1}.
\end{aligned}
\]
Taking square roots proves
\eqref{eq:moving-family-energy-contract}.
\end{proof}

\subsection{Endpoint remainders and temporal refinement}

We next show that the fixed-length and increment remainders remain negligible after summation over the successive time grids.
Fix integers \(r\geq6\) and \(J\geq3\) such that
\begin{equation}
q=2r,
\qquad \frac1q<\alpha_-.
\label{eq:moving-local-moment-choice}
\end{equation}
Use common constants \(c_q,C_q,K_q,c_{q,J},K_{q,J}\) in
\eqref{eq:direct-endpoint-remainder} and
\eqref{eq:direct-complete-refinement-budget}.
Choose
\begin{equation}
0<\delta<\min\left\{
\frac1q,\frac{c_{q,J}}2,\frac{c_q}{C_q}
\right\}.
\label{eq:moving-direct-jet-margins}
\end{equation}

Let \(c_0>0\) satisfy the uniform Hoeffding bound
\begin{equation}
\sup_{\substack{n\geq1\\\theta\in\Torus}}
\Prob\bigl(S_n(\theta)>x\bigr)
\leq4e^{-c_0x},
\qquad x\geq0.
\label{eq:moving-hoeffding-bound}
\end{equation}
Choose \(\varepsilon\) and \(L\) so that
\begin{equation}
0<\varepsilon<\min\left\{\frac14,\frac{c_0}{4q}\right\},
\qquad \frac{c_0L}{2q}>1+\delta,
\label{eq:moving-tail-parameter-choice}
\end{equation}
with \(L\) also sufficiently large for
Lemma~\ref{lem:moving-family-energy-contract} at orders
\(1\leq s\leq2J+3\).
Finally, choose \(\gamma_0>0\) sufficiently small that
Proposition~\ref{prop:moving-upper-boundary-growth} holds with the chosen \(\varepsilon,\delta\) and derivative orders
\(0,\ldots,2J+4\), and
Proposition~\ref{prop:moving-boundary-finite-jets} holds with
\[
K=\max\{4,K_q,K_{q,J}+2\},
\qquad M=K_{q,J}+2,
\]
and exponent \(\delta\) in the derivative bound.

These parameters remain fixed as \(N\to\infty\), with choices
made in the order
\begin{equation}
\mathfrak A,T\longrightarrow(q,J)
\longrightarrow\delta\longrightarrow\varepsilon
\longrightarrow L\longrightarrow\gamma_0.
\label{eq:moving-parameter-order}
\end{equation}
In the remainder of this appendix, \(C_{\mathfrak A,T}\) may depend
on \(\chi\) and these fixed choices, but is uniform in \(N\),
\(\alpha\in\mathfrak A\), and \(i\in\{0,1\}\).
Constants \(C_T\) in refinement and cardinality bounds depend only on \(T\).

\begin{lemma}
\label{lem:moving-compactified-endpoint}
There exists \(c>0\) such that, for all sufficiently large \(N\),
\begin{equation}
\sup_{\substack{\alpha\in\mathfrak A\\i\in\{0,1\}}}
\mathcal R_{N,q,J}^{\mathrm{end}}(u_{N,\alpha}^{(i)})
\leq N^{-2c},\qquad i\in\{0,1\}.
\label{eq:moving-endpoint-remainder-bound}
\end{equation}

\end{lemma}

\begin{proof}
Substituting \eqref{eq:moving-family-energy-contract} and
\eqref{eq:moving-finite-jet-bounds} into
\eqref{eq:direct-endpoint-remainder}, and using
\(\dot{\mathfrak J}_K(f)\leq2\mathfrak J_K(f)\), yields
\[
\begin{aligned}
\sup_{\substack{\alpha\in\mathfrak A\\i\in\{0,1\}}}
\mathcal R_{N,q,J}^{\mathrm{end}}(u_{N,\alpha}^{(i)})
\leq C\left[
N^{-1+o(1)}+N^{1-J+2\delta}+N^{-c_{q,J}+2\delta}
\right]\leq N^{-2c}
\end{aligned}
\]
for some \(c>0\), since \(\Lambda_N,\ell_N=N^{o(1)}\)
and \eqref{eq:moving-direct-jet-margins} holds.

\end{proof}

Use the refinement constructed in
Lemma~\ref{lem:active-scale-A-refinement-tree}
with \(\Gamma=\Lambda_N\), and write
\[
\mathcal E_{N,j}=\mathcal E_{N,\Lambda_N,j},
\qquad \mathcal C_N=\mathcal C_{N,\Lambda_N},
\qquad M_j=\Lambda_N2^j.
\]
For \(u=u_{N,\alpha}^{(i)}\), evaluate the quantities in
\eqref{eq:direct-complete-refinement-budget} and
\eqref{eq:block-refinement-max-error} with
\begin{equation}
\Gamma=\Lambda=\Lambda_N,
\qquad \mathcal A=\ell_N^i\Lambda_N^{1/2-\alpha_-}.
\label{eq:moving-direct-two-scales}
\end{equation}

\begin{proposition}
\label{prop:moving-error-absorption}
\begin{equation}
\sup_{\substack{\alpha\in\mathfrak A\\i\in\{0,1\}}}
\operatorname{Err}^{\mathrm{max}}(u_{N,\alpha}^{(i)})
\longrightarrow0.
\label{eq:moving-direct-max-error}
\end{equation}
\end{proposition}

\begin{proof}
All estimates below are uniform in \(\alpha\in\mathfrak A\) and
\(i\in\{0,1\}\).
By \eqref{eq:moving-endpoint-remainder-bound},
\[
\Lambda_N^{1/q}
\bigl[\mathcal R_{N,q,J}^{\mathrm{end}}(u)\bigr]^{1/2}
\leq\Lambda_N^{1/q}N^{-c}\longrightarrow0.
\]

Using \eqref{eq:moving-family-energy-contract} and
\eqref{eq:moving-finite-jet-bounds} to estimate \(\mathcal S_N^{\mathrm{dir}}\) defined in
\eqref{eq:direct-complete-refinement-budget},
and applying \eqref{eq:direct-nonprojected-refinement-bound},
we obtain
\[
\begin{aligned}
\sum_{\substack{j\geq0\\
\mathcal E_{N,j}\neq\varnothing}}
M_j^{1/q}
\max_{(m,n)\in\mathcal E_{N,j}}
\left[
\operatorname{Err}_{m,n}^{\mathrm{nonlin,dir}}(u)
-
\operatorname{Err}_{m,n}^{\mathrm{diff,proj}}(u)
\right]
\leq
C\log(N+2)\mathcal S_N^{\mathrm{dir}}
\longrightarrow0.
\end{aligned}
\]
The convergence follows from
\(\Lambda_N,\ell_N=N^{o(1)}\) and
\eqref{eq:moving-direct-jet-margins}.

By \eqref{eq:moving-family-energy-contract} with \(s=2\)
and \eqref{eq:moving-finite-jet-bounds},
using \(\Lambda_N,\ell_N=N^{o(1)}\) and \(K\geq4\), we have
\[
\mathfrak E_2(u)\leq N^{o(1)},
\qquad \mathfrak J_4(u)\leq N^\delta.
\]
For every nonempty \(\mathbf k\in\mathsf C_r\),
\[
2K(\mathbf k)\delta\leq q\delta<1
\leq a_r\bigl(L(\mathbf k),K(\mathbf k)\bigr).
\]
Proposition~\ref{prop:direct-projected-composition-absorption},
with \((\delta_E,\delta_J)=(0,\delta)\), therefore shows that
the projected sum in \eqref{eq:direct-projected-refinement-decay}
tends to zero.
Combining these bounds with \eqref{eq:block-direct-local-error}
and \eqref{eq:block-refinement-max-error} proves
\eqref{eq:moving-direct-max-error}.
\end{proof}

\subsection{Maximal estimates above the cutoff}

Now we combine the energy and remainder bounds to control the smoothed Laguerre tail above the moving cutoff uniformly over nested lengths.

\begin{proposition}
\label{cor:moving-upper-boundary-input}
With the parameter choices in \eqref{eq:moving-parameter-order},
\begin{equation}
\max_{i\in\{0,1\}}
\mathfrak M_{N,T;\mathfrak A}
\bigl(B_{\Lambda_N,\bullet}^{+,i}\bigr)
\longrightarrow0.
\label{eq:moving-upper-boundary-input}
\end{equation}
\end{proposition}

\begin{proof}
All estimates below are uniform in \(\alpha\in\mathfrak A\)
and \(i\in\{0,1\}\).

Apply Theorem~\ref{thm:refinement-safe-maximal-interface}
to \(u_{N,\alpha}^{(i)}\), with the parameters in
\eqref{eq:moving-direct-two-scales} and the energy bounds of
Lemma~\ref{lem:moving-family-energy-contract}.
Taking the \(L^2(\mathfrak A)\)-norm in \(\alpha\) and using
Proposition~\ref{prop:moving-error-absorption}, we obtain
\begin{equation}
\mathfrak M_{N,T;\mathfrak A}
\bigl(u_{N,\bullet}^{(i)}\bigr)
\leq C_{\mathfrak A,T}\ell_N^i
\Lambda_N^{1/q-\alpha_-}+o(1).
\label{eq:moving-compactified-maximum}
\end{equation}

For the remainder outside the intensity cutoff, \eqref{eq:moving-remainder-growth}
at derivative order zero and \eqref{eq:moving-hoeffding-bound}
imply, uniformly in \(n\geq1\) and \(\theta\in\Torus\),
\[
\begin{aligned}
\E\left|v_{N,\alpha}^{(i)}(S_n(\theta))\right|^q
&\leq CN^{q\delta}
\sum_{\substack{k\in\mathbb N_0\\
k\geq\lfloor\mathsf R_N\rfloor}}
(2+k)^{\iota q}e^{q\varepsilon(k+1)}
\Prob\bigl(S_n(\theta)\geq k\bigr)\\
&\leq CN^{q\delta}(1+\mathsf R_N)^{\iota q}
e^{-c_0\mathsf R_N/2}.
\end{aligned}
\]
Since \(\mathsf R_N=L\log(N+2)\),
\eqref{eq:moving-tail-parameter-choice} implies, for sufficiently
large \(N\),
\[
\sup_{\substack{\alpha\in\mathfrak A\\
n\geq1,\ \theta\in\Torus}}
\left\|v_{N,\alpha}^{(i)}(S_n(\theta))\right\|_{L^q}
\leq CN^\delta(1+\mathsf R_N)^\iota
e^{-c_0\mathsf R_N/(2q)}
\leq N^{-1}.
\]
By \eqref{eq:model-In}, Minkowski's inequality, and
\eqref{eq:model-Xn},
\[
\begin{aligned}
\left\|I_n(v_{N,\alpha}^{(i)})\right\|_{L^q}
&\leq\int_{\Torus}
\left\|v_{N,\alpha}^{(i)}(S_n(\theta))\right\|_{L^q}
\,\dd\mu(\theta)
\leq N^{-1},\\
\left\|X_n(v_{N,\alpha}^{(i)})\right\|_{L^q}
&\leq2\sqrt n\,
\left\|I_n(v_{N,\alpha}^{(i)})\right\|_{L^q}
\leq C_TN^{-1/2},
\qquad n_-\leq n\leq n_+.
\end{aligned}
\]
Since \(n_+-n_-+1\leq C_TN\),
Lemma~\ref{lem:block-finite-family-maximum} yields
\[
\left\|\max_{n_-\leq n\leq n_+}
\left|X_n(v_{N,\alpha}^{(i)})\right|\right\|_{L^2}
\leq C_{\mathfrak A,T}N^{-1/2+1/q}.
\]
Taking the \(L^2(\mathfrak A)\) norm gives
\begin{equation}
\mathfrak M_{N,T;\mathfrak A}
\bigl(v_{N,\bullet}^{(i)}\bigr)
\leq C_{\mathfrak A,T}N^{-1/2+1/q}.
\label{eq:moving-upper-spatial-remainder}
\end{equation}

Combining \eqref{eq:moving-compactified-maximum} and
\eqref{eq:moving-upper-spatial-remainder} with
\eqref{eq:moving-boundary-split} and subadditivity yields
\[
\begin{aligned}
\mathfrak M_{N,T;\mathfrak A}
\bigl(B_{\Lambda_N,\bullet}^{+,i}\bigr)
&\leq C_{\mathfrak A,T}
\left[
\ell_N^i\Lambda_N^{1/q-\alpha_-}
+N^{-1/2+1/q}
\right]+o(1)\longrightarrow0,
\end{aligned}
\]
because \(1/q<\alpha_-\), \(q>2\), and
\(\ell_N=O(\log\Lambda_N)\).
\end{proof}

\section{Uniform expectation estimates for deterministic centering}
\label{app:expectation-envelope}

The fixed-exponent estimate
\(\E T_n^{(p)}=m_p+O_p(n^{-1})\) and the associated
small-ball estimates, Gaussian smoothing, Stein comparison, and
correction for covariance anisotropy appear in the companion paper
\cite[Sections~3--4 and Theorem~1.2]{ChangChenDuanFangZhangAS}.
We reproduce that method here, keeping its constants uniform for
\(p\) in a fixed compact subset of \((0,\infty)\).
In particular, a single derivative bound that is integrable under the
distributions used in the comparison is needed when the exponent interval
contains \(2\).  The proof below establishes
this uniformity directly; the fixed-exponent statement alone does not
supply it.  We also bound the measure of angles at which the intensity
is small and prove estimates for the Abel blocks, as needed for
Proposition~\ref{prop:global-Abel-endpoint}.

\subsection{Angular covariance and small-ball estimates}
\label{subsec:expectation-one-angle-gaussian-comparison}

Now we estimate the covariance at each angle and the probability of small intensities, preparing the uniform comparison of positive powers.
For
\[
j\geq0,
\qquad
\theta\in\Torus,
\]
put
\begin{equation}
v_j(\theta)
=
        \begin{pmatrix}
\cos(j\theta)\\
\sin(j\theta)
        \end{pmatrix}
\in\R^2.       \label{eq:expectation-vj}
\end{equation}
Define
\begin{equation}
\begin{aligned}
X_{n,\theta}
= \frac1{\sqrt n}
\sum_{j=0}^{n-1}
\varepsilon_jv_j(\theta)                                 =
\frac1{\sqrt n}
\begin{pmatrix}
\operatorname{Re}q_n(e^{i\theta})\\
\operatorname{Im}q_n(e^{i\theta})
\end{pmatrix}.
\end{aligned}
\label{eq:expectation-Xntheta}
\end{equation}
Then
\begin{equation}
S_n(\theta)=|X_{n,\theta}|^2.
\label{eq:expectation-Sn-vector}
\end{equation}
Recall from \eqref{eq:model-In} that, for a measurable function
\(f:[0,\infty)\to\R\),
\[
I_n(f)
=
\int_{\Torus}
f\!\left(S_n(\theta)\right)\,\dd\mu(\theta).
\]

The covariance matrix at one angle is
\begin{equation}
\Sigma_n(\theta)
=
\E
\left[
X_{n,\theta}X_{n,\theta}^{\mathsf T}
\right]
=
\frac1n
\sum_{j=0}^{n-1}
v_j(\theta)v_j(\theta)^{\mathsf T}.
\label{eq:expectation-An}
\end{equation}
Recall that the one-sided Dirichlet kernel is
\begin{equation}
D_n(u)
=
\sum_{j=0}^{n-1}e^{iju}.
\label{eq:expectation-Dirichlet-kernel}
\end{equation}
Directly from the definition of
\(\Sigma_n(\theta)\),
\begin{equation}
\begin{aligned}
\Sigma_n(\theta)
&= \frac12I_2
+
\frac1{2n}
\begin{pmatrix}
\operatorname{Re}D_n(2\theta) & \operatorname{Im}D_n(2\theta)\\
\operatorname{Im}D_n(2\theta) & -\operatorname{Re}D_n(2\theta)
\end{pmatrix}.
\end{aligned}
\label{eq:expectation-covariance-formula}
\end{equation}
Consequently,
\begin{equation}
\left\|
\Sigma_n(\theta)-\frac12I_2
\right\|_{\mathrm{op}}
=
\frac{|D_n(2\theta)|}{2n}.
\label{eq:expectation-covariance-deviation}
\end{equation}
Here, \(\|\cdot\|_{\mathrm{op}}\) denotes the operator norm of a matrix.
The Dirichlet kernel satisfies
\begin{equation}
|D_n(u)|
\leq
\min
\left\{
n,\,
\frac{C}
{\operatorname{dist}(u,2\pi\Z)}
\right\}.
\label{eq:expectation-Dirichlet-pointwise}
\end{equation}

Recall the Rademacher field \(Z_n=Z_n^{\mathrm R}\) from
the coordinate-path notation and \eqref{eq:direct-three-fields}:
\begin{equation}\label{eq:sharp-Sn-h}
Z_n(\theta)=\frac{q_n(e^{i\theta})}{\sqrt n}
=\frac1{\sqrt n}\sum_{j=0}^{n-1}\varepsilon_je^{ij\theta}.
\end{equation}
Thus
\[
X_{n,\theta}
=\bigl(\operatorname{Re}Z_n(\theta),
\operatorname{Im}Z_n(\theta)\bigr)^{\mathsf T},
\qquad S_n(\theta)=|Z_n(\theta)|^2.
\]
For real matrices \(A\) and \(B\) of the same size, we write
\(A:B=\operatorname{tr}(A^{\mathsf T}B)\).

The following bound is a consequence of Khintchine's inequality
and the second-moment estimate; see
\cite[Theorem~9.21]{ODonnell2014}.

\begin{lemma}\label{lem:sharp-positive}
For every \(r>0\),
\begin{equation}\label{eq:sharp-positive}
 \sup_{n\ge1}\sup_{\theta\in\Torus}\E|Z_n(\theta)|^r<\infty.
\end{equation}
\end{lemma}

The following small-ball estimate at the scale of the coefficients
is a standard consequence of the multidimensional Ess\'een
concentration inequality; see
\cite[Lemma~3.4]{RudelsonVershynin2009}.
The same coefficient-scale estimate is recorded in
\cite[Lemma~3.2]{ChangChenDuanFangZhangAS}.
We include the Fourier proof, keeping the dependence on the dimension
and the size of the coefficients explicit for the later application
to several angles.

\begin{lemma}\label{lem:sharp-smallball}
Let \(d\ge1\) be an integer, fix \(L,\lambda>0\), and let
\(0<h\le1\). Suppose \(X=\sum_j\varepsilon_ja_j\in\R^d\) is a
finite sum with deterministic vectors satisfying
\[
 |a_j|\leq Lh,
 \qquad \sum_j a_ja_j^{\mathsf T}\geq\lambda I_d.
\]
Then, for every \(r\geq h\),
\begin{equation}\label{eq:sharp-smallball}
 \sup_{x\in\R^d}\Prob(|X-x|\leq r)
 \leq C_{d,L,\lambda}\min(1,r^d).
\end{equation}
\end{lemma}

\begin{proof}

Let
\[
 \kappa_{\mathrm{sb}}(u)=\prod_{\ell=1}^d
 \left(\frac{\sin u_\ell}{u_\ell}\right)^2,
\]
where \(\sin t/t\) is extended continuously at \(t=0\).
Since \(\kappa_{\mathrm{sb}}\) is nonnegative and bounded below near zero,
and its Fourier transform is bounded and supported in
\([-2,2]^d\), Fourier inversion yields, for sufficiently
small \(b=b(d,L)>0\),
\[
 \Prob(|X-x|\leq r)
 \leq C\E\kappa_{\mathrm{sb}}\left(\frac{b}r(X-x)\right)
 \leq Cr^d\int_{[-2b/r,2b/r]^d}
       |\phi_X(\xi)|\dd\xi,
\]
where \(\phi_X(\xi)=\E e^{i\langle\xi,X\rangle}\).
 On this integration box,
\(|\langle a_j,\xi\rangle|
 \leq 2\sqrt d\,Lb\le1\). Using \(|\cos t|\leq e^{-ct^2}\) for \(|t|\le1\), we obtain
\[
 \begin{aligned}
 |\phi_X(\xi)|
 =\prod_j|\cos\langle a_j,\xi\rangle|
 \leq e^{-c\lambda|\xi|^2}.
 \end{aligned}
\]
Combining these
bounds proves \eqref{eq:sharp-smallball}.

\end{proof}

\begin{corollary}\label{cor:sharp-negative}
Under the assumptions of Lemma~\ref{lem:sharp-smallball},
if \(d=2\), then for every \(0<\beta<2\),
\begin{equation}\label{eq:sharp-negative}
 \sup_{z\in\R^2}\E(|X-z|+h)^{-\beta}
 \leq C_{\beta,L,\lambda}.
\end{equation}
\end{corollary}

\begin{proof}
We integrate the small-ball bound as in
\cite[Corollary~3.3]{ChangChenDuanFangZhangAS}.
For \(a\geq0\),
\[
 (a+h)^{-\beta}
 \leq 1+\beta\int_h^1
 t^{-\beta-1}\ind_{\{a\leq t\}}\dd t.
\]
By Tonelli's theorem and \eqref{eq:sharp-smallball} with \(d=2\),
\[
 \sup_{z\in\R^2}\E(|X-z|+h)^{-\beta}
 \leq 1+C_{L,\lambda}\beta\int_h^1t^{1-\beta}\dd t
 \leq C_{\beta,L,\lambda},
\]
since \(\beta<2\).
\end{proof}

Define the bad-angle and good-angle sets by
\begin{equation}\label{eq:sharp-angular-sets}
\begin{aligned}
\Xi_n^{\mathrm{bad}}
&=\left\{\theta\in\Torus:|D_n(2\theta)|>\frac n{16}\right\},\\
\Xi_n^{\mathrm{good}}
&=\left\{\theta\in\Torus:|D_n(2\theta)|\leq\frac n{16}\right\}
=\Torus\setminus\Xi_n^{\mathrm{bad}}.
\end{aligned}
\end{equation}

\begin{lemma}\label{lem:sharp-angular}
The bad-angle set satisfies
\begin{equation}\label{eq:sharp-angular-size}
 \mu(\Xi_n^{\mathrm{bad}})\leq Cn^{-1}.
\end{equation}
For \(\theta\in\Xi_n^{\mathrm{good}}\),
\[
 \Sigma_n(\theta)\geq\frac{15}{32}I_2.
\]
Moreover,
\begin{equation}\label{eq:sharp-rho-square}
 \int_{\Torus}\frac{|D_n(2\theta)|^2}{n^2}\dd\mu(\theta)
 =\frac1n.
\end{equation}
\end{lemma}

\begin{proof}
By \eqref{eq:expectation-Dirichlet-pointwise},
\[
 \mu\{t\in\Torus:|D_n(t)|>n/16\}\leq Cn^{-1}.
\]
Since the map \(\theta\mapsto2\theta\) preserves normalized Haar
measure, this proves \eqref{eq:sharp-angular-size}.
By \eqref{eq:expectation-covariance-formula}, the eigenvalues of
\(\Sigma_n(\theta)\) are
\[
 \frac12\left(1\pm\frac{|D_n(2\theta)|}{n}\right),
\]
which gives the covariance lower bound on \(\Xi_n^{\mathrm{good}}\).
Finally, \eqref{eq:sharp-rho-square} follows from Fourier orthogonality.
\end{proof}

\subsection{Uniform Gaussian comparison}
\label{subsec:expectation-uniform-Stein}

Now we establish Gaussian comparison bounds for positive powers,
uniformly over compact exponent intervals.
Throughout this subsection, fix a compact interval
\(K=[p_-,p_+]\Subset(0,\infty)\).
All Gaussian vectors below are real two-dimensional vectors.
Lemmas~\ref{lem:sharp-gaussian-smoothing}
and~\ref{lem:sharp-stein-derivative}, and
Proposition~\ref{prop:sharp-stein-comparison} make the estimates in
\cite[Lemmas~4.1--4.2 and Proposition~4.3]{ChangChenDuanFangZhangAS}
uniform over compact exponent intervals.
The Ornstein--Uhlenbeck resolvent is the standard solution of the
multivariate Gaussian Stein equation; see, for example,
\cite[Eqs.~(1.3)--(1.4)]{Gaunt2020}.
We give the details because the nonsmooth power at the origin and the
uniformity in \(p\) require estimates beyond a smooth-test bound.

\begin{lemma}\label{lem:sharp-gaussian-smoothing}
Fix \(0<c\leq1\). Suppose \(cI_2\leq\Sigma\leq I_2\),
and let \(Z_\Sigma\sim N(0,\Sigma)\). For \(\sigma>0\) and \(p\in K\), define
\(\mathsf H_{\sigma,p}^{\Sigma}(x)=\E|x+\sigma Z_\Sigma|^p\).
Then, uniformly for \(p\in K\),
\begin{equation}\label{eq:sharp-smoothing-error}
 \bigl|\mathsf H_{\sigma,p}^{\Sigma}(x)-|x|^p\bigr|
 \leq C_{K,c}\sigma^2(|x|+\sigma)^{p-2},
\end{equation}
and
\begin{equation}\label{eq:sharp-heat-fourth}
 \|\mathrm D^4\mathsf H_{\sigma,p}^{\Sigma}(x)\|
 \leq C_{K,c}(|x|+\sigma)^{p-4}.
\end{equation}
\end{lemma}

\begin{proof}
By scaling, it suffices to take \(\sigma=1\).
Write \(f(x)=|x|^p\), and let \(\varphi_\Sigma\) be
the density of \(Z_\Sigma\). Uniformly for
\(cI_2\leq\Sigma\leq I_2\),
\[
 |\partial^\nu\varphi_\Sigma(y)|
 \leq C_c(1+|y|^4)e^{-|y|^2/4},
 \qquad |\nu|\le4.
\]
For \(|x|\leq2\), the integrands are bounded by
\(C_K(1+|y|^{p_+})\) times the displayed Gaussian derivative bound.
This gives both estimates uniformly for \(p\in K\).

For \(R=|x|>2\), let \(P_{3,x}\) be the third-order
Taylor polynomial of \(y\mapsto f(x+y)\) at zero.
The first four derivatives of \(x\mapsto|x|^p\) on the unit circle
are bounded uniformly for \(p\in K\), since their coefficients
depend continuously on \(p\).  Homogeneity and Taylor's theorem
therefore give, for \(|y|\leq R/2\),
\[
 \begin{aligned}
 |f(x+y)-f(x)-\mathrm Df(x)[y]|
 &\leq C_KR^{p-2}|y|^2,\\
 |f(x+y)-P_{3,x}(y)|
 &\leq C_KR^{p-4}|y|^4.
 \end{aligned}
\]
Since \(\E Z_\Sigma=0\) and
\[
 \int_{\R^2}P_{3,x}(y)
 \partial^\nu\varphi_\Sigma(y)\dd y=0,
 \qquad |\nu|=4,
\]
we have
\[
 \begin{aligned}
 \mathsf H_{1,p}^{\Sigma}(x)-f(x)
 &=\int_{\R^2}
 [f(x+y)-f(x)-\mathrm Df(x)[y]]
 \varphi_\Sigma(y)\dd y,\\
 \partial^\nu \mathsf H_{1,p}^{\Sigma}(x)
 &=\int_{\R^2}
 [f(x+y)-P_{3,x}(y)]
 \partial^\nu\varphi_\Sigma(y)\dd y,
 \qquad |\nu|=4.
 \end{aligned}
\]
The integrals over \(|y|\leq R/2\) are bounded by
\(C_{K,c}R^{p-2}\) and \(C_{K,c}R^{p-4}\),
respectively. 
On \(|y|>R/2\), homogeneity gives
\[
\begin{aligned}
|f(x+y)-f(x)-\mathrm Df(x)[y]|
 +|f(x+y)-P_{3,x}(y)|
&\leq C_K\left(
 R^p+|y|^p+
 \sum_{j=1}^3R^{p-j}|y|^j
 \right)\\
&\leq C_KR^p(1+|y|)^{p_++3}.
\end{aligned}
\]
The Gaussian derivative bound therefore controls
both tail integrals by
\(C_{K,c}R^pe^{-R^2/32}\).

Thus
\[
 |\mathsf H_{1,p}^{\Sigma}(x)-f(x)|\leq C_{K,c}R^{p-2},
 \qquad
 \|\mathrm D^4\mathsf H_{1,p}^{\Sigma}(x)\|\leq C_{K,c}R^{p-4}.
\]
Rescaling proves \eqref{eq:sharp-smoothing-error}
and \eqref{eq:sharp-heat-fourth}.
\end{proof}

Next, we construct the Stein solution and bound its fourth derivatives, including their possible growth near the origin.
\begin{lemma}\label{lem:sharp-stein-derivative}
Under the covariance assumptions of Lemma~\ref{lem:sharp-gaussian-smoothing},
let \(0<h\leq1\), put
\[
 \beta_0=2-\min\{p_-/2,1\}\in[1,2),
\]
and, for \(p\in K\), define
\begin{equation}\label{eq:sharp-stein-solution}
 \mathsf U_{h,p}^{\Sigma}(x)=-\int_0^\infty
 \left[\mathsf P_t^{\mathrm{OU},\Sigma}\mathsf H_{h,p}^{\Sigma}(x)
 -\E\mathsf H_{h,p}^{\Sigma}(Z_\Sigma)\right]\dd t,
\end{equation}
where
\(\mathsf P_t^{\mathrm{OU},\Sigma}f(x)
=\E f(e^{-t}x+\sqrt{1-e^{-2t}}Z_\Sigma)\).
Then \(\mathsf U_{h,p}^{\Sigma}\in C^4(\R^2)\) and
\begin{equation}\label{eq:sharp-stein-equation}
 \Sigma:\mathrm D^2\mathsf U_{h,p}^{\Sigma}(x)
 -x\cdot\nabla\mathsf U_{h,p}^{\Sigma}(x)
 =\mathsf H_{h,p}^{\Sigma}(x)-\E\mathsf H_{h,p}^{\Sigma}(Z_\Sigma).
\end{equation}
Moreover, uniformly for \(p\in K\) and \(0<h\leq1\),
\begin{equation}\label{eq:sharp-uniform-Stein-envelope}
 \|\mathrm D^4\mathsf U_{h,p}^{\Sigma}(x)\|
 \leq C_{K,c}\left[(|x|+h)^{-\beta_0}+(1+|x|)^{p_+}\right].
\end{equation}
\end{lemma}

\begin{proof}
Fix \(h>0\). Gaussian convolution gives \(\mathsf H_{h,p}^{\Sigma}\in C^\infty\),
with \(\mathsf H_{h,p}^{\Sigma}\) and its first four derivatives of polynomial
growth, with bounds uniform for \(p\in K\). Put
\[
 Y_t=e^{-t}x+\sqrt{1-e^{-2t}}Z_\Sigma.
\]
Then
\[
 \begin{aligned}
 |Y_t-Z_\Sigma|
 &\leq e^{-t}(|x|+|Z_\Sigma|),\\
 \mathrm D^r\mathsf P_t^{\mathrm{OU},\Sigma}\mathsf H_{h,p}^{\Sigma}(x)
 &=e^{-rt}\E[\mathrm D^r\mathsf H_{h,p}^{\Sigma}(Y_t)],
 \qquad 1\leq r\le4.
 \end{aligned}
\]
The mean value theorem and Gaussian moment bounds imply that,
for every compact \(K_{\mathrm{sp}}\subset\R^2\),
\[
 \begin{aligned}
 &\sup_{x\in K_{\mathrm{sp}}}
 |\mathsf P_t^{\mathrm{OU},\Sigma}\mathsf H_{h,p}^{\Sigma}(x)-\E \mathsf H_{h,p}^{\Sigma}(Z_\Sigma)|
 \leq C_{h,K,c,K_{\mathrm{sp}}}e^{-t},\\
 &\sup_{x\in K_{\mathrm{sp}}}\|\mathrm D^r\mathsf P_t^{\mathrm{OU},\Sigma}\mathsf H_{h,p}^{\Sigma}(x)\|
 \leq C_{h,K,c,K_{\mathrm{sp}}}e^{-rt},
 \qquad 1\leq r\le4.
 \end{aligned}
\]
These bounds justify convergence and differentiation under
the integral in \eqref{eq:sharp-stein-solution}, giving \(\mathsf U_{h,p}^{\Sigma}\in C^4\).

Let \(\mathcal L_{\mathrm{OU},\Sigma}\) denote the Ornstein--Uhlenbeck generator
\[
 \mathcal L_{\mathrm{OU},\Sigma} f(x)
 =\Sigma:\mathrm D^2f(x)-x\cdot\nabla f(x).
\]
Gaussian integration by parts
(\cite[Theorem~5.1.8]{Bogachev1998}) yields
\[
 \begin{aligned}
 \partial_t\mathsf P_t^{\mathrm{OU},\Sigma}\mathsf H_{h,p}^{\Sigma}(x)
 =-e^{-t}x\cdot\E[\nabla\mathsf H_{h,p}^{\Sigma}(Y_t)]+e^{-2t}\E[\Sigma:\mathrm D^2\mathsf H_{h,p}^{\Sigma}(Y_t)]=\mathcal L_{\mathrm{OU},\Sigma}\mathsf P_t^{\mathrm{OU},\Sigma}\mathsf H_{h,p}^{\Sigma}(x).
 \end{aligned}
\]
Consequently,
\[
 \begin{aligned}
 \mathcal L_{\mathrm{OU},\Sigma}\mathsf U_{h,p}^{\Sigma}(x)
 =-\int_0^\infty\partial_t\mathsf P_t^{\mathrm{OU},\Sigma}\mathsf H_{h,p}^{\Sigma}(x)\dd t=\mathsf H_{h,p}^{\Sigma}(x)-\E\mathsf H_{h,p}^{\Sigma}(Z_\Sigma),
 \end{aligned}
\]
which proves \eqref{eq:sharp-stein-equation}.

Gaussian convolution also gives
\[
 \mathsf P_t^{\mathrm{OU},\Sigma}\mathsf H_{h,p}^{\Sigma}(x)=\mathsf H_{\sigma_t,p}^{\Sigma}(e^{-t}x),
 \qquad
 \sigma_t^2=h^2+1-e^{-2t}.
\]
By \eqref{eq:sharp-heat-fourth},
\begin{equation}\label{eq:sharp-stein-integral}
 \|\mathrm D^4\mathsf U_{h,p}^{\Sigma}(x)\|
 \leq C_{K,c}\int_0^\infty e^{-4t}
 (e^{-t}|x|+\sigma_t)^{p-4}\dd t.
\end{equation}

For every \(a>0\) and \(p\in K\),
\(a^{p-4}\leq a^{-\beta_0-2}+a^{p_+}\).
Put
\[
b_t(x,h)=e^{-t}|x|+\sqrt{h^2+1-e^{-2t}}.
\]
Using \eqref{eq:sharp-stein-integral},
we obtain, for \(0<h\leq1\),
\[
\|\mathrm D^4\mathsf U_{h,p}^{\Sigma}(x)\|
\leq C_{K,c}\int_0^\infty e^{-4t}
\left[b_t(x,h)^{-\beta_0-2}+b_t(x,h)^{p_+}\right]\,\dd t.
\]
For \(0\leq t\leq1\),
\(b_t(x,h)\geq c_0\sqrt{(|x|+h)^2+t}\), with an absolute
\(c_0>0\).  Hence
\begin{align*}
\int_0^1 e^{-4t}b_t(x,h)^{-\beta_0-2}\,\dd t
&\leq C_K\int_0^\infty
\bigl((|x|+h)^2+t\bigr)^{-(\beta_0+2)/2}\,\dd t\\
&=\frac{2C_K}{\beta_0}(|x|+h)^{-\beta_0}.
\end{align*}
For \(t\geq1\), \(b_t(x,h)\) is bounded below by a positive
absolute constant.  Also \(b_t(x,h)\leq |x|+\sqrt2\) for all
\(t\geq0\).  These bounds prove
\eqref{eq:sharp-uniform-Stein-envelope}.
\end{proof}

\begin{proposition}\label{prop:sharp-stein-comparison}
Let \(W=\sum_{j=1}^n\varepsilon_ja_j\in\R^2\), where
\(|a_j|=n^{-1/2}\), and suppose
\(\Sigma=\sum_j a_ja_j^{\mathsf T}\geq cI_2\) for some fixed \(c>0\).
Then, with \(Z_\Sigma\sim N(0,\Sigma)\),
\begin{equation}\label{eq:sharp-stein-comparison}
 \sup_{p\in K}\left|\E|W|^p-\E|Z_\Sigma|^p\right|
 \leq\frac{C_{K,c}}n.
\end{equation}
\end{proposition}

\begin{proof}
Since \(\operatorname{tr}\Sigma=1\), Khintchine's inequality and
Gaussian moment bounds give
\[
 \sup_{p\in K}\left(\E|W|^p+\E|Z_\Sigma|^p\right)\leq C_K.
\]
Thus it suffices to consider \(n\geq2/c\).

Put \(h=n^{-1/2}\) and \(W_j=W-\varepsilon_ja_j\). Then
\[
 \Cov(W_j)=\Sigma-a_ja_j^{\mathsf T}
 \geq(c-h^2)I_2\geq\frac c2I_2.
\]
With \(\beta_0\) as in Lemma~\ref{lem:sharp-stein-derivative},
Corollary~\ref{cor:sharp-negative}, applied with \(L=1\) and
\(\lambda=c/2\), and Khintchine's inequality give
\[
 \sup_j\sup_{|u|\leq1}
 \E\left[(|W_j+ua_j|+h)^{-\beta_0}
 +(1+|W_j+ua_j|)^{p_+}\right]\leq C_{K,c}.
\]
The same bound holds with \(W_j+ua_j\) replaced by \(W\).
For \(0<h\leq1\) and \(p\in K\),
\[
 (|x|+h)^{p-2}
 \leq (|x|+h)^{-\beta_0}+C_K(1+|x|)^{p_+}.
\]
Consequently, \eqref{eq:sharp-uniform-Stein-envelope} and
\eqref{eq:sharp-smoothing-error} yield
\begin{align*}
 \sup_{p\in K}\sup_j\sup_{|u|\leq1}
 \E\|\mathrm D^4\mathsf U_{h,p}^{\Sigma}(W_j+ua_j)\|
 &\leq C_{K,c},\\
 \sup_{p\in K}\E\left|\mathsf H_{h,p}^{\Sigma}(W)-|W|^p\right|
 &\leq C_{K,c}h^2.
\end{align*}

For \(g\in C^4([-1,1])\), integration by parts gives
\begin{equation}\label{eq:sharp-trapezoid-identity}
 \frac{g''(1)+g''(-1)-g'(1)+g'(-1)}2
 =\frac14\int_{-1}^1(1-u^2)g^{(4)}(u)\dd u.
\end{equation}
For each \(p\in K\), set \(g_j(u)=\mathsf U_{h,p}^{\Sigma}(W_j+ua_j)\).
Since \(W_j\) and \(\varepsilon_j\) are independent,
\eqref{eq:sharp-stein-equation} and
\eqref{eq:sharp-trapezoid-identity} give
\begin{equation}\label{eq:sharp-smoothed-comparison}
\begin{aligned}
\left|\E \mathsf H_{h,p}^{\Sigma}(W)-\E \mathsf H_{h,p}^{\Sigma}(Z_\Sigma)\right|&=
\left|\frac14\sum_{j=1}^n
 \int_{-1}^1(1-u^2)
 \E\bigl[\mathrm D^4\mathsf U_{h,p}^{\Sigma}(W_j+ua_j)
 [a_j,a_j,a_j,a_j]\bigr]\dd u\right|\\
&\leq C_{K,c}\sum_{j=1}^n|a_j|^4
 =C_{K,c}n^{-1}.
\end{aligned}
\end{equation}

By Gaussian convolution,
\begin{align*}
 \sup_{p\in K}
 \left|\E\mathsf H_{h,p}^{\Sigma}(Z_\Sigma)-\E|Z_\Sigma|^p\right|=
 \sup_{p\in K}\left((1+h^2)^{p/2}-1\right)\E|Z_\Sigma|^p
 \leq C_Kh^2.
\end{align*}
Taking the supremum over \(p\in K\) in
\eqref{eq:sharp-smoothed-comparison} and combining the two
smoothing errors with \(h^2=n^{-1}\) proves
\eqref{eq:sharp-stein-comparison}.
\end{proof}

The following estimate for the effect of covariance anisotropy makes
\cite[Lemma~4.4]{ChangChenDuanFangZhangAS} uniform over compact
exponent intervals.  Its quadratic dependence on the anisotropy
comes from symmetry.

\begin{lemma}\label{lem:sharp-elliptic-gaussian}
Let \(\Sigma\) be a positive semidefinite \(2\times2\) matrix
with eigenvalues \((1+r)/2,(1-r)/2\), where \(0\leq r\leq1\),
and let \(Z_\Sigma\sim N(0,\Sigma)\). Then
\begin{equation}\label{eq:sharp-elliptic-gaussian}
 \sup_{p\in K}\left|\E|Z_\Sigma|^p-m_p\right|\leq C_Kr^2.
\end{equation}
\end{lemma}

\begin{proof}
Let \(U,V\) be independent standard normal variables.
For \(-1\leq r\leq1\), set
\[
 Q_r=\frac{1+r}{2}U^2+\frac{1-r}{2}V^2,
 \qquad F_p(r)=\E Q_r^{p/2}.
\]
Orthogonal diagonalization gives \(\E|Z_\Sigma|^p=F_p(r)\).
Since \((U^2+V^2)/2\) has the exponential distribution with mean
one, \(F_p(0)=m_p\). Interchanging \(U\) and \(V\) shows that
\(F_p\) is even.

For \(|r|\leq1/2\),
\[
 \frac14(U^2+V^2)\leq Q_r\leq\frac34(U^2+V^2).
\]
The integrand defining \(F_p\) and its first two derivatives in
\(r\) are bounded in absolute value, uniformly for \(p\in K\), by
\[
 C_K\left[1+(U^2+V^2)^{p_+/2}\right].
\]
This integrable bound justifies differentiation under the expectation,
and hence
\begin{align*}
 F_p''(r)
 =\frac{p(p-2)}{16}
   \E\left[(U^2-V^2)^2Q_r^{p/2-2}\right],\qquad
 \sup_{\substack{p\in K\\|r|\leq1/2}}|F_p''(r)|\leq C_K.
\end{align*}
Evenness gives \(F_p'(0)=0\), so Taylor's theorem yields
\[
 \sup_{p\in K}|F_p(r)-m_p|\leq C_Kr^2,
 \qquad 0\leq r\leq1/2.
\]
For \(1/2\leq r\leq1\), use \(Q_r\leq U^2+V^2\),
the uniform moment bound over \(p\in K\), and \(r^2\geq1/4\).
\end{proof}

\subsection{Uniform deterministic centering}
\label{subsec:expectation-fractional-transfer-centering}

We next show that the error in replacing the exact expectations by \(m_p\) vanishes uniformly in the exponent after multiplication by \(\sqrt n\).
First we integrate the uniform comparison estimates over the angle.
This is the argument for the fixed-\(p\) expectation estimate in
\cite[Theorem~1.2]{ChangChenDuanFangZhangAS}; here every bound is uniform
on \(K\).

\begin{proposition}
\label{prop:uniform-deterministic-expectation}
For every compact interval
\(K\Subset(0,\infty)\), there exists \(C_K<\infty\) such that, for every
\(n\geq2\),
\begin{equation}
\sup_{p\in K}
\left|\E T_n^{(p)}-m_p\right|
\leq\frac{C_K}{n},
\label{eq:uniform-expectation-envelope}
\end{equation}
where \(m_p\) is defined in \eqref{eq:intro-first-order-mean}, i.e.,
\[
m_p=\Gamma\left(1+\frac p2\right).
\]
\end{proposition}

\begin{proof}
Fix \(K=[p_-,p_+]\Subset(0,\infty)\).
For \(\theta\in\Xi_n^{\mathrm{good}}\), Lemma~\ref{lem:sharp-angular}
gives \(\Sigma_n(\theta)\geq15I_2/32\).
Proposition~\ref{prop:sharp-stein-comparison} and
Lemma~\ref{lem:sharp-elliptic-gaussian}, with
\(r=|D_n(2\theta)|/n\), therefore give
\[
 \sup_{p\in K}\left|\E|X_{n,\theta}|^p-m_p\right|
 \leq C_K\left(\frac1n+\frac{|D_n(2\theta)|^2}{n^2}\right).
\]
On \(\Xi_n^{\mathrm{bad}}\), Lemma~\ref{lem:sharp-positive},
\(|x|^p\leq1+|x|^{p_+}\), and boundedness of \(m_p\) on \(K\)
give a uniform bound for the same difference.
Integrating over the angle and using
\eqref{eq:sharp-angular-size} and \eqref{eq:sharp-rho-square},
we obtain
\[
\begin{aligned}
 \sup_{p\in K}|\E T_n^{(p)}-m_p|
 \leq C_K\left[
 \frac1n+\mu(\Xi_n^{\mathrm{bad}})
 +\int_{\Torus}\frac{|D_n(2\theta)|^2}{n^2}\,\dd\mu(\theta)
 \right]\leq\frac{C_K}{n}.
\end{aligned}
\]
\end{proof}

For the fixed logarithmic-time window, \(n\asymp_T N\), so
\eqref{eq:uniform-expectation-envelope} also gives
\[
 \max_{n_-\leq n\leq n_+}\sup_{p\in K}
 \sqrt n\,|\E T_n^{(p)}-m_p|
 \leq C_{K,T}N^{-1/2}\longrightarrow0.
\]
This is the deterministic-centering error used in
Subsection~\ref{subsec:assembly-deterministic-centering}.

\subsection{Averaged small-ball estimates and moment bounds}
\label{subsec:expectation-small-ball-sublevel}

We first average the small-ball estimate over the angle.  We then bound
moments of the measure of angles at which the intensity is small, and of
angular integrals whose integrands have polynomial growth in the intensity.

\begin{lemma}
\label{lem:expectation-averaged-small-ball}
There exists \(C<\infty\) such that, for every \(n\geq2\) and
\(r\geq0\),
\begin{equation}
\int_{\Torus}
\Prob\!\left(
|X_{n,\theta}|\leq r
\right)
\,\dd\mu(\theta)
\leq
C
\left(
r^2+\frac1n
\right).
\label{eq:expectation-averaged-small-ball}
\end{equation}
Equivalently, for every \(u\geq0\),
\begin{equation}
\E I_n\!\left(\ind_{[0,u]}\right)
=
\int_{\Torus}
\Prob\!\left(
S_n(\theta)\leq u
\right)
\,\dd\mu(\theta)
\leq
C
\left(
u+\frac1n
\right).
\label{eq:expectation-averaged-intensity-small-ball}
\end{equation}
\end{lemma}

\begin{proof}
For \(\theta\in\Xi_n^{\mathrm{good}}\),
Lemma~\ref{lem:sharp-angular} gives a fixed positive lower bound for
\(\Sigma_n(\theta)\).  Apply Lemma~\ref{lem:sharp-smallball} with
\(d=2\), \(L=1\), and \(h=n^{-1/2}\).  For \(r\geq h\),
\[
\Prob(|X_{n,\theta}|\leq r)\leq Cr^2,
\qquad \theta\in\Xi_n^{\mathrm{good}}.
\]
The probability on \(\Xi_n^{\mathrm{bad}}\) is at most one,
and \(\mu(\Xi_n^{\mathrm{bad}})\leq C/n\) by
\eqref{eq:sharp-angular-size}.  Averaging proves
\eqref{eq:expectation-averaged-small-ball} for \(r\geq h\).
For \(0\leq r<h\), use monotonicity and the bound at \(r=h\).
The relation \(S_n(\theta)=|X_{n,\theta}|^2\) and Tonelli's theorem
give \eqref{eq:expectation-averaged-intensity-small-ball}.
\end{proof}

\begin{proposition}
\label{prop:gae-occupation}
For every \(q\geq1\),
\begin{equation}
\left\|
I_n\!\left(\ind_{[0,u]}\right)
\right\|_{L^q}
\leq
C_q
\left(
u+\frac1n
\right),
\qquad
n\geq2,\quad u\geq0.
\label{eq:gae-occupation}
\end{equation}
For every \(q\geq1\) and \(M\geq0\),
\begin{equation}
\left\|
I_n\!\left((1+x)^M\right)
\right\|_{L^q}
\leq
C_{q,M},
\qquad n\geq2.
\label{eq:gae-polynomial-moments}
\end{equation}
\end{proposition}

\begin{proof}
For \eqref{eq:gae-occupation}, the case \(q=1\) follows from
\eqref{eq:expectation-averaged-intensity-small-ball}.
Assume \(q>1\) and set \(m=\lceil q\rceil\).
Since \(S_n(-\theta)=S_n(\theta)\), it is enough in the moment
calculation below to work with angular variables in \([0,\pi]\).

Choose \(L=L_m\) sufficiently large.
For \(1\leq k\leq m\), consider angles
\[
\theta_1,\ldots,\theta_k\in[0,\pi]
\]
satisfying
\begin{equation}
\begin{gathered}
\frac Ln
\leq
\theta_j
\leq
\pi-\frac Ln,
\qquad
1\leq j\leq k,
\\
|\theta_j-\theta_\ell|
\geq
\frac Ln,
\qquad
j\neq\ell.
\end{gathered}
\label{eq:gae-separated-angles}
\end{equation}
Set
\[
V=(X_{n,\theta_1},\ldots,X_{n,\theta_k})\in\R^{2k},
\qquad
\Sigma=\E[VV^{\mathsf T}].
\]
The diagonal blocks of \(\Sigma\) are \(\Sigma_n(\theta_j)\), while the
entries of every off-diagonal block are linear combinations of the
real and imaginary parts of
\(n^{-1}D_n(\theta_j-\theta_\ell)\) and
\(n^{-1}D_n(\theta_j+\theta_\ell)\).
Thus
\eqref{eq:expectation-covariance-deviation},
\eqref{eq:expectation-Dirichlet-pointwise}, and
\eqref{eq:gae-separated-angles} give
\[
\left\|
\Sigma-\frac12I_{2k}
\right\|_{\mathrm{op}}
\leq
\frac{C_m}{L}.
\]
Choosing \(L=L_m\) sufficiently large, we obtain
\begin{equation}
\Sigma\geq\frac14I_{2k}.
\label{eq:gae-covariance}
\end{equation}
Write
\[
 V=\sum_{a=0}^{n-1}\varepsilon_a A_a,
 \qquad
 A_a=n^{-1/2}
 \bigl(v_a(\theta_1),\ldots,v_a(\theta_k)\bigr)
 \in\R^{2k}.
\]
Then \(|A_a|=\sqrt{k/n}\).
Apply Lemma~\ref{lem:sharp-smallball} in dimension \(2k\), with
\[
 h=n^{-1/2},\qquad L=\sqrt k,\qquad \lambda=1/4,
 \qquad r=\max\{\sqrt{ku},n^{-1/2}\}.
\]
The covariance hypothesis is \eqref{eq:gae-covariance}, and the
coefficient hypothesis follows from the displayed norm of \(A_a\).
Moreover,
\[
 \{S_n(\theta_j)\leq u,\ 1\leq j\leq k\}
 \subseteq\{|V|\leq\sqrt{ku}\}\subseteq\{|V|\leq r\}.
\]
Since \(k\leq m\) is fixed, Lemma~\ref{lem:sharp-smallball} gives
\[
 \Prob(|V|\leq r)
 \leq C_m\min(1,r^{2k})
 \leq C_m(u+n^{-1})^k.
\]
Consequently,
\begin{equation}
\Prob\!\left(
S_n(\theta_j)\leq u,\ 1\leq j\leq k
\right)
\leq
C_m
\left(
u+\frac1n
\right)^k.
\label{eq:gae-separated-small-ball}
\end{equation}

We now expand the \(m\)-th moment.
By the symmetry \(S_n(-\theta)=S_n(\theta)\),
\[
\begin{aligned}
\E
I_n\!\left(\ind_{[0,u]}\right)^m
=
\frac1{\pi^m}
\int_{[0,\pi]^m}
\Prob\!\left(
S_n(\theta_j)\leq u,
\quad
1\leq j\leq m
\right)
\,\dd\theta_1\cdots\dd\theta_m.
\end{aligned}
\]
For a given
\[
(\theta_1,\ldots,\theta_m)\in[0,\pi]^m,
\]
choose a maximal subcollection
\[
\theta_{i_1},\ldots,\theta_{i_k}
\]
satisfying
\eqref{eq:gae-separated-angles}.
By maximality, every remaining angle lies either within distance
\(L/n\) of \(0\) or \(\pi\), or within distance \(L/n\) of one of the
selected angles.
Consequently, once the \(k\) selected angles are fixed, the admissible
set for each of the remaining \(m-k\) angles has measure
\(O_m(n^{-1})\).

For each \(k\), we sum over all choices of the \(k\) selected
indices and, for each remaining angle, over the
\(L/n\)-neighbourhoods of \(0\), \(\pi\), and the selected
angles. There are at most
\(\binom{m}{k}(k+2)^{m-k}\) such choices, so we absorb this factor into \(C_m\).
Using
\eqref{eq:gae-separated-small-ball}
for the selected angles, with the probability factor interpreted
as \(1\) when \(k=0\), we therefore obtain
\[
\E
I_n\!\left(\ind_{[0,u]}\right)^m
\leq
C_m
\sum_{k=0}^{m}
n^{-(m-k)}
\left(
u+\frac1n
\right)^k
\leq
C_m
\left(
u+\frac1n
\right)^m.
\]
After enlarging \(C_m\) to cover the finitely many smaller values of
\(n\), we obtain
\[
\left\|
I_n(\ind_{[0,u]})
\right\|_{L^m}
\leq
C_m\left(u+\frac1n\right).
\]
Since \(m\geq q\), monotonicity of \(L^p\)-norms proves
\eqref{eq:gae-occupation}.

For \eqref{eq:gae-polynomial-moments}, let \(q\geq1\)
and \(M\geq0\). Each coordinate of
\(X_{n,\theta}\) is uniformly subgaussian, and hence
\[
\Prob(S_n(\theta)>x)
\leq
4e^{-x/4},
\qquad x\geq0,
\]
uniformly in \(n\) and \(\theta\).
Therefore
\[
\sup_{n\geq2}\sup_{\theta\in\Torus}
\E(1+S_n(\theta))^R<\infty
\]
for every fixed \(R\geq0\).
Jensen's inequality, applied with \(R=Mq\), then gives
\eqref{eq:gae-polynomial-moments}.
\end{proof}

\subsection{Abel-block estimates}
\label{subsec:expectation-one-time-Abel-blocks}

Now we combine the occupation bounds, sign flips, and derivative estimates to
control the ordinary and logarithmic Abel blocks at each fixed length.
Recall from \eqref{eq:calculus-Abel-blocks} the Abel blocks
\(f_{A,\alpha}\) and \(\widetilde f_{A,\alpha}\).
The bounds \eqref{eq:gae-occupation} and
\eqref{eq:gae-polynomial-moments}, together with an inequality for
coordinate sign flips and derivative estimates valid for all intensities,
will prove Proposition~\ref{prop:global-Abel-endpoint}.

\begin{lemma}
\label{lem:gae-weighted-occupation}
Let \(r\geq1\), let \(s\) vary in a fixed compact interval of
\(\mathbb R\), and let \(k\in\{0,1,2\}\).
Uniformly for \(n\geq2\) and
\(
n^{-1}\leq u\leq 1/2,
\)
one has
\begin{equation}
\left\|
I_n\!\left(
(x+u)^s
\bigl(1+|\log(x+u)|\bigr)^k
\ind_{\{x\leq1\}}
\right)
\right\|_{L^r}
\leq
C
\int_0^1
(x+u)^s
\bigl(1+|\log(x+u)|\bigr)^k
\,\dd x.
\label{eq:gae-weighted-occupation}
\end{equation}
The constant is uniform over the stated compact interval of \(s\).
\end{lemma}

\begin{proof}
Choose a dyadic number
\[
u_*\in[u,2u]\cap(0,1/2].
\]
Partition \([0,1]\) into \([0,u_*]\) and the intervals
\((R,2R]\), where \(R\) ranges over the dyadic numbers satisfying
\[
u_*\leq R\leq\frac12.
\]
On the first interval,
\[
\sup_{0\leq x\leq u_*}
(x+u)^s
\bigl(1+|\log(x+u)|\bigr)^k
\leq
C
u^s
\bigl(1+\log(u^{-1})\bigr)^k,
\]
while \eqref{eq:gae-occupation} gives
\[
\left\|
I_n(\ind_{[0,u_*]})
\right\|_{L^r}
\leq
Cu.
\]
On every interval \((R,2R]\),
\[
\sup_{R\leq x\leq2R}
(x+u)^s
\bigl(1+|\log(x+u)|\bigr)^k
\leq
C
R^s
\bigl(1+|\log R|\bigr)^k,
\]
and again by \eqref{eq:gae-occupation},
\[
\left\|
I_n(\ind_{[0,2R]})
\right\|_{L^r}
\leq
CR.
\]
On each of these intervals the supremum of the weight is bounded by a
fixed multiple of its infimum, uniformly for \(s\) in the prescribed
compact interval and \(k\in\{0,1,2\}\).
Minkowski's inequality therefore gives
\eqref{eq:gae-weighted-occupation}.
\end{proof}

\begin{lemma}
\label{lem:gae-sign-flip}
Let \(f\in C^2([0,\infty))\) be real-valued, and let
\(H_{n,f}\) be defined by
\eqref{eq:gae-local-jacobian-envelope}, i.e., 
\[
H_{n,f}(x)
=
\sup_{\substack{r\geq0\\|r-\sqrt x|\leq2n^{-1/2}}}
\bigl(|f'(r^2)|+2r^2|f''(r^2)|\bigr).
\]
Then, for every \(q\geq2\),
\begin{equation}
\|X_n(f)\|_{L^q}
\leq C_q\left(
\|I_n(x|f'(x)|^2)\|_{L^{q/2}}^{1/2}
+\|I_n(H_{n,f})\|_{L^q}
\right).
\label{eq:gae-sign-flip}
\end{equation}
\end{lemma}

\begin{proof}
We use the dimension-free scalar \(L^q\)-Poincar\'e inequality
\eqref{eq:endpoint-product-space-square-function}, with the
half-difference convention there.  It holds for every finite
\(q\geq2\); see \cite[Eqs.~(2), (3), and (5)]{CorderoErausquinEskenazis2024}.
This inequality controls the required higher moments, as well as the
variance.  The argument using sign flips and Bessel's inequality also
appears in the variance estimates of
\cite[Section~5]{ChangChenDuanFangZhangAS}.
We verify the remainder needed here for general \(f\).
Set \(F_n(f)=\sqrt n\,I_n(f)\) and, for \(y\in\R^2\), put
\[
\Phi_f(y)=yf'(|y|^2).
\]
For a sign vector \(\varepsilon\), let \(\varepsilon^{(k)}\) denote the
vector with only its \(k\)th sign reversed, and set
\[
\Delta_kF_n(f)
=\tfrac12\bigl(F_n(f)(\varepsilon)
-F_n(f)(\varepsilon^{(k)})\bigr).
\]
We use the real-vector form of the coordinate path in
Lemma~\ref{lem:exact-endpoint-sign-flip}.
Along the path
\[
X_{n;k,t}(\theta)
=X_{n,\theta}+(t-\varepsilon_k)n^{-1/2}v_k(\theta),
\qquad -1\leq t\leq1,
\]
the chain rule yields
\[
\Delta_kF_n(f)
=\varepsilon_k\int_{-1}^1\int_{\Torus}
\Phi_f(X_{n;k,t}(\theta))\cdot v_k(\theta)
\,\dd\mu(\theta)\,\dd t.
\]
The real differential of \(\Phi_f\) satisfies
\[
\|\mathrm D\Phi_f(y)\|_{\mathrm{op}}
\leq |f'(|y|^2)|+2|y|^2|f''(|y|^2)|.
\]
Since \(|X_{n;k,t}(\theta)-X_{n,\theta}|\leq2n^{-1/2}\), the mean-value
formula yields
\begin{equation}
\begin{split}
&\Delta_kF_n(f)
=2\varepsilon_k c_k+E_k,\\
&c_k
=
\int_{\Torus}
\Phi_f(X_{n,\theta})\cdot v_k(\theta)
\,\dd\mu(\theta),\\
&|E_k|\leq4n^{-1/2}I_n(H_{n,f}).
\end{split}
\label{eq:gae-flip-remainder}
\end{equation}
The functions \(v_k\) are orthonormal in \(L^2(\Torus;\R^2)\), since
\[
\int_{\Torus}v_k(\theta)\cdot v_\ell(\theta)\,\dd\mu(\theta)
=\int_{\Torus}\cos((k-\ell)\theta)\,\dd\mu(\theta)
=\ind_{\{k=\ell\}}.
\]
Bessel's inequality therefore yields
\[
\sum_{k=0}^{n-1}|c_k|^2
\leq I_n(x|f'(x)|^2).
\]
Consequently, by \eqref{eq:gae-flip-remainder},
\[
\begin{aligned}
\sum_{k=0}^{n-1}|\Delta_kF_n(f)|^2
&\leq8\sum_{k=0}^{n-1}|c_k|^2
+32I_n(H_{n,f})^2\\
&\leq8I_n(x|f'(x)|^2)+32I_n(H_{n,f})^2.
\end{aligned}
\]
Apply \eqref{eq:endpoint-product-space-square-function} to \(F_n(f)\) and take
square roots.  Since \(F_n(f)-\E F_n(f)=X_n(f)\), this proves
\eqref{eq:gae-sign-flip}.
\end{proof}

\begin{lemma}
\label{lem:gae-abel-derivatives}
Let
\(
\mathfrak A=[\alpha_-,\alpha_+]\Subset(0,\infty).
\)
There exists \(C_{\mathfrak A}<\infty\) such that, for every
\(A\geq2\), uniformly for
\(\alpha\in\mathfrak A\) and \(x\geq0\), the Abel blocks
\(f_{A,\alpha}\) and \(\widetilde f_{A,\alpha}\) in
\eqref{eq:calculus-Abel-blocks} satisfy
\begin{align}
|f_{A,\alpha}(x)|
&\leq
\frac{C_{\mathfrak A}}{A}
(1+x)(x+A^{-1})^{\alpha-1},
\label{eq:gae-value-ordinary}\\
|\widetilde f_{A,\alpha}(x)|
&\leq
\frac{C_{\mathfrak A}}{A}
(1+x)(x+A^{-1})^{\alpha-1}
\bigl(1+|\log(x+A^{-1})|\bigr).
\label{eq:gae-value-log}
\end{align}
Moreover, for \(j\in\{1,2\}\),
\begin{align}
|\partial_x^j f_{A,\alpha}(x)|
&\leq
\frac{C_{\mathfrak A}}{A}
\left[
(x+A^{-1})^{\alpha-j}
+
|\alpha-1|(x+A^{-1})^{\alpha-j-1}
\right],
\label{eq:gae-derivative-ordinary}\\
|\partial_x^j\widetilde f_{A,\alpha}(x)|
&\leq
\frac{C_{\mathfrak A}}{A}
\left[
(x+A^{-1})^{\alpha-j}
\bigl(1+|\log(x+A^{-1})|\bigr)
\right.
\nonumber\\
&\qquad\qquad\left.
+
(x+A^{-1})^{\alpha-j-1}
\left(
1+
|\alpha-1|
\bigl(1+|\log(x+A^{-1})|\bigr)
\right)
\right].
\label{eq:gae-derivative-log}
\end{align}
\end{lemma}

\begin{proof}
Put
\[
F_{\tau,\alpha}=P_\tau h_\alpha,
\qquad
d_\tau=1-e^{-\tau},
\qquad
y=\frac{e^{-\tau}x}{d_\tau},
\]
and
\[
\pi_j(\alpha)
=
\prod_{\ell=0}^{j-1}(\alpha-\ell),
\qquad
\pi_0=1.
\]
By \eqref{eq:kummer-semigroup-power}, the classical differentiation
formula for \({}_1F_1\) and Kummer's transformation
\cite[Eqs.~13.3.16 and 13.2.39]{NISTDLMF}, for \(0\leq j\leq4\),
\begin{equation}
\partial_x^jF_{\tau,\alpha}(x)
=
e^{-j\tau}
d_\tau^{\alpha-j}
\pi_j(\alpha)
\E\left[
\frac{\Gamma(K_y+\alpha+1)}
{\Gamma(K_y+j+1)}
\right],
\qquad
K_y\sim\operatorname{Poisson}(y).
\label{eq:gae-poisson-formula}
\end{equation}
The identity is valid also at the integral values of \(\alpha\), since
\eqref{eq:gae-poisson-formula} involves no division by
\(\pi_j(\alpha)\).
Uniformly for
\(k\geq0,\)
\(\alpha\in\mathfrak A,\)
\(0\leq j\leq4,\)
the product formula for the gamma ratio gives
\[
\frac{\Gamma(k+\alpha+1)}
{\Gamma(k+j+1)}
\asymp_{\mathfrak A}
(k+1)^{\alpha-j},
\]
and
\[
\left|
\partial_\alpha
\log
\frac{\Gamma(k+\alpha+1)}
{\Gamma(k+j+1)}
-
\log(k+1)
\right|
\leq
C_{\mathfrak A}.
\]
The factorial moments of \(K_y\), together with the Poisson lower-tail
bound
\[
\Prob(K_y<y/2)\leq e^{-cy},
\qquad y\geq2,
\]
imply uniform bounds for every fixed positive and negative moment of
\((K_y+1)/(1+y)\).  Hence, uniformly for \(s\) in a fixed compact
interval,
\[
\E\!\left[
\left(\frac{K_y+1}{1+y}\right)^s
\left(
1+
\left|
\log\frac{K_y+1}{1+y}
\right|
\right)
\right]
\leq C,
\qquad y\geq0.
\]
Since
\[
d_\tau(1+y)
\asymp
x+\tau,
\qquad
0<\tau\leq\frac12,
\]
\eqref{eq:gae-poisson-formula} implies, uniformly for
\(\alpha\in\mathfrak A\) and \(0\leq j\leq4\),
\begin{align}
|\partial_x^jF_{\tau,\alpha}(x)|
&\leq
C_{\mathfrak A}
|\pi_j(\alpha)|
(x+\tau)^{\alpha-j},
\label{eq:gae-semigroup-derivatives}\\
|\partial_\alpha\partial_x^jF_{\tau,\alpha}(x)|
&\leq
C_{\mathfrak A}
(x+\tau)^{\alpha-j}
\left[
|\pi_j'(\alpha)|
+
|\pi_j(\alpha)|
\bigl(1+|\log(x+\tau)|\bigr)
\right].
\label{eq:gae-semigroup-log-derivatives}
\end{align}
Differentiation in \(\alpha\) in
\eqref{eq:gae-poisson-formula}
is justified by local uniform convergence of the differentiated
Poisson series on compact exponent intervals.
By
\eqref{eq:spectral-analytic-semigroup},
\eqref{eq:spectral-Laguerre-eigenvalue},
and
\eqref{eq:spectral-Laguerre-generator},
\[
\partial_\tau F_{\tau,\alpha}
=
\LagGen F_{\tau,\alpha}
=
x\partial_x^2F_{\tau,\alpha}
+
(1-x)\partial_xF_{\tau,\alpha}.
\]
Therefore, by \eqref{eq:calculus-Abel-blocks}, for \(j=0,1,2\),
\begin{equation}
\partial_x^j f_{A,\alpha}(x)
=
-\int_{1/(2A)}^{1/A}
\Bigl[
x\partial_x^{j+2}F_{\tau,\alpha}(x)
+
(j+1-x)\partial_x^{j+1}F_{\tau,\alpha}(x)
-
j\partial_x^jF_{\tau,\alpha}(x)
\Bigr]
\,\dd\tau.
\label{eq:gae-generator-integral}
\end{equation}
Moreover,
\eqref{eq:spectral-semigroup-p-derivative}
and
\eqref{eq:calculus-Abel-blocks}
give
\[
\widetilde f_{A,\alpha}
=
\partial_\alpha f_{A,\alpha}.
\]
For \(2\leq j\leq4\),
\[
|\pi_j(\alpha)|
\leq
C_{\mathfrak A}|\alpha-1|,
\qquad
|\pi_j'(\alpha)|
\leq
C_{\mathfrak A}.
\]
On
\[
\frac1{2A}\leq \tau\leq\frac1A,
\]
one has
\[
x+\tau\asymp x+A^{-1}.
\]
Substituting
\eqref{eq:gae-semigroup-derivatives}
and
\eqref{eq:gae-semigroup-log-derivatives}
into
\eqref{eq:gae-generator-integral}
and its \(\alpha\)-derivative gives
\eqref{eq:gae-value-ordinary}--\eqref{eq:gae-derivative-log}.
\end{proof}

\section*{Acknowledgments}

\noindent 
\textit{Funding.}
X. F. was supported by the National Science and Technology
Council of Taiwan (Grant No.~114-2115-M-A49-003-MY3).

\bigskip

\noindent 
\textit{AI Statement.}
The authors used OpenAI ChatGPT (GPT-6 Astra) for language editing, \LaTeX\ formatting, and assistance with local mathematical reasoning and numerical verification. The overall mathematical strategy and framework were developed by the authors. All mathematical content was written and checked by the authors who take full responsibility.

\clearpage
\begingroup
\small
\setlength{\tabcolsep}{4pt}
\renewcommand{\arraystretch}{1.12}
\begin{center}{\large\textsc{Index of notation}}\end{center}
\begin{longtable}{@{}>{\raggedright\arraybackslash}p{.26\textwidth}>{\raggedright\arraybackslash}p{.59\textwidth}>{\raggedright\arraybackslash}p{.11\textwidth}@{}}
\toprule
\textbf{Notation} & \textbf{Meaning} & \textbf{Location}\\
\midrule
\endfirsthead
\toprule
\textbf{Notation} & \textbf{Meaning} & \textbf{Location}\\
\midrule
\endhead
\midrule
\endfoot
\bottomrule
\endlastfoot
\addlinespace[5pt]\multicolumn{3}{@{}l}{\textit{A. Model, parameters, and logarithmic time}}\\*
$q_n$, $(\varepsilon_j)_{j\ge0}$ & Littlewood polynomial $q_n(z)=\sum_{j=0}^{n-1}\varepsilon_jz^j$, using one shared sequence of independent Rademacher signs. & \S\,\ref{sec:model-power-statistics}\\
$\mathbb T$, $\mu$ & Circle $\mathbb R/(2\pi\mathbb Z)$ with normalized measure $\mathrm d\mu=\mathrm d\theta/(2\pi)$. & \S\,\ref{sec:model-power-statistics}\\
$S_n(\theta)$ & Normalized intensity $|q_n(e^{i\theta})|^2/n$. & \S\,\ref{sec:model-power-statistics}\\
$I_n(f)$, $X_n(f)$ & $I_n(f)=\int_{\mathbb T}f(S_n)\,\mathrm d\mu$; $X_n(f)=\sqrt n\,[I_n(f)-\mathbb E I_n(f)]$. & \S\,\ref{sec:model-power-statistics}\\
$K$, $\alpha$, $\mathfrak A$ & Exponent intervals $K=[p_-,p_+]\Subset(0,\infty)$ and $\mathfrak A=[p_-/2,p_+/2]$, with rescaled exponent $\alpha=p/2$. & \S\,\ref{sec:model-power-statistics}, \S\,\ref{subsec:spectral-power-coefficients}\\
$h_p$, $h_\alpha$, $\widetilde h_\alpha$ & Power function $h_p(x)=x^{p/2}$, its form $h_\alpha(x)=x^\alpha$ for $\alpha=p/2$, and derivative with respect to $\alpha$, $\widetilde h_\alpha(x)=x^\alpha\log x$, with $\widetilde h_\alpha(0)=0$. & \S\,\ref{sec:model-power-statistics}, \S\,\ref{subsec:spectral-power-coefficients}\\
$T_n^{(p)}$, $m_p$ & Normalized power $T_n^{(p)}=I_n(h_p)=\|q_n\|_p^p/n^{p/2}$ and Gaussian mean $m_p=\Gamma(1+p/2)$. & \S\,\ref{sec:model-power-statistics}, \S\,\ref{sec:model-prelimit-fields}\\
$I_T$, $n_N(t)$, $n_-$, $n_+$ & Time window $I_T=[-T,T]$, length $n_N(t)=\lfloor Ne^t\rfloor$, and endpoints $n_\pm=\lfloor Ne^{\pm T}\rfloor$. & \S\,\ref{sec:model-log-time-interpolation}\\
$\vartheta_{N,n}$, $\operatorname{Poly}_{N,T}$ & Interpolation times ($\vartheta_{N,n_-}=-T$, otherwise $\log(n/N)$) and polygonal interpolation, constant after the last knot. & \S\,\ref{sec:model-log-time-interpolation}\\
$C(I_T;C(K))$, $D(I_T;C(K))$ & Continuous paths with the uniform norm, and right-continuous paths with left limits under the Skorokhod $J_1$ topology. & \S\,\ref{sec:model-log-time-interpolation}\\
\addlinespace[5pt]\multicolumn{3}{@{}l}{\textit{B. Prelimit fields and Gaussian limits}}\\*
$\mathrm{ec}$, $\mathrm{dc}$, $\mathrm{poly}$, $\mathrm{fl}$ & Expectation centering, deterministic centering, polygonal interpolation, and floor sampling. & \S\,\ref{sec:model-prelimit-fields}\\
$Y_n^{\mathrm{ec}}$, $Y_n^{\mathrm{dc}}$ & Centered power statistics at integer lengths: $\sqrt n\,(T_n^{(p)}-\mathbb E T_n^{(p)})$ and $\sqrt n\,(T_n^{(p)}-m_p)$. & \S\,\ref{sec:model-prelimit-fields}\\
$\mathcal Z_N^{\mathrm{ec},\mathrm{poly}}$, $\ldots$ & The four fields obtained from the two centerings $Y_n^{\mathrm{ec}}$, $Y_n^{\mathrm{dc}}$ and the two path constructions: polygonal interpolation and floor sampling. & \S\,\ref{sec:model-prelimit-fields}\\
$f_a$, $V_a$ & Convolution density $f_a=\mathbf1_{[0,1]}^{*a}$ and mode constant $V_a=\int_0^a f_a(u)^2\,\mathrm du$. & \S\,\ref{subsec:intro-main-results}\\
$c_a(p)$, $\gamma_a(\alpha)$ & Power coefficients $\gamma_a(\alpha)=\Gamma(1+\alpha)(-1)^a\binom\alpha a$, with $c_a(p)=\gamma_a(p/2)$ for $a\ge2$. & \S\,\ref{subsec:spectral-power-coefficients}\\
$G^{(a)}$ & Independent centered stationary Gaussian modes: $\operatorname{Cov}(G^{(a)}(s),G^{(a)}(t))=2V_ae^{-(a-1/2)|t-s|}$, $a\ge2$. & \S\,\ref{subsec:intro-main-results}\\
$G_d(\alpha,t)$, $G$ & $G_d=\sum_{a=2}^d\gamma_a(\alpha)G^{(a)}$; its limit is denoted $G(\alpha,t)$, or $G(p,t)$ in the coordinate $p=2\alpha$. & \S\,\ref{subsec:fixed-power-truncations}, \S\,\ref{sec:full-gaussian-field}\\
$r_p$, $\mathfrak a_p$ & Norm profile $r_p=m_p^{1/p}$ and derivative $\mathfrak a_p=p^{-1}m_p^{1/p-1}$ of $x\mapsto x^{1/p}$ with respect to $x$, at $x=m_p$. & \S\,\ref{subsec:intro-main-results}\\
$W_n$, $\mathcal W_N^{\mathrm{poly}}$, $\mathcal W_N^{\mathrm{fl}}$ & $W_n(p)=\sqrt n\,[\|q_n\|_p/\sqrt n-r_p]$ and the fields obtained by polygonal interpolation and floor sampling. & \S\,\ref{subsec:assembly-normalized-transfer}\\
$H(p,t)$ & Gaussian limit for normalized norms: $H(p,t)=\mathfrak a_pG(p,t)$. & \S\,\ref{subsec:intro-main-results}\\
\noalign{\penalty-10000}
\addlinespace[5pt]\multicolumn{3}{@{}l}{\textit{C. Laguerre expansions and semigroup}}\\*
$\nu$, $L_a$ & Exponential measure $\mathrm d\nu=e^{-x}\mathrm dx$ on $[0,\infty)$ and its orthonormal Laguerre basis, with $L_0=1$, $L_1=1-x$. & \S\,\ref{subsec:spectral-power-coefficients}\\
$L_j^{(\eta)}$ & Generalized Laguerre polynomial with parameter $\eta$; $L_j^{(0)}=L_j$. & \S\,\ref{subsec:weighted-derivative-bounds}\\
$I_j$ & Modified Bessel function of the first kind of order $j$. & \S\,\ref{subsec:block-conditional-mean}\\
$\widehat f_a$, $\widetilde\gamma_a$ & Laguerre coefficient $\widehat f_a=\int_0^\infty fL_a\,\mathrm d\nu$; logarithmic coefficient $\widetilde\gamma_a=\partial_\alpha\gamma_a$. & \S\,\ref{subsec:spectral-power-coefficients}\\
$\Pi_{\le d}^{\mathrm{Lag}}$, $h_{\alpha,\le d}$ & Projection onto Laguerre degrees $0,\ldots,d$ and truncation $h_{\alpha,\le d}=\Pi_{\le d}^{\mathrm{Lag}}h_\alpha$; $d$ is fixed before $N\to\infty$. & \S\,\ref{subsec:spectral-sharp-truncations}\\
$R_{\alpha,d}$, $\widetilde R_{\alpha,d}$ & Sharp tails $(I-\Pi_{\le d}^{\mathrm{Lag}})h_\alpha$ and $(I-\Pi_{\le d}^{\mathrm{Lag}})\widetilde h_\alpha$; $R_{p,d}$ is the introductory parametrization. & \S\,\ref{subsec:spectral-sharp-truncations}\\
$\mathfrak E_s(f)$ & Laguerre energy $\sum_{a\ge1}a^s|\widehat f_a|^2$, excluding the constant mode. & \S\,\ref{subsec:spectral-laguerre-semigroup}\\
$\mathcal L_{\mathrm{Lag}}$, $\mathcal N_{\mathrm{Lag}}$ & Generator $\mathcal L_{\mathrm{Lag}}f=xf^{\prime\prime}+(1-x)f^{\prime}$ and number operator $\mathcal N_{\mathrm{Lag}}=-\mathcal L_{\mathrm{Lag}}$. & \S\,\ref{subsec:spectral-laguerre-semigroup}\\
$P_\tau$ & Laguerre semigroup $P_\tau f=\sum_{a\ge0}e^{-a\tau}\widehat f_aL_a$. & \S\,\ref{subsec:spectral-laguerre-semigroup}\\
\addlinespace[5pt]\multicolumn{3}{@{}l}{\textit{D. Abel blocks, scales, and boundaries}}\\*
$A$, $m_A(a)$ & Abel scale $A\ge2$ and multiplier $m_A(a)=e^{-a/(2A)}-e^{-a/A}$; $a$ is the Laguerre degree. & \S\,\ref{subsec:abel-blocks-dyadic-telescoping}\\
$f_{A,\alpha}$, $\widetilde f_{A,\alpha}$ & Ordinary and logarithmic Abel blocks $(P_{1/(2A)}-P_{1/A})h_\alpha$ and $(P_{1/(2A)}-P_{1/A})\widetilde h_\alpha$. & \S\,\ref{subsec:abel-blocks-dyadic-telescoping}\\
$f_{A,\alpha}^{\circ}$, $f_{A,\alpha}^{\mathrm{ren}}$ & $L_1$ component removed, then value at zero subtracted: $f^{\mathrm{ren}}_{A,\alpha}=f^\circ_{A,\alpha}-f^\circ_{A,\alpha}(0)L_0$; tildes give logarithmic analogues. & \S\,\ref{subsec:process-equivalent-normalization}\\
$F_{A,\alpha}$, $\Delta F_{A,\alpha}^{\circ}$ & Scaled ordinary block $F_{A,\alpha}(y)=A^\alpha f_{A,\alpha}(y/A)$ and the function $\Delta F^\circ$ after subtracting $q_{A,\alpha}L_1(y/A)$ and the resulting value at zero. & \S\,\ref{subsec:scaled-Abel-functions}\\
$\mathscr D_A$, $\mathcal P_{A,\alpha}$ & Shifted derivative $\mathscr D_A=\partial_\alpha-\log A$ and scaled logarithmic function $\mathcal P_{A,\alpha}=\mathscr D_A\Delta F^\circ_{A,\alpha}$. & \S\,\ref{subsec:scaled-Abel-functions}\\
$\Lambda_N$ & Least dyadic $A\ge2$ with $A>\gamma_0\log(N+2)/\log\log(N+3)$, for a fixed sufficiently small $\gamma_0>0$. & \S\,\ref{sec:active-abel-range}\\
$B_{\Lambda,\alpha}^{-,i}$, $B_{\Lambda,\alpha}^{+,i}$ & Abel--sharp comparison terms with degrees $a\leq\Lambda$ and $a>\Lambda$, respectively; $i=0,1$ select ordinary and logarithmic functions. & \S\,\ref{subsec:Abel-sharp-boundary-comparison}\\
$g_{A,\mathbf u}$, $\mathcal E_A(\mathbf u)$ & Sharp dyadic block $\sum_{A<a\le2A}u_aL_a$ and coefficient energy $\sum_{A<a\le2A}a|u_a|^2$. & \S\,\ref{subsec:coefficient-uniform-low-blocks}\\
$Q(y)$ & Limiting logarithmic function at $\alpha=1$: $\mathcal P_{A,1}(y)\to Q(y)$ as $A\to\infty$, locally with fixed-order $y$-derivatives, directly from the Kummer formula. & \S\,\ref{subsec:kummer-mehler-heine-critical}\\
\noalign{\penalty-10000}
\addlinespace[5pt]\multicolumn{3}{@{}l}{\textit{E. Intensity localization and maximal estimates}}\\*
$\lambda_{\mathrm{low}}$, $\lambda_{\mathrm{mid}}$, $\lambda_{\mathrm{up}}$ & Smooth partition of unity separating lower, middle, and upper intensities. & \S\,\ref{subsec:spatial-localization}\\
$\chi_0$, $\omega_R$, $\varpi_R$ & Lowest-scale and lower-shell cutoffs in $y=Ax$, and upper-shell cutoffs in $x$; $R$ labels the corresponding dyadic scale. & \S\,\ref{subsec:spatial-localization}\\
$r_{A,0,\alpha}$, $r_{A,R,\alpha}$ & Lowest-scale and lower-shell pieces of $f^{\mathrm{ren}}_{A,\alpha}$; tildes denote logarithmic counterparts. & \S\,\ref{subsec:spatial-localization}\\
$f_{A,\alpha}^{\mathrm{low}}$, $f_{A,\alpha}^{\mathrm{mid}}$, $f_{A,\alpha}^{\mathrm{up}}$ & Regional parts of $f^{\mathrm{ren}}_{A,\alpha}$ obtained by multiplication by the intensity partition; tildes give logarithmic parts. & \S\,\ref{subsec:spatial-localization}, \S\,\ref{sec:active-abel-range}\\
$f_{A,R,\alpha}^{\mathrm{up}}$, $\widetilde f_{A,R,\alpha}^{\mathrm{up}}$ & Individual upper shells $\varpi_R\lambda_{\mathrm{up}}f^{\mathrm{ren}}_{A,\alpha}$ and their logarithmic analogues. & \S\,\ref{subsec:spatial-localization}\\
$\mathfrak M_{N,T;\mathfrak A}(f_\bullet)$ & Integrated maximal quantity $\bigl(\int_{\mathfrak A}\mathbb E\max_{n_-\le n\le n_+}|X_n(f_\alpha)|^2\,\mathrm d\alpha\bigr)^{1/2}$. & \S\,\ref{sec:active-abel-range}\\
$\sigma_*$, $\beta$ & Decay exponent $\sigma_*=\min\{\alpha_-,1\}$ and fixed small range exponent $\beta>0$, for $\Lambda_N\le A\le N^\beta$. & \S\,\ref{sec:active-abel-range}\\
\addlinespace[5pt]\multicolumn{3}{@{}l}{\textit{F. Chaos, comparison, and temporal refinement}}\\*
$\varepsilon_E$, $\mathsf P_r^{(n)}$ & Walsh monomial $\varepsilon_E=\prod_{j\in E}\varepsilon_j$ and projection onto Walsh degree $r$. & \S\,\ref{subs:app-A-Walsh-top-chaos}\\
$\mathcal G_{n,a}$ & Top homogeneous chaos $\mathsf P_{2a}^{(n)}X_n(L_a)$ of Laguerre mode $a$. & \S\,\ref{subs:app-A-Walsh-top-chaos}\\
$\mathcal F_n$, $M_{n,a}$ & Filtration $\sigma(\varepsilon_0,\ldots,\varepsilon_{n-1})$ and martingale $M_{n,a}=n^{a-1/2}\mathcal G_{n,a}$. & \S\,\ref{subs:app-A-balanced-covariance}\\
$\Phi_\psi$ & Radial field multiplier $\Phi_\psi(z)=\psi^{\prime}(|z|^2)z$. & \S\,\ref{subsec:weighted-derivative-bounds}\\
$\mathscr R$ & Lowering operator $\mathscr R\psi=\psi^{\prime}+x\psi^{\prime\prime}$, with $\mathscr R L_a=-aL_{a-1}$. & \S\,\ref{subsec:common-field-expansion}\\
$\rho_{m,n}$, $\tau_{m,n}$ & Relative added length $(n-m)/n$ and Laguerre semigroup time $-\log(m/n)$. & \S\,\ref{app:block-addition-refinement}, \S\,\ref{subsec:block-conditional-mean}\\
$U_{m,n}(u)$, $V_{m,n}(u)$ & Conditionally centered and centered conditional-mean terms of $X_n(u)-X_m(u)$, conditional on $\mathcal F_m$. & \S\,\ref{app:block-addition-refinement}\\
$\Pi_I$, $D_I$ & Fourier projection onto a finite frequency set $I\subset\mathbb Z$ and kernel $D_I(x)=\sum_{j\in I}e^{ijx}$. & \S\,\ref{subsec:block-fourier-projection}\\
$\mathcal C_{N,\Gamma}$, $\mathcal E_{N,\Gamma,\ell}$ & Coarse lengths and level-$\ell$ refinement edges at scale $\Gamma$; the Abel argument uses $\Gamma=A$. & \S\,\ref{subsec:block-temporal-refinement}\\
\end{longtable}\endgroup

\bibliographystyle{amsalpha}
\bibliography{CLT__p_t__process}

\end{document}